\documentclass[11pt,reqno]{amsart}
\usepackage[letterpaper,margin=1in]{geometry}
\usepackage{amsmath,amssymb,amsthm,mathtools}
\usepackage{booktabs}
\usepackage{enumitem}
\usepackage{tikz}
\usepackage{pgfplots}
\usetikzlibrary{decorations.pathreplacing}
\pgfplotsset{compat=1.18}
\definecolor{figblue}{RGB}{38,70,115}
\definecolor{figorange}{RGB}{204,122,40}
\definecolor{figgray}{RGB}{221,228,238}
\definecolor{figred}{RGB}{150,36,42}
\definecolor{figdark}{RGB}{62,66,74}
\usepackage[hidelinks]{hyperref}
\usepackage{microtype}

\newtheorem{theorem}{Theorem}[section]
\newtheorem{proposition}[theorem]{Proposition}
\newtheorem{lemma}[theorem]{Lemma}
\newtheorem{corollary}[theorem]{Corollary}

\theoremstyle{definition}
\newtheorem{definition}[theorem]{Definition}
\theoremstyle{remark}
\newtheorem{remark}[theorem]{Remark}

\newcommand{\E}{\mathbb{E}}
\renewcommand{\P}{\mathbb{P}}
\newcommand{\R}{\mathbb{R}}

\newcommand{\dd}{\mathrm{d}}
\newcommand{\1}{\mathbf{1}}
\newcommand{\cut}{\square}
\newcommand{\tr}{\operatorname{tr}}
\newcommand{\Vbal}{V_{\mathrm{SBM}}}

\pgfplotstableread{
r V m
0.020073 0.499933 0.499882
0.024141 0.499888 0.499805
0.028153 0.499830 0.499704
0.032327 0.499752 0.499571
0.036551 0.499653 0.499403
0.040692 0.499536 0.499205
0.044951 0.499392 0.498963
0.049271 0.499221 0.498677
0.053586 0.499023 0.498347
0.057827 0.498801 0.497979
0.061922 0.498560 0.497579
0.066295 0.498271 0.497103
0.070429 0.497967 0.496604
0.074806 0.497613 0.496022
0.078846 0.497253 0.495434
0.083088 0.496843 0.494763
0.087539 0.496373 0.493997
0.092207 0.495837 0.493124
0.096384 0.495318 0.492279
0.100727 0.494739 0.491335
0.105239 0.494092 0.490281
0.109922 0.493371 0.489106
0.113956 0.492709 0.488025
0.118109 0.491985 0.486842
0.122381 0.491197 0.485549
0.126771 0.490338 0.484137
0.131277 0.489404 0.482596
0.135897 0.488389 0.480916
0.140626 0.487290 0.479086
0.145459 0.486101 0.477093
0.150392 0.484817 0.474925
0.154404 0.483718 0.473056
0.158470 0.482554 0.471060
0.162584 0.481323 0.468929
0.166742 0.480022 0.466656
0.170935 0.478651 0.464234
0.175156 0.477209 0.461653
0.179397 0.475696 0.458906
0.183648 0.474113 0.455983
0.187898 0.472460 0.452877
0.192135 0.470740 0.449578
0.196347 0.468955 0.446078
0.200518 0.467112 0.442365
0.204632 0.465215 0.438432
0.208673 0.463272 0.434269
0.209670 0.462780 0.433190
}\Vbaldata
\pgfplotstableread{
r h
0.020000 0.901760
0.022750 0.891182
0.025500 0.880939
0.028250 0.870992
0.031000 0.861313
0.033750 0.851879
0.036500 0.842669
0.039250 0.833666
0.042000 0.824856
0.044750 0.816227
0.047500 0.807766
0.050250 0.799465
0.053000 0.791314
0.055750 0.783307
0.058500 0.775434
0.061250 0.767692
0.064000 0.760072
0.066750 0.752571
0.069500 0.745183
0.072250 0.737904
0.075000 0.730730
0.077750 0.723657
0.080500 0.716680
0.083250 0.709798
0.086000 0.703007
0.088750 0.696304
0.091500 0.689685
0.094250 0.683150
0.097000 0.676695
0.099750 0.670317
0.102500 0.664016
0.105250 0.657788
0.108000 0.651633
0.110750 0.645547
0.113500 0.639529
0.116250 0.633579
0.119000 0.627693
0.121750 0.621871
0.124500 0.616111
0.127250 0.610411
0.130000 0.604771
0.132750 0.599190
0.135500 0.593665
0.138250 0.588196
0.141000 0.582782
0.143750 0.577421
0.146500 0.572113
0.149250 0.566857
0.152000 0.561651
0.154750 0.556495
0.157500 0.551388
0.160250 0.546329
0.163000 0.541317
0.165750 0.536352
0.168500 0.531432
0.171250 0.526557
0.174000 0.521726
0.176750 0.516939
0.179500 0.512194
0.182250 0.507492
0.185000 0.502831
0.187750 0.498211
0.190500 0.493631
0.193250 0.489091
0.196000 0.484590
0.198750 0.480128
0.201500 0.475704
0.204250 0.471317
0.207000 0.466967
0.209750 0.462654
0.212500 0.458377
0.215250 0.454136
0.218000 0.449929
0.220750 0.445758
0.223500 0.441620
0.226250 0.437517
0.229000 0.433446
0.231750 0.429409
0.234500 0.425405
0.237250 0.421432
0.240000 0.417492
0.242750 0.413583
0.245500 0.409705
0.248250 0.405859
0.251000 0.402042
0.253750 0.398256
0.256500 0.394499
0.259250 0.390772
0.262000 0.387074
0.264750 0.383405
0.267500 0.379765
0.270250 0.376153
0.273000 0.372569
0.275750 0.369012
0.278500 0.365483
0.281250 0.361981
0.284000 0.358506
0.286750 0.355058
0.289500 0.351636
0.292250 0.348240
0.295000 0.344870
0.297750 0.341526
0.300500 0.338207
0.303250 0.334913
0.306000 0.331644
0.308750 0.328400
0.311500 0.325180
0.314250 0.321985
0.317000 0.318813
0.319750 0.315666
0.322500 0.312542
0.325250 0.309442
0.328000 0.306365
0.330750 0.303311
0.333500 0.300280
0.336250 0.297271
0.339000 0.294285
0.341750 0.291321
0.344500 0.288380
0.347250 0.285460
0.350000 0.282562
}\hdata

\begin{document}

\title[The sparse triangle lower tail]{The lower tail of triangle counts in sparse random graphs}

\author{Dong Huang}
\address{Department of Statistics and Data Science, Tsinghua University}
\email{hd23@mails.tsinghua.edu.cn}

\author{Shuangping Li}
\address{Department of Statistics and Data Science, Yale University}
\email{shuangping.li@yale.edu}


\begin{abstract}
Let $T_n$ be the number of triangles in the Erdős–Rényi random graph
$G(n,p)$ with $n^{-1/2}\ll p\ll1$ and fix $\eta\in(0,1)$. We study the
large deviation rate of the lower tail event $T_n\le\eta\E[T_n]$, and the
typical structure of $G(n,p)$ conditioned on this event. By earlier
work~\cite{Chi25,CD20,CDP,KS23}, both the rate and the structure are governed
by a variational problem over graphons. In this work we solve this variational
problem for every $\eta\in(0,1)$. We show that every minimizer is either a
constant graphon or a step graphon with two blocks of equal size.
Correspondingly, the conditioned graph is close to an Erdős–Rényi graph or
to a symmetric stochastic block model with two blocks.
These two regimes are separated by a first order phase transition at the threshold $\eta_*=0.009217\ldots$.
This resolves a conjecture of Zhao~\cite{Zha17} and an
open problem of Janson and Warnke~\cite{JW16} on triangles. It thereby gives a
complete description of the lower tail of triangle counts in this range of
$p$. 
\end{abstract}

\maketitle

\setcounter{tocdepth}{1}
\tableofcontents

\section{Introduction}\label{sec:intro}

\subsection{Lower-tail large deviations for triangle counts}\label{s1:sub:lt}

Let $T_n$ be the number of triangles in the Erdős–Rényi random graph
$G(n,p)$, and let $A$ be its adjacency matrix. Since $T_n$ is concentrated
around its mean, the events that $T_n$ is much larger or much smaller than
$\E[T_n]$ have exponentially small probability. These events, the upper and
lower tails of $T_n$, are basic examples of large deviation events in random
graphs. Large deviation theory studies the probabilities of rare events and
the behavior of a random system conditioned on them. It goes back to
Cramér~\cite{Cra38}, and Varadhan~\cite{Var66} developed it into a general
theory; see~\cite{DZ10} for a textbook account.

The two tails behave quite differently. The upper tail, known as the
infamous upper tail problem in random graph theory~\cite{JR02}, is now well
understood in a wide range of $p$, as we discuss in
Section~\ref{s1:sub:related}. There the conditioned graph looks like $G(n,p)$
with a small planted structure, such as a clique or a hub, on a vanishing
fraction of the vertices. The lower tail is less well understood, and it is the subject
of this paper. Fix $\eta\in(0,1)$ and consider the lower tail event
$\mathcal L=\{T_n\le\eta\E[T_n]\}$. We ask how unlikely this event is and what
a typical graph in it looks like. The first question asks for the leading
order asymptotics of the lower-tail large deviation rate
\[
\mathrm{LT}_{n,p}(\eta)=-\log\P[T_n\le\eta\E[T_n]].
\]
The second asks for a description of $A$ that holds with probability
$1-o(1)$ conditionally on $\mathcal L$, up to an error
that is small compared with $p$ in the cut norm \eqref{eq:cutnorm}.

The standard prediction for the lower tail comes from the naive mean field
approximation, which has been developed for both tails; see,
e.g.,~\cite{CD16,Eld18,Aug20,CD20,HMS22,KS23}. It
replaces the conditioned graph by a random graph with independent edges,
whose edge probabilities may differ from pair to pair, and it measures the
cost of such a random graph by its relative entropy with respect to
$G(n,p)$. The predicted rate is the smallest cost of a random graph with
independent edges that has at most $\eta\E[T_n]$ triangles in expectation,
and the minimizing random graph is the predicted description of the
conditioned graph. The edge probabilities of a minimizer are at most $p$, since raising an
edge probability above $p$ increases the cost and adds triangles.
After division by $p$ they form a \emph{graphon}, a symmetric measurable
function from $[0,1]^2$ to $[0,1]$, and as $p\to0$ the prediction becomes a
variational problem over graphons. See Section~\ref{s1:sub:previous} for
more discussion.

For the upper tail, the analogous prediction is correct to leading order,
and the resulting variational problem has been solved~\cite{LZ17,BGLZ17}.
For the lower tail and $p\gg n^{-1/2}$, Kozma and Samotij~\cite{KS23} proved
that the prediction is also correct to leading order, extending earlier
results in smaller ranges of $p$~\cite{CD20,CDP}. Together with the work of
Chin~\cite{Chi25,Chi26} on the structure of the conditioned graph, this
reduces both of our questions to the variational problem over graphons, as
we explain in Section~\ref{s1:sub:reduction}. Unlike for the upper tail, this
problem has remained open since the work of Zhao~\cite{Zha17}, as also noted
by Cook and Dembo~\cite{CD24} and by Chin~\cite{Chi25}.

The difficulty is that the optimizer need not be homogeneous. Deleting
edges uniformly, which leads to $G(n,\eta^{1/3}p)$, is optimal when $\eta$
is close to $1$. Zhao showed that it is not optimal for small $\eta$, where
random graphs close to a stochastic block model with two blocks have
smaller relative entropy. He conjectured that a single threshold separates
the two behaviors, and for triangles he predicted its value. Very recently,
Chin~\cite{Chi26} proved that such a threshold exists for every nonbipartite
graph in place of the triangle, and that the transition is of first order.
These results, though, do not determine the value of the threshold or
the minimizers below it. In this paper we solve the variational problem completely, identifying the
threshold and the minimizers on both sides of it. This answers both
questions for every $\eta\in(0,1)$ when $n^{-1/2}\ll p\ll1$. It settles questions of Janson and Warnke~\cite{JW16} and of
Zhao~\cite{Zha17} on triangles that have been open for about a decade. The condition
$p\gg n^{-1/2}$ cannot be dropped, since for smaller $p$ a different
mechanism governs the lower tail, as we discuss in
Section~\ref{s1:sub:related}.

\subsection{Main results}\label{s1:sub:results}

We now state our main results. Let $h(u)=u\log u-u+1$ for $u>0$ and
$h(0)=1$. For $r=\eta^{1/3}$, the Erdős–Rényi graph $G(n,rp)$ has
$\eta\E[T_n]$ triangles in expectation, and its relative entropy with respect
to $G(n,p)$ is $(1+o(1))\frac{n^2p}2h(r)$. We compare this graph with
stochastic block models with two blocks. For $a,b\in[0,1]$, the \emph{two
block SBM} is the random graph on $[n]$ with two blocks $S$ and
$[n]\setminus S$, where $|S|=\lfloor n/2\rfloor$, in which each pair is an
edge independently with probability $pa$ inside the blocks and $pb$ across
them. It has $(1+o(1))\frac{a^3+3ab^2}4\E[T_n]$ triangles in expectation,
and its relative entropy with respect to $G(n,p)$ is
$(1+o(1))\frac{n^2p}2\cdot\frac{h(a)+h(b)}2$. Hence, among two block SBMs
with at most $\eta\E[T_n]$ triangles in expectation, the smallest relative
entropy is $(1+o(1))\frac{n^2p}2\Vbal(r)$, where
\begin{equation}\label{eq:Vbal}
\Vbal(r)=\min\Big\{\frac{h(a)+h(b)}2:a,b\in[0,1]\text{ and
}\frac{a^3+3ab^2}4\le r^3\Big\}\quad(0<r<1).
\end{equation}
Let
\[
r_*=\sup\{r\in(0,1):\Vbal(r)<h(r)\},\qquad\eta_*=r_*^3.
\]
Thus $\eta_*$ is the transition point below which the optimal two block SBM
has smaller relative entropy than $G(n,rp)$. Our first theorem shows that the
lower-tail large deviation rate is the smaller of these two relative
entropies.

\begin{theorem}\label{thm:rate}
Let $n^{-1/2}\ll p\ll1$, $\eta\in(0,1)$ and $r=\eta^{1/3}$. Then
\[
\mathrm{LT}_{n,p}(\eta)=(1+o(1))\frac{n^2p}2\times
\begin{cases}
h(r) & \text{if }\eta\ge\eta_*,\\
\Vbal(r) & \text{if }\eta<\eta_*,
\end{cases}
\]
where $\eta_*=0.009217\ldots$.
\end{theorem}

Our second theorem describes the conditioned graph. We use the
normalized cut norm of a real $n\times n$ matrix $M$, introduced by Frieze
and Kannan~\cite{FK99},
\begin{equation}\label{eq:cutnorm}
\|M\|_\cut=\max_{x,y\in[0,1]^n}\frac1{n^2}\big|x^{\top}My\big|.
\end{equation}
The cut norm underlies
the theory of graph limits~\cite{LS06,BCLSV08,Lov12} and the large
deviation principle of Chatterjee and Varadhan~\cite{CV11}, and versions of it rescaled by the edge
density are used for sparse graphs~\cite{BCCZ19,BCCZ18,CDP}. We state Theorem~\ref{thm:structure} through bounds of the form
$\|A-pQ\|_\cut\le\varepsilon p$ for a symmetric matrix $Q$. Such a bound means that, the graph has
about as many edges as the random graph with edge probabilities $pQ_{ij}$ between any two vertex sets.
More precisely, taking $x$ and $y$ in \eqref{eq:cutnorm} to be the
indicator vectors of sets $X,Y\subset[n]$ shows that the number of pairs
$(i,j)\in X\times Y$ with $A_{ij}=1$ differs from
$p\sum_{i\in X,j\in Y}Q_{ij}$ by at most $\varepsilon pn^2$. This error is small compared with the roughly
$pn^2$ such pairs in $G(n,p)$.

Let $\1\in\R^n$ be the all ones vector, so that the edge
probabilities of $G(n,rp)$ form the matrix $rp\1\1^{\top}$. For $S\subset[n]$
and $a,b\in[0,1]$, let $M^S_{a,b}$ be the $n\times n$ matrix with entries $a$
at pairs $(i,j)$ on the same side of the partition $\{S,[n]\setminus S\}$
and $b$ at the other pairs, so that the edge probabilities of the two block
SBM with blocks $S$ and $[n]\setminus S$ form the matrix $pM^S_{a,b}$. For $r\in(0,r_*]$,
exactly one pair $(a,b)$ with $a<b$ attains the minimum in \eqref{eq:Vbal} (Theorem~\ref{thm:main}(b) and
Corollary~\ref{cor:Vbal}(b)), and we denote it by $(a_r,b_r)$.

\begin{theorem}\label{thm:structure}
Let $n^{-1/2}\ll p\ll1$, $\eta\in(0,1)$, $r=\eta^{1/3}$ and
$\varepsilon>0$. Conditionally on
$T_n\le\eta\E[T_n]$, the following holds with probability $1-o(1)$.
\begin{enumerate}[label=(\alph*)]
\item If $\eta>\eta_*$, then $\|A-rp\1\1^{\top}\|_\cut\le\varepsilon p$.
\item If $\eta<\eta_*$, then there is $S\subset[n]$ with
$\big||S|-\frac n2\big|\le\varepsilon n$ and
$\|A-pM^S_{a_r,b_r}\|_\cut\le\varepsilon p$.
\item If $\eta=\eta_*$, then the conclusion of (a) or the conclusion of (b)
holds.
\end{enumerate}
\end{theorem}

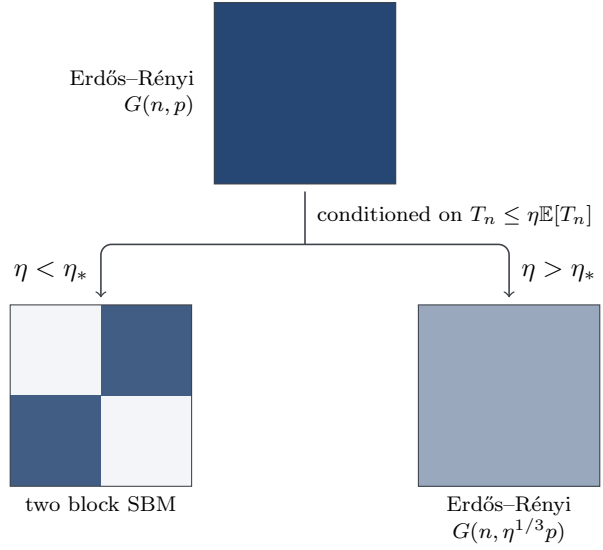
\begin{figure}[t]
\centering
\begin{tikzpicture}[x=0.2cm,y=0.2cm,
  frame/.style={figdark,line width=0.4pt},
  flow/.style={figdark,line width=0.6pt,rounded corners=4pt},
  lab/.style={font=\scriptsize,inner sep=1pt}]
%
\begin{scope}[shift={(-6,32)}]
\fill[figblue] (0,0) rectangle (12,-12);
\draw[frame] (0,0) rectangle (12,-12);
\end{scope}
\node[lab,anchor=east,align=right] at (-7.1,26) {Erdős–Rényi\\$G(n,p)$};
\draw[flow] (0,19.5) -- (0,16);
\draw[flow,->] (0,16) -| (-13.5,12.5);
\draw[flow,->] (0,16) -| (13.5,12.5);
\node[lab,anchor=west] at (0.6,17.9) {conditioned on $T_n\le\eta\E[T_n]$};
\node[anchor=east,font=\small,inner sep=1pt] at (-14.1,14.2) {$\eta<\eta_*$};
\node[anchor=west,font=\small,inner sep=1pt] at (14.1,14.2) {$\eta>\eta_*$};
\begin{scope}[shift={(-19.5,12)}]
\fill[figblue!5.347!white] (0,0) rectangle (6,-6);
\fill[figblue!5.347!white] (6,-6) rectangle (12,-12);
\fill[figblue!87.67!white] (6,0) rectangle (12,-6);
\fill[figblue!87.67!white] (0,-6) rectangle (6,-12);
\draw[frame] (0,0) rectangle (12,-12);
\end{scope}
\node[lab,anchor=north] at (-13.5,-0.45) {two block SBM};
\begin{scope}[shift={(7.5,12)}]
\fill[figblue!47.2!white] (0,0) rectangle (12,-12);
\draw[frame] (0,0) rectangle (12,-12);
\end{scope}
\node[lab,anchor=north,align=center] at (13.5,-0.45) {Erdős–Rényi\\$G(n,\eta^{1/3}p)$};
\end{tikzpicture}
\caption{The graphons of $G(n,p)$ (top) and of $G(n,p)$ conditioned on
$T_n\le\eta\E[T_n]$ (bottom), with edge probabilities divided by $p$. By
Theorem~\ref{thm:structure}, the conditioned graph is close to a two block
SBM for $\eta<\eta_*$ (left) and to $G(n,\eta^{1/3}p)$ for $\eta>\eta_*$
(right). Darker colors indicate larger values.}
\label{fig:schematic}
\end{figure}

Intuitively, there are two ways to have few triangles
(Figure~\ref{fig:schematic}). The first is to delete edges uniformly.
Keeping a fraction $r=\eta^{1/3}$ of the edges leads to $G(n,rp)$, which is
natural when $\eta$ is close to $1$ and the constraint is weak. The second
is to move towards a bipartite graph, which has no triangles at all. This
leads to the two block SBM, which is nearly bipartite since $a_r$ is small,
and it is natural when $\eta$ is close to $0$. Theorem~\ref{thm:structure}
shows that for every $\eta$ the conditioned graph takes one of these two
structures, and for $\eta\ne\eta_*$ it takes the more likely one.

The transition at $\eta_*$ is of first order, as the jump in
Figure~\ref{fig:phase} shows. Indeed, as $\eta\uparrow\eta_*$, the numbers
$a_r$ and $b_r$ tend to $0.01703\ldots$ and $0.8493\ldots$, so the two block
SBM does not approach $G(n,r_*p)$, and the structure of the conditioned
graph changes discontinuously at $\eta_*$. For instance, its edge density
divided by $p$ is close to $(a_r+b_r)/2$ for $\eta<\eta_*$, which tends to
$0.4331\ldots$, and close to $r$ for $\eta>\eta_*$, which tends to
$r_*=0.2096\ldots$. Thus, just below $\eta_*$, the conditioned graph has
about twice as many edges as just above it. Chin~\cite{Chi26}
proved that the edge density jumps at the threshold for every nonbipartite
graph in place of the triangle, and for triangles the values above give the
size of this jump. Moreover, $\mathrm{LT}_{n,p}(\eta)/(n^2p/2)$ converges to a function of
$\eta$ with a corner at $\eta_*$ (see the discussion after
Theorem~\ref{thm:main}). Section~\ref{s1:sub:solution} explains the
transition through a function of one variable with two competing local minima.

Theorem~\ref{thm:rate} answers an open problem of Janson
and Warnke~\cite{JW16} for triangles, which we describe in Section~\ref{s1:sub:related},
and its constant agrees with theirs as $\eta\to1$ and as $\eta\to0$. For the
structure, Chin~\cite{Chi25} showed that the conditioned
graph is close to $G(n,rp)$ when $r>0.4664\ldots$, and
Theorem~\ref{thm:structure}(a) extends this to every $\eta>\eta_*$. More
recently, Chin~\cite{Chi26} showed that, for every $\eta$, the
conditioned graph is close to some minimizer of the variational problem of
Section~\ref{s1:sub:previous}, without identifying the minimizers.
Theorem~\ref{thm:main} below identifies them, and this gives the explicit
description in Theorem~\ref{thm:structure}(b) and~(c) for $\eta\le\eta_*$.

We derive Theorems~\ref{thm:rate} and~\ref{thm:structure} from
Theorem~\ref{thm:main}, which solves the variational problem and is the
main contribution of this paper. The reduction from the conditioned graph
to this problem is due to Kozma and Samotij~\cite{KS23} and
Chin~\cite{Chi25,Chi26}, and we explain it in Section~\ref{s1:sub:reduction}.

\subsection{The variational problem}\label{s1:sub:previous}\label{s1:sub:reduction}

The variational problem comes from the naive mean field approximation. Let
$D(\cdot\|\cdot)$ denote Kullback–Leibler (KL) divergence, and let $\P$ also denote the law
of $G(n,p)$. The law $\P[\cdot\mid\mathcal L]$ of the conditioned graph, that is, of
$G(n,p)$ conditioned on $\mathcal L$, has density
$\1_{\mathcal L}/\P[\mathcal L]$ with respect to $\P$. Hence every law
$\nu$ of a random graph on $[n]$ with $\nu(\mathcal L)=1$ satisfies
\[
D(\nu\|\P)=D(\nu\|\P[\cdot\mid\mathcal L])-\log\P[\mathcal L]\ge-\log\P[\mathcal L],
\]
with equality if and only if $\nu=\P[\cdot\mid\mathcal L]$. Thus
$-\log\P[\mathcal L]$ is the smallest relative entropy with respect to $\P$
of a random graph that lies in $\mathcal L$, and the law of the conditioned
graph is the unique minimizer. The naive mean field approximation makes two changes to
this minimization. First, it restricts $\nu$ to the
laws of random graphs with independent edges. For a vector $q=(q_e)$
indexed by the pairs $e$ of vertices, let $G_q$ be the random graph on $[n]$
in which each pair $e$ is an edge with probability $q_e$, independently. Second, it replaces $\mathcal L$ by
the condition that $G_q$ has at most $\eta\E[T_n]$ triangles in
expectation. For such $\nu$ we have $D(\nu\|\P)=\sum_ei_p(q_e)$. Here
$i_p(x)=x\log\frac xp+(1-x)\log\frac{1-x}{1-p}$ is the relative entropy of
the Bernoulli law with mean $x$ with respect to that with mean $p$. This
gives the minimization
\[
\Phi_{n,p}(\eta)=\min\Big\{\sum_ei_p(q_e):G_q\text{ has at most
}\eta\E[T_n]\text{ triangles in expectation}\Big\}.
\]
The minimizing $G_q$ is the predicted description of the conditioned
graph. Cook and Dembo~\cite{CD24} explain this approximation
through the Gibbs variational principle.

Recall that a graphon is a symmetric measurable function
$U:[0,1]^2\to[0,1]$. If we write the edge probabilities as $q_e=pU(x,y)$
for a graphon $U$, then as $p\to0$ the minimization defining $\Phi_{n,p}$
becomes the following problem. For a graphon $U$ let
\[
F(U)=\iint h(U(x,y))\dd x\dd y,\quad
t(U)=\iiint U(x,y)U(y,z)U(x,z)\dd x\dd y\dd z ,
\]
and for a \emph{target} $r\in(0,1)$ define
\begin{equation}\label{eq:V}
V(r)=\min\{F(U):U\text{ a graphon with }t(U)\le r^3\}.
\end{equation}
A \emph{minimizer at $r$} is a graphon $U$ with $t(U)\le r^3$ and
$F(U)=V(r)$. Since $h(u)=\lim_{p\to0}p^{-1}i_p(pu)$, the
quantity $F(U)$ is the normalized cost of the random graph with edge probabilities
$pU$. For $r=\eta^{1/3}$, the constraint $t(U)\le r^3$ means that this
random graph has at most $\eta\E[T_n]$ triangles in expectation.

These predictions are correct when $n^{-1/2}\ll p\ll1$. After earlier
results of Cook and Dembo~\cite{CD20} and of Cook, Dembo and
Pham~\cite{CDP} in smaller ranges of $p$, Kozma and
Samotij~\cite{KS23} proved that $\Phi_{n,p}(\eta)$ agrees with
$-\log\P[\mathcal L]$ to leading order, up to an arbitrarily small change of
$\eta$, whenever $p\gg n^{-1/2}$, and they observed that this range of $p$
is optimal up to constants. Combined with the limit above, this gives
$\mathrm{LT}_{n,p}(\eta)=(1+o(1))\frac{n^2p}2V(r)$ for $r=\eta^{1/3}$ when
$p\to0$~\cite{Chi26}. For the structure, we compare graphons up to relabeling. For a measure
preserving map $\varphi:[0,1]\to[0,1]$ let
$U^\varphi(x,y)=U(\varphi(x),\varphi(y))$, and call a graphon that equals
$U^\varphi$ a.e.\ a \emph{relabeling} of $U$. The \emph{cut distance} of
graphons $U$ and $W$ is $\delta_\cut(U,W)=\inf_\varphi\|U^\varphi-W\|_\cut$,
where the infimum is over measure preserving bijections $\varphi$ of
$[0,1]$. Here
$\|D\|_\cut=\sup_{S,T}\big|\int_{S\times T}D(x,y)\dd x\dd y\big|$, with the
supremum over measurable $S,T\subset[0,1]$, is the analogue of
\eqref{eq:cutnorm} for functions on $[0,1]^2$. Chin~\cite{Chi26}
showed that, conditionally on $\mathcal L$, the adjacency matrix divided by $p$, viewed as a step function on
$[0,1]^2$, is close in cut distance to the set of minimizers at $r$ with
probability $1-o(1)$. Both results hold for every nonbipartite graph in
place of the triangle. Thus $V(r)$ determines the probability of the lower
tail, and the minimizers determine the structure of the conditioned graph.
We obtained these two statements for triangles independently
of~\cite{Chi26}, and we give our self-contained proof in Appendix~\ref{sec:reduction}
(Theorem~\ref{thm:reduction}).

The two ways to have few triangles correspond to two kinds of graphons.
Uniform deletion corresponds to the constant graphon $r$, which costs
$h(r)$. Moving towards a bipartite graph corresponds to the graphon
$B_{a,b}$ of the two block SBM, which equals $a$ on
$[0,\frac12)^2\cup[\frac12,1]^2$ and $b$ elsewhere. Restricting
\eqref{eq:V} to the graphons $B_{a,b}$ gives \eqref{eq:Vbal}.

\subsection{Solution of the variational problem}\label{s1:sub:solution}

Relabeling preserves $t$ and $F$, so minimizers of \eqref{eq:V} are
determined at best up to relabeling. In the proof we also use the graphons of
general SBMs with two blocks. A graphon $U$ is a \emph{two block SBM
graphon} if there are a measurable set $P\subset[0,1]$ and numbers
$a,d,b\in[0,1]$ with the following property. With $P^c=[0,1]\setminus P$,
the graphon $U$ equals $a$ on $P\times P$, $d$ on $P^c\times P^c$ and $b$
on $(P\times P^c)\cup(P^c\times P)$ a.e. Such a graphon describes an SBM
with two blocks of proportions $|P|$ and $1-|P|$. In our scaling, its edge
probabilities are $pa$ inside the first block, $pd$ inside the second block
and $pb$ across them. For example, $B_{a,b}$ is the two block SBM graphon with $P=[0,\frac12)$ and $d=a$. In the literature on graphons of maximal entropy with given edge and triangle densities, such graphons are called bipodal~\cite{KRRS17,NRS24}.

\begin{figure}[t]
\centering
\pgfplotsset{phase/.style={scale only axis,width=0.6\textwidth,height=0.34\textwidth,
  xmin=0,xmax=0.5,xtick={0,0.1,0.2,0.3,0.4,0.5},
  xlabel={$r$},
  tick label style={font=\scriptsize,/pgf/number format/fixed},
  label style={font=\small},
  axis on top,tick align=outside,tick pos=left}}
\tikzset{blockline/.style={figblue,line width=0.9pt,line join=round},
  constline/.style={figorange,line width=0.9pt,dashed},
  curvelabel/.style={font=\scriptsize,inner sep=1pt}}
\begin{tikzpicture}
\begin{axis}[phase,
  ymin=-0.06,ymax=1.08,ytick={0,0.25,0.5,0.75,1},ylabel={values of the minimizer}]
\draw[figdark,dotted,line width=0.6pt] (axis cs:0.20967,-0.06) -- (axis cs:0.20967,1.08);
\node[curvelabel,anchor=north west] at (axis cs:0.2125,1.06) {$r_*$};
\addplot[blockline] coordinates {
  (0,1)(0.00419,1)(0.00839,1)(0.01258,0.9999)(0.01677,0.9998)
  (0.02097,0.9997)(0.02516,0.9995)(0.02935,0.9993)(0.03355,0.999)
  (0.03774,0.9986)(0.04193,0.9982)(0.04613,0.9976)(0.05032,0.997)
  (0.05451,0.9963)(0.05871,0.9955)(0.0629,0.9946)(0.06709,0.9936)
  (0.07129,0.9925)(0.07548,0.9913)(0.07967,0.9899)(0.08387,0.9885)
  (0.08806,0.9869)(0.09226,0.9851)(0.09645,0.9833)(0.10064,0.9813)
  (0.10483,0.9792)(0.10903,0.9769)(0.11322,0.9744)(0.11742,0.9718)
  (0.12161,0.969)(0.1258,0.9661)(0.12999,0.9629)(0.13419,0.9596)
  (0.13838,0.9561)(0.14258,0.9523)(0.14677,0.9484)(0.15096,0.9442)
  (0.15516,0.9398)(0.15935,0.9351)(0.16354,0.9301)(0.16774,0.9248)
  (0.17193,0.9193)(0.17612,0.9133)(0.18032,0.9071)(0.18451,0.9004)
  (0.1887,0.8933)(0.1929,0.8857)(0.19709,0.8776)(0.20128,0.8689)
  (0.20548,0.8595)(0.20967,0.8493)};
\addplot[blockline] coordinates {
  (0,0)(0.00419,0)(0.00839,0)(0.01258,0)(0.01677,0)(0.02097,0)(0.02516,0)
  (0.02935,0)(0.03355,0.0001)(0.03774,0.0001)(0.04193,0.0001)
  (0.04613,0.0001)(0.05032,0.0002)(0.05451,0.0002)(0.05871,0.0003)
  (0.0629,0.0003)(0.06709,0.0004)(0.07129,0.0005)(0.07548,0.0006)
  (0.07967,0.0007)(0.08387,0.0008)(0.08806,0.0009)(0.09226,0.0011)
  (0.09645,0.0012)(0.10064,0.0014)(0.10483,0.0016)(0.10903,0.0018)
  (0.11322,0.002)(0.11742,0.0023)(0.12161,0.0026)(0.1258,0.0028)
  (0.12999,0.0032)(0.13419,0.0035)(0.13838,0.0039)(0.14258,0.0043)
  (0.14677,0.0047)(0.15096,0.0051)(0.15516,0.0056)(0.15935,0.0062)
  (0.16354,0.0067)(0.16774,0.0074)(0.17193,0.008)(0.17612,0.0087)
  (0.18032,0.0095)(0.18451,0.0103)(0.1887,0.0112)(0.1929,0.0122)
  (0.19709,0.0133)(0.20128,0.0144)(0.20548,0.0157)(0.20967,0.017)};
\addplot[constline,domain=0.20967:0.5,samples=2] {x};
\addplot[figred,only marks,mark=*,mark size=1.4pt] coordinates {(0.20967,0.01703) (0.20967,0.84935)};
\addplot[figred,only marks,mark=square*,mark size=1.3pt] coordinates {(0.20967,0.20967)};
\node[curvelabel,text=figblue,anchor=north] at (axis cs:0.08,0.97) {$b_r$};
\node[curvelabel,text=figblue,anchor=south] at (axis cs:0.08,0.03) {$a_r$};
\node[curvelabel,text=figorange,anchor=south east] at (axis cs:0.42,0.43) {$r$};
\node[curvelabel,anchor=west] at (axis cs:0.218,0.84935) {$b_*$};
\node[curvelabel,anchor=west] at (axis cs:0.218,0.01703) {$a_*$};
\end{axis}
\end{tikzpicture}
\caption{The values of the minimizer at target $r$. For $r<r_*$, the
minimizer is $B_{a_r,b_r}$ up to relabeling. It equals $a_r$ inside the
blocks and $b_r$ across them (solid). For $r>r_*$, it is the constant
graphon $r$ (dashed). At $r_*$, both the constant graphon $r_*$ and the
graphon $B_*=B_{a_*,b_*}$ of Theorem~\ref{thm:main}(c) are minimizers
(marks), and the minimizer jumps.}
\label{fig:phase}
\end{figure}
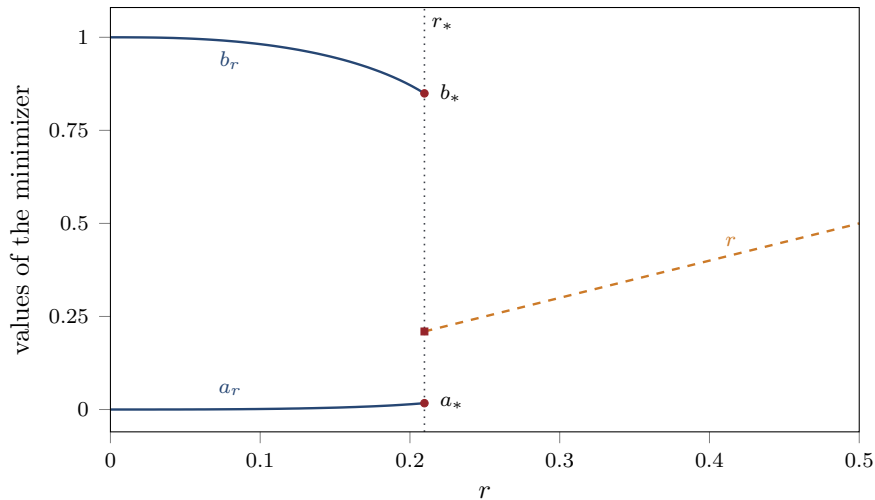

\begin{theorem}\label{thm:main}
Up to relabeling and null sets, the minimizers at $r\in(0,1)$ are as
follows.
\begin{enumerate}[label=(\alph*)]
\item If $r>r_*$, the unique minimizer is the constant graphon $r$.
\item If $r<r_*$, the unique minimizer is $B_r=B_{a_r,b_r}$. Here
$(a_r,b_r)$ is the only pair attaining the minimum in \eqref{eq:Vbal}, and
$0<a_r<r<b_r<1$.
\item If $r=r_*=0.2096\ldots$, the minimizers are exactly the constant
graphon $r_*$ and $B_*=B_{a_*,b_*}$. Here $(a_*,b_*)$ is the only pair
with $a_*<b_*$ attaining the minimum in \eqref{eq:Vbal} at $r_*$, and
$a_*=0.01703\ldots$ and $b_*=0.8493\ldots$.
\end{enumerate}
\end{theorem}

Zhao~\cite{Zha17} showed that the constant graphon is the unique minimizer
for $r>0.4664\ldots$ and is not a minimizer for $r<0.2096\ldots$, and he
conjectured that $0.2096\ldots$ is the threshold.
Proposition~\ref{prop:zhao} shows that $r_*$ equals this number, so
Theorem~\ref{thm:main} proves his conjecture.

We now explain why the transition is of first order. Fix $r$ close to $r_*$ and parametrize the graphons $B_{a,b}$ with
$t(B_{a,b})=r^3$ by $a$. Their cost
equals $h(r)$ at $a=r$, which gives the constant graphon, and has a second
local minimum at a small value of $a$. This minimum is above $h(r)$ for $r>r_*$ and below
$h(r)$ for $r<r_*$. So, as $r$ decreases
through $r_*$, the best graphon $B_{a,b}$ jumps from the constant graphon to
one with $a$ close to $a_*$ (Figure~\ref{fig:phase}), and $V$ has a corner at
$r_*$, with left derivative $-0.4961\ldots$ (Lemma~\ref{lem:Bderivative})
and right derivative $\log r_*=-1.562\ldots$. This picture goes back to
Zhao~\cite{Zha17}, and Theorem~\ref{thm:main} shows that it describes the
minimizers among all graphons.

\subsection{Related work}\label{s1:sub:related}
\noindent\textbf{The upper tail.} The upper tail of subgraph counts has been
studied extensively; see~\cite{Cha16,Cha17} for surveys. For triangles and
$\delta>0$, let $\mathrm{UT}_{n,p}(\delta)=-\log\P[T_n\ge(1+\delta)\E[T_n]]$. A long series of works~\cite{JR02,KV04,JOR04} determined the
order of $\mathrm{UT}_{n,p}(\delta)$ up to a factor $\log(1/p)$.
Chatterjee~\cite{Cha12} and DeMarco and Kahn~\cite{DK12a} closed this gap.
They showed that $\mathrm{UT}_{n,p}(\delta)$ is of order
$n^2p^2\log(1/p)$ for $p\ge n^{-1}\log n$, and DeMarco and
Kahn~\cite{DK12b} proved the analogous result for cliques. The leading
constant, which these results leave open, is the one predicted by the naive
mean field approximation. For constant $p$,
this follows from the large deviation principle of Chatterjee and
Varadhan~\cite{CV11}, which rests on the theory of graph
limits~\cite{LS06,BCLSV08,Lov12}; see also~\cite{CDey10}. In this dense
setting, Lubetzky and Zhao~\cite{LZ15} determined exactly when the constant
graphon is optimal. For $p\to0$,
Chatterjee and Dembo~\cite{CD16} developed a theory of nonlinear large
deviations, which establishes the approximation for $p\ge n^{-c}$ with a
small constant $c>0$. Subsequent works extended the range of
$p$~\cite{Eld18,CD20,Aug20,CDP}, up to essentially optimal ranges for
cliques~\cite{HMS22} and for regular graphs~\cite{BB23}. The resulting
variational problem was solved by Lubetzky and Zhao~\cite{LZ17} for cliques
and by Bhattacharya, Ganguly, Lubetzky and Zhao~\cite{BGLZ17} for general
graphs. The optimal strategy is to plant a clique, a hub, or both on a vanishing
proportion of the vertices. Harel, Mousset and Samotij~\cite{HMS22} and Cook and
Dembo~\cite{CD24} determined the structure of the conditioned graph. Thus the upper tail is now well understood.

\medskip
\noindent\textbf{The lower tail for fixed $p$.} The large deviation principle of
Chatterjee and Varadhan also reduces both of our questions to the dense
problem
\[
V_p(q)=\min\Big\{\iint i_p(W(x,y))\dd x\dd y:W\text{ a graphon with
}t(W)\le q^3\Big\},\qquad q=\eta^{1/3}p,
\]
with $-\log\P[\mathcal L]=(1+o(1))\frac{n^2}2V_p(q)$~\cite{Zha17}. Zhao also
found two curves in the $(p,q)$ plane such that the constant graphon is the unique
minimizer above the first and is not a minimizer below the second, and the
problem remains open between them. He also obtained \eqref{eq:V} as the
limit of this problem as $p\to0$, in the sense that $p^{-1}V_p(rp)\to V(r)$.
For the probability that $G(n,p)$ is claw free, Perkins and van der
Poel~\cite{PvdP25} found a first order transition at $p=(3-\sqrt5)/2$.
Deviations of smaller order than the mean, called moderate deviations, were
studied, mainly in the random graph $G(n,m)$, by Goldschmidt, Griffiths and
Scott~\cite{GGS20}, by Neeman, Radin and Sadun~\cite{NRSm23}, and by
Alvarado, Dias do Couto and Griffiths~\cite{ADG25}.

\medskip
\noindent\textbf{The lower tail for $n^{-1/2}\ll p\ll1$.} Janson and
Warnke~\cite{JW16} determined $\mathrm{LT}_{n,p}(\eta)$ to leading order in
this range of $p$ when $\eta=\eta(n)$ tends to $1$ or to $0$, and left the
determination of the constant for fixed $\eta\in(0,1)$ as an open problem.
The naive mean field prediction for $\mathrm{LT}_{n,p}(\eta)$ was proved
correct to leading order by Cook and Dembo~\cite{CD20} and by Cook, Dembo
and Pham~\cite{CDP} in smaller ranges of $p$, and by Kozma and
Samotij~\cite{KS23} in this full range. Chin~\cite{Chi25} determined
the structure of the conditioned graph for the values of $\eta$ at which the
constant graphon is known to be the unique minimizer. Section~\ref{s1:sub:lt} describes what was known about the
variational problem, including the recent work of Chin~\cite{Chi26}.

\medskip
\noindent\textbf{The critical window.} By results of Janson, Łuczak and
Ruciński~\cite{JLR90} and Janson~\cite{Jan90}, $\mathrm{LT}_{n,p}(\eta)$ is of
order $\min\{n^3p^3,n^2p\}$ (see~\cite{JW16}). The first term is the cost of
deleting edges from individual triangles, and the second is the cost of
changing the edge probabilities of a positive fraction of the pairs. The two
cross at $p\asymp n^{-1/2}$, and the second is smaller in the range
$p\gg n^{-1/2}$ of this paper. In the window $p=cn^{-1/2}$, Jenssen, Perkins,
Potukuchi and Simkin~\cite{JPPS24} determined the lower tail for every $c$
when $\eta\ge0.4993$ and for small $c$ otherwise, by an approach from
statistical physics, and they extended this approach to other
subgraphs~\cite{JPPS26a}. They also gave efficient algorithms that sample
$G(n,p)$ conditioned to be triangle free on both sides of the
window~\cite{JPPS24S}. Kwan and Pham~\cite{KP26} proved a variational formula
for the lower tail in the whole window, which answers a question of Jenssen,
Perkins, Potukuchi and Simkin. For the probability of being triangle free,
Jenssen, Perkins, Potukuchi and Simkin obtained this formula independently,
as a special case of their results on hypergraph independent
sets~\cite{JPPS26b}. Kwan and Pham also solved a penalized version of
\eqref{eq:V}, which is a key input to our proof for small targets
(Section~\ref{sec:smalltail}).

\medskip
\noindent\textbf{The lower tail for $p\ll n^{-1/2}$.} Here deleting edges from
individual triangles is cheaper, and the lower tail is of Poisson type. Indeed, for $np\to\infty$, Janson and
Warnke~\cite{JW16} showed that $\mathrm{LT}_{n,p}(\eta)=(1+o(1))(1-\eta+\eta\log\eta)\E[T_n]$.
For $\eta=0$, finer asymptotics of this probability were obtained by
Wormald~\cite{Wor96}, by Stark and Wormald~\cite{SW18}, and by Mousset,
Noever, Panagiotou and Samotij~\cite{MNPS20}. Deleting triangles locally is then asymptotically
optimal~\cite{JPPS24}, and the conditioned graph stays close
to $G(n,p)$ in cut norm, since
$\P[\mathcal L]\ge\P[T_n=0]\ge(1-p^3)^{\binom n3}=e^{-o(n^2p)}$ by the
Harris inequality (see, e.g., \cite{Har60}), whereas a cut norm deviation of
order $p$ has probability $e^{-\Omega(n^2p)}$ by Lemma~\ref{lem:cutconc}.
Lintott, Perkins and Yap~\cite{LPY26} refined this picture for
$p=o(n^{-2/3})$, finding the first order asymptotics of $\P[T_n=k]$ for every
$k\le(1-\varepsilon)\E[T_n]$ and the law of the triangles under this
conditioning.

\medskip
\noindent\textbf{Triangle-free graphs.} As $\eta\to0$, the optimal two block
SBM tends to the complete bipartite graphon $B_{0,1}$, in line with the
bipartite structure of typical triangle-free graphs. Erdős, Kleitman and
Rothschild~\cite{EKR76} showed that almost all triangle-free graphs are
bipartite, and Kolaitis, Prömel and Rothschild~\cite{KPR87} extended this to
$K_{\ell+1}$-free graphs. For a uniformly random triangle-free graph with $m$
edges, Prömel and Steger~\cite{PS96}, Łuczak~\cite{Luc00} and
Steger~\cite{Ste05} studied when it is bipartite or close to bipartite, and
Osthus, Prömel and Taraz~\cite{OPT03} showed that the threshold for being
bipartite is $m\approx\frac{\sqrt3}4n^{3/2}\sqrt{\log n}$. Balogh, Morris,
Samotij and Warnke~\cite{BMSW16} found the corresponding threshold for
$K_{r+1}$-free graphs. For $m$ of the same order, Jenssen, Perkins and
Potukuchi~\cite{JPP25} described how the edges inside the two parts of a
maximum cut evolve, and Mizgerd, Perkins and Wang~\cite{MPW26} studied the
chromatic number. For $G(n,p)$ itself, DeMarco and Kahn~\cite{DK15} showed that its largest triangle-free subgraph has as many edges as its
largest bipartite subgraph when
$p>C(\log n/n)^{1/2}$, and Conlon and Gowers~\cite{CG16} and
Schacht~\cite{Sch16} proved sparse random versions of Turán's theorem for
$p\gg n^{-1/2}$, the range of this paper. The hypergraph container method of
Balogh, Morris and Samotij~\cite{BMS15} and of Saxton and
Thomason~\cite{ST15} gives another approach to such results.

\medskip
\noindent\textbf{Entropy maximization with constraints.} Structurally
similar questions
arise for graphons of maximal entropy with prescribed edge and triangle
densities. There, optimizers with two blocks appear~\cite{KRRS17,NRS24,RS13}, and
numerical experiments of DiCarlo and Sadun~\cite{DS25} suggest optimizers with three blocks when the edge density is below
$(3-\sqrt3)/6$ and the triangle density is slightly below the cube of the
edge density.
The objective there is the entropy of the graphon under two equality
constraints. Here the function $h$ and the one sided triangle constraint
lead to different first order conditions and a different phase diagram.
The certificates of Sections~\ref{sec:threshold}
and~\ref{sec:intermediate} are close in spirit to the sum of squares
certificates of flag algebras \cite{Raz07}. 

\subsection{Overview of the proof}
\begin{figure}[t]
\centering
\begin{tikzpicture}[x=1cm,y=1cm,
  sbmaxis/.style={figblue,line width=0.9pt},
  constaxis/.style={figorange,line width=0.9pt,dashed},
  tick/.style={figdark,line width=0.6pt},
  move/.style={line width=0.6pt,->},
  brace/.style={figdark,line width=0.4pt,
    decorate,decoration={brace,mirror,amplitude=4pt,raise=20pt}},
  ticklab/.style={font=\small,anchor=base,inner sep=1pt},
  secname/.style={font=\small,anchor=base,inner sep=1pt},
  lab/.style={font=\scriptsize,text=figdark,anchor=base,inner sep=1pt}]
%
\draw[sbmaxis] (0,0) -- (10.4,0);
\draw[constaxis] (10.4,0) -- (14,0);
\foreach \x in {0,3,7.6,14} \draw[tick] (\x,-0.12) -- (\x,0.12);
\foreach \x/\lbl in {0/$0$,3/$r_-$,7.6/$r_*-\delta$,14/$1$}
  \node[ticklab] at (\x,-0.45) {\lbl};
\node[ticklab,text=figred] at (10.4,-0.45) {$r_*$};
\draw[figred,line width=0.6pt] (10.4,0) -- (10.4,1.1);
\fill[figred] (10.4,0) circle (2.2pt);
\node[lab] at (10.4,1.25) {certificates at $r_*$};
\node[secname,text=figred] at (10.4,1.6) {Section~\ref{sec:threshold}};
\draw[move,figorange,dashed] (10.4,0.45) -- (13.8,0.45);
\node[lab] at (12.1,0.62) {scaling};
\draw[move,figblue] (7.4,0.45) -- (3.2,0.45);
\node[lab] at (5.3,0.62) {global continuation};
\draw[brace] (0,0) -- (3,0);
\draw[brace] (3,0) -- (7.6,0);
\draw[brace] (7.6,0) -- (10.4,0);
\foreach \x/\key/\desc in {1.5/sec:smalltail/penalized problem,
    5.3/sec:intermediate/left endpoint,9/sec:rigidity/local continuation}{
  \node[secname,text=figblue] at (\x,-1.25) {Section~\ref{\key}};
  \node[lab] at (\x,-1.6) {\desc};}
\end{tikzpicture}
\caption{How the proof covers all targets $r\in(0,1)$. Below $r_*$ the
minimizers are the relabelings of $B_r$ (solid line), above $r_*$ the
minimizer is the constant graphon $r$ (dashed line), and at $r_*$ (dot) both
are minimizers. Section~\ref{sec:threshold} treats $r_*$ and, by scaling,
all $r>r_*$. Section~\ref{sec:rigidity} extends the result from $r_*$ to a
short interval below it, Section~\ref{sec:smalltail} treats the small
targets $r\le r_-$, and Section~\ref{sec:intermediate} closes the gap
between them. Sections~\ref{sec:threshold} and~\ref{sec:intermediate} use
computer assisted certificates. Lengths on the axis are not proportional to
the values of $r$, since $r_-\approx0.040$ and $r_*\approx0.210$.}
\label{fig:roadmap}
\end{figure}
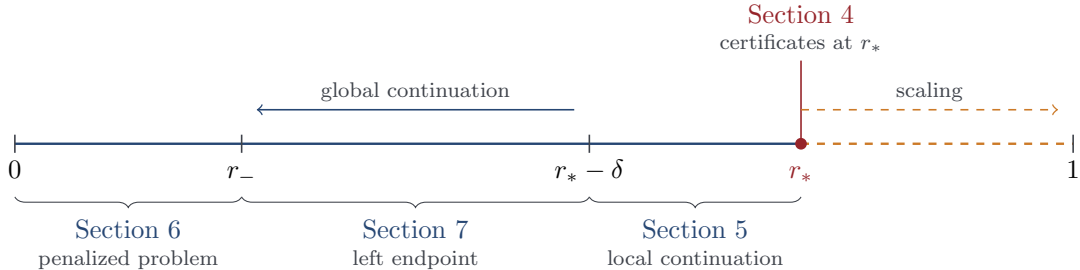

As explained in Section~\ref{s1:sub:reduction}, both questions about the
conditioned graph can be reduced to the variational problem \eqref{eq:V}, which
we solve in the rest of the paper. The problem is nonconvex and infinite
dimensional, and at the threshold two different graphons are optimal, so
first order conditions alone cannot identify the minimizers. Our plan is as
follows.
We first compare the two block graphons $B_{a,b}$ with the constant graphon. This comparison gives the threshold $r_*$, at which $B_*$ and the constant graphon $r_*$ have the same cost. We then show that these two graphons are the only minimizers at $r_*$. Next, a scaling argument shows that the constant graphon is the unique minimizer at every $r>r_*$. Finally, below $r_*$, we continue the minimizer $B_*$ downwards to the minimizers $B_r$, treat very small targets separately, and show that the two ranges meet. Figure~\ref{fig:roadmap} shows the range of targets covered by each section.

\begin{itemize}
\item \textbf{Two block graphons (Section~\ref{sec:balanced}).} We first ask
which conditions every minimizer $U$ must satisfy. Changing the values of $U$
on a small set of pairs gives the Euler--Lagrange equation
$\log U=-\mu_UC_U$, where $C_U(x,y)=\int_0^1U(x,w)U(y,w)\dd w$ is the
codegree function of $U$ and $\mu_U>0$ is a multiplier. For the graphons
$B_{a,b}$, this equation becomes two explicit equations in $a$, $b$ and $\mu$,
and we draw three consequences from them. First, if we also require that
$B_{a,b}$ has the same cost as the constant graphon with the same triangle
density, there is exactly one solution $(a_*,b_*,\mu_*)$ with $a_*\ne b_*$,
and this gives the threshold $r_*=((a_*^3+3a_*b_*^2)/4)^{1/3}$. Second,
$\Vbal(r)<h(r)$ exactly when $r<r_*$. Third, using a second first order
condition, which comes from reweighting the vertices, we show that a
minimizer at $r<r_*$ that is a two block SBM graphon is a relabeling of
$B_r$. The rest of the proof therefore has to show that minimizers are two
block SBM graphons.

\item \textbf{The threshold and above (Section~\ref{sec:threshold}).} At
the threshold $r_*$, the constant graphon $r_*$ and the two block graphon
$B_*=B_{a_*,b_*}$ from Section~\ref{sec:balanced} have the same cost. We
show that they are the only minimizers there. This is the heart of the
proof, because any argument at $r_*$ must be tight at two different
graphons. Testing the first order conditions of a minimizer against
finitely many functions gives linear constraints on the joint law of the
three edge values of a random triangle. These constraints depend only on
two numbers, the edge density and the cost of the minimizer, and elementary
bounds confine this pair to a compact region of the plane. We cover the
region by twelve polygons (Figure~\ref{fig:cover}). For each polygon we give an explicit
inequality, a certificate, which shows that no minimizer has its pair in
the polygon, except at the pair of $B_*$, where it forces the minimizer to
be a relabeling of $B_*$ (Proposition~\ref{prop:cover}). Targets above
$r_*$ then follow by a contraction argument. If a nonconstant graphon $U$
did at least as well as the constant graphon at some $r>r_*$, then the
scaled graphon $(r_*/r)U$ would do strictly better than the constant
graphon at $r_*$, which is impossible.

\item \textbf{Continuation below the threshold
(Section~\ref{sec:rigidity}).} Below $r_*$ we prove that the minimizer is
$B_r$ by continuation in $r$. Suppose that at some target the minimizers
are exactly the relabelings of $B_r$. We show that the same holds at all
nearby targets. By compactness, a minimizer at a nearby target is close to
$B_r$ in cut distance. Theorem~\ref{thm:rigidity} shows that a minimizer
close to a two block graphon satisfying the Euler--Lagrange equation is
itself a two block SBM graphon. Its proof works row by row. Every row of a
minimizer minimizes an explicit functional, and for $B_{a,b}$ this
functional has exactly two minimizers, the two rows of $B_{a,b}$, both
nondegenerate. Hence every row of a nearby minimizer lies close to one of
these two rows, and nondegeneracy forces rows in the same block to
coincide. By Section~\ref{sec:balanced}, the minimizer is then a
relabeling of $B_r$. Thus the set $G$ of targets $r\in(0,r_*)$ at which the
minimizers are exactly the relabelings of $B_r$ is open
(Corollary~\ref{cor:open}). Starting the same argument from $B_*$ at $r_*$
shows that $G$ contains an interval $(r_*-\delta,r_*)$
(Corollary~\ref{cor:seed}).

\item \textbf{Small targets (Section~\ref{sec:smalltail}).} For very small targets, $r\le r_-=0.03965\ldots$, we use a different
argument. Kwan and Pham~\cite{KP26} showed that the minimum of $F+\gamma t$
over all graphons equals its minimum over the graphons $B_{a,b}$, and this
result is the key input here. We take the penalty $\gamma=\mu_r/3$, where
$\mu_r$ is the multiplier of $B_r$. With this choice, the first order
condition for minimizing $F+\gamma t$ is $\log U=-3\gamma C_U=-\mu_rC_U$.
This is the Euler--Lagrange equation of the first step, with the multiplier
of $B_r$. The graphon $B_r$ attains this penalized minimum. Hence $B_r$ is a
minimizer at $r$, and every minimizer at $r$ also minimizes $F+\gamma t$. A
second variation argument, which uses that the triangle density is small,
shows that such a minimizer is a two block SBM graphon, hence a relabeling
of $B_r$. Thus $(0,r_-]\subset G$.

\item \textbf{Joining the two ends (Section~\ref{sec:intermediate}).} It
remains to show that $G=(0,r_*)$, that is, that the continuation from $r_*$
downwards does not stop before it reaches $r_-$. Suppose that it stops at
some $r_0\ge r_-$, so that $(r_0,r_*)\subset G$ but $r_0\notin G$. A
minimizer $U$ at $r_0$ then has the cost of $B_{r_0}$ by continuity, and
comparing scaled copies of $U$ with $B_s$ for $s$ slightly above $r_0$ shows
that $U$ has at least the edge density of $B_{r_0}$. An explicit family of
certificates, now depending on $r_0$, shows that these two facts force $U$
to be a relabeling of $B_{r_0}$ (Proposition~\ref{prop:family}). Hence
$r_0\in G$, a contradiction. Together with Section~\ref{sec:threshold}, this
completes the proof of Theorem~\ref{thm:main}.
\end{itemize}

\medskip
\noindent\textbf{Organization of the rest of the paper.} Section~\ref{sec:prelim}
collects notation and preliminary facts about graphons, and
Sections~\ref{sec:balanced} to~\ref{sec:intermediate} carry out the five
steps above. Section~\ref{sec:concluding} discusses open problems.
Appendix~\ref{sec:reduction} contains our self-contained proof of the
reduction for triangles and derives Theorems~\ref{thm:rate}
and~\ref{thm:structure} from Theorem~\ref{thm:main}, and
Appendix~\ref{app:numerics} contains the validated numerics behind the
computer-assisted steps.

\medskip
\noindent\textbf{Computer assistance.} We reduce the proofs of
Propositions~\ref{prop:cover} and~\ref{prop:family} to explicit inequalities
for functions of at most three variables and one parameter, and we give all the coefficients needed to verify them in
Appendix~\ref{app:numerics}.
The verification programs and their outputs are available at
\url{https://github.com/fifalsp/sparse-triangle-lower-tail-verification}.

\medskip
\noindent\textbf{Formal verification.} We have also formalized the arguments
of this paper in Lean~4 with mathlib. The formalization covers the
reduction of Appendix~\ref{sec:reduction}, the solution of \eqref{eq:V} in
Sections~\ref{sec:prelim} to~\ref{sec:intermediate}, a replay of the
certificates of Propositions~\ref{prop:cover} and~\ref{prop:family}, and the
deduction of Theorems~\ref{thm:rate} and~\ref{thm:structure}. It assumes the
results of Kozma and Samotij~\cite{KS23} and of Kwan and Pham~\cite{KP26}
that we use, together with classical facts such as the Harris inequality,
the compactness of the space of graphons in the cut distance and the
counting lemma. The formalization is available in the same repository.

\medskip
\noindent\textbf{Statement of AI use.}
The authors used GPT-6 Astra and Claude Opus 5.5 during the research and
preparation of this paper. From the work of Zhao~\cite{Zha17}, the authors expected the minimizers
to be constant graphons above the threshold $r_*$ and balanced two block
graphons below it. The difficulty was to prove this, and in particular to
find inequalities that hold for every graphon and are tight exactly at the
two block graphons. Discussions with GPT-6 Astra
helped the authors identify how the penalized value theorem of Kwan and
Pham~\cite{KP26}, applied with penalty coefficient
$\gamma=\mu_r/3$, could determine the optimal value for
$0<r\le r_-=0.03965\ldots$. These discussions also helped the authors find the scaling argument with which they reduce the constant
regime $r>r_*$ to the critical target $r_*$. With assistance from Claude Opus 5.5, the authors developed the
computer-assisted proof at $r_*$. Further
discussions with GPT-6 Astra helped develop the continuation
argument, together with a family of computer-assisted certificates,
that connects the small-target and critical regimes and covers
the remaining interval $(r_-,r_*)$. Claude Opus 5.5 was also used for language editing and figure preparation. Parts of the verification programs and of the Lean~4 formalization were
written by GPT-6 Astra and GPT-6.1 Sol. The authors take full
responsibility for the mathematical results, proofs, verification
programs, and presentation of this paper.

\section{Preliminaries}\label{sec:prelim}


\subsection{Notation}\label{s2:sub:notation}

For positive functions $f$ and $g$ of $n$, we write $f\ll g$ and $g\gg f$ if
$f/g\to0$ as $n\to\infty$, and $f\asymp g$ if $f/g$ is bounded above and
below by positive constants. For a measurable set $S\subset[0,1]^k$ let $|S|$ be its Lebesgue measure.
For an integrable function $W$ on $[0,1]^2$ we write
$\iint W=\int_0^1\int_0^1W(x,y)\dd x\dd y$. Similarly, $\iiint$ denotes
the integral over $[0,1]^3$. Let $\|W\|_p$ denote the norm of $W$ in
$L^p([0,1]^2)$. A \emph{kernel} is a bounded measurable function
$W:[0,1]^2\to\R$ with $W(x,y)=W(y,x)$. A \emph{graphon} is a kernel
with values in $[0,1]$. Relations between kernels hold almost everywhere
(a.e.) unless they are stated pointwise. For kernels define the edge
density, the degree, the codegree, the triangle degree and the triangle
density by
\begin{align*}
e(W)&=\iint W(x,y)\dd x\dd y,\quad d_W(x)=\int_0^1W(x,y)\dd y,\\
C_W(x,y)&=\int_0^1W(x,z)W(y,z)\dd z,\quad
\tau_W(x)=\int_0^1W(x,y)C_W(x,y)\dd y,\\
t(W_1,W_2,W_3)&=\iiint W_1(x,y)W_2(y,z)W_3(x,z)\dd x\dd y\dd z,
\end{align*}
and $t(W)=t(W,W,W)$. Note that $\int_0^1\tau_W=t(W)$.

Let $h(u)=u\log u-u+1$ for $u\ge0$, with $h(0)=1$. Then $h$ is
continuous on $[0,\infty)$, with $h'(u)=\log u$ and $h''(u)=1/u$ for $u>0$.
Hence $h$ is strictly convex, strictly decreasing on $[0,1]$ and strictly
increasing on $[1,\infty)$, with minimum $h(1)=0$. In particular $0\le h\le1$
on $[0,1]$. For a nonnegative kernel $W$ define
\[
H_W(x)=\int_0^1h(W(x,y))\dd y,\quad F(W)=\iint h(W(x,y))\dd x\dd y=\int_0^1H_W(x)\dd x.
\]
Let $\langle f,g\rangle=\int_0^1fg$ be the inner product on $L^2([0,1])$,
with norm $\|f\|_2$. The integral operator
$(T_Wf)(x)=\int_0^1W(x,y)f(y)\dd y$ of a kernel $W$ is self-adjoint on
$L^2([0,1])$. By the Cauchy--Schwarz inequality, its operator norm is at
most $\|W\|_2$. Note that $(T_WW(x,\cdot))(y)=C_W(x,y)$ and
$\langle W(x,\cdot),T_WW(x,\cdot)\rangle=\tau_W(x)$.

Let $D$ be a bounded measurable function on $[0,1]^2$. Define its \emph{cut norm} by
\[
\|D\|_\cut=\sup_{S,T}\Big|\int_{S\times T}D(x,y)\dd x\dd y\Big|,
\]
where the supremum is over measurable $S,T\subset[0,1]$. Thus
$\|D\|_\cut\le\|D\|_1$. A map $\varphi:[0,1]\to[0,1]$ is \emph{measure
preserving} if it is Borel measurable and $|\varphi^{-1}(B)|=|B|$ for every
Borel set $B\subset[0,1]$. A \emph{measure preserving bijection} is a
bijection $\varphi$ of $[0,1]$ such that $\varphi$ and $\varphi^{-1}$ are
measure preserving. For a kernel $W$ let
$W^\varphi(x,y)=W(\varphi(x),\varphi(y))$. A kernel $U$ is a
\emph{relabeling} of $W$ if $U=W^\varphi$ a.e.\ for some measure preserving
$\varphi$. The \emph{cut distance} of kernels $U$ and $W$ is
$\delta_\cut(U,W)=\inf_\varphi\|U^\varphi-W\|_\cut$, where the infimum is
over all measure preserving bijections $\varphi$ of $[0,1]$.

Every kernel agrees a.e.\ with a Borel measurable kernel. Changing a
kernel on a null set does not change $e$, $t$ or $F$, and it changes $d_W$,
$C_W$, $\tau_W$ and $H_W$ only on null sets. We always use Borel
representatives, so that these functions are defined everywhere and every
row $W(x,\cdot)$ is Borel measurable. If $\varphi$ is measure preserving,
then so is $(x_1,\dots,x_k)\mapsto(\varphi(x_1),\dots,\varphi(x_k))$ on
$[0,1]^k$ for $k\in\{2,3\}$, since Lebesgue measure and its image under this map agree on
products of Borel sets. Hence $W^\varphi$ is a kernel with the same $e$, $t$ and $F$ as $W$.
If $W$ is a graphon, then so is $W^\varphi$. For this reason, uniqueness of
minimizers is understood up to relabeling and null sets.

The formula defining $B_{a,b}$ in Section~\ref{sec:intro} gives a
nonnegative kernel for all $a,b\ge0$. Since three uniform points of $[0,1]$
lie in the same half with probability $\frac14$, these kernels satisfy
\begin{equation}\label{eq:Bdens}
e(B_{a,b})=\frac{a+b}2,\quad t(B_{a,b})=\frac{a^3+3ab^2}4,\quad
F(B_{a,b})=\frac{h(a)+h(b)}2.
\end{equation}
For a measurable $A\subset[0,1]$ let
$A^c=[0,1]\setminus A$. If $|A|=\frac12$, let $B^A_{a,b}$ equal $a$ on
$(A\times A)\cup(A^c\times A^c)$ and $b$ elsewhere. Then
$B_{a,b}=B^{[0,1/2)}_{a,b}$. For $a,b\in[0,1]$, Lemma~\ref{lem:relabelB}
shows that the graphons $B^A_{a,b}$ are exactly the relabelings of
$B_{a,b}$. Hence they also satisfy \eqref{eq:Bdens}.

Finally, $V(r)$ denotes the infimum in \eqref{eq:V}, which is attained by
Lemma~\ref{lem:existence}(a).

\subsection{Cut norm, scaling and existence of minimizers}\label{s2:sub:cut}

The trilinear form $t$ satisfies the following identities and bounds.

\begin{lemma}\label{lem:trilinear}
Let $W,Y,W_1,W_2,W_3$ be kernels.
\begin{enumerate}[label=(\alph*)]
\item The form $t(W_1,W_2,W_3)$ is invariant under permutations of its
arguments. Moreover, $t(Y,W,W)=\iint YC_W$. In particular
$t(W)=\iint WC_W$. Finally,
\begin{equation}\label{eq:texp}
t(W+Y)=t(W)+3\iint YC_W+3t(Y,Y,W)+t(Y).
\end{equation}
\item If $|W_2|\le1$ and $|W_3|\le1$, then $|t(W_1,W_2,W_3)|\le\|W_1\|_1$.
\item $|t(Y)|\le\|Y\|_2^3$.
\end{enumerate}
\end{lemma}

\begin{proof}
Renaming the integration variables and using the symmetry of the kernels
gives the permutation invariance. Integrating first in $z$ gives
$t(Y,W,W)=\iint YC_W$ by Fubini's theorem, and \eqref{eq:texp} follows by
trilinearity and permutation invariance. In (b), the integrand is at most $|W_1(x,y)|$ in absolute
value. For (c), the Cauchy--Schwarz inequality in $z$ gives
$C_Y(x,y)^2\le\int_0^1Y(x,z)^2\dd z\int_0^1Y(y,z)^2\dd z$, so
$\|C_Y\|_2\le\|Y\|_2^2$. Hence $|t(Y)|=|\iint YC_Y|\le\|Y\|_2\|C_Y\|_2\le
\|Y\|_2^3$.
\end{proof}

The cut norm controls the edge and triangle densities, graphons are compact
in the cut distance, and the cut distance is a pseudometric.

\begin{lemma}\label{lem:cut}
\begin{enumerate}[label=(\alph*)]
\item Let $D$ be a bounded measurable function on $[0,1]^2$. For all
measurable $f,g:[0,1]\to[0,1]$,
\[
\Big|\iint D(x,y)f(x)g(y)\dd x\dd y\Big|\le\|D\|_\cut.
\]
\item For all graphons $U$ and $W$, $|e(U)-e(W)|\le\|U-W\|_\cut$ and
$|t(U)-t(W)|\le3\|U-W\|_\cut$.
\item Every sequence of graphons has a subsequence $(U_n)$ for which there
are measure preserving bijections $\varphi_n$ of $[0,1]$ and a graphon $U$
with $\|U_n^{\varphi_n}-U\|_\cut\to0$.
\item For all kernels $U$, $Y$ and $W$ and every measure preserving
$\varphi$, we have $\delta_\cut(U,W)=\delta_\cut(W,U)\le\delta_\cut(U,Y)
+\delta_\cut(Y,W)$ and $\delta_\cut(W^\varphi,W)=0$.
\end{enumerate}
\end{lemma}

\begin{proof}
(a) Write $f(x)=\int_0^1\1\{f(x)>s\}\dd s$ and
$g(y)=\int_0^1\1\{g(y)>s'\}\dd s'$. By Fubini's theorem, the integral
equals $\int_0^1\int_0^1\big(\int_{S_s\times T_{s'}}D\big)\dd s\dd s'$ with
$S_s=\{f>s\}$ and $T_{s'}=\{g>s'\}$. Each inner integral is at most
$\|D\|_\cut$ in absolute value.

(b) The first bound is the case $S=T=[0,1]$ of the definition. The second
is the triangle case of the counting lemma (see, e.g.,
\cite[Lemma~10.23]{Lov12}), since the cut distance in that lemma is at most
$\|U-W\|_\cut$.

(c) and (d). We compare $\delta_\cut$ with the cut distance
$\delta'_\cut$ of \cite{BCLSV08}, which is defined for bounded symmetric
measurable functions through couplings of Lebesgue measure on $[0,1]$ with
itself. A Lebesgue measurable set differs from a Borel set by a subset of a
Borel null set. Hence our measure preserving maps and bijections are also
measure preserving in the sense of \cite{BCLSV08}. For $k\ge1$ and a
permutation $\pi$ of $[k]$, let $\tilde\pi$ translate each interval
$[\frac{i-1}k,\frac ik)$ onto $[\frac{\pi(i)-1}k,\frac{\pi(i)}k)$ and fix
$1$. This is a measure preserving bijection in our sense. It differs from
the interval permutation of \cite{BCLSV08} at finitely many points, which
does not affect cut norms. By \cite[Lemma~3.5]{BCLSV08}, $\delta'_\cut(W,U)$
equals the infimum of $\|W-U^\psi\|_\cut$ over the measure preserving
bijections $\psi$ in the sense of \cite{BCLSV08}, and also equals
$\lim_{k\to\infty}\min_\pi\|W-U^{\tilde\pi}\|_\cut$. The infimum defining
$\delta_\cut(U,W)$ lies between these two expressions, so
$\delta_\cut(U,W)=\delta'_\cut(W,U)$. Since $\delta'_\cut$ is symmetric
by its definition and satisfies the triangle inequality (see
\cite[Section~3.4]{BCLSV08}), so does $\delta_\cut$. By
\cite[Lemma~3.5]{BCLSV08} again, $\delta'_\cut(W,W^\varphi)$ is the infimum of
$\|W^{\varphi_1}-(W^\varphi)^{\psi_1}\|_\cut$ over measure preserving maps
$\varphi_1$ and $\psi_1$. Taking $\varphi_1=\varphi$ and
$\psi_1=\mathrm{id}$ gives $\delta_\cut(W^\varphi,W)=0$, which completes
(d). For (c), by \cite[Proposition~3.6]{BCLSV08} with $I=[0,1]$, the
space of graphons is compact under $\delta'_\cut$ when graphons at distance
zero are identified. See also \cite[Theorem~5.1]{LS07} and
\cite[Theorem~9.23]{Lov12}. Hence some subsequence $(U_n)$ and some graphon $U$
satisfy $\delta_\cut(U_n,U)\to0$. We may take $U$ Borel, since changing $U$
on a null set does not change cut norms. It remains to choose $\varphi_n$ with
$\|U_n^{\varphi_n}-U\|_\cut\le\delta_\cut(U_n,U)+1/n$.
\end{proof}

We deduce the lower semicontinuity of $F$ from block averaging, which we
also use in Appendix~\ref{sec:reduction} to discretize minimizers. For an integer $k\ge0$ let $K=2^k$ and
$J_a=[\frac{a-1}K,\frac aK)$ for $a\in[K]=\{1,\dots,K\}$, with $J_K$ closed.
For an integrable function $g$ on $[0,1]^2$ let $A_kg$ be the function equal
to $K^2\iint_{J_a\times J_b}g$ on $J_a\times J_b$, for all $a,b\in[K]$.

\begin{lemma}\label{lem:lsc}
\begin{enumerate}[label=(\alph*)]
\item For every graphon $W$ and $k\ge0$, $A_kW$ is a graphon with
$F(A_kW)\le F(W)$, and $\|A_kW-W\|_1\to0$ as $k\to\infty$.
\item If $U_n$ and $U$ are graphons with $\|U_n-U\|_\cut\to0$, then
$F(U)\le\liminf_nF(U_n)$.
\end{enumerate}
\end{lemma}

\begin{proof}
(a) The function $A_kW$ is symmetric with values in $[0,1]$, and Jensen's
inequality for $h$ on each square $J_a\times J_b$ gives
$F(A_kW)\le F(W)$. The map $A_k$ is linear with $\|A_kg\|_1\le\|g\|_1$,
and $\|A_kg-g\|_1\to0$ for continuous $g$ by uniform continuity. Since
continuous functions are dense in $L^1([0,1]^2)$, we get
$\|A_kW-W\|_1\to0$.

(b) Fix $k$. The functions $A_kU_n$ and $A_kU$ are constant on each of
the $K^2$ squares $J_a\times J_b$, where their values differ by
$K^2|\iint_{J_a\times J_b}(U_n-U)|\le K^2\|U_n-U\|_\cut\to0$. Since $h$
is continuous, part (a) gives $F(A_kU)=\lim_nF(A_kU_n)\le
\liminf_nF(U_n)$. By part (a), $A_kU\to U$ in $L^1$, so some subsequence
$A_{k_j}U$ converges to $U$ a.e. The function $h$ is bounded and
continuous on $[0,1]$, so bounded convergence gives
$F(U)=\lim_jF(A_{k_j}U)\le\liminf_nF(U_n)$.
\end{proof}

The cost $F$ satisfies an exact scaling identity. We use it to restore the
triangle density after a perturbation and to compare different targets.

\begin{lemma}\label{lem:scaling}
Let $a>0$, $s\ge0$, and let $W$ be a nonnegative kernel. Then
$t(aW)=a^3t(W)$,
\begin{align}
F(aW)&=aF(W)+(a\log a)e(W)+1-a,\label{eq:scale1}\\
F(aW)-h(as)&=a\big(F(W)-h(s)\big)+(a\log a)\big(e(W)-s\big).\label{eq:scale2}
\end{align}
\end{lemma}

\begin{proof}
The first identity follows from trilinearity. For $u\ge0$, both
$h(au)$ and $ah(u)+(a\log a)u+1-a$ equal $au\log u+au\log a-au+1$, with
$0\log0=0$. Integrating this identity with $u=W(x,y)$ gives
\eqref{eq:scale1}. Subtracting the identity with $u=s$ gives
\eqref{eq:scale2}.
\end{proof}

The kernels obtained by perturbing minimizers in
Section~\ref{s3:sub:euler} and later sections may take values larger than $1$. We next compare the cost of such kernels with $V(r)$.

\begin{lemma}\label{lem:truncation}
Let $W$ be a nonnegative kernel and $W'=\min(W,1)$. Then $W'$ is a graphon,
$t(W')\le t(W)$, and $F(W')\le F(W)$ with strict inequality if $|\{W>1\}|>0$. Hence,
for $r\in(0,1)$, every nonnegative kernel $W$ with $t(W)\le r^3$ satisfies
$F(W)\ge V(r)$, with equality only if $W\le1$ a.e. Thus $V(r)$ is also
the infimum of $F$ over such kernels.
\end{lemma}

\begin{proof}
Since $0\le W'\le W$ and the integrand of $t$ is nondecreasing in each
nonnegative factor, $t(W')\le t(W)$. Moreover $h(W')=h(W)$ on $\{W\le1\}$
and $h(W')=0<h(W)$ on $\{W>1\}$. This gives the claim on $F$. If
$t(W)\le r^3$, then $W'$ is feasible in \eqref{eq:V}, so $F(W)\ge F(W')\ge
V(r)$. By the strict inequality for $F$, equality $F(W)=V(r)$ forces
$|\{W>1\}|=0$.
\end{proof}

We next show that minimizers exist and that $V$ is continuous.

\begin{lemma}\label{lem:existence}
\begin{enumerate}[label=(\alph*)]
\item Let $r\in(0,1)$. Every graphon $U$ with $t(U)\le r^3$ satisfies
$F(U)>0$. The infimum in \eqref{eq:V} is attained. Moreover,
$0<V(r)\le h(r)<1$.
\item For $0<s<r<1$ we have $0\le V(s)-V(r)\le1-s/r$. In particular $V$ is
nonincreasing and continuous on $(0,1)$.
\end{enumerate}
\end{lemma}

\begin{proof}
(a) If $F(U)=0$, then $h(U)=0$ a.e., so $U=1$ a.e.\ and $t(U)=1>r^3$. Let
$(U_j)$ be graphons with $t(U_j)\le r^3$ and $F(U_j)\to V(r)$. By
Lemma~\ref{lem:cut}(c), after passing to a subsequence, there are measure
preserving bijections $\varphi_j$ and a graphon $U$ such that
$W_j=U_j^{\varphi_j}$ satisfies $\|W_j-U\|_\cut\to0$. Since relabeling
preserves $t$ and $F$, we also have $t(W_j)\le r^3$ and $F(W_j)\to V(r)$.
By Lemma~\ref{lem:cut}(b), $t(U)=\lim_jt(W_j)\le r^3$. By
Lemma~\ref{lem:lsc}(b), $F(U)\le\lim_jF(W_j)=V(r)$. Thus $U$ is
a minimizer at $r$, and $V(r)=F(U)>0$. Finally, the constant graphon $r$
is feasible, so $V(r)\le h(r)<h(0)=1$.

(b) Every graphon feasible at $s$ is feasible at $r$, so $V(s)\ge V(r)$.
Let $\lambda=s/r$ and let $U$ be a graphon with $t(U)\le r^3$. Then
$\lambda U$ is a graphon with $t(\lambda U)=\lambda^3t(U)\le s^3$. By
\eqref{eq:scale1},
\[
F(\lambda U)=\lambda F(U)+(\lambda\log\lambda)e(U)+1-\lambda\le F(U)+1-\lambda,
\]
where the inequality holds because $F(U)\ge0$, $e(U)\ge0$ and
$0<\lambda<1$. Taking the infimum over
$U$ gives $V(s)\le V(r)+1-s/r$. Continuity follows, since $1-s/r\to0$ as
$s\uparrow r$ with $r$ fixed and as $r\downarrow s$ with $s$ fixed.
\end{proof}

\subsection{Elementary bounds that isolate the constant graphon}\label{s2:sub:bounds}

We use the bounds of this subsection at the threshold $r_*$ in
Section~\ref{sec:threshold}, where they confine the point $(e(U),F(U))$ to a compact
region for every nonconstant graphon $U$ with $t(U)=r_*^3$ and
$F(U)\le h(r_*)$. The bounds rest on a quadratic lower bound for the
Bregman divergence of $h$, whose constant is the ratio of the divergence to
$(u-s)^2$ at $u=1$. For $s\in(0,1)$ define
\begin{equation}\label{eq:beta}
\beta(s)=\frac{h(1)-h(s)-(1-s)h'(s)}{(1-s)^2}=\frac{s-1-\log s}{(1-s)^2}.
\end{equation}

\begin{lemma}\label{lem:bregman}
Let $s\in(0,1)$. Then $\beta(s)\ge\frac12$, and every $u\in[0,1]$
satisfies
\begin{equation}\label{s2:bregman}
h(u)-h(s)-(u-s)\log s\ge\beta(s)(u-s)^2.
\end{equation}
\end{lemma}

\begin{proof}
For $u\in[0,1]$ let $R(u)=\int_0^1(1-\theta)((1-\theta)s+\theta
u)^{-1}\dd\theta$. This function is finite and nonincreasing in $u$. For
$u\in(0,1]$, Taylor's formula with integral remainder and $h''(v)=1/v$ show
that the left side of \eqref{s2:bregman} equals $(u-s)^2R(u)$. For $u=0$,
the left side and $s^2R(0)$ both equal $s$. For $u=1$, the identity for
$u\in(0,1]$ and \eqref{eq:beta} give $R(1)=\beta(s)$. Since $R$ is nonincreasing,
$R(u)\ge R(1)$, and this gives \eqref{s2:bregman}.
Finally $(1-\theta)s+\theta\le1$, so $\beta(s)=R(1)\ge\int_0^1(1-\theta)
\dd\theta=\frac12$.
\end{proof}

A first consequence is that a nonconstant graphon with cost at most $h(r)$
has edge density above $r$. Moreover, the only feasible constant graphon with such a
cost is $r$.

\begin{lemma}\label{lem:meangap}
Let $r\in(0,1)$, and let $U$ be a graphon with $F(U)\le h(r)$. If $U=c$
a.e.\ for a constant $c$ and $t(U)\le r^3$, then $c=r$. If $U$ is not a.e.\
constant, then $m=e(U)$ lies in $(0,1)$,
$F(U)\ge h(m)+\beta(m)\|U-m\|_2^2>h(m)$, and $m>r$.
\end{lemma}

\begin{proof}
If $U=c$ a.e.\ and $t(U)\le r^3$, then $c^3\le r^3$ and $h(c)\le h(r)$, so
$c=r$ because $h$ is strictly decreasing on $[0,1]$. Now suppose that $U$
is not a.e.\ constant. If $m=0$ or $m=1$, then $U=m$ a.e.\ because
$0\le U\le1$. Hence $m\in(0,1)$. Integrating \eqref{s2:bregman} with $s=m$
and $u=U(x,y)$ and using $\iint(U-m)=0$, we get
$F(U)-h(m)\ge\beta(m)\|U-m\|_2^2$. The right side is positive since $U$ is
not a.e.\ constant. Hence $h(m)<F(U)\le h(r)$. Since $h$ is strictly
decreasing on $[0,1]$, this gives $m>r$.
\end{proof}

Under the constraint $t(U)\le r^3$ the gap of Lemma~\ref{lem:meangap}
becomes quantitative. For $r\in(0,1)$ define
\begin{equation}\label{eq:Mr}
M(r)=r+\frac{9\beta(r)^3r^4}{(-\log r)^3},\quad M_1=M(r_*),
\end{equation}
where $r_*$ is the constant of Theorem~\ref{thm:main}. The function $M$ is
continuous on $(0,1)$ and satisfies $M(r)>r$.

\begin{lemma}\label{lem:isolation}
Let $r\in(0,1)$ and $L=-\log r$, and let $U$ be a graphon with
$t(U)\le r^3$, $F(U)\le h(r)$ and $|\{U\ne r\}|>0$. Then
\[
\|U-r\|_2\ge\frac{3\beta(r)r^2}L\quad\text{and}\quad e(U)\ge M(r).
\]
\end{lemma}

\begin{proof}
Let $\beta=\beta(r)$, $A=U-r$ and $\zeta=\|A\|_2>0$. Then
\[
3r^2\big(e(U)-r\big)\overset{\mathrm{(a)}}{=}t(U)-r^3-3r\|d_A\|_2^2-t(A)
\overset{\mathrm{(b)}}{\le}|t(A)|\overset{\mathrm{(c)}}{\le}\zeta^3,
\]
where $\mathrm{(a)}$ is \eqref{eq:texp} with $W=r$ and $Y=A$, together with
$C_r=r^2$ and $t(A,A,r)=r\|d_A\|_2^2$; $\mathrm{(b)}$ uses $t(U)\le r^3$
and $\|d_A\|_2^2\ge0$; and $\mathrm{(c)}$ is Lemma~\ref{lem:trilinear}(c). Integrating
\eqref{s2:bregman} with $s=r$ and $u=U(x,y)$ gives
$F(U)-h(r)\ge-L(e(U)-r)+\beta\zeta^2$. Since $F(U)\le h(r)$, this yields
$e(U)-r\ge\beta\zeta^2/L$. Combining the two bounds on $e(U)-r$ gives
$\beta\zeta^2/L\le\zeta^3/(3r^2)$, that is, $\zeta\ge3\beta r^2/L$. Then
$e(U)-r\ge\beta\zeta^2/L\ge9\beta^3r^4/L^3$.
\end{proof}

The next bound is the graphon form of Goodman's inequality~\cite{Goo59} on
the number of triangles.

\begin{lemma}\label{lem:goodman}
Every graphon $U$ satisfies $t(U)\ge e(U)(2e(U)-1)$, and hence
$e(U)\le\frac14\big(1+\sqrt{1+8t(U)}\big)$.
\end{lemma}

\begin{proof}
Integrating $(1-U(x,z))(1-U(y,z))\ge0$ in $z$ gives $C_U(x,y)\ge
d_U(x)+d_U(y)-1$. We multiply this bound by $U(x,y)$ and integrate. By
Lemma~\ref{lem:trilinear}(a) and the Cauchy--Schwarz inequality, we obtain
\[
t(U)=\iint UC_U\ge2\int_0^1d_U(x)^2\dd x-e(U)\ge2e(U)^2-e(U).
\]
Thus $e(U)$ is at most the larger root of $2y^2-y-t(U)$.
\end{proof}

The last bound controls the cost at given edge and triangle densities.

\begin{lemma}\label{lem:FB1}
Let $U$ be a graphon with $m=e(U)\in(0,1)$. Then
$\|U-m\|_2^2\ge(m^3-t(U))_+^{2/3}$ and
\[
F(U)\ge h(m)+\beta(m)\big(m^3-t(U)\big)_+^{2/3}.
\]
\end{lemma}

\begin{proof}
Let $A=U-m$, so that $\iint A=0$. By \eqref{eq:texp} with $W=m$ and $Y=A$
and by Lemma~\ref{lem:trilinear}(c), we have
$t(U)=m^3+3m\|d_A\|_2^2+t(A)\ge m^3-\|A\|_2^3$. Hence
$\|A\|_2^2\ge(m^3-t(U))_+^{2/3}$.
Integrating \eqref{s2:bregman} with $s=m$ and $u=U(x,y)$ gives
$F(U)-h(m)\ge\beta(m)\|A\|_2^2$.
\end{proof}

\subsection{Graphons with two values}\label{s2:sub:parity}

Several of our arguments end with a graphon that takes two values $a\ne b$
and in which almost every triangle has an even number of edges of value
$b$. In the next two lemmas we show that such a graphon is a two block SBM
graphon with the value $a$ inside both blocks and $b$ across them. If
moreover its edge density is $(a+b)/2$, then it is a relabeling of
$B_{a,b}$.

\begin{lemma}\label{lem:relabelB}
Let $a,b\in[0,1]$, and let $U$ be a graphon. Then $U$ is a relabeling of
$B_{a,b}$ if and only if $U=B^A_{a,b}$ a.e.\ for some measurable
$A\subset[0,1]$ with $|A|=\frac12$.
\end{lemma}

\begin{proof}
Let $\varphi$ be measure preserving and $A=\varphi^{-1}([0,\frac12))$.
Then $|A|=\frac12$. Moreover, $\varphi(x)$ and $\varphi(y)$ lie in the same
block $[0,\frac12)$ or $[\frac12,1]$ exactly when $\1_A(x)=\1_A(y)$, so
$B_{a,b}^\varphi=B^A_{a,b}$ everywhere.

Conversely, let $U=B^A_{a,b}$ a.e.\ with $|A|=\frac12$, where we may take
$A$ Borel. For
a Borel set $E\subset[0,1]$, the function $G_E(x)=|E\cap[0,x]|$ is
continuous and nondecreasing from $0$ to $|E|$. Hence, for $0\le c<|E|$,
we have $\{G_E\le c\}=[0,x_c]$ with $G_E(x_c)=c$, so
$|\{x\in E:G_E(x)\le c\}|=c$. Thus $G_E$ maps Lebesgue measure on $E$ to
Lebesgue measure on $[0,|E|]$. Let $\varphi=G_A$ on $A$ and
$\varphi=\frac12+G_{A^c}$ on $A^c$. This Borel map satisfies
$|\varphi^{-1}(B)|=|B\cap[0,\frac12]|+|(B-\frac12)\cap[0,\frac12]|=|B|$
for every Borel $B\subset[0,1]$, so it is measure preserving. In
particular, $P=\varphi^{-1}(\{\frac12\})$ is a null set. Since $\varphi$
maps $A\setminus P$ into $[0,\frac12)$ and $A^c\setminus P$ into
$(\frac12,1]$, we have $B_{a,b}^\varphi(x,y)=B^A_{a,b}(x,y)$ for
$x,y\notin P$. Hence $U=B_{a,b}^\varphi$ a.e.
\end{proof}

\begin{lemma}\label{lem:parity}
Let $a,b\in[0,1]$ with $a\ne b$, and let $U$ be a graphon with
$U\in\{a,b\}$ a.e. Let $S=\{(x,y):U(x,y)=b\}$, and suppose that
$\1_S(x,y)+\1_S(y,z)+\1_S(x,z)$ is even for a.e.\ $(x,y,z)\in[0,1]^3$.
Then there is a measurable $A\subset[0,1]$ such that $U=a$ a.e.\ on
$(A\times A)\cup(A^c\times A^c)$ and $U=b$ a.e.\ elsewhere. If moreover
$e(U)=(a+b)/2$, then $|A|=\frac12$. In this case $U=B^A_{a,b}$ a.e., and
$U$ is a relabeling of $B_{a,b}$.
\end{lemma}

\begin{proof}
Changing $U$ on a null set of $[0,1]^2$ changes the parity condition only
on a null set of $[0,1]^3$. Hence we may assume that $U$ is Borel and takes
only the values $a$ and $b$. Let $N\subset[0,1]^3$ be the null set on which
the parity condition fails. By Fubini's theorem there is $z_0\in[0,1]$
such that $(x,y,z_0)\notin N$ for a.e.\ $(x,y)$. Let
$A=\{x:(x,z_0)\in S\}$. For a.e.\ $(x,y)$,
\[
\1_S(x,y)\equiv\1_S(x,z_0)+\1_S(y,z_0)=\1_A(x)+\1_A(y)\pmod2.
\]
Hence $(x,y)\in S$ exactly when one of $x,y$ lies in $A$ and the other does
not. This is the first claim. Let $\alpha=|A|$. Then
\[
e(U)=a\big(\alpha^2+(1-\alpha)^2\big)+2b\alpha(1-\alpha)
=a+2(b-a)\alpha(1-\alpha).
\]
If $e(U)=(a+b)/2$, then $\alpha(1-\alpha)=\frac14$ because $a\ne b$.
Hence $\alpha=\frac12$. Now Lemma~\ref{lem:relabelB} gives the last claim.
\end{proof}
\section{The Euler--Lagrange equation and two block SBM graphons}\label{sec:balanced}

In this section, we explain where the threshold $r_*$ comes from. We obtain
$r_*$ by comparing the constant graphon with the graphons $B_{a,b}$. We also
show that a minimizer that is a two block SBM graphon must be a relabeling of
some $B_{a,b}$. Thus its two blocks have measure $\frac12$, and it takes the
same value inside both blocks.

Perturbing the edge values of a minimizer $U$ gives the Euler--Lagrange
equation $\log U=-\mu_UC_U$, with an explicit multiplier $\mu_U>0$
(Section~\ref{s3:sub:euler}). For the graphons $B_{a,b}$, this equation
becomes the two equations \eqref{s4:eq:EL} in $a$, $b$ and $\mu$. Equating
the cost of $B_{a,b}$ with that of the constant graphon with the same
triangle density gives a third equation. The resulting system has exactly
one solution $(a,b,\mu)$ with $a,b\in[0,1]$ and $a\ne b$, namely
$(a_*,b_*,\mu_*)$. The threshold $r_*$ is the cube root of the triangle
density of $B_{a_*,b_*}$. Along the solutions of \eqref{s4:eq:EL}, the
difference between the cost of $B_{a,b}$ and that of the constant graphon
changes sign exactly once (Theorem~\ref{thm:crossing}). This sign change
also gives $\Vbal(r)$ at every target $r$ (Corollary~\ref{cor:Vbal}).
Reweighting the vertices gives a second first order condition, the vertex
identity (Section~\ref{s3:sub:vertex}). Together, the two conditions force
a nonconstant two block SBM graphon that satisfies them to have blocks of
measure $\frac12$ and the same value inside both blocks
(Proposition~\ref{prop:bipodal}). Hence a minimizer at $r<r_*$ that is a two
block SBM graphon is a relabeling of $B_r$ (Corollary~\ref{cor:bipodalmin}).

\subsection{The Euler--Lagrange equation of a minimizer}\label{s3:sub:euler}


The first order conditions of minimizers are the main tool of this paper.
We first show that the triangle constraint is active.

\begin{lemma}\label{lem:active}
Let $r\in(0,1)$. Every minimizer $U$ at $r$ satisfies $t(U)=r^3$.
\end{lemma}

\begin{proof}
Suppose that $t(U)<r^3$. For $s\in(0,1)$ let $U_s=(1-s)U+s$, which is a
graphon. Since $t(U_s)$ is a polynomial in $s$ with value $t(U)<r^3$ at
$s=0$, we have $t(U_s)<r^3$ for all small $s>0$. Convexity of $h$ and
$h(1)=0$ give $h(U_s)\le(1-s)h(U)$. Since $F(U)>0$ by
Lemma~\ref{lem:existence}(a), it follows that $F(U_s)\le(1-s)F(U)<F(U)$.
This contradicts the minimality of $U$.
\end{proof}

For a minimizer $U$ at $r$ define its \emph{multiplier}
\begin{equation}\label{eq:mult}
\mu_U=\frac{1-F(U)-e(U)}{r^3}.
\end{equation}
Zhao~\cite[Lemma~5.1]{Zha17} obtained the Euler--Lagrange equation by
varying the edge values of a minimizer. In the next proposition we identify
its multiplier with $\mu_U$ and derive two sided bounds on $U$.

\begin{proposition}\label{prop:euler}
Let $U$ be a minimizer at $r\in(0,1)$, with multiplier $\mu_U$.
\begin{enumerate}[label=(\alph*)]
\item $U>0$ a.e.
\item $U$ satisfies the \emph{Euler--Lagrange equation}
\begin{equation}\label{eq:euler}
\log U=-\mu_UC_U\quad\text{a.e.}
\end{equation}
\item $\mu_U>0$ and $\mu_Ur^3=-\iint U\log U$. Moreover, $\mu_U$ is the only
real number $\mu$ for which $\log U=-\mu C_U$ a.e.
\item $e^{-\mu_U}\le U\le\exp(-\mu_Ue^{-2\mu_U})<1$ and
$C_U\ge e^{-2\mu_U}$ a.e.
\end{enumerate}
\end{proposition}

\begin{proof}
(a), (b) The problem \eqref{eq:V} is the problem $\mathrm{LT}(K_3,r)$ of
\cite[(2.2)]{Zha17}. Both minimize $\iint h(W)$, with the same function
$h$, over graphons $W$ with values in $[0,1]$ and $t(W)\le r^3$. Thus $U$
minimizes $\mathrm{LT}(K_3,r)$. Applying \cite[Lemma~5.2]{Zha17} to $K_3$,
which has three edges, gives $U\ge r^{3r^{-3}}>0$ a.e. This proves (a).
For the triangle $K_3$, the functional derivative of $t$ at $U$ in
\cite[(5.1)]{Zha17} is $3C_U$. Hence \cite[Lemma~5.1]{Zha17} gives
$\lambda\ge0$ with $\log U+3\lambda C_U=0$ a.e. We multiply this equation by $U$ and integrate. Recall that
$\iint UC_U=t(U)$ (Lemma~\ref{lem:trilinear}(a)) and $t(U)=r^3$
(Lemma~\ref{lem:active}). Since $u\log u=h(u)+u-1$, we get
\[
3\lambda r^3=-\iint U\log U=1-F(U)-e(U)=\mu_Ur^3,
\]
where the last equality is \eqref{eq:mult}. Hence $3\lambda=\mu_U$, which
proves (b).

(c) The proof of (b) gives $\mu_Ur^3=-\iint U\log U$. Recall that $U>0$
a.e.\ by (a). Moreover, $U=1$ a.e.\ is excluded by $t(U)=r^3<1$. Since
$-u\log u>0$ for $u\in(0,1)$, the function $-U\log U\ge0$ is positive on a
set of positive measure. Hence $\mu_U>0$. We prove uniqueness at the end
of (d).

(d) By \eqref{eq:euler}, $U=\exp(-\mu_UC_U)$ a.e. Since $0\le C_U\le1$ and
$\mu_U>0$, this gives $U\ge e^{-\mu_U}$ a.e. By Fubini's theorem, for
a.e.\ $x$ we have $U(x,z)\ge e^{-\mu_U}$ for a.e.\ $z$. For two such points
$x,y$ we get $C_U(x,y)\ge e^{-2\mu_U}$. Hence $C_U\ge e^{-2\mu_U}$ a.e.
Thus $U=\exp(-\mu_UC_U)\le\exp(-\mu_Ue^{-2\mu_U})<1$ a.e. Finally, let
$\mu$ be a real number with $\log U=-\mu C_U$ a.e. Then $(\mu_U-\mu)C_U=0$
a.e. Since $C_U>0$ a.e., this gives $\mu=\mu_U$.
\end{proof}

In Sections~\ref{s3:sub:vertex} and~\ref{sec:rigidity} we perturb a
minimizer $U$ by reweighting its vertices and by inserting new vertices. We
next derive the first order condition for such a perturbation, followed by
a rescaling that restores the triangle density. Let $c>0$, and let $I$ be
$[0,c]$ or $[-c,c]$. For nonnegative kernels $(W_s)_{s\in I}$ such that
$F(W_s)$ and $t(W_s)$ are differentiable at $s=0$, from the right if
$I=[0,c]$, let $F'_0$ and $t'_0$ be these derivatives.

\begin{lemma}\label{s3:rescale}
Let $U$ be a minimizer at $r\in(0,1)$, and let $(W_s)_{s\in I}$ be
nonnegative kernels with $W_0=U$. Suppose that $F(W_s)$, $e(W_s)$ and
$t(W_s)$ are differentiable at $s=0$, from the right if $I=[0,c]$. Then
$F'_0+\frac{\mu_U}3t'_0\ge0$, with equality if $I=[-c,c]$.
\end{lemma}

\begin{proof}
By differentiability and Lemma~\ref{lem:active}, $t(W_s)\to t(U)=r^3>0$
as $s\to0$. Hence $\lambda(s)=(r^3/t(W_s))^{1/3}$ is defined for $s\in I$
near $0$. It is differentiable at $0$, with $\lambda(0)=1$ and
$\lambda'(0)=-t'_0/(3r^3)$. Let $g(s)=F(\lambda(s)W_s)$. The nonnegative
kernel $\lambda(s)W_s$ has triangle density $r^3$, so
Lemma~\ref{lem:truncation} gives $g(s)\ge V(r)=g(0)$. By \eqref{eq:scale1},
$g(s)=\lambda(s)F(W_s)+\lambda(s)\log\lambda(s)e(W_s)+1-\lambda(s)$. Since
the derivative of $\lambda\log\lambda$ at $0$ is $\lambda'(0)$, we get
\[
g'(0)=F'_0+\lambda'(0)\big(F(U)+e(U)-1\big)=F'_0+\frac{\mu_U}3t'_0,
\]
where the last equality uses \eqref{eq:mult}. Since $g$ has a minimum at
$0$, we have $g'(0)\ge0$, and $g'(0)=0$ if $I=[-c,c]$.
\end{proof}

\subsection{Solving the Euler--Lagrange equation for \texorpdfstring{$B_{a,b}$}{B\_\{a,b\}}}\label{s4:sec:euler}

Recall from Section~\ref{sec:prelim} that $B_{a,b}$ is defined for all
$a,b\ge0$ and satisfies \eqref{eq:Bdens}. We call
$[0,\frac12)^2\cup[\frac12,1]^2$ the diagonal blocks and its complement the
off diagonal blocks. Let $a,b\ge0$. If $x$ and $y$ lie in the same one of
the blocks $[0,\frac12)$ and $[\frac12,1]$, then the function
$w\mapsto B_{a,b}(x,w)B_{a,b}(y,w)$ equals $a^2$ on that block and $b^2$ on
the other. Otherwise it equals $ab$. Hence $C_{B_{a,b}}$ equals
$(a^2+b^2)/2$ on the diagonal blocks and $ab$ on the off diagonal blocks.
We also use the identity
\begin{equation}\label{s4:eq:cube}
\Big(\frac{a+b}2\Big)^3-\frac{a^3+3ab^2}4=\frac{(b-a)^3}8,
\end{equation}
which holds for all real $a$ and $b$ because
$(a+b)^3-2(a^3+3ab^2)=(b-a)^3$. By \eqref{eq:Bdens}, its left side equals
$e(B_{a,b})^3-t(B_{a,b})$.

Let $\mu\in\R$. Both $B_{a,b}$ and $C_{B_{a,b}}$ are constant on the
diagonal blocks and on the off diagonal blocks. Hence the Euler--Lagrange
equation $\log B_{a,b}=-\mu C_{B_{a,b}}$ of Proposition~\ref{prop:euler}
holds a.e.\ if and only if it holds at every point of $[0,1]^2$. Both
conditions are equivalent to
\begin{equation}\label{s4:eq:EL}
\log a=-\mu\frac{a^2+b^2}2,\quad \log b=-\mu ab.
\end{equation}
If $a=0$ or $b=0$, the equation holds for no $\mu$, since
$\log B_{a,b}=-\infty$ on a set of positive measure.

The threshold $r_*$ comes from equating the cost of $B_{a,b}$ with that of
the constant graphon with the same triangle density. By \eqref{eq:Bdens},
this gives the third equation
\begin{equation}\label{s4:eq:third}
h\Big(\Big(\frac{a^3+3ab^2}4\Big)^{1/3}\Big)=\frac{h(a)+h(b)}2.
\end{equation}
In Section~\ref{s4:sec:crossing} we show that the system of \eqref{s4:eq:EL}
and \eqref{s4:eq:third} has exactly one solution $(a,b,\mu)$ with
$a,b\in[0,1]$ and $a\ne b$ (Lemma~\ref{lem:branch} and
Theorem~\ref{thm:crossing}(c)). This solution defines $a_*$, $b_*$ and
$\mu_*$. We also show there that the threshold $r_*$ defined before
Theorem~\ref{thm:rate} equals $((a_*^3+3a_*b_*^2)/4)^{1/3}$
(Corollary~\ref{cor:Vbal}).

We first solve \eqref{s4:eq:EL}, whose nonconstant solutions with
$a,b\le1$ are parametrized by $z=a/b$ (Lemma~\ref{lem:branch}). For
$z\in(0,1)$, let $L=\log(1/z)$ and define
\begin{equation}\label{eq:branch}
b(z)=\exp\Big(-\frac{2zL}{(1-z)^2}\Big),\quad a(z)=zb(z),\quad
\mu(z)=\frac{2L}{(1-z)^2b(z)^2},
\end{equation}
and
\[
B_z=B_{a(z),b(z)},\quad \kappa(z)=\Big(\frac{z^3+3z}4\Big)^{1/3},\quad
\rho(z)=\kappa(z)b(z).
\]
By \eqref{eq:Bdens}, $t(B_z)=b(z)^3(z^3+3z)/4=\rho(z)^3$. Thus $B_z$ has
the same triangle density as the constant graphon $\rho(z)$. Define
$D(z)=h(\rho(z))-F(B_z)$, the difference between the costs of these two
graphons. By \eqref{eq:Bdens}, $D(z)=0$ if and only if $(a(z),b(z))$
satisfies \eqref{s4:eq:third}.

\begin{lemma}\label{lem:branch}
\begin{enumerate}[label=(\alph*)]
\item Let $z\in(0,1)$. Then $0<a(z)<b(z)<1$ and $\mu(z)>0$. Moreover,
$(a(z),b(z),\mu(z))$ satisfies \eqref{s4:eq:EL}. Equivalently,
$\log B_z=-\mu(z)C_{B_z}$ at every point of $[0,1]^2$.
\item Let $a,b>0$ with $a\ne b$ and $(a^3+3ab^2)/4<1$, and let $\mu\in\R$ be
such that $(a,b,\mu)$ satisfies \eqref{s4:eq:EL}. Then $a<b$ and
$(a,b,\mu)=(a(z),b(z),\mu(z))$ with $z=a/b\in(0,1)$.
\end{enumerate}
\end{lemma}

\begin{proof}
(a) Write $a,b,\mu$ for $a(z),b(z),\mu(z)$. Since $2zL/(1-z)^2>0$, we have
$b\in(0,1)$. Hence $0<a<b<1$ and $\mu>0$. By \eqref{eq:branch},
$\mu b^2(1-z)^2/2=L$ and $\log b=-z\mu b^2=-\mu ab$. Since
$(1+z^2)/2=(1-z)^2/2+z$,
\[
-\mu\frac{a^2+b^2}2=-\mu b^2\frac{(1-z)^2}2-z\mu b^2=\log z+\log b=\log a.
\]

(b) If $\mu\le0$, the right sides of \eqref{s4:eq:EL} are nonnegative. Then
$a,b\ge1$ and $(a^3+3ab^2)/4\ge1$, which is excluded. Hence $\mu>0$.
Subtracting the second equation from the first gives
$\log(a/b)=-\mu(b-a)^2/2<0$. Thus $a<b$. With $z=a/b\in(0,1)$, this
identity reads $\mu b^2=2L/(1-z)^2$. The second equation now reads
$\log b=-z\mu b^2=-2zL/(1-z)^2$, so $b=b(z)$, $a=zb=a(z)$ and $\mu=\mu(z)$.
\end{proof}

Note that $t(B_{a,b})<1$ for $a,b\in[0,1]$ with $a\ne b$. Hence, by
Lemma~\ref{lem:branch}, the nonconstant graphons $B_{a,b}$ that satisfy an
Euler--Lagrange equation are exactly the graphons $B_z$ with $z\in(0,1)$.
We next show that $\rho$ is unimodal.

\begin{lemma}\label{lem:unimodal}
The function $\rho$ is continuously differentiable on $(0,1)$. Moreover,
$\rho(z)\to0$ as $z\to0$ and as $z\to1$. There is $z_{\mathrm{sp}}\in(0,1)$
such that $\rho'>0$ on $(0,z_{\mathrm{sp}})$ and $\rho'<0$ on
$(z_{\mathrm{sp}},1)$.
\end{lemma}

\begin{proof}
As $z\to0$, $z\log z\to0$, so $b(z)\to1$ and $\rho(z)\le\kappa(z)\to0$. As
$z\to1$, we have $\kappa(z)\to1$. Moreover, $\log z\le z-1$ gives
$\log b(z)\le-2z/(1-z)\to-\infty$, so $\rho(z)\to0$. Since
\[
(\log b)'(z)=\frac{2(1+z)\log z+2(1-z)}{(1-z)^3},\quad
(\log\kappa)'(z)=\frac{z^2+1}{z(z^2+3)},
\]
we obtain $(\log\rho)'(z)=2(1+z)(1-z)^{-3}\vartheta(z)$, where
\[
\vartheta(z)=\log z+\frac{1-z}{1+z}+\frac{(z^2+1)(1-z)^3}{2z(z^2+3)(1+z)}.
\]
Thus $\rho'$ has the sign of $\vartheta$. A direct computation gives
$\vartheta'(z)=p(z)/(2z^2(z^2+3)^2)$ with
$p(z)=-z^6+2z^5-5z^4+20z^3-15z^2+18z-3$. We have
$p'(z)=18-30z+40z^2-2z^2(3z^3-5z^2+10z-10)$. For $z\in[0,1]$, the
polynomial $3z^3-5z^2+10z-10=z^2(3z-5)-10(1-z)$ is negative. Moreover,
$18-30z+40z^2>0$ because its discriminant $30^2-4\cdot18\cdot40$ is negative.
Hence $p'>0$ on $[0,1]$. Since $p(0)=-3$ and $p(1)=16$, the polynomial $p$ has
exactly one zero $z_0\in(0,1)$. Thus $\vartheta$ is strictly decreasing on
$(0,z_0]$ and strictly increasing on $[z_0,1)$. Since $\vartheta(z)\to0$ as
$z\to1$, we get $\vartheta<0$ on $[z_0,1)$. Since $z\vartheta(z)\to\frac16$
as $z\to0$, we get $\vartheta(z)\to+\infty$ as $z\to0$. Therefore $\vartheta$ has
exactly one zero $z_{\mathrm{sp}}\in(0,z_0)$, with $\vartheta>0$ on
$(0,z_{\mathrm{sp}})$ and $\vartheta<0$ on $(z_{\mathrm{sp}},1)$.
\end{proof}

Let $r_{\mathrm{sp}}=\rho(z_{\mathrm{sp}})$. By Lemma~\ref{lem:unimodal}, for
each $r\in(0,r_{\mathrm{sp}})$ the equation $\rho(z)=r$ has exactly two
solutions $z_1(r)<z_{\mathrm{sp}}<z_2(r)$. Moreover, for each such $r$,
we have $\rho(z)>r$ exactly for $z\in(z_1(r),z_2(r))$. We also set $z_1(r_{\mathrm{sp}})=z_{\mathrm{sp}}$. Then
the map $z_1:(0,r_{\mathrm{sp}}]\to(0,z_{\mathrm{sp}}]$ is the inverse of
$\rho$ restricted to $(0,z_{\mathrm{sp}}]$. Hence $z_1$ is continuous and
strictly increasing. By the inverse function theorem, $z_1$ is continuously
differentiable on $(0,r_{\mathrm{sp}})$.

\subsection{The threshold \texorpdfstring{$r_*$}{r*} among the graphons \texorpdfstring{$B_{a,b}$}{B\_\{a,b\}}}\label{s4:sec:crossing}

For $z\in(0,1]$, let $Y_z=B_{z/\kappa(z),1/\kappa(z)}$. This kernel has
triangle density one. Every kernel $B_{a,b}$ with $0<a<b$ is
a positive multiple of $Y_z$ with $z=a/b$. For $z\in(0,1]$, let
$A(z)=e(Y_z)-1=(1+z)/(2\kappa(z))-1$. For $r\in(0,1)$, let
$g_r(z)=A(z)\log r+F(Y_z)$. Using the scaling identity \eqref{eq:scale1}, we reduce the comparison of
$rY_z$ with the constant graphon $r$ to the sign of $g_r(z)$.

\begin{lemma}\label{lem:ray}
The following hold.
\begin{enumerate}[label=(\alph*)]
\item For $z\in(0,1]$ and $r\in(0,1)$, we have $t(rY_z)=r^3$ and
$F(rY_z)=h(r)+rg_r(z)$. Conversely, if $0\le a<b$ and $t(B_{a,b})=r^3$ with
$r\in(0,1)$, then $a>0$ and $B_{a,b}=rY_{a/b}$.
\item We have $A(1)=F(Y_1)=0$. Moreover, $A>0$ and $A'<0$ on $(0,1)$. For
$r\in(0,1)$, the function $g_r$ is continuous on $(0,1]$ and continuously
differentiable on $(0,1)$, with $g_r(1)=0$.
\item For $r\in(0,1)$ and $z\in(0,1)$,
$g_r'(z)=A'(z)(\log r-\log\rho(z))$.
\item For $z\in(0,1)$, we have $\rho(z)\in(0,1)$, $B_z=\rho(z)Y_z$ and
$D(z)=-\rho(z)g_{\rho(z)}(z)$.
\end{enumerate}
\end{lemma}

\begin{proof}
(a) By \eqref{eq:Bdens}, $t(Y_z)=(z^3+3z)/(4\kappa(z)^3)=1$. By
\eqref{eq:scale1}, $F(rY_z)=rF(Y_z)+(r\log r)e(Y_z)+1-r$. Subtracting
$h(r)=r\log r-r+1$ gives $F(rY_z)-h(r)=r(F(Y_z)+A(z)\log r)=rg_r(z)$.
Conversely, $a>0$ because $t(B_{0,b})=0$. By \eqref{eq:Bdens}, $z=a/b$
satisfies $b^3\kappa(z)^3=t(B_{a,b})=r^3$. Hence $b=r/\kappa(z)$ and
$a=rz/\kappa(z)$.

(b) Since $\kappa(1)=1$, $Y_1$ is the constant $1$, so $A(1)=F(Y_1)=0$.
Applying \eqref{s4:eq:cube} to the two values of $Y_z$ and using
$t(Y_z)=1$, we obtain $e(Y_z)^3=1+(1-z)^3/(8\kappa(z)^3)=1+(1-z)^3/(2(z^3+3z))$. Hence $A>0$
on $(0,1)$. The derivative of $(1-z)^3/(z^3+3z)$ is
$-3(1-z)^2(1+z)^2/(z^2(z^2+3)^2)<0$, so $A'<0$ on $(0,1)$. The functions $A$
and $z\mapsto F(Y_z)=\frac12(h(z/\kappa(z))+h(1/\kappa(z)))$ are continuous
on $(0,1]$ and continuously differentiable on $(0,1)$. Hence so is $g_r$.

(c) Let $u=z/\kappa(z)$ and $v=1/\kappa(z)$. Then $u^3+3uv^2=4$ by (a), so
$(u^2+v^2)u'+2uvv'=0$. Moreover
$F(Y_z)=(h(u)+h(v))/2$ and $A=(u+v)/2-1$. Since $\rho=\kappa b$, we have
$a(z)=\rho(z)u$ and $b(z)=\rho(z)v$. Hence Lemma~\ref{lem:branch}(a) gives
$\log u=-\log\rho(z)-\mu(z)\rho(z)^2(u^2+v^2)/2$ and
$\log v=-\log\rho(z)-\mu(z)\rho(z)^2uv$. Therefore
\begin{align*}
\frac{\dd}{\dd z}F(Y_z)&=\frac{u'\log u+v'\log v}2\\
&=-\frac{u'+v'}2\log\rho(z)-\frac{\mu(z)\rho(z)^2}4\big((u^2+v^2)u'+2uvv'\big)\\
&=-A'(z)\log\rho(z).
\end{align*}
Since $g_r'(z)=A'(z)\log r+\frac{\dd}{\dd z}F(Y_z)$, this proves the claim.

(d) Since $z^3<z$, we have $\kappa(z)<1$. Together with $b(z)<1$, this gives
$\rho(z)\in(0,1)$.
As $a(z)<b(z)$ and $t(B_z)=\rho(z)^3$, part (a) gives $B_z=\rho(z)Y_z$ and
$F(B_z)=h(\rho(z))+\rho(z)g_{\rho(z)}(z)$.
\end{proof}

For $r\in(0,r_{\mathrm{sp}}]$ let $\delta(r)=g_r(z_1(r))$.

\begin{theorem}\label{thm:crossing}
The following hold.
\begin{enumerate}[label=(\alph*)]
\item The function $\delta$ is continuous and strictly increasing on
$(0,r_{\mathrm{sp}}]$. Moreover, $\delta(r)\to-\infty$ as $r\to0$,
$\delta(r_{\mathrm{sp}})>0$, and $\delta'(r)=A(z_1(r))/r$ for
$r\in(0,r_{\mathrm{sp}})$. Let $r_*\in(0,r_{\mathrm{sp}})$ be the unique zero
of $\delta$, and let $z_*=z_1(r_*)$.
\item For $r\in(0,r_*)$, the minimum of $g_r$ over $(0,1)$ equals
$\delta(r)<0$ and is attained only at $z_1(r)$. The function $g_{r_*}$ is
nonnegative on $(0,1)$ and vanishes only at $z_*$. For $r\in(r_*,1)$, we have
$g_r>0$ on $(0,1)$.
\item $D>0$ on $(0,z_*)$, $D(z_*)=0$, and $D<0$ on $(z_*,1)$.
\end{enumerate}
\end{theorem}

\begin{proof}
\noindent\textbf{Shape of $g_r$.} By Lemma~\ref{lem:ray}(b),(c), $g_r'(z)$
has the sign of $\rho(z)-r$ for $z\in(0,1)$. First let $r<r_{\mathrm{sp}}$.
Then $\rho(z)>r$ exactly for $z\in(z_1(r),z_2(r))$. Hence $g_r$ is strictly
decreasing on $(0,z_1(r)]$, strictly increasing on $[z_1(r),z_2(r)]$ and
strictly decreasing on $[z_2(r),1]$. In particular, $g_r>g_r(1)=0$ on
$[z_2(r),1)$. Now let $r\ge r_{\mathrm{sp}}$. Then
$\rho(z)\le r_{\mathrm{sp}}\le r$, with equality only if $z=z_{\mathrm{sp}}$
and $r=r_{\mathrm{sp}}$. Hence $g_r$ is strictly decreasing on $(0,1]$, so
$g_r>0$ on $(0,1)$.

(a) The map $(r,z)\mapsto g_r(z)$ is continuous on $(0,1)\times(0,1]$ and
continuously differentiable on $(0,1)\times(0,1)$. Since $z_1$ is continuous,
$\delta$ is continuous. Let $r\in(0,r_{\mathrm{sp}})$. Since
$\rho(z_1(r))=r$, Lemma~\ref{lem:ray}(c) gives $g_r'(z_1(r))=0$. Since also
$\partial_rg_r(z)=A(z)/r$, the chain rule gives
$\delta'(r)=A(z_1(r))/r>0$. Hence $\delta$ is strictly increasing on
$(0,r_{\mathrm{sp}}]$. By the shape of $g_r$ with $r=r_{\mathrm{sp}}$, we have
$\delta(r_{\mathrm{sp}})=g_{r_{\mathrm{sp}}}(z_{\mathrm{sp}})>0$. Fix
$w\in(0,z_{\mathrm{sp}})$ and let $r\in(0,\rho(w))$. Then
$z_1(r)<w<z_2(r)$. Since $g_r$ is increasing on $[z_1(r),z_2(r)]$, we get
$\delta(r)\le g_r(w)=A(w)\log r+F(Y_w)$. Since $A(w)>0$, the right side
tends to $-\infty$ as $r\to0$. Hence $\delta$ has a unique zero
$r_*\in(0,r_{\mathrm{sp}})$, with $\delta<0$ on $(0,r_*)$ and $\delta>0$ on
$(r_*,r_{\mathrm{sp}}]$.

(b) Let $r\in(0,r_{\mathrm{sp}})$. By the shape of $g_r$, we have
$g_r\ge\delta(r)$ on $(0,z_2(r)]$, with equality only at $z_1(r)$. Moreover,
$g_r>0$ on $[z_2(r),1)$. Together with the signs of $\delta$ from (a), this
proves all claims for $r<r_{\mathrm{sp}}$. For $r\ge r_{\mathrm{sp}}$, the
shape of $g_r$ gives $g_r>0$ on $(0,1)$.

(c) By Lemma~\ref{lem:ray}(d), $D(z)=-\rho(z)g_{\rho(z)}(z)$. If
$z\in(0,z_{\mathrm{sp}}]$, then $z=z_1(\rho(z))$, so
$D(z)=-\rho(z)\delta(\rho(z))$. By (a), this is positive, zero or negative
according as $\rho(z)<r_*$, $\rho(z)=r_*$ or $\rho(z)>r_*$. Since $\rho$ is
strictly increasing on $(0,z_{\mathrm{sp}}]$ with $\rho(z_*)=r_*$, these are
the cases $z<z_*$, $z=z_*$ and $z>z_*$. If $z\in(z_{\mathrm{sp}},1)$, then
$r=\rho(z)<r_{\mathrm{sp}}$ and $z=z_2(r)$. Hence $g_r(z)>0$ by the shape of
$g_r$, and $D(z)<0$.
\end{proof}

By Lemma~\ref{lem:branch}(b) and Theorem~\ref{thm:crossing}(c), $B_{z_*}$ is
the only nonconstant graphon $B_{a,b}$ that satisfies an Euler--Lagrange
equation and has the cost of the constant graphon with the same triangle
density. We define $a_*=a(z_*)$, $b_*=b(z_*)$, $\mu_*=\mu(z_*)$ and
$B_*=B_{z_*}=B_{a_*,b_*}$. Thus $(a_*,b_*,\mu_*)$ is the only solution of
the system of \eqref{s4:eq:EL} and \eqref{s4:eq:third} with $a,b\in[0,1]$
and $a\ne b$. For $r\in(0,r_*]$ let $z_r=z_1(r)$,
$a_r=a(z_r)$, $b_r=b(z_r)$, $\mu_r=\mu(z_r)$, $m_r=(a_r+b_r)/2$ and
$B_r=B_{z_r}=B_{a_r,b_r}$. In this notation, the subscripts $z_*$ and $z_r$
of $B$ are values of the branch parameter $z$, whereas the subscript $r$ of
$B_r$ is a target. Note that $B_{r_*}=B_*$. Since $t(B_r)=\rho(z_r)^3=r^3$ and
$a_r^3<(a_r^3+3a_rb_r^2)/4<b_r^3$, Lemma~\ref{lem:branch}(a) gives
\begin{equation}\label{s4:eq:arbr}
0<a_r<r<b_r<1\quad\text{for every }r\in(0,r_*].
\end{equation}

We now solve the problem \eqref{eq:Vbal}. For $r\in(0,1)$, we call
$(a,b)\in[0,1]^2$ a \emph{minimizing pair} at $r$ if $t(B_{a,b})\le r^3$
and $F(B_{a,b})=\Vbal(r)$.

\begin{corollary}\label{cor:Vbal}
Let $r\in(0,1)$.
\begin{enumerate}[label=(\alph*)]
\item If $r<r_*$, then $\Vbal(r)=F(B_r)=h(r)+r\delta(r)<h(r)$, and
$(a_r,b_r)$ is the only minimizing pair at $r$.
\item $\Vbal(r_*)=h(r_*)$, and the minimizing pairs at $r_*$ are exactly
$(r_*,r_*)$ and $(a_*,b_*)$.
\item If $r>r_*$, then $\Vbal(r)=h(r)$, and $(r,r)$ is the only minimizing
pair at $r$.
\end{enumerate}
\end{corollary}

\begin{proof}
Let $(a,b)\in[0,1]^2\setminus\{(r,r)\}$ with $t(B_{a,b})\le r^3$.

Suppose first that $a\ge b$, and let $m=(a+b)/2$. By \eqref{s4:eq:cube},
$m^3\le t(B_{a,b})\le r^3$. Since $h$ is strictly convex and strictly
decreasing on $[0,1]$, we get $F(B_{a,b})=(h(a)+h(b))/2\ge h(m)\ge h(r)$,
with equality only if $a=b=r$. Hence $F(B_{a,b})>h(r)$.

Now let $a<b$. The function $s\mapsto t(B_{(1-s)a+s,(1-s)b+s})$ is
continuous on $[0,1]$. It is at most $r^3$ at $s=0$ and equal to $1$ at
$s=1$. Let $s\in[0,1)$ be the smallest point where it equals $r^3$, and let
$a'=(1-s)a+s$ and $b'=(1-s)b+s$. Then $a'<b'$ and $t(B_{a',b'})=r^3$. Since
$h$ is strictly decreasing on $[0,1]$, we have $F(B_{a,b})\ge F(B_{a',b'})$,
with equality only if $s=0$. By Lemma~\ref{lem:ray}(a), $B_{a',b'}=rY_{z'}$
and $F(B_{a',b'})=h(r)+rg_r(z')$ with $z'=a'/b'$.

Let $r<r_*$. For $a<b$, Theorem~\ref{thm:crossing}(b) now gives
$F(B_{a,b})\ge h(r)+r\delta(r)$, with equality only if $B_{a,b}=rY_{z_r}$.
For $a\ge b$, we have shown that $F(B_{a,b})>h(r)$. Moreover,
$F(B_{r,r})=h(r)$. By Lemma~\ref{lem:ray}(d), we have $B_r=rY_{z_r}$ and
$t(B_r)=r^3$. Since $\delta(r)<0$ by Theorem~\ref{thm:crossing}(a),
Lemma~\ref{lem:ray}(d) also gives $F(B_r)=h(r)+r\delta(r)<h(r)$. This proves
(a). Now let
$r\ge r_*$. For $a<b$, Theorem~\ref{thm:crossing}(b) gives
$F(B_{a,b})\ge h(r)$, with equality only if $r=r_*$ and
$B_{a,b}=r_*Y_{z_*}=B_*$. Together with the case $a\ge b$ and
$F(B_{r,r})=h(r)$, this proves (b) and (c).
\end{proof}

In particular, the zero $r_*$ of $\delta$ equals
$\sup\{r\in(0,1):\Vbal(r)<h(r)\}$. This is the number $r_*$ defined before
Theorem~\ref{thm:rate}. Since $r_*=\rho(z_*)$ and $\rho(z_*)^3=t(B_*)$,
\eqref{eq:Bdens} gives $r_*=((a_*^3+3a_*b_*^2)/4)^{1/3}$. Moreover,
Corollary~\ref{cor:Vbal} proves the conjecture of Zhao described in
Section~\ref{s1:sub:solution}, restricted to the graphons $B_{a,b}$. Zhao~\cite[Section~2.1]{Zha17} denotes
the graphon $B_{a,b}$ by $\mathrm{BIP}_{a,b}$. The number $0.2096\ldots$
of \cite[Theorem~2.6]{Zha17} is defined there as the largest $r<1$ for which
$t(B_{a,b})=r^3$ and $F(B_{a,b})=h(r)$ have a solution with $0<a<r$ and
$b>0$. We next show that this number is $r_*$,
and we write $r_*$ as the maximum of an explicit function of one variable.
For $z\in(0,1)$ let
\[
R(z)=\exp\Big(-\frac{F(Y_z)}{A(z)}\Big)
=\exp\Big(-\frac{h(z/\kappa(z))+h(1/\kappa(z))}{(1+z)/\kappa(z)-2}\Big).
\]

\begin{proposition}\label{prop:zhao}
Let $r\in(0,1)$. There are $a\in(0,r)$ and $b>0$ with $t(B_{a,b})=r^3$ and
$F(B_{a,b})=h(r)$ if and only if $r=R(z)$ for some $z\in(0,1)$. Moreover
$R(z)\le r_*$ for every $z\in(0,1)$, with equality if and only if $z=z_*$.
Consequently $r_*=\max_{z\in(0,1)}R(z)$ is the largest such $r$.
\end{proposition}

\begin{proof}
Let $z\in(0,1)$. The two values of $Y_z$ differ, so $F(Y_z)>0$. Moreover,
$A(z)>0$ by Lemma~\ref{lem:ray}(b). Hence $R(z)\in(0,1)$ and
\[
g_r(z)=A(z)\log\frac r{R(z)}\quad\text{for every }r\in(0,1).
\]
By Theorem~\ref{thm:crossing}(b), $g_{r_*}(z)\ge0$, with equality if and
only if $z=z_*$. Thus $R(z)\le r_*$, with equality if and only if $z=z_*$.

Let $a\in(0,r)$ and $b>0$ with $t(B_{a,b})=r^3$. Then $a<b$, since $b\le a$
would give $t(B_{a,b})\le a^3<r^3$. By Lemma~\ref{lem:ray}(a),
$F(B_{a,b})=h(r)+rg_r(a/b)$. Hence $F(B_{a,b})=h(r)$ if and only if
$r=R(a/b)$. Conversely, if $r=R(z)$, then $g_r(z)=0$, and
Lemma~\ref{lem:ray}(a) shows that $(a,b)=(rz/\kappa(z),r/\kappa(z))$ is a
solution. Here $a<r$ because
$\kappa(z)^3-z^3=3z(1-z^2)/4>0$.
\end{proof}

In the last result of this subsection we identify the multiplier of $B_r$
and compute the derivative of $r\mapsto F(B_{z_1(r)})$. By
Corollary~\ref{cor:Vbal}, this function equals $\Vbal$ on $(0,r_*]$. We
use this result in Sections~\ref{sec:rigidity} and~\ref{sec:intermediate}.

\begin{lemma}\label{lem:Bderivative}
The following hold.
\begin{enumerate}[label=(\alph*)]
\item For $z\in(0,1)$, we have
$\mu(z)\rho(z)^3=-\frac12(a(z)\log a(z)+b(z)\log b(z))=1-F(B_z)-e(B_z)$. In
particular, $\mu_r=(1-F(B_r)-m_r)/r^3$ for $r\in(0,r_*]$. Hence, if $B_r$ is
a minimizer at $r$, then its multiplier \eqref{eq:mult} is $\mu_{B_r}=\mu_r$.
\item The function $r\mapsto F(B_{z_1(r)})$ is continuously differentiable
on $(0,r_{\mathrm{sp}})$, and its derivative is $-\mu(z_1(r))r^2$. In
particular, $\frac{\dd}{\dd r}F(B_r)=-\mu_rr^2$ for $r\in(0,r_*]$.
\end{enumerate}
\end{lemma}

\begin{proof}
(a) Write $a,b,\mu,\rho$ for $a(z),b(z),\mu(z),\rho(z)$. By \eqref{s4:eq:EL},
$a\log a+b\log b=-\mu(a(a^2+b^2)/2+ab^2)=-2\mu\rho^3$. Since
$u\log u=h(u)+u-1$, \eqref{eq:Bdens} gives
$(a\log a+b\log b)/2=F(B_z)+e(B_z)-1$. At $z=z_r$ we have $\rho=r$ and
$e(B_r)=m_r$.

(b) Let $r\in(0,r_{\mathrm{sp}})$, $z=z_1(r)$ and $B=B_z$. By
Lemma~\ref{lem:ray}(d), $B=rY_z$ and $F(B)=h(r)+r\delta(r)$. By
Theorem~\ref{thm:crossing}(a), this is continuously differentiable in $r$,
with derivative $\log r+\delta(r)+A(z)$. Since $e(B)=re(Y_z)=r(1+A(z))$ and
$h(r)=r\log r-r+1$, we have
$r(\log r+\delta(r)+A(z))=r\log r+F(B)-h(r)+rA(z)=F(B)+e(B)-1$. By (a), the
right side equals $-\mu(z)r^3$. Since $r_*<r_{\mathrm{sp}}$, this covers
$(0,r_*]$.
\end{proof}

\subsection{Numerical values and properties of \texorpdfstring{$B_*$}{B*}}\label{s4:sec:Bstar}

By Theorem~\ref{thm:crossing}(c), $D$ changes sign exactly once. Hence we can
determine the digits of $z_*$, $a_*$, $b_*$, $\mu_*$ and $r_*$ by evaluating
explicit functions of one variable.

\begin{lemma}\label{lem:digits}
Let $z_-=0.02005543105488$ and $z_+=0.02005543105489$. Then
$D(z_-)>0>D(z_+)$, and hence $z_*\in(z_-,z_+)$. Moreover
\begin{align*}
a_*&\in[0.01703401908,0.01703401909], & b_*&\in[0.8493469443,0.8493469444],\\
\mu_*&\in[11.28627816,11.28627817], & r_*&\in[0.2096696762,0.2096696763].
\end{align*}
\end{lemma}

\begin{proof}
Let $y_-=0.0200554310548823255291088$ and $y_+=0.0200554310548823255291090$,
so that $z_-<y_-<y_+<z_+$. Using ball arithmetic (Appendix~\ref{app:constants}), we find
$D(y_-)>0>D(y_+)$, so
$z_*\in(y_-,y_+)$ and $D(z_-)>0>D(z_+)$ by Theorem~\ref{thm:crossing}(c).
The same computation encloses the ranges of $a$, $b$, $\mu$ and $\rho$ over $[y_-,y_+]$
in the four stated intervals. These ranges contain $a_*=a(z_*)$,
$b_*=b(z_*)$, $\mu_*=\mu(z_*)$ and $r_*=\rho(z_*)$.
\end{proof}

Let $m_*=(a_*+b_*)/2=e(B_*)$. Define $F_*=h(r_*)$ and
$t_*=r_*^3$. For $u\ge0$ define
\begin{equation}\label{eq:k}
k(u)=2h(u)-u\log u=u\log u-2u+2.
\end{equation}
In Section~\ref{s3:sub:vertex} we show that almost every row of a minimizer
has the same integral of $k$. We next record the properties of $B_*$ used in
Section~\ref{sec:threshold}.

\begin{corollary}\label{cor:Bstar}
The graphon $B_*$ has the following properties.
\begin{enumerate}[label=(\alph*)]
\item $0<a_*<r_*<b_*<1$, $\mu_*>0$, $t(B_*)=t_*$, $F(B_*)=F_*$ and $m_*>r_*$.
\item $\log a_*=-\mu_*(a_*^2+b_*^2)/2$ and $\log b_*=-\mu_*a_*b_*$. Moreover,
$\log B_*(x,y)=-\mu_*C_{B_*}(x,y)$ for every $(x,y)\in[0,1]^2$.
\item $\mu_*t_*=1-F_*-m_*$. Moreover,
$\int_0^1k(B_*(x,y))\dd y=2F_*+\mu_*t_*=1+F_*-m_*$ for every $x\in[0,1]$.
\end{enumerate}
\end{corollary}

\begin{proof}
In (a), the inequalities $0<a_*<r_*<b_*<1$ are \eqref{s4:eq:arbr} with
$r=r_*$, and $\mu_*>0$ by Lemma~\ref{lem:branch}(a). Since
$t(B_z)=\rho(z)^3$ and $\rho(z_*)=r_*$, we have $t(B_*)=t_*$. Moreover,
$F(B_*)=h(r_*)-D(z_*)=F_*$ by Theorem~\ref{thm:crossing}(c). Finally,
\eqref{s4:eq:cube} gives $m_*^3-t_*=(b_*-a_*)^3/8>0$, so $m_*>r_*$. Part (b) is
Lemma~\ref{lem:branch}(a) with $z=z_*$. In (c), the first identity and
$(a_*\log a_*+b_*\log b_*)/2=-\mu_*t_*$ are Lemma~\ref{lem:Bderivative}(a)
with $z=z_*$. Each row of $B_*$ equals $a_*$ on a set of measure $\frac12$
and $b_*$ on its complement. Since $k(u)=2h(u)-u\log u$, the integral of
$k$ over each row is
$h(a_*)+h(b_*)-(a_*\log a_*+b_*\log b_*)/2=2F_*+\mu_*t_*$, which equals
$1+F_*-m_*$ by the first identity.
\end{proof}

\subsection{The vertex identity of a minimizer}\label{s3:sub:vertex}

Consider a minimizer that is a two block SBM graphon with blocks of
arbitrary measures (Section~\ref{s4:sec:bipodal}). The Euler--Lagrange
equation gives three equations relating its three values, the measure of one
block and its multiplier. Reweighting the vertices gives one more condition,
the vertex identity, which for such a graphon is the fourth equation of
\eqref{s4:eq:bipsys}. The vertex identity holds for every minimizer, and
we use it again in Sections~\ref{sec:threshold} to~\ref{sec:smalltail}.

Let $U$ be a minimizer at $r\in(0,1)$. Let $\psi:[0,1]\to[-1,1]$ be
measurable with $\int_0^1\psi=0$, and define
$(\psi\otimes\psi)(x,y)=\psi(x)\psi(y)$. For $|\varepsilon|\le\frac12$ let
$w_\varepsilon=1+\varepsilon\psi\ge\frac12$, and define
$G_\varepsilon(x)=\int_0^xw_\varepsilon$. Then $G_\varepsilon$ is an
increasing continuous bijection of $[0,1]$. Let $Q_\varepsilon$ be its
inverse, and let $U^\varepsilon(x,y)=U(Q_\varepsilon(x),Q_\varepsilon(y))$.
Since $|\{x:Q_\varepsilon(x)\le c\}|=G_\varepsilon(c)=\int_0^cw_\varepsilon$,
the map $Q_\varepsilon$ carries Lebesgue measure to
$w_\varepsilon(x)\dd x$. Hence, for $j\in\{2,3\}$, the product of $j$ copies
of $Q_\varepsilon$ carries Lebesgue measure on $[0,1]^j$ to the measure
with density $w_\varepsilon(x_1)\cdots w_\varepsilon(x_j)$. It follows that
$U^\varepsilon$ is a graphon. Moreover, $F(U^\varepsilon)$,
$e(U^\varepsilon)$ and $t(U^\varepsilon)$ are the integrals defining $F(U)$,
$e(U)$ and $t(U)$ with these densities inserted. Expanding the products of
weights and using the symmetry of the integrands and $t(U)=r^3$
(Lemma~\ref{lem:active}), we obtain
\begin{align}
F(U^\varepsilon)&=F(U)+2\varepsilon\int_0^1\psi H_U
+\varepsilon^2\iint h(U)(\psi\otimes\psi),\notag\\
e(U^\varepsilon)&=e(U)+2\varepsilon\int_0^1\psi d_U
+\varepsilon^2\iint U(\psi\otimes\psi),\label{s3:expand}\\
t(U^\varepsilon)&=r^3+3\varepsilon\int_0^1\psi\tau_U
+3\varepsilon^2\iint UC_U(\psi\otimes\psi)+\varepsilon^3c_\psi,\notag
\end{align}
where $c_\psi=\iiint U(x,y)U(y,z)U(x,z)\psi(x)\psi(y)\psi(z)\dd x\dd y\dd z$.
For a minimizer $U$ at $r$ define $\bar k_U=2F(U)+\mu_Ur^3$. We now state
the vertex identity in terms of $\bar k_U$ and the function $k$ of
\eqref{eq:k}.

\begin{proposition}\label{prop:vertex}
Let $U$ be a minimizer at $r\in(0,1)$. Then $\bar k_U=1+F(U)-e(U)$. For
a.e.\ $x\in[0,1]$,
\begin{equation}\label{eq:vertex}
2H_U(x)+\mu_U\tau_U(x)=\bar k_U\quad\text{and}\quad
\int_0^1k(U(x,y))\dd y=\bar k_U.
\end{equation}
Moreover, for a.e.\ $x\in[0,1]$,
\begin{equation}\label{s3:degree}
\mu_U\tau_U(x)+2d_U(x)=\mu_Ur^3+2e(U).
\end{equation}
\end{proposition}

\begin{proof}
The first claim follows from \eqref{eq:mult}. Fix $\psi$ as in the
construction of $U^\varepsilon$, and let $g=2H_U+\mu_U\tau_U$. By
\eqref{s3:expand}, Lemma~\ref{s3:rescale} applies to
$W_\varepsilon=U^\varepsilon$ with $I=[-\frac12,\frac12]$,
$F'_0=2\int_0^1\psi H_U$ and $t'_0=3\int_0^1\psi\tau_U$. It gives
$\int_0^1\psi g=0$. The mean of $g$ is
$\int_0^1g=2F(U)+\mu_Ut(U)=\bar k_U$. Since $0\le H_U,\tau_U\le1$ and
$\mu_U>0$ (Proposition~\ref{prop:euler}(c)), both $g$ and $\bar k_U$ lie in
$[0,2+\mu_U]$. Hence $\psi=(g-\bar k_U)/(2+\mu_U)$ is admissible. For this
$\psi$ we get $0=\int_0^1\psi g=\int_0^1(g-\bar k_U)^2/(2+\mu_U)$. Thus
$g=\bar k_U$ a.e., which is the first identity in \eqref{eq:vertex}.

Multiplying \eqref{eq:euler} by $U$ and integrating in $y$ gives
\[
\mu_U\tau_U(x)=-\int_0^1U(x,y)\log U(x,y)\dd y\quad\text{for a.e.\ }x.
\]
Substituting this into the first identity gives the second, because
$2h(u)-u\log u=k(u)$. Since also $k(u)=h(u)-u+1$, the second identity
reads $H_U=\bar k_U-1+d_U$ a.e. Substituting this into the first identity
gives $\mu_U\tau_U+2d_U=2-\bar k_U=1-F(U)+e(U)$. By \eqref{eq:mult}, the
right side equals $\mu_Ur^3+2e(U)$.
\end{proof}

\subsection{Classification of two block SBM minimizers}\label{s4:sec:bipodal}

Recall from Section~\ref{sec:intro} that a graphon $U$ is a two block SBM
graphon if there are a measurable set $P\subset[0,1]$ and numbers
$a,d,b\in[0,1]$ such that
\begin{equation}\label{s4:eq:bipodal}
U=a\1_{P\times P}+d\1_{P^c\times P^c}+b\1_{(P\times P^c)\cup(P^c\times P)}
\quad\text{a.e.},
\end{equation}
where $P^c=[0,1]\setminus P$. Let $c=|P|$. Then $C_U$ equals
$ca^2+(1-c)b^2$ on $P\times P$, $cb^2+(1-c)d^2$ on $P^c\times P^c$, and
$b(ca+(1-c)d)$ on the rest.

Let $U$ be a minimizer of this form with $c\in(0,1)$, and let $\mu=\mu_U$.
By Proposition~\ref{prop:euler}, $\mu>0$, $U>0$ and $\log U=-\mu C_U$. The
sets $P\times P$, $P^c\times P^c$ and $(P\times P^c)\cup(P^c\times P)$ have
positive measures $c^2$, $(1-c)^2$ and $2c(1-c)$. Hence $a,d,b>0$, and the
Euler--Lagrange equation on these three sets gives the first three
equations of \eqref{s4:eq:bipsys}. Their right sides are negative, so
$a,d,b\in(0,1)$. By Proposition~\ref{prop:vertex}, the integral
$\int_0^1k(U(x,y))\dd y$ is the same for a.e.\ $x$. It equals
$ck(a)+(1-c)k(b)$ for $x\in P$ and $ck(b)+(1-c)k(d)$ for $x\in P^c$. Since
$P$ and $P^c$ have positive measure, this gives the last equation of the
system
\begin{align}
\log a&=-\mu\big(ca^2+(1-c)b^2\big), & \log d&=-\mu\big(cb^2+(1-c)d^2\big),\notag\\
\log b&=-\mu b\big(ca+(1-c)d\big), & ck(a)+(1-c)k(b)&=ck(b)+(1-c)k(d).
\label{s4:eq:bipsys}
\end{align}
A solution with $a\ne d$ would force an explicit function of two variables to
vanish. We next show that this function is positive. For $0<x<y<1$
let $\alpha=1-x$, $\beta=1-y$, $\eta=y-x$, $\xi=-\log x$, $\zeta=-\log y$ and
\begin{equation}\label{s4:eq:J}
J(x,y)=\beta y(1+x)\xi^2-\alpha x(1+y)\zeta^2-2\eta\xi\zeta
+2\alpha\beta(\beta\xi-\alpha\zeta).
\end{equation}

\begin{lemma}\label{lem:J}
For all $0<x<y<1$, $J(x,y)>0$.
\end{lemma}

\begin{proof}
Fix $0<x<y<1$. We write $J(x,y)$ as an integral of a function that is
positive on an open square.

\noindent\textbf{Step 1: an integral representation.} For $v,w\in\R^4$ define
the symmetric bilinear form
\begin{align*}
\mathcal J(v,w)&=\beta y(1+x)v_1w_1-\alpha x(1+y)v_2w_2
-\eta(v_1w_2+v_2w_1)\\
&\quad-\beta\big((yv_1+v_2)w_3+(yw_1+w_2)v_3\big)\\
&\quad+\alpha\big((v_1+xv_2)w_4+(w_1+xw_2)v_4\big).
\end{align*}
Let $\hat\phi=(\xi,\zeta,\alpha,\beta)$. Then
$\mathcal J(\hat\phi,\hat\phi)=J(x,y)$, since
$2\alpha\beta(\xi+x\zeta-y\xi-\zeta)=2\alpha\beta(\beta\xi-\alpha\zeta)$. For
$s\in[x,1]$ let
\[
\phi(s)=\Big(s-x,(s-y)_+,\frac{s^2-x^2}2,\frac{(s^2-y^2)_+}2\Big).
\]
For $g(u)=1/u=h''(u)$, the coordinates of $\hat\phi$ are
$\int_x^1g(u)\dd u$, $\int_y^1g(u)\dd u$, $\int_x^1ug(u)\dd u$ and
$\int_y^1ug(u)\dd u$. The coordinates of $\phi(s)$ are the same integrals for
$g=\1_{[x,s]}$. For $u\in[x,1]$, we have
$1/u=1+\int_u^1s^{-2}\dd s=\1_{[x,1]}(u)+\int_x^1\1_{[x,s]}(u)s^{-2}\dd s$.
Hence Fubini's theorem gives $\hat\phi=\phi(1)+\int_x^1\phi(s)s^{-2}\dd s$.
Let $\omega(s,t)=\mathcal J(\phi(s),\phi(t))$. This is a continuous symmetric
function on $[x,1]^2$. By bilinearity,
\begin{equation}\label{s4:eq:Jrep}
J(x,y)=\omega(1,1)+2\int_x^1\omega(s,1)s^{-2}\dd s
+\int_x^1\int_x^1\omega(s,t)s^{-2}t^{-2}\dd s\dd t.
\end{equation}

\noindent\textbf{Step 2: positivity of $\omega$.} We claim that $\omega\ge0$
on $[x,1]^2$ and $\omega>0$ on $(x,1)^2$. By symmetry, we may assume
$s\le t$. In each of the following three regions, each positive part in
$\phi$ equals either $0$ or its argument. The stated formula for
$\omega(s,t)$ in each region is then a polynomial identity.

If $x\le s\le t\le y$, then $\phi(s)=(s-x,0,(s^2-x^2)/2,0)$ and
$\omega(s,t)=\beta y(s-x)(t-x)(1-(s+t)/2)$. Since $s+t<2$, this is positive
for $x<s\le t\le y$.

If $y\le s\le t\le1$, let $p=1-s$ and $q=1-t$, so that $0\le q\le p\le\beta$.
Then
\[
\omega(s,t)=\frac\eta2\Big((1-xy)(p-q)^2+pq\big(2(x+y-2xy)-(x+y)(p+q)\big)\Big).
\]
Since $p+q\le2\beta$, the factor $2(x+y-2xy)-(x+y)(p+q)$ is at least
$2(x+y-2xy)-2(x+y)(1-y)=2y\eta>0$. Hence $\omega(s,t)\ge0$, with equality only
if $p=q$ and $pq=0$, that is, only if $s=t=1$.

If $x\le s\le y\le t\le1$, let $p=s-x\in[0,\eta]$ and $q=1-t\in[0,\beta]$.
Then $\omega(s,t)=pQ(p)/2$, where
\[
Q(p)=\beta y(\alpha-q)(\alpha-p+q)+(\beta-q)\big(\alpha(\beta-q)-\beta p\big).
\]
The function $Q$ is affine in $p$. Its slope is
$-\beta(y(\alpha-q)+\beta-q)\le-\beta y\eta<0$, because
$\alpha-q\ge\alpha-\beta=\eta$ and $\beta-q\ge0$. At $p=\eta$, that is, at
$s=y$, the point $(y,t)$ also lies in the second region. Since $y<1$, the
second region gives $\eta Q(\eta)=2\omega(y,t)>0$. Hence
$Q(p)\ge Q(\eta)>0$ for $p\in[0,\eta]$. Thus $\omega(s,t)>0$ whenever $s>x$.

The three regions cover $\{(s,t):x\le s\le t\le1\}$. This proves the claim.
Hence the first two terms in \eqref{s4:eq:Jrep} are nonnegative. The last
term is positive, because its integrand is continuous and positive on the
open square $(x,1)^2$. Thus $J(x,y)>0$.
\end{proof}

We now return to the system \eqref{s4:eq:bipsys}. In terms of
\eqref{s4:eq:bipodal}, the conclusion of the next proposition says that
$|P|=\frac12$ and $d=a$. By Lemma~\ref{lem:relabelB}, such a two block SBM
graphon is a relabeling of $B_{a,b}$.

\begin{proposition}\label{prop:bipodal}
Let $c\in(0,1)$, $\mu>0$ and $a,d,b\in(0,1)$ satisfy \eqref{s4:eq:bipsys}, and
suppose that $a,d,b$ are not all equal. Then $a=d<b$ and $c=\frac12$.
\end{proposition}

\begin{proof}
The system \eqref{s4:eq:bipsys} and the conclusion are invariant under the
substitution $(c,a,d)\mapsto(1-c,d,a)$. Hence we may assume $a\le d$. Adding
the first two equations and subtracting twice the third gives
$\log(a/b)+\log(d/b)=-\mu(c(a-b)^2+(1-c)(d-b)^2)$. If $b\le a$, then the left
side is nonnegative and the right side is nonpositive. Hence both sides
vanish and $a=b=d$, which is excluded. Thus $a<b$. Now suppose that $b\le d$.
The second equation minus the third gives
$\log(d/b)=-\mu(cb(b-a)+(1-c)d(d-b))$. Here the left side is nonnegative and
the right side is negative, a contradiction. Hence $a\le d<b$.

Suppose that $a<d$. Let $x=a/b$ and $y=d/b$, so that $0<x<y<1$. Let
$\alpha,\beta,\eta,\xi,\zeta$ be as in \eqref{s4:eq:J}. Let $\sigma=-\log b>0$,
$\lambda=\mu b^2>0$, $S=1-xy>0$, $\theta_1=c\alpha$, $\theta_2=(1-c)\beta$,
$N_1=\zeta+y\xi$ and $N_2=\xi+x\zeta$. Since $\log a=-\sigma-\xi$ and
$\log d=-\sigma-\zeta$, the first three equations of \eqref{s4:eq:bipsys} read
\[
\sigma+\xi=\lambda(cx^2+1-c),\quad \sigma+\zeta=\lambda\big(c+(1-c)y^2\big),
\quad \sigma=\lambda\big(cx+(1-c)y\big).
\]
Subtracting the third equation from the first two gives
$\xi=\lambda(\theta_2-x\theta_1)$ and $\zeta=\lambda(\theta_1-y\theta_2)$.
This is a linear system for $(\lambda\theta_1,\lambda\theta_2)$ with
determinant $xy-1=-S\ne0$. Since $N_2-xN_1=S\xi$ and $N_1-yN_2=S\zeta$, its
solution is $\lambda S\theta_1=N_1$, $\lambda S\theta_2=N_2$. Thus
$\lambda Sc=N_1/\alpha$ and $\lambda S(1-c)=N_2/\beta$. Since
$x\beta+\alpha=\beta+y\alpha=S$, we have
$x\beta N_1+y\alpha N_2=S(y\xi+x\zeta)$. Inserting $\lambda Sc$ and
$\lambda S(1-c)$ into the third equation, we obtain
\begin{equation}\label{s4:eq:sigma}
\sigma=\frac{x\beta N_1+y\alpha N_2}{\alpha\beta S}=\frac{y\xi+x\zeta}{\alpha\beta}.
\end{equation}
Next, $k(u)=2-u(2-\log u)$, so $k(a)=2-bx(2+\sigma+\xi)$,
$k(b)=2-b(2+\sigma)$ and $k(d)=2-by(2+\sigma+\zeta)$. Dividing the last
equation of \eqref{s4:eq:bipsys} by $-b$ gives
\begin{align*}
0&=cx(2+\sigma+\xi)+(1-c)(2+\sigma)-c(2+\sigma)-(1-c)y(2+\sigma+\zeta)\\
&=(2+\sigma)(\theta_2-\theta_1)+cx\xi-(1-c)y\zeta.
\end{align*}
We multiply by $\alpha\beta\lambda S$ and use
$\lambda S(\theta_2-\theta_1)=N_2-N_1=\beta\xi-\alpha\zeta$,
$\alpha\beta\lambda Sc=\beta N_1$, $\alpha\beta\lambda S(1-c)=\alpha N_2$ and
\eqref{s4:eq:sigma}. This gives
\begin{align*}
0&=\alpha\beta(2+\sigma)(\beta\xi-\alpha\zeta)+\beta x\xi N_1-\alpha y\zeta N_2\\
&=2\alpha\beta(\beta\xi-\alpha\zeta)+(y\xi+x\zeta)(\beta\xi-\alpha\zeta)
+\beta x\xi(\zeta+y\xi)-\alpha y\zeta(\xi+x\zeta)\\
&=2\alpha\beta(\beta\xi-\alpha\zeta)+\beta y(1+x)\xi^2-\alpha x(1+y)\zeta^2
+2(\beta x-\alpha y)\xi\zeta.
\end{align*}
Since $\beta x-\alpha y=x-y=-\eta$, the right side is $J(x,y)$. By
Lemma~\ref{lem:J}, $J(x,y)>0$. This contradiction shows that $a=d$. Then the
first equation minus the second gives $0=-\mu(2c-1)(a^2-b^2)$. Since $\mu>0$
and $a<b$, we get $c=\frac12$.
\end{proof}

In Sections~\ref{sec:rigidity} and~\ref{sec:smalltail} we prove that certain
minimizers are two block SBM graphons. We next show that these minimizers are relabelings of $B_r$.

\begin{corollary}\label{cor:bipodalmin}
Let $r\in(0,r_*)$, and let $U$ be a minimizer at $r$ that is a two block SBM
graphon. Then $U$ is a relabeling of $B_r$, and $V(r)=\Vbal(r)$.
\end{corollary}

\begin{proof}
Let $P$, $c=|P|$ and $a,d,b$ be as in \eqref{s4:eq:bipodal}. By
Corollary~\ref{cor:Vbal}(a), $F(U)=V(r)\le\Vbal(r)<h(r)$. By
Lemma~\ref{lem:meangap}, an a.e.\ constant minimizer would equal $r$ a.e.,
and its cost would be $h(r)$. Hence $U$ is not a.e.\ constant. In
particular, $c\in(0,1)$ and $a,d,b$ are not all equal. By the derivation
of \eqref{s4:eq:bipsys}, $\mu_U>0$, $a,d,b\in(0,1)$, and
\eqref{s4:eq:bipsys} holds with $\mu=\mu_U$.
Proposition~\ref{prop:bipodal} now gives $a=d<b$ and $c=\frac12$. Thus
$U=B^P_{a,b}$. By Lemma~\ref{lem:relabelB}, $U$ is a relabeling of
$B_{a,b}$. Relabeling preserves $t$ and $F$, so
$t(B_{a,b})\le r^3$ and $F(B_{a,b})=V(r)\le\Vbal(r)$. Hence $(a,b)$ is a
minimizing pair at $r$. Corollary~\ref{cor:Vbal}(a) then gives
$(a,b)=(a_r,b_r)$.
Therefore $U$ is a relabeling of $B_r$, and $V(r)=F(B_r)=\Vbal(r)$.
\end{proof}

In Sections~\ref{sec:rigidity} and~\ref{sec:smalltail} we use the following
criterion to show that a graphon is a two block SBM graphon.

\begin{lemma}\label{lem:rowsbipodal}
Let $U$ be a graphon. Let $P_1\subset[0,1]$ be measurable with
$0<|P_1|<1$, and let $P_2=[0,1]\setminus P_1$. Let $f_1,f_2$ be measurable
functions on $[0,1]$ such that $U(x,y)=f_i(y)$ for a.e.\
$(x,y)\in P_i\times[0,1]$, for $i=1,2$. Then $U$ is a two block SBM graphon.
\end{lemma}

\begin{proof}
Let $i,j\in\{1,2\}$. By the symmetry of $U$, $f_i(y)=U(x,y)=U(y,x)=f_j(x)$ for
a.e.\ $(x,y)\in P_i\times P_j$. By Fubini's theorem, some $x_0\in P_i$
satisfies $f_i(y)=f_j(x_0)$ for a.e.\ $y\in P_j$. Hence $f_i$ equals a
constant $c_{ij}$ a.e.\ on $P_j$, and then $f_j=c_{ij}$ a.e.\ on $P_i$.
Thus $U=c_{ij}$ a.e.\ on $P_i\times P_j$. Moreover, $c_{12}=c_{21}$, because
both equal $f_2$ a.e.\ on $P_1$. Hence \eqref{s4:eq:bipodal} holds with
$P=P_1$, $a=c_{11}$, $d=c_{22}$ and $b=c_{12}$.
\end{proof}
\section{Minimizers at and above the threshold}\label{sec:threshold}

In this section we prove parts~(a) and~(c) of Theorem~\ref{thm:main}.
We first classify the minimizers at $r_*$, and then use scaling to treat
every $r>r_*$ (Section~\ref{s5:reduction}).
For any nonconstant minimizer $U$ at $r_*$, the pair
$(e(U),F(U))$ lies in a compact region $\Omega$
(Section~\ref{s5:region}). 
For each fixed pair, the first order conditions give an explicit
positive lower bound on the edge values $U(x,y)$ almost everywhere.
If $X,Y,Z$ are independent uniform points in $[0,1]$, then the same
conditions impose linear constraints on the joint law of
$(U(X,Y),U(Y,Z),U(Z,X))$, and the lower bound restricts this law to a
compact cube bounded away from zero. We construct dual certificates to exclude every parameter
pair except $(m_*,F_*)$. At this remaining pair, the certificates force
the law to be supported on $(a_*,a_*,a_*)$ and the permutations of
$(a_*,b_*,b_*)$ (Section~\ref{s5:relax}).
This identifies every nonconstant minimizer as a relabeling of $B_*$. 
We cover $\Omega$ by twelve convex polygons and assign one certificate
to each polygon (Section~\ref{s5:cover}).
In Appendix~\ref{app:cover} we verify the required nonnegativity and
equality cases. Throughout, we use the constants $t_*$, $m_*$, and $F_*$
from Section~\ref{s4:sec:Bstar} and the properties of $B_*$
in Corollary~\ref{cor:Bstar}.

\subsection{Reduction to the target \texorpdfstring{$r_*$}{r*} by scaling}\label{s5:reduction}

Suppose that a nonconstant graphon does at least as well as the constant at
some target. Scaling this graphon down gives a strict improvement at every
smaller target.

\begin{lemma}\label{lem:push}
Let $0<r'<r<1$, and let $U$ be a graphon that is not a.e. constant, with
$t(U)\le r^3$ and $F(U)\le h(r)$. Then $W=(r'/r)U$ is a graphon with
$t(W)\le r'^3$ and $F(W)<h(r')$. In particular $V(r')<h(r')$.
\end{lemma}

\begin{proof}
Let $a=r'/r\in(0,1)$. By Lemma~\ref{lem:scaling}, $t(W)=a^3t(U)\le r'^3$.
Moreover, \eqref{eq:scale2} with $s=r$ gives
\[
F(W)-h(r')=a\big(F(U)-h(r)\big)+(a\log a)\big(e(U)-r\big).
\]
The first term is at most $0$. Since $\log a<0$ and $e(U)>r$ by
Lemma~\ref{lem:meangap}, the second term is negative.
\end{proof}

Once $V(r_*)=h(r_*)$ is established, Lemma~\ref{lem:push} implies
that the constant graphon $r$ is the unique minimizer for every $r>r_*$.
Indeed, a nonconstant minimizer at such an $r$ would scale to a feasible
graphon at $r_*$ with cost strictly below $h(r_*)$.
In the next proposition we establish this equality and classify all
minimizers at $r_*$. Recall that $B_r=B_{a_r,b_r}$ is the balanced two block graphon attaining $\Vbal(r)$ for $0<r<r_*$, and $B_*=B_{a_*,b_*}$ is its limit as $r\uparrow r_*$.

\begin{proposition}\label{prop:critical}
$V(r_*)=h(r_*)$, and every minimizer at $r_*$ either equals $r_*$ a.e. or is
a relabeling of $B_*$.
\end{proposition}

We defer the proof of Proposition~\ref{prop:critical} to
Section~\ref{s5:proofcritical}.
Together with Lemma~\ref{lem:push}, the proposition gives parts~(i) and~(ii) of
Theorem~\ref{thm:threshold}.
Part~(iii) follows from Corollary~\ref{cor:Vbal}, which shows that
$B_r$ has strictly smaller cost than the constant graphon whenever
$r<r_*$.

\begin{theorem}\label{thm:threshold}
The following hold.
\begin{enumerate}[label=(\roman*)]
\item For every $r\in(r_*,1)$, $V(r)=h(r)$, and every minimizer at $r$
equals $r$ a.e.
\item $V(r_*)=h(r_*)$, and a graphon is a minimizer at $r_*$ if and only if
it equals $r_*$ a.e. or is a relabeling of $B_*$.
\item For every $r\in(0,r_*)$, $V(r)\le F(B_r)<h(r)$, and no minimizer at
$r$ is a.e. constant.
\end{enumerate}
Consequently, $V(r)<h(r)$ if and only if $r<r_*$. Moreover, the constant
graphon $r$ is a minimizer at $r$ if and only if $r\ge r_*$.
\end{theorem}

\begin{proof}
Let $r\in(0,1)$. The constant $r$ is feasible at $r$ with cost $h(r)$. By
Lemma~\ref{lem:meangap}, every a.e.\ constant graphon $U$ with
$t(U)\le r^3$ and $F(U)\le h(r)$ equals $r$ a.e.

If $r<r_*$, then $(a_r,b_r)$ is a minimizing pair by
Corollary~\ref{cor:Vbal}(a), so $t(B_r)\le r^3$ and $F(B_r)=\Vbal(r)<h(r)$. Hence $V(r)\le F(B_r)<h(r)$, so no minimizer at $r$
is a.e.\ constant. This is (iii). If $r>r_*$, then every graphon $U$ with
$t(U)\le r^3$ and $F(U)\le h(r)$ is a.e.\ constant. Otherwise
Lemma~\ref{lem:push} with $r'=r_*$ would give $V(r_*)<h(r_*)$, contrary to
Proposition~\ref{prop:critical}. Hence every such $U$ equals $r$ a.e.\ and
costs $h(r)$, which gives (i). Part (ii) is Proposition~\ref{prop:critical},
since relabeling preserves $t$ and $F$, and $t(B_*)=t_*$ and $F(B_*)=F_*$ by
Corollary~\ref{cor:Bstar}(a). The last two assertions follow from (i), (ii)
and (iii).
\end{proof}

In Section~\ref{sec:rigidity} we continue the minimizers $B_r$ downward from
$r_*$, starting from limits in cut distance of minimizers at targets
increasing to $r_*$. We first show that such a limit is again a minimizer.
We also need this fact at other targets in Section~\ref{sec:rigidity},
together with the convergence of the multipliers.

\begin{lemma}\label{lem:limits}
Let $r_n,r\in(0,1)$ with $r_n\to r$. For every $n$ let $U_n$ be a minimizer
at $r_n$, and let $U$ be a graphon with $\delta_\cut(U_n,U)\to0$. Then $U$
is a minimizer at $r$, $e(U_n)\to e(U)$, $F(U_n)\to F(U)$ and
$\mu_{U_n}\to\mu_U$.
\end{lemma}

\begin{proof}
Choose measure preserving bijections $\varphi_n$ with
$\|U_n^{\varphi_n}-U\|_\cut\to0$. The graphon $U_n^{\varphi_n}$ has the same
edge density, triangle density and cost as $U_n$, so it is a minimizer at
$r_n$. Hence Lemmas~\ref{lem:active} and~\ref{lem:cut}(b) give
$e(U_n)\to e(U)$ and $t(U)=\lim_nr_n^3=r^3$. By Lemma~\ref{lem:lsc}(b) and
the continuity of $V$ (Lemma~\ref{lem:existence}(b)),
$F(U)\le\liminf_nF(U_n)=\lim_nV(r_n)=V(r)$. Hence $U$ is a minimizer at $r$
and $F(U_n)=V(r_n)\to V(r)=F(U)$. Finally, \eqref{eq:mult} gives
$\mu_{U_n}\to\mu_U$.
\end{proof}

We can now identify these limits.

\begin{lemma}\label{lem:leftlimits}
Let $r_n\in(0,r_*)$ satisfy $r_n\to r_*$, and for every $n$ let $U_n$ be a
minimizer at $r_n$. If a graphon $U$ satisfies $\delta_\cut(U_{n_k},U)\to0$
along a subsequence, then $U$ is a relabeling of $B_*$.
\end{lemma}

\begin{proof}
By Lemma~\ref{lem:limits}, $U$ is a minimizer at $r_*$ and
$e(U_{n_k})\to e(U)$. By Theorem~\ref{thm:threshold}(iii), $U_{n_k}$ is not
a.e.\ constant and $F(U_{n_k})=V(r_{n_k})<h(r_{n_k})$. Hence
Lemma~\ref{lem:isolation} gives $e(U_{n_k})\ge M(r_{n_k})$, where $M$ is the
continuous function of \eqref{eq:Mr}. Letting $k\to\infty$ gives
$e(U)\ge M(r_*)=M_1>r_*$. Thus $U$ is not a.e.\ equal to $r_*$, and
Theorem~\ref{thm:threshold}(ii) shows that $U$ is a relabeling of $B_*$.
\end{proof}

\subsection{A compact region for edge density and cost}\label{s5:region}

The bounds of Section~\ref{s2:sub:bounds} confine the edge density and the
cost of a nonconstant minimizer at $r_*$ to a compact region $\Omega$. In the
rest of the section we use minimality to exclude every point of $\Omega$ other
than $(m_*,F_*)$, which is the point of $B_*$.

Let $\beta$ and $M_1=M(r_*)$ be as in \eqref{eq:beta} and \eqref{eq:Mr}. Let
$m_G=(1+\sqrt{1+8t_*})/4$ be the bound of Lemma~\ref{lem:goodman} at triangle
density $t_*$. For $m\in(0,1)$, let
$F_{\mathrm{B1}}(m)=h(m)+\beta(m)(m^3-t_*)_+^{2/3}$ be the bound of
Lemma~\ref{lem:FB1}. We write a point of the plane as $\pi=(m,F)$, where $m$
stands for an edge density and $F$ for a cost. Define
\begin{align}
\Omega&=\big\{(m,F)\in[M_1,m_G]\times\R:F_{\mathrm{B1}}(m)\le F\le F_*\big\},\label{s5:Omega}\\
\Pi&=\big\{(m,F)\in\R^2:m>0\text{ and }m+F<1\big\}.\notag
\end{align}
Since $0<r_*<M_1$ and $m_G<1$, the function $F_{\mathrm{B1}}$ is continuous on
$[M_1,m_G]$. Hence $\Omega$ is compact. Let $\pi_*=(m_*,F_*)$. The relaxation of
Section~\ref{s5:relax} is defined for parameters in $\Pi$.

\begin{lemma}\label{lem:region}
The following hold.
\begin{enumerate}[label=(\alph*)]
\item Let $U$ be a graphon with $t(U)=t_*$, $F(U)\le F_*$ and
$|\{U\ne r_*\}|>0$. Then $(e(U),F(U))\in\Omega$.
\item The conclusion of (a) holds for every minimizer $U$ at $r_*$ that is
not a.e. constant.
\item $\pi_*\in\Omega$ and $\Omega\subset\Pi$.
\end{enumerate}
\end{lemma}

\begin{proof}
(a) Let $m=e(U)$. Lemma~\ref{lem:isolation} with $r=r_*$ gives $m\ge M_1$, and
Lemma~\ref{lem:goodman} gives $m\le m_G$. Hence $m\in(0,1)$, so
Lemma~\ref{lem:FB1} gives $F(U)\ge F_{\mathrm{B1}}(m)$.

(b) By Lemma~\ref{lem:active}, $t(U)=t_*$. The constant $r_*$ is feasible, so
$F(U)=V(r_*)\le F_*$. Since $U$ is not a.e. constant, $|\{U\ne r_*\}|>0$.

(c) Since $t(B_*)=t_*$, $F(B_*)=F_*$ and $B_*$ takes only the values
$a_*,b_*\ne r_*$, part (a) applies to $B_*$ and gives $\pi_*\in\Omega$.
Every $(m,F)\in\Omega$ satisfies $m\ge M_1>0$ and $m+F\le m_G+F_*$, so it
remains to show $m_G+F_*<1$. By Lemma~\ref{lem:digits}, $0.205<r_*<0.215$, so
$t_*<0.01$. Hence $m_G<(1+\sqrt{1.08})/4<0.51$, since $1.04^2>1.08$.
Note that $e^{1.5}<9/2<1/0.205$, because $e<2.72$ and $2.72^3<81/4$. Since
$u\mapsto u(1-\log u)$ is increasing on $(0,1)$, we get
$1-F_*=r_*(1-\log r_*)>0.205(1-\log0.205)>0.205\cdot2.5>0.51>m_G$.
\end{proof}

\subsection{A relaxation by the law of a random triangle}\label{s5:relax}

Let $U$ be a minimizer at $r_*$. The relaxation of this subsection uses the
multiplier $\mu_U$ of \eqref{eq:mult}, the constant $\bar k_U$ of
Proposition~\ref{prop:vertex} and a lower bound for $U$. All three are
explicit functions of the point $(e(U),F(U))$. The lower bound comes from
the vertex identity \eqref{eq:vertex}, which bounds the rows of $U$ in
$L^2$, and from the Euler--Lagrange equation \eqref{eq:euler}. We prove this bound
at every target, since we use it again in Section~\ref{sec:intermediate}. For
$t>0$ and real numbers $m$ and $y$, define
\begin{equation}\label{eq:Lambda}
\Lambda_t(m,y)=\exp\Big(-\frac{(1-y)^2-m^2}{2t}\Big).
\end{equation}

\begin{lemma}\label{lem:support}
Let $U$ be a minimizer at $r\in(0,1)$, and let $m=e(U)$ and $F=F(U)$. Then
$m\ge r$, $\mu_Ur^3=1-F-m>0$ and $\bar k_U=1+F-m$. Moreover, $U\ge\ell_U$
a.e., where $\ell_U=\Lambda_{r^3}(m,F)$.
\end{lemma}

\begin{proof}
The identities are \eqref{eq:mult} and Proposition~\ref{prop:vertex}, and
$\mu_U>0$ by Proposition~\ref{prop:euler}(c). The constant $r$ is feasible,
so $F\le h(r)$, and Lemma~\ref{lem:meangap} gives $m\ge r$. Next, the
function $k$ of \eqref{eq:k} satisfies $k(u)\le2-2u^2$ on $[0,1]$. Indeed,
both sides equal $2$ at $u=0$, and for $0<u\le1$ the ratio
$(2-k(u))/u^2=(2-\log u)/u$ has the negative derivative $(\log u-3)/u^2$, so
it is at least its value $2$ at $u=1$. Hence, by \eqref{eq:vertex}, for
a.e.\ $x$,
\[
\bar k_U=\int_0^1k(U(x,y))\dd y\le2-2\int_0^1U(x,y)^2\dd y .
\]
For two such points $x,y$ the Cauchy--Schwarz inequality gives
$C_U(x,y)\le\|U(x,\cdot)\|_2\|U(y,\cdot)\|_2\le(2-\bar k_U)/2=(1-F+m)/2$.
Finally, \eqref{eq:euler} and $\mu_U>0$ give a.e.
\[
U=\exp(-\mu_UC_U)\ge\exp\Big(-\frac{(1-F-m)(1-F+m)}{2r^3}\Big)=\Lambda_{r^3}(m,F) .
\qedhere
\]
\end{proof}

At the target $r_*$, we regard $\mu_U$, $\bar k_U$ and $\ell_U$ as functions
of a point of the plane. For $\pi=(m,F)\in\R^2$ define
\begin{align}
\mu(\pi)&=\frac{1-F-m}{t_*},\quad \bar k(\pi)=1+F-m,\notag\\
\ell(\pi)&=\Lambda_{t_*}(m,F)=\exp\Big(-\frac{(1-F)^2-m^2}{2t_*}\Big).\label{s5:params}
\end{align}
If $\pi\in\Pi$, then $1-F>m>0$, so $\mu(\pi)>0$ and $0<\ell(\pi)<1$.

The relaxation of Definition~\ref{def:relaxed} keeps finitely many
consequences of the Euler--Lagrange equation \eqref{eq:euler} and the vertex
identity \eqref{eq:vertex}. It uses the following test functions of an edge
value $u$:
\[
(\varphi_1,\dots,\varphi_5)=(1,u,u^2,u^3,u\log u),\quad
(\chi_1,\dots,\chi_4)=(u,u^3,u^4,u\log u).
\]
In the relaxation we multiply the residual of \eqref{eq:euler} by each
$\varphi_j$ and the residual of \eqref{eq:vertex} by each $\chi_j$, and we
integrate the products. We also add one inequality built from the function
\[
g_*(u)=u+\omega_*h(u)\quad(u\ge0),\quad\text{where}\quad
\omega_*=\frac{b_*-a_*}{h(a_*)-h(b_*)} .
\]
Here $\omega_*>0$ because $h$ is strictly decreasing on $[0,1]$. We have chosen
$\omega_*$ so that $g_*(a_*)=g_*(b_*)$. Averaging these two values
and using $F_*=(h(a_*)+h(b_*))/2$ gives
\begin{equation}\label{s5:gstar}
g_*(a_*)=g_*(b_*)=m_*+\omega_*F_*.
\end{equation}
Hence $g_*\circ B_*$ is constant.

\begin{definition}\label{def:trilaw}
Let $U$ be a graphon and let $X_1,X_2,X_3$ be independent and uniform on
$[0,1]$. Let $\nu_U$ be the law on $[0,1]^3$ of
\[
(A,B,C)=\big(U(X_1,X_2),U(X_1,X_3),U(X_2,X_3)\big).
\]
We call $\nu_U$ the \emph{law of the edge values of a random triangle} in
$U$.
\end{definition}

Since $U$ is symmetric, exchanging $X_1$ and $X_2$ maps $(A,B,C)$ to
$(A,C,B)$, and exchanging $X_2$ and $X_3$ maps it to $(B,A,C)$. These
exchanges preserve the law of $(X_1,X_2,X_3)$, and the two transpositions
generate all permutations. Hence $\nu_U$ is invariant under permutations of
the coordinates.

\begin{definition}\label{def:relaxed}
Let $k$ be the function of \eqref{eq:k}. For $\pi=(m,F)\in\Pi$, let
$\mathcal T(\pi)$ be the set of Borel probability measures $\nu$ on
$[0,1]^3$ that are invariant under permutations of the coordinates and
satisfy the following conditions, in which $(A,B,C)$ has law
$\nu$.
\begin{enumerate}[label=\textup{(T\arabic*)}]
\item\label{s5:T1} $\E[A]=m$, $\E[ABC]=t_*$ and $\E[h(A)]=F$.
\item\label{s5:T2} $\E[\varphi_j(A)(\log A+\mu(\pi)BC)]=0$ for
$j=1,\dots,5$.
\item\label{s5:T3} $\E[\chi_j(A)(k(B)+k(C)-2\bar k(\pi))]=0$ for
$j=1,\dots,4$.
\item\label{s5:T4} $\E[g_*(A)^2]-2\E[g_*(A)g_*(B)]+(m+\omega_*F)^2\ge0$.
\item\label{s5:T5} $\nu([\ell(\pi),1]^3)=1$.
\end{enumerate}
\end{definition}

By \ref{s5:T5}, all integrands are bounded on a set that carries $\nu$. For
fixed $\pi$, every condition is linear in $\nu$. If $\nu=\nu_U$ for a graphon
$U$ with edge density $m$ and cost $F$, then $\E[g_*(A)]=m+\omega_*F$, so
\ref{s5:T4} holds by the next lemma.

\begin{lemma}\label{lem:hoeffding}
Let $U$ be a graphon, and let $(A,B,C)$ have law $\nu_U$. Let $H=g_*\circ U$,
$d(x)=\int_0^1H(x,y)\dd y$ and $\eta=\iint H$. Then
\[
\E[g_*(A)^2]-2\E[g_*(A)g_*(B)]+\E[g_*(A)]^2
=\iint\big(H(x,y)-d(x)-d(y)+\eta\big)^2\dd x\dd y\ge0 .
\]
In particular, both sides vanish if $g_*\circ U$ is a.e.\ constant.
\end{lemma}

\begin{proof}
Since $g_*$ is continuous on $[0,1]$, $H$ is a bounded symmetric kernel.
Hence $\int_0^1H(x,y)\dd x=d(y)$, $\int_0^1d=\eta$, and
$\iint H(x,y)d(x)\dd x\dd y=\iint H(x,y)d(y)\dd x\dd y=\int_0^1d^2$.
Expanding the square with these identities shows that the right side equals
$\iint H^2-2\int_0^1d^2+\eta^2$. By the Fubini theorem,
$\E[g_*(A)^2]=\iint H^2$ and $\E[g_*(A)]=\eta$. Since $A$ and $B$ share
only the vertex $X_1$,
\[
\E[g_*(A)g_*(B)]=\int_0^1\Big(\int_0^1H(x,y)\dd y\Big)\Big(\int_0^1H(x,z)\dd z\Big)\dd x
=\int_0^1d^2 .
\]
If $H$ equals a constant $c$ a.e., then $d=\eta=c$, and the integrand on the
right vanishes.
\end{proof}

We next show that the law $\nu_U$ of a minimizer $U$ at $r_*$ is
feasible for the relaxation at the parameter $(e(U),F(U))$ of $U$.

\begin{lemma}\label{lem:valid}
Let $U$ be a minimizer at $r_*$, and let $\pi=(e(U),F(U))$. Then $\pi\in\Pi$,
$\mu_U=\mu(\pi)$, and $\nu_U\in\mathcal T(\pi)$.
\end{lemma}

\begin{proof}
By Lemma~\ref{lem:support} and \eqref{s5:params}, $e(U)\ge r_*>0$,
$1-F(U)-e(U)=\mu_Ut_*>0$, $\mu_U=\mu(\pi)$, $\bar k_U=\bar k(\pi)$ and
$U\ge\ell(\pi)$ a.e. In particular $\pi\in\Pi$. Let $(A,B,C)$ be as in
Definition~\ref{def:trilaw}.

\ref{s5:T5}: Each of $A,B,C$ has the law of $U(X_1,X_2)$, and
$|\{U<\ell(\pi)\}|=0$.

\ref{s5:T1}: $\E[A]=e(U)$, $\E[h(A)]=F(U)$, and $\E[ABC]=t(U)=t_*$ by
Lemma~\ref{lem:active}.

\ref{s5:T2}: Let $\varphi$ be bounded and Borel on $[\ell(\pi),1]$. Since
$\E[BC\mid X_1,X_2]=C_U(X_1,X_2)$,
\[
\E[\varphi(A)(\log A+\mu_UBC)]=\iint\varphi(U)\big(\log U+\mu_UC_U\big)=0
\]
by \eqref{eq:euler}. The integrals are finite since $\ell(\pi)\le U\le1$ a.e.

\ref{s5:T3}: Let $\chi$ be bounded and Borel on $[\ell(\pi),1]$. Let
$K(x)=\int_0^1k(U(x,z))\dd z$. By \eqref{eq:vertex}, $K(x)$ equals
$\bar k(\pi)$ for a.e. $x$. By the Fubini theorem,
\[
\E[\chi(A)k(B)]=\iint\chi(U(x,y))K(x)\dd x\dd y=\bar k(\pi)\E[\chi(A)].
\]
In the same way, with $K(y)$ in place of $K(x)$, we get
$\E[\chi(A)k(C)]=\bar k(\pi)\E[\chi(A)]$.

\ref{s5:T4}: This is Lemma~\ref{lem:hoeffding}, since
$\E[g_*(A)]=e(U)+\omega_*F(U)$.
\end{proof}


A certificate is a combination of the integrands of
Definition~\ref{def:relaxed} that has nonpositive mean under every
$\nu\in\mathcal T(\pi)$. If it is positive on $[\ell(\pi),1]^3$, then
$\mathcal T(\pi)$ is empty. For
$u=(u_1,u_2,u_3)\in(0,\infty)^3$ and $\pi=(m,F)\in\R^2$, define the functions
$\mathcal E_0,\dots,\mathcal E_{11}$ by
\begin{align}
\mathcal E_0(u;\pi)&=\frac13\sum_{i=1}^3u_i-m,\quad
\mathcal E_1(u;\pi)=u_1u_2u_3-t_*,\notag\\
\mathcal E_2(u;\pi)&=\frac13\sum_{i=1}^3h(u_i)-F,\notag\\
\mathcal E_{2+j}(u;\pi)&=\frac13\sum_{i=1}^3\varphi_j(u_i)\Big(\log u_i+\mu(\pi)\frac{u_1u_2u_3}{u_i}\Big)
\quad(j=1,\dots,5),\label{eq:features}\\
\mathcal E_{7+j}(u;\pi)&=\frac13\sum_{i=1}^3\chi_j(u_i)\Big(\sum_{l\ne i}k(u_l)-2\bar k(\pi)\Big)
\quad(j=1,\dots,4),\notag
\end{align}
and define the function
\[
\mathcal H_*(u;\pi)=\frac13\sum_{i=1}^3g_*(u_i)^2-\frac23\sum_{1\le i<l\le3}g_*(u_i)g_*(u_l)+(m+\omega_*F)^2 .
\]
Each of these functions averages an integrand of Definition~\ref{def:relaxed}
over the permutations of the coordinates.

\begin{definition}\label{def:certificate}
Let $\theta=(\theta_0,\dots,\theta_{11})\in\R^{12}$ and $\lambda\ge0$. The
\emph{certificate function} is
\begin{equation}\label{s5:Psi}
\Psi_{\theta,\lambda}(u;\pi)=-\sum_{i=0}^{11}\theta_i\mathcal E_i(u;\pi)-\lambda\mathcal H_*(u;\pi),
\quad u\in(0,\infty)^3,\quad\pi\in\R^2 .
\end{equation}
\end{definition}

The function $\Psi_{\theta,\lambda}(\cdot;\pi)$ is real analytic and
invariant under permutations of the coordinates. We next prove that a certificate has nonpositive mean under every $\nu\in\mathcal T(\pi)$
(weak duality). Part~(d) of the lemma, the concavity in the parameter,
lets a single certificate serve a whole polygon of parameters
(Section~\ref{s5:cover}).

\begin{lemma}\label{lem:weakdual}
Let $\theta\in\R^{12}$, $\lambda\ge0$, $\pi\in\Pi$ and $0<\ell_0\le\ell(\pi)$.
\begin{enumerate}[label=(\alph*)]
\item Every $\nu\in\mathcal T(\pi)$ satisfies
$\E_\nu[\Psi_{\theta,\lambda}(\cdot;\pi)]=-\lambda\E_\nu[\mathcal H_*(\cdot;\pi)]\le0$.
\item If $\Psi_{\theta,\lambda}(\cdot;\pi)>0$ on $[\ell_0,1]^3$, then
$\mathcal T(\pi)=\emptyset$.
\item If $\Psi_{\theta,\lambda}(\cdot;\pi)\ge0$ on $[\ell_0,1]^3$, then every
$\nu\in\mathcal T(\pi)$ is carried by the zero set of
$\Psi_{\theta,\lambda}(\cdot;\pi)$ in $[\ell_0,1]^3$.
\item For every $u\in(0,\infty)^3$, the function $\Psi_{\theta,\lambda}(u;\cdot)$
is concave on $\R^2$.
\end{enumerate}
\end{lemma}

\begin{proof}
(a) Let $(u_1,u_2,u_3)$ have law $\nu\in\mathcal T(\pi)$. By \ref{s5:T5},
$\nu$ is carried by $[\ell(\pi),1]^3$, where all functions involved are
bounded. For distinct $i,l,l'$, permutation invariance shows that
$(u_i,u_l,u_{l'})$ has the law of $(A,B,C)$. Hence the expectation of each
summand in \eqref{eq:features} equals that of an integrand of \ref{s5:T1},
\ref{s5:T2} or \ref{s5:T3}, so $\E_\nu[\mathcal E_i(\cdot;\pi)]=0$ for every
$i$. In the same way, $\E_\nu[\mathcal H_*(\cdot;\pi)]$ is the left side of
\ref{s5:T4}, which is nonnegative. Since $\lambda\ge0$, \eqref{s5:Psi}
gives (a).

(b) and (c). The function $\Psi_{\theta,\lambda}(\cdot;\pi)$ is continuous
on the compact cube $[\ell_0,1]^3$, which carries every
$\nu\in\mathcal T(\pi)$. In (b), its minimum $\varepsilon$ there is positive,
so every $\nu\in\mathcal T(\pi)$ would satisfy
$\E_\nu[\Psi_{\theta,\lambda}(\cdot;\pi)]\ge\varepsilon>0$, contrary to (a).
In (c), it is nonnegative with expectation at most $0$ by (a), so it vanishes
$\nu$ almost surely.

(d) By \eqref{s5:params}, $\mu(\pi)$ and $\bar k(\pi)$ are affine in $\pi$,
so every $\mathcal E_i(u;\pi)$ is affine in $\pi$. In $\mathcal H_*(u;\pi)$,
only the term $(m+\omega_*F)^2$ depends on $\pi$. Hence
$\Psi_{\theta,\lambda}(u;\pi)$ is the sum of an affine function of $\pi$ and
$-\lambda(m+\omega_*F)^2$, which is concave since $\lambda\ge0$.
\end{proof}

We will choose the coefficients $\theta_i$ and $\lambda\ge0$
separately for each parameter polygon (see Figure~\ref{fig:cover}).
Concavity allows us to verify nonnegativity on a common cube by
checking the polygon vertices.
For each $\pi\in\Omega\setminus\{\pi_*\}$, we will obtain a certificate
that is strictly positive on $[\ell(\pi),1]^3$.
Weak duality then excludes any minimizer with this parameter pair.
At $\pi_*=(m_*,F_*)$, we require a nonnegative certificate whose zeros
on this cube are exactly $(a_*,a_*,a_*)$ and the permutations of
$(a_*,b_*,b_*)$.
The law $\nu_U$ of a minimizer $U$ with parameter pair $\pi_*$ is then
supported on these triples. Since also $e(U)=m_*=(a_*+b_*)/2$,
Lemma~\ref{lem:patterns} identifies $U$ as a relabeling of $B_*$.

Lemma~\ref{lem:patterns} also records the identities needed to
construct these certificates. We state it for general $B_{a,b}$,
since we use it at $B_*$ in this section and at $B_z$ in
Section~\ref{sec:intermediate}.
For $P\in\R^3$, let $\mathcal O(P)$ be the set of permutations of $P$.
For $a,b\in[0,1]$, let $P_1=(a,a,a)$ and $P_2=(a,b,b)$.

\begin{lemma}\label{lem:patterns}
Let $a,b\in[0,1]$ with $a\ne b$.
\begin{enumerate}[label=(\alph*)]
\item The law $\nu_{B_{a,b}}$ is
\begin{equation}\label{s5:nustar}
\nu_{B_{a,b}}=\frac14\big(\delta_{(a,a,a)}+\delta_{(a,b,b)}+\delta_{(b,a,b)}+\delta_{(b,b,a)}\big).
\end{equation}
Hence $\E_{\nu_{B_{a,b}}}[\Psi]=\frac14\Psi(P_1)+\frac34\Psi(P_2)$ for every
function $\Psi$ on $\{a,b\}^3$ that is invariant under permutations of the
coordinates.
\item Let $G\subset\R^3$ be open, and let $\Psi:G\to\R$ be differentiable at
$P\in G$. Assume that $G$ and $\Psi$ are invariant under permutations of the
coordinates. If $P=(p,p,p)$, then
$\partial_1\Psi(P)=\partial_2\Psi(P)=\partial_3\Psi(P)$. If $P=(p,p',p')$,
then $\partial_2\Psi(P)=\partial_3\Psi(P)$.
\item If a graphon $U$ satisfies
$\nu_U(\mathcal O(P_1)\cup\mathcal O(P_2))=1$ and $e(U)=(a+b)/2$, then $U$
is a relabeling of $B_{a,b}$.
\end{enumerate}
\end{lemma}

\begin{proof}
(a) With probability $\frac14$, the points $X_1,X_2,X_3$ of
Definition~\ref{def:trilaw} lie in the same block of $B_{a,b}$, and all three
edge values equal $a$. Otherwise exactly one $X_i$ lies in the other block,
each with probability $\frac14$, and the two edges at $X_i$ take the value
$b$. This gives \eqref{s5:nustar}. The formula for $\E_{\nu_{B_{a,b}}}[\Psi]$ then follows from the permutation
invariance of $\Psi$.

(b) Differentiating $\Psi(u_{\sigma(1)},u_{\sigma(2)},u_{\sigma(3)})=\Psi(u)$
at a point $P$ fixed by a permutation $\sigma$ gives
$\partial_i\Psi(P)=\partial_{\sigma^{-1}(i)}\Psi(P)$. Every permutation
fixes $(p,p,p)$, and the transposition of the last two coordinates fixes
$(p,p',p')$.

(c) For a.e.\ $(x,y,z)$, the triple $(U(x,y),U(x,z),U(y,z))$ lies in
$\mathcal O(P_1)\cup\mathcal O(P_2)\subset\{a,b\}^3$, so it has an even
number of coordinates equal to $b$. By the Fubini theorem, $U\in\{a,b\}$
a.e. Since $e(U)=(a+b)/2$, Lemma~\ref{lem:parity} gives the claim.
\end{proof}

We now turn to the parameter $\pi_*$ of $B_*$. Let $\bar k_*=\bar k(\pi_*)$
and $\ell_*=\ell(\pi_*)$. From now on, let $P_1=(a_*,a_*,a_*)$ and
$P_2=(a_*,b_*,b_*)$, and let $\nu_*=\nu_{B_*}$.

In Appendix~\ref{app:cover} we check positivity in the coordinates $q=-\log u$,
where the logarithm is taken in each coordinate. Let $q_a=-\log a_*$ and
$q_b=-\log b_*$. Then $q_a>q_b>0$. Let $\hat P_1=(q_a,q_a,q_a)$ and
$\hat P_2=(q_a,q_b,q_b)$, and define
\[
\Phi_{\theta,\lambda,\pi}(q)=\Psi_{\theta,\lambda}\big((e^{-q_1},e^{-q_2},e^{-q_3});\pi\big),\quad q\in\R^3 .
\]
For $Q>0$, the map $q\mapsto(e^{-q_1},e^{-q_2},e^{-q_3})$ is a bijection from
$[0,Q]^3$ onto $[e^{-Q},1]^3$ that commutes with permutations. Hence it
carries the positivity, the nonnegativity and the zeros of
$\Phi_{\theta,\lambda,\pi}$ on $[0,Q]^3$ to those of
$\Psi_{\theta,\lambda}(\cdot;\pi)$ on $[e^{-Q},1]^3$. By symmetry, it
suffices to consider the ordered set
\[
\mathcal D_Q=\{q\in\R^3:Q\ge q_1\ge q_2\ge q_3\ge0\}.
\]
In Lemma~\ref{lem:apex} we collect the facts at $\pi_*$ used in Section~\ref{s5:cover}.

\begin{lemma}\label{lem:apex}
The following hold.
\begin{enumerate}[label=(\alph*)]
\item $\pi_*\in\Pi$, $\mu(\pi_*)=\mu_*$ and $\ell_*<a_*$. Moreover
$\nu_*\in\mathcal T(\pi_*)$ and $\E_{\nu_*}[\mathcal H_*(\cdot;\pi_*)]=0$.
\item For every $\theta\in\R^{12}$ and $\lambda\ge0$,
$\Psi_{\theta,\lambda}(P_1;\pi_*)+3\Psi_{\theta,\lambda}(P_2;\pi_*)=0$.
\end{enumerate}
\end{lemma}

\begin{proof}
(a) By Corollary~\ref{cor:Bstar}(c), $1-F_*-m_*=\mu_*t_*>0$, so $\pi_*\in\Pi$
and $\mu(\pi_*)=\mu_*$. Moreover, $2-\bar k_*=1-F_*+m_*=\mu_*t_*+a_*+b_*$, so
$-\log\ell_*=\mu_*(\mu_*t_*+a_*+b_*)/2$. Since
$\log a_*=-\mu_*(a_*^2+b_*^2)/2$ by Corollary~\ref{cor:Bstar}(b), we get
\[
\log\frac{a_*}{\ell_*}=\frac{\mu_*}2\big(\mu_*t_*+a_*(1-a_*)+b_*(1-b_*)\big)>0 .
\]
By \eqref{s5:nustar}, $\nu_*$ is permutation invariant and carried by
$\{a_*,b_*\}^3\subset[\ell_*,1]^3$, which is \ref{s5:T5}. Condition
\ref{s5:T1} holds because $e(B_*)=m_*$, $t(B_*)=t_*$ and $F(B_*)=F_*$. By
Corollary~\ref{cor:Bstar}(b),(c), $\log B_*=-\mu_*C_{B_*}$ everywhere and
$\int_0^1k(B_*(x,y))\dd y=\bar k_*$ for every $x$. Hence the computations in
the proof of Lemma~\ref{lem:valid} give \ref{s5:T2} and \ref{s5:T3}.
Finally, by permutation invariance,
$\E_{\nu_*}[\mathcal H_*(\cdot;\pi_*)]$ is the left side of \ref{s5:T4} for
$\nu_*$. Since $g_*\circ B_*$ is the constant $m_*+\omega_*F_*$ by
\eqref{s5:gstar}, we have $\E_{\nu_*}[g_*(A)]=m_*+\omega_*F_*$, and the
last claim of Lemma~\ref{lem:hoeffding} shows that this left side is $0$.
This gives \ref{s5:T4} and the last claim.

(b) Let $\Psi=\Psi_{\theta,\lambda}(\cdot;\pi_*)$. By
Lemma~\ref{lem:patterns}(a), $\E_{\nu_*}[\Psi]=\frac14\Psi(P_1)+\frac34\Psi(P_2)$.
This vanishes by (a) and Lemma~\ref{lem:weakdual}(a).
\end{proof}

Let $\Psi=\Psi_{\theta,\lambda}(\cdot;\pi_*)$, and suppose that $\Psi$ is
nonnegative on a cube $[\ell_0,1]^3$ with $0<\ell_0<a_*$. Then
Lemma~\ref{lem:apex}(b) forces $\Psi(P_1)=\Psi(P_2)=0$. Since $P_1$ and $P_2$
are then interior minimum points of $\Psi$, we get
$\nabla\Psi(P_1)=\nabla\Psi(P_2)=0$. Conversely, by
Lemmas~\ref{lem:apex}(b) and~\ref{lem:patterns}(b), these zeros and critical
points follow from the four linear \emph{tangency equations}
\begin{equation}\label{s5:tangency}
\Psi(P_1)=0,\quad \partial_{u_1}\Psi(P_1)=0,\quad \partial_{u_1}\Psi(P_2)=0,\quad
\partial_{u_2}\Psi(P_2)=0 .
\end{equation}
Positivity near $P_1$ and $P_2$ then follows from the third order Taylor
bound of Lemma~\ref{lem:taylor}, which we apply in Appendix~\ref{app:cover}
at $\hat P_1$ and $\hat P_2$. For $y\in\R^3$, write $|y|_2$ for the Euclidean
norm and $|y|_\infty=\max_i|y_i|$.

\begin{lemma}\label{lem:taylor}
Let $\Phi:\R^3\to\R$ be three times continuously differentiable, $P\in\R^3$,
$\delta>0$ and $\mathcal C_\delta=\{q\in\R^3:|q-P|_\infty\le\delta\}$. Assume that
$\Phi(P)=0$ and $\nabla\Phi(P)=0$. Assume also that
$\nabla^2\Phi(P)-\lambda_0I$ is positive semidefinite for some
$\lambda_0>0$, and that
\[
S_3\ge\sum_{i,j,k=1}^3\max_{q\in \mathcal C_\delta}|\partial_i\partial_j\partial_k\Phi(q)|,\quad
\delta S_3<3\lambda_0 .
\]
Then for every $y\in\R^3$ with $0<|y|_\infty\le\delta$,
\[
\Phi(P+y)\ge\frac{\lambda_0}2|y|_2^2-\frac{S_3}6|y|_\infty^3
\ge\Big(\frac{\lambda_0}2-\frac{\delta S_3}6\Big)|y|_\infty^2>0 .
\]
\end{lemma}

\begin{proof}
Let $f(s)=\Phi(P+sy)$ for $s\in[0,1]$. Then $f(0)=f'(0)=0$ and
$f''(0)=y^{\top}\nabla^2\Phi(P)y\ge\lambda_0|y|_2^2$. Since
$P+sy\in\mathcal C_\delta$,
$|f'''(s)|\le\sum_{i,j,k}|\partial_i\partial_j\partial_k\Phi(P+sy)||y_iy_jy_k|\le S_3|y|_\infty^3$.
Hence Taylor's theorem with the Lagrange remainder,
$f(1)=f(0)+f'(0)+\frac12f''(0)+\frac16f'''(\xi)$ for some $\xi\in(0,1)$,
gives the first inequality. The second follows from $|y|_2\ge|y|_\infty$
and $|y|_\infty\le\delta$, and the last from $\delta S_3<3\lambda_0$.
\end{proof}

\subsection{Covering the region by twelve cells}\label{s5:cover}

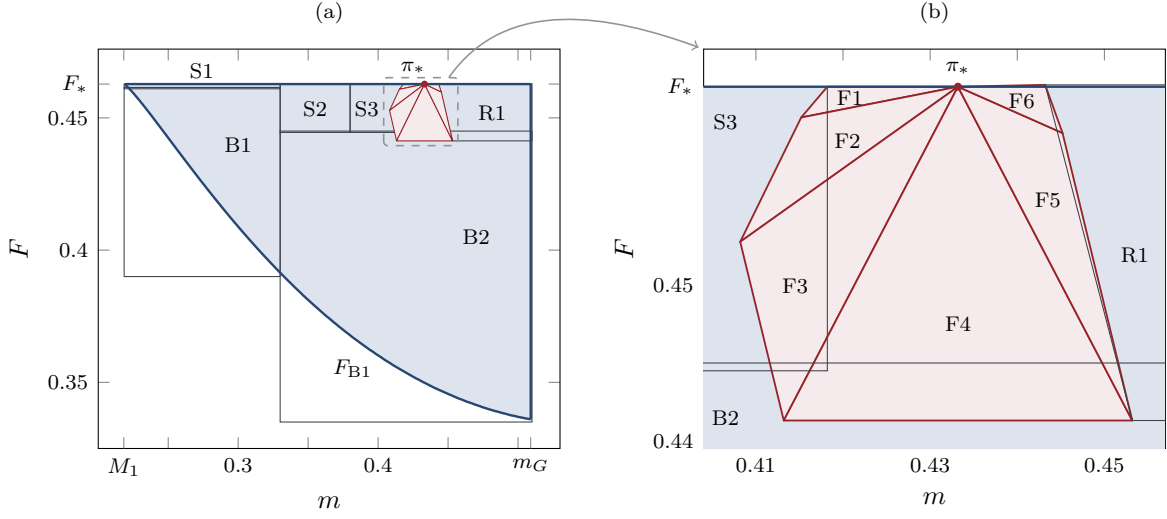
\begin{figure}[t]
\centering
\pgfplotsset{cover/.style={scale only axis,width=0.37\textwidth,height=0.32\textwidth,
  xlabel={$m$},ylabel={$F$},
  tick label style={font=\scriptsize,/pgf/number format/fixed},
  label style={font=\small},title style={font=\scriptsize},
  extra y ticks={0.46277975},extra y tick labels={$F_*$},
  extra tick style={grid=none}}}
\tikzset{box/.style={figdark,line width=0.4pt},
  fan/.style={figred,line width=0.7pt,line join=round,fill=figred!9},
  celllabel/.style={font=\scriptsize,inner sep=1pt},
  zoom/.style={gray!85,line width=0.6pt,dashed,rounded corners=1.5pt}}
\begin{tikzpicture}
\begin{axis}[cover,title={(a)},name=left,
  xmin=0.2,xmax=0.53,ymin=0.325,ymax=0.476,
  xtick={0.25,0.3,0.35,0.4,0.45,0.5},xticklabels={,$0.3$,,$0.4$,,},
  ytick={0.35,0.4,0.45},
  extra x ticks={0.21827910,0.50905344},extra x tick labels={$M_1$,$m_G$},
  declare function={hh(\u)=\u*ln(\u)-\u+1;
    bb(\u)=(\u-1-ln(\u))/(1-\u)^2;
    FB(\u)=hh(\u)+bb(\u)*(\u^3-0.00921736687)^(2/3);}]
\addplot[draw=none,fill=figgray,domain=0.2188573:0.5090534,samples=160]
  {FB(x)} -- (axis cs:0.5090534,0.46277975) -- (axis cs:0.2188573,0.46277975) -- cycle;
\draw[box] (0.2185,0.39) rectangle (0.33,0.4615);
\draw[box] (0.2185,0.461) rectangle (0.33,0.46278);
\draw[box] (0.33,0.335) rectangle (0.51,0.445);
\draw[box] (0.33,0.4445) rectangle (0.38,0.46278);
\draw[box] (0.38,0.4445) rectangle (0.4182,0.46278);
\draw[box] (0.4532,0.4413) -- (0.51,0.4413) -- (0.51,0.4629) -- (0.4432,0.4629) -- cycle;
\draw[fan,line width=0.4pt] (0.43319048,0.46277975) -- (0.4182,0.46278) -- (0.4152,0.4608) -- cycle;
\draw[fan,line width=0.4pt] (0.43319048,0.46277975) -- (0.4152,0.4608) -- (0.4082,0.4528) -- cycle;
\draw[fan,line width=0.4pt] (0.43319048,0.46277975) -- (0.4082,0.4528) -- (0.4132,0.4413) -- cycle;
\draw[fan,line width=0.4pt] (0.43319048,0.46277975) -- (0.4132,0.4413) -- (0.4532,0.4413) -- cycle;
\draw[fan,line width=0.4pt] (0.43319048,0.46277975) -- (0.4532,0.4413) -- (0.4452,0.4598) -- cycle;
\draw[fan,line width=0.4pt] (0.43319048,0.46277975) -- (0.4452,0.4598) -- (0.4432,0.4629) -- cycle;
\addplot[figblue,line width=0.9pt,domain=0.2188573:0.5090534,samples=160] {FB(x)};
\draw[figblue,line width=0.9pt] (0.5090534,0.3359913) -- (0.5090534,0.46277975);
\draw[figblue,line width=0.9pt] (0.2188573,0.46277975) -- (0.5090534,0.46277975);
\fill[figred] (0.43319048,0.46277975) circle[radius=1.3pt];
\draw[zoom] (0.404,0.4395) rectangle (0.457,0.4652);
\coordinate (zoomtop) at (axis cs:0.4505,0.4652);
\node[celllabel] at (0.300,0.440) {B1};
\node[celllabel,anchor=south] at (0.274,0.4637) {S1};
\node[celllabel] at (0.355,0.4535) {S2};
\node[celllabel] at (0.392,0.4535) {S3};
\node[celllabel] at (0.4815,0.452) {R1};
\node[celllabel] at (0.470,0.405) {B2};
\node[celllabel,anchor=south] at (0.4250,0.4656) {$\pi_*$};
\node[celllabel,anchor=north east] at (0.398,0.360) {$F_{\mathrm{B1}}$};
\end{axis}
\begin{axis}[cover,title={(b)},name=right,
  at={(left.south east)},anchor=south west,xshift=1.9cm,
  xmin=0.404,xmax=0.457,ymin=0.4395,ymax=0.4652,
  xtick={0.41,0.42,0.43,0.44,0.45},xticklabels={$0.41$,,$0.43$,,$0.45$},
  ytick={0.44,0.45},
  yticklabel style={/pgf/number format/precision=2}]
\fill[figgray] (0.404,0.4385) rectangle (0.457,0.46277975);
\draw[fan] (0.43319048,0.46277975) -- (0.4182,0.46278) -- (0.4152,0.4608) -- cycle;
\draw[fan] (0.43319048,0.46277975) -- (0.4152,0.4608) -- (0.4082,0.4528) -- cycle;
\draw[fan] (0.43319048,0.46277975) -- (0.4082,0.4528) -- (0.4132,0.4413) -- cycle;
\draw[fan] (0.43319048,0.46277975) -- (0.4132,0.4413) -- (0.4532,0.4413) -- cycle;
\draw[fan] (0.43319048,0.46277975) -- (0.4532,0.4413) -- (0.4452,0.4598) -- cycle;
\draw[fan] (0.43319048,0.46277975) -- (0.4452,0.4598) -- (0.4432,0.4629) -- cycle;
\draw[box] (0.38,0.4445) rectangle (0.4182,0.46278);
\draw[box] (0.33,0.335) rectangle (0.51,0.445);
\draw[box] (0.4532,0.4413) -- (0.51,0.4413) -- (0.51,0.4629) -- (0.4432,0.4629) -- cycle;
\draw[figblue,line width=0.9pt] (0.404,0.46277975) -- (0.457,0.46277975);
\fill[figred] (0.43319048,0.46277975) circle[radius=1.5pt];
\node[celllabel] at (0.4208,0.46208) {F1};
\node[celllabel] at (0.4205,0.4593) {F2};
\node[celllabel] at (0.4145,0.4500) {F3};
\node[celllabel] at (0.4332,0.4475) {F4};
\node[celllabel] at (0.4435,0.4555) {F5};
\node[celllabel] at (0.4405,0.4619) {F6};
\node[celllabel] at (0.4065,0.4605) {S3};
\node[celllabel] at (0.4065,0.4415) {B2};
\node[celllabel] at (0.4535,0.4520) {R1};
\node[celllabel,anchor=south] at (0.43319048,0.4634) {$\pi_*$};
\end{axis}
\draw[gray!85,line width=0.6pt,->,shorten >=2pt] (zoomtop) to[out=55,in=160] (right.north west);
\end{tikzpicture}
\caption{The cover of Proposition~\ref{prop:cover}. (a) The region $\Omega$
(shaded) and the twelve cells that cover it. The region is bounded by the
graph of $F_{\mathrm{B1}}$, the line $F=F_*$ and the line $m=m_G$. (b) An
enlargement of the framed part of (a), showing the triangles
$\mathrm{F1},\dots,\mathrm{F6}$ with the common vertex $\pi_*$.}
\label{fig:cover}
\end{figure}

A \emph{cell} is a convex polygon $K$ that is the convex hull of a finite set
$\mathcal V(K)$ of vertices. Each vertex is either a point with dyadic
rational coordinates or the point $\pi_*$. Let $m_K$ and $F_K$ be the minima
of the first and of the second coordinates over $\mathcal V(K)$. Our cover of $\Omega$
consists of the twelve cells of Table~\ref{tab:cells}. Six of them are the
triangles $\mathrm{F1},\dots,\mathrm{F6}$, which have the common vertex
$\pi_*$. The other six are five rectangles and one quadrilateral. Let
$w_0,\dots,w_6$ be the vertices of $\mathrm{F1},\dots,\mathrm{F6}$ other than
$\pi_*$, in the order of Table~\ref{tab:cells}. Then
$\mathrm{F}j=\mathrm{conv}\{\pi_*,w_{j-1},w_j\}$ for $j=1,\dots,6$.
Figure~\ref{fig:cover}(a) shows the twelve cells and the region $\Omega$.
Figure~\ref{fig:cover}(b) enlarges the triangles $\mathrm{F1},\dots,\mathrm{F6}$
around $\pi_*$.

By Lemma~\ref{lem:weakdual}(d), a certificate that is positive at the
vertices of a cell is positive on the whole cell. The point $\pi_*$ is a
vertex of every cell that contains it, because no certificate is positive
at $\pi_*$, as $\nu_*\in\mathcal T(\pi_*)$ by Lemma~\ref{lem:apex}(a). At $\pi_*$, the certificates of
$\mathrm{F1},\dots,\mathrm{F6}$ are nonnegative and vanish exactly on
$\mathcal O(P_1)\cup\mathcal O(P_2)$. Lemma~\ref{lem:weakdual}(d) still
gives positivity at the other points of these triangles. The zero set at
$\pi_*$ gives the classification in Step~2 of
Section~\ref{s5:proofcritical}.

\begin{table}[ht]
\centering
\small
\caption{The cells of Proposition~\ref{prop:cover}. Each
decimal denotes the nearest binary64 number. The triangles
$\mathrm{F1},\dots,\mathrm{F6}$ are listed counterclockwise around $\pi_*$. A
rectangle $[m',m'']\times[F',F'']$ has the vertices $(m',F')$, $(m'',F')$,
$(m'',F'')$, $(m',F'')$.}
\label{tab:cells}
\begin{tabular}{lll}
\toprule
Cell & Type & Vertices\\
\midrule
F1 & triangle & $\pi_*$, $(0.4182,0.46278)$, $(0.4152,0.4608)$\\
F2 & triangle & $\pi_*$, $(0.4152,0.4608)$, $(0.4082,0.4528)$\\
F3 & triangle & $\pi_*$, $(0.4082,0.4528)$, $(0.4132,0.4413)$\\
F4 & triangle & $\pi_*$, $(0.4132,0.4413)$, $(0.4532,0.4413)$\\
F5 & triangle & $\pi_*$, $(0.4532,0.4413)$, $(0.4452,0.4598)$\\
F6 & triangle & $\pi_*$, $(0.4452,0.4598)$, $(0.4432,0.4629)$\\
S1 & rectangle & $[0.2185,0.33]\times[0.461,0.46278]$\\
S2 & rectangle & $[0.33,0.38]\times[0.4445,0.46278]$\\
S3 & rectangle & $[0.38,0.4182]\times[0.4445,0.46278]$\\
R1 & quadrilateral & \begin{tabular}[t]{@{}l@{}}$(0.4532,0.4413)$, $(0.51,0.4413)$,\\ $(0.51,0.4629)$, $(0.4432,0.4629)$\end{tabular}\\
B1 & rectangle & $[0.2185,0.33]\times[0.39,0.4615]$\\
B2 & rectangle & $[0.33,0.51]\times[0.335,0.445]$\\
\bottomrule
\end{tabular}
\end{table}

Let $\mathcal K$ be the family of the twelve cells of Table~\ref{tab:cells}.
Each cell $K\in\mathcal K$ carries a certificate $\Psi_{\theta_K,\lambda_K}$
with multipliers $\theta_K\in\R^{12}$ and a weight $\lambda_K\ge0$. It also
carries a number $Q_K>0$. Tables~\ref{tab:cellparams} to~\ref{tab:coef3} in
Appendix~\ref{app:cover} list $\lambda_K$, $Q_K$ and the entries of
$\theta_K$, with one exception. For each of the triangles
$K=\mathrm{F1},\dots,\mathrm{F6}$, the tables omit the entries
$\theta_{K,0},\theta_{K,2},\theta_{K,3},\theta_{K,4}$. For such $K$, the
tangency equations \eqref{s5:tangency} for
$\Psi_{\theta_K,\lambda_K}(\cdot;\pi_*)$ form a linear system in these four
entries. By Proposition~\ref{prop:cover}, it has exactly one solution, which
completes $\theta_K$. For every $K\in\mathcal K$, let
$\Psi_K=\Psi_{\theta_K,\lambda_K}$.

\begin{proposition}[The cover, computer assisted]\label{prop:cover}
For each of the triangles $K=\mathrm{F1},\dots,\mathrm{F6}$, the tangency
equations \eqref{s5:tangency} for $\Psi_{\theta_K,\lambda_K}(\cdot;\pi_*)$
have exactly one solution
$(\theta_{K,0},\theta_{K,2},\theta_{K,3},\theta_{K,4})$. Moreover, the
following hold.
\begin{enumerate}[label=(\roman*)]
\item $\Omega\subset\bigcup_{K\in\mathcal K}K$.
\item For every $K\in\mathcal K$, $m_K>0$, $F_K<1$ and
$((1-F_K)^2-m_K^2)/(2t_*)\le Q_K$.
\item For every $K\in\mathcal K$ and every $v\in\mathcal V(K)$ with
$v\ne\pi_*$, $\Psi_K(u;v)>0$ for all $u\in[e^{-Q_K},1]^3$.
\item The cells $K\in\mathcal K$ with $\pi_*\in\mathcal V(K)$ are the six
triangles $\mathrm{F1},\dots,\mathrm{F6}$. For each of them, $\Psi_K(u;\pi_*)\ge0$ for all
$u\in[e^{-Q_K},1]^3$, with equality if and only if
$u\in\mathcal O(P_1)\cup\mathcal O(P_2)$.
\end{enumerate}
\end{proposition}

We verify Proposition~\ref{prop:cover} in Appendix~\ref{app:cover}. There we
solve the tangency equations by validated linear algebra. Item~(i) follows
from lower bounds for $F_{\mathrm{B1}}$ in ball arithmetic and from exact
rational geometry, and we check item~(ii) in ball arithmetic. For
items~(iii) and~(iv) we use a branch and bound search over boxes, except
that near $\mathcal O(P_1)\cup\mathcal O(P_2)$ we treat item~(iv) by
interval Cholesky factorizations and Lemma~\ref{lem:taylor}.

Items~(iii) and~(iv) concern the cube $[e^{-Q_K},1]^3$. Let
$K\in\mathcal K$, and let $\pi=(m,F)\in K$ with $F\le1$. Since
$K=\mathrm{conv}(\mathcal V(K))$, we have $m\ge m_K$ and $F_K\le F\le1$.
Moreover, $m_K>0$ by item~(ii), so $(1-F)^2-m^2\le(1-F_K)^2-m_K^2$.
Hence \eqref{s5:params} and item~(ii) give
\begin{equation}\label{s5:ellK}
\ell(\pi)\ge\ell(m_K,F_K)\ge e^{-Q_K} .
\end{equation}
By \ref{s5:T5}, the cube therefore carries every $\nu\in\mathcal T(\pi)$ if
$\pi\in K\cap\Pi$. Now let $K$ be one of the triangles
$\mathrm{F1},\dots,\mathrm{F6}$. Then \eqref{s5:ellK} with $\pi=\pi_*$ gives
$e^{-Q_K}\le\ell_*$. Moreover, $\ell_*<a_*<b_*<1$ by
Lemma~\ref{lem:apex}(a) and Corollary~\ref{cor:Bstar}(a). Hence
\begin{equation}\label{s5:patterns}
Q_K\ge-\log\ell_*>q_a>q_b>0 .
\end{equation}
In particular, $P_1,P_2\in[e^{-Q_K},1]^3$, so the zero set in item~(iv) is
the whole of $\mathcal O(P_1)\cup\mathcal O(P_2)$.

\begin{corollary}\label{cor:cover}
For every $\pi\in\Omega\setminus\{\pi_*\}$ the set $\mathcal T(\pi)$ is empty.
\end{corollary}

\begin{proof}
Let $\pi\in\Omega\setminus\{\pi_*\}$. Then $\pi\in\Pi$ by
Lemma~\ref{lem:region}(c), and $\pi\in K$ for some $K\in\mathcal K$ by
Proposition~\ref{prop:cover}(i). Since $F\le F_*<1$ on $\Omega$ by
\eqref{s5:Omega} and Lemma~\ref{lem:region}(c), \eqref{s5:ellK} gives
$\ell(\pi)\ge e^{-Q_K}$. Write $\pi=\sum_{v\in\mathcal V(K)}\alpha_vv$ with
nonnegative weights summing to one. Let $\alpha_*=\alpha_{\pi_*}$ if
$\pi_*\in\mathcal V(K)$, and let $\alpha_*=0$ otherwise. Since $\pi\ne\pi_*$,
we have $\alpha_*<1$. By Proposition~\ref{prop:cover}(iii), the continuous
functions $\Psi_K(\cdot;v)$ with $v\in\mathcal V(K)\setminus\{\pi_*\}$ are
positive on the compact cube $[e^{-Q_K},1]^3$. Hence they are bounded below
on this cube by some $\varepsilon>0$. If $\alpha_*>0$, then $K$ is one of the
triangles $\mathrm{F1},\dots,\mathrm{F6}$, and $\Psi_K(\cdot;\pi_*)\ge0$ on
the cube by Proposition~\ref{prop:cover}(iv). Now let
$u\in[e^{-Q_K},1]^3$. Then the concavity of Lemma~\ref{lem:weakdual}(d) gives
\[
\Psi_K(u;\pi)\ge\alpha_*\Psi_K(u;\pi_*)+\sum_{v\in\mathcal V(K)\setminus\{\pi_*\}}\alpha_v\Psi_K(u;v)
\ge(1-\alpha_*)\varepsilon>0 .
\]
Hence Lemma~\ref{lem:weakdual}(b) with $\ell_0=e^{-Q_K}$ shows that
$\mathcal T(\pi)$ is empty.
\end{proof}

\subsection{Proof of Proposition~\ref{prop:critical}}\label{s5:proofcritical}

\begin{proof}[Proof of Proposition~\ref{prop:critical}]
Minimizers at $r_*$ exist by Lemma~\ref{lem:existence}; let $U$ be one of them. The
constant $r_*$ is feasible, so $F(U)=V(r_*)\le F_*$.

\noindent\textbf{Step 1: the value.} If $U$ is a.e.\ constant, then $U=r_*$
a.e.\ by Lemma~\ref{lem:meangap}, and $F(U)=h(r_*)$. Otherwise, let
$\pi=(e(U),F(U))$. By Lemmas~\ref{lem:region}(b) and~\ref{lem:valid},
$\pi\in\Omega$ and $\nu_U\in\mathcal T(\pi)$. Hence
Corollary~\ref{cor:cover} forces $\pi=\pi_*$, so $F(U)=F_*=h(r_*)$. Thus
every minimizer at $r_*$ costs $h(r_*)$, and $V(r_*)=h(r_*)$.

\noindent\textbf{Step 2: classification.} Let $U$ be a minimizer that is not
a.e.\ constant. By Step~1 and Lemma~\ref{lem:valid}, $e(U)=m_*$ and
$\nu_U\in\mathcal T(\pi_*)$. Let $K$ be one of the triangles
$\mathrm{F1},\dots,\mathrm{F6}$. By \eqref{s5:patterns}, $e^{-Q_K}\le\ell_*$,
and by Proposition~\ref{prop:cover}(iv), $\Psi_K(\cdot;\pi_*)$ is nonnegative
on $[e^{-Q_K},1]^3$ with zero set $\mathcal O(P_1)\cup\mathcal O(P_2)$. Hence
Lemma~\ref{lem:weakdual}(c) with $\ell_0=e^{-Q_K}$ gives
$\nu_U(\mathcal O(P_1)\cup\mathcal O(P_2))=1$. Since $e(U)=m_*=(a_*+b_*)/2$,
Lemma~\ref{lem:patterns}(c) shows that $U$ is a relabeling of $B_*$.
\end{proof}

Step~2 needs only one certificate tangent at $\pi_*$, and the member
$z=z_*$ of the family \eqref{eq:certpsi} of Section~\ref{sec:intermediate}
provides one in closed form. Indeed, expanding \eqref{eq:certpsi} with
$\mu_*t_*=1-F_*-m_*$ shows that this member equals $\Psi_{\theta,0}(\cdot;\pi_*)$ for
a vector $\theta$ whose only nonzero entries are
$\theta_0,\theta_1,\theta_2,\theta_3,\theta_7$. By
Proposition~\ref{prop:family}, it is nonnegative on $(0,1]^3$ with zero set
$\mathcal O(P_1)\cup\mathcal O(P_2)$. Hence it can replace
$\Psi_K(\cdot;\pi_*)$ in Step~2, with $\ell_0=\ell_*$ in
Lemma~\ref{lem:weakdual}(c).

\section{Local continuation of the minimizer below the threshold}\label{sec:rigidity}

We have proved that the constant graphon is the unique minimizer
for $r>r_*$. At $r=r_*$, the minimizers are the constant graphon and
the relabelings of $B_*$. We now turn to the case $0<r<r_*$.
In Sections~\ref{sec:rigidity}--\ref{sec:intermediate}, we prove that every minimizer of the variational
problem in this range is a relabeling of the balanced two block SBM
graphon $B_r$.

In this section we study minimizers near the graphons $B_z$ with $z\in(0,1)$
(Section~\ref{s4:sec:euler}). In Theorem~\ref{thm:rigidity} we show that a
minimizer close to some $B_z$ in cut distance, with multiplier close to
$\mu(z)$, is a two block SBM graphon. Together with
Corollary~\ref{cor:bipodalmin}, this theorem gives two facts that we use in
Section~\ref{sec:intermediate}, namely that the set of $r\in(0,r_*)$ at which
the minimizers are exactly the relabelings of $B_r$ is open
(Corollary~\ref{cor:open}) and contains an interval $(r_*-\delta,r_*)$
(Corollary~\ref{cor:seed}).

We prove Theorem~\ref{thm:rigidity} by controlling individual rows. We first show that almost every row of a
minimizer $U$ minimizes a functional $J_U$ built from $U$
(Lemma~\ref{lem:insertion}). If $U$ is close to $B=B_z$ and $\mu_U$ is close
to $\mu(z)$, then $J_U$ is uniformly close to the corresponding functional
$J_B$. The functional $J_B$ has exactly two minimizers, namely the two rows of $B$,
and both are nondegenerate (Theorem~\ref{thm:twowell}). Hence every row of a
nearby minimizer is close to
a row of $B$. The Euler--Lagrange equation then shows that the minimizer is uniformly
close to a two block SBM kernel. Finally, we compare the row equations of two
rows in the same block, and the nondegeneracy of the two minimizers of $J_B$
forces these rows to coincide. Kwan and Pham~\cite[Section~7.4]{KP26} use a similar argument for a
penalized problem with a large penalty. They upgrade closeness in cut distance
to uniform closeness by local optimality and then conclude by convexity.

For a minimizer $U$ at $r\in(0,1)$ and a measurable function
$v:[0,1]\to[0,1]$, define
\begin{equation}\label{s3:J}
J_U(v)=2\int_0^1h(v(y))\dd y+\mu_U\langle v,T_Uv\rangle.
\end{equation}
Recall that $\bar k_U=2F(U)+\mu_Ur^3$ (Section~\ref{s3:sub:vertex}). In the
proof of the next lemma we insert new vertices attached to $U$ according to
$v$ and then apply Lemma~\ref{s3:rescale}, as in the proof of
Proposition~\ref{prop:vertex}.

\begin{lemma}\label{lem:insertion}
Let $U$ be a minimizer at $r\in(0,1)$. Then $J_U(v)\ge\bar k_U$ for every
measurable $v:[0,1]\to[0,1]$. Moreover, for a.e.\ $x\in[0,1]$ the row
$u_x=U(x,\cdot)$ satisfies $J_U(u_x)=\bar k_U$ and
$\log u_x+\mu_UT_Uu_x=0$ a.e.\ on $[0,1]$.
\end{lemma}

\begin{proof}
Fix $v$. For $\varepsilon\in[0,\frac12]$ let $\alpha=1-\varepsilon$, and
define the graphon $W_\varepsilon$ as follows. For $x,y\in[0,\alpha]$ let
$W_\varepsilon(x,y)=U(x/\alpha,y/\alpha)$. For $x\in(\alpha,1]$ and
$y\in[0,\alpha]$ let $W_\varepsilon(x,y)=v(y/\alpha)$, and define
$W_\varepsilon$ symmetrically for $x\in[0,\alpha]$ and $y\in(\alpha,1]$.
On $(\alpha,1]^2$ let $W_\varepsilon=0$. Thus $W_\varepsilon$ compresses
$U$ to $[0,\alpha]$ and adds a set of new vertices of measure
$\varepsilon$, each attached to the old vertices according to $v$. We split pairs and triangles by the number of their vertices in
$(\alpha,1]$. Since a triangle with two vertices in $(\alpha,1]$ contains
an edge of value $0$, this gives
\[
F(W_\varepsilon)=\alpha^2F(U)+2\alpha\varepsilon\int_0^1h(v)
+\varepsilon^2,\quad
t(W_\varepsilon)=\alpha^3r^3+3\alpha^2\varepsilon\langle v,T_Uv\rangle,
\]
and $e(W_\varepsilon)=\alpha^2e(U)+2\alpha\varepsilon\int_0^1v$. We apply
Lemma~\ref{s3:rescale} with $I=[0,\frac12]$,
$F'_0=2\int_0^1h(v)-2F(U)$ and $t'_0=3\langle v,T_Uv\rangle-3r^3$, and obtain
$F'_0+\frac{\mu_U}3t'_0=J_U(v)-2F(U)-\mu_Ur^3\ge0$, that is,
$J_U(v)\ge\bar k_U$.

Each row $u_x$ is measurable with values in $[0,1]$. Moreover,
$\int_0^1h(u_x)=H_U(x)$ and $\langle u_x,T_Uu_x\rangle=\tau_U(x)$. Hence
\eqref{eq:vertex} gives $J_U(u_x)=2H_U(x)+\mu_U\tau_U(x)=\bar k_U$ for
a.e.\ $x$. Finally, $(T_Uu_x)(y)=C_U(x,y)$. By \eqref{eq:euler} and
Fubini's theorem, for a.e.\ $x$ we have $\log U(x,y)+\mu_UC_U(x,y)=0$ for
a.e.\ $y$.
\end{proof}

\subsection{Two nondegenerate minimizers of \texorpdfstring{$J_B$}{J\_B}}

Throughout this subsection, $0<a<b<1$ and $\mu>0$ satisfy the system
\eqref{s4:eq:EL}, and $B=B_{a,b}$. Hence $\log B=-\mu C_B$. Let
$A_1=[0,\frac12)$ and $A_2=[\frac12,1]$, and let $f_1=a\1_{A_1}+b\1_{A_2}$ and
$f_2=b\1_{A_1}+a\1_{A_2}$ be the two rows of $B$. For measurable
$v:[0,1]\to[0,1]$, let $J_B(v)=2\int_0^1h(v)+\mu\langle v,T_Bv\rangle$.
This is \eqref{s3:J} with $U=B$ and with $\mu$ in place of $\mu_U$.
Finally, let
\[
M_B=\begin{pmatrix}1/a+\mu a/2&\mu b/2\\ \mu b/2&1/b+\mu a/2\end{pmatrix}.
\]

\begin{theorem}\label{thm:twowell}
Let $0<a<b<1$ and $\mu>0$ satisfy \eqref{s4:eq:EL}. Then the following
hold.
\begin{enumerate}[label=(\alph*)]
\item The matrix $M_B$ is positive definite.
\item We have $J_B(f_1)=J_B(f_2)=\min J_B$, where the minimum is over all
measurable $v:[0,1]\to[0,1]$. Every minimizer of $J_B$ equals $f_1$ or $f_2$
a.e.
\item There is $c_B>0$ such that for every measurable $v:[0,1]\to[0,1]$,
\[
J_B(v)-J_B(f_1)\ge c_B\min_{i=1,2}\|v-f_i\|_2^2 .
\]
\end{enumerate}
\end{theorem}

In the proof we use Jensen's inequality to reduce $J_B$ to a function $q$ of the
means of $v$ on $A_1$ and $A_2$. The Hessian of $q$ at the minimum point
$(a,b)$ is $M_B$.

\begin{proof}
\noindent\textbf{Step 1: reduction to two variables.} Let $v:[0,1]\to[0,1]$
be measurable, let $u$ and $w$ be the means of $v$ on $A_1$ and $A_2$, and
let $\bar v=u\1_{A_1}+w\1_{A_2}$. Since $h''(s)=1/s\ge1$ on $(0,1]$, the
function $h(s)-s^2/2$ is convex and continuous on $[0,1]$. Jensen's
inequality for this function on $A_1$ and on $A_2$ gives
\[
\int_0^1h(v)\ge\frac{h(u)+h(w)}2+\frac12\Big(\int_0^1v^2-\int_0^1\bar v^2\Big)
=\frac{h(u)+h(w)}2+\frac12\|v-\bar v\|_2^2 .
\]
Since $B$ is constant on each product $A_i\times A_j$, we have
$\langle v,T_Bv\rangle=\langle\bar v,T_B\bar v\rangle=(a(u^2+w^2)+2buw)/4$.
Let $q(u,w)=h(u)+h(w)+\frac\mu4(a(u^2+w^2)+2buw)$. The Jensen bound and the
formula for $\langle v,T_Bv\rangle$ give
\begin{equation}\label{s6:jensen}
J_B(v)\ge q(u,w)+\|v-\bar v\|_2^2,
\end{equation}
with equality if $v=\bar v$. In particular, $J_B(f_1)=q(a,b)=q(b,a)=J_B(f_2)$.

\noindent\textbf{Step 2: the minimizers of $q$.} Let $\alpha=\mu a/2$ and
$\beta=\mu b/2$. Then $\partial_uq(u,w)=\log u+\alpha u+\beta w$. The
continuous function $q$ attains its minimum on $[0,1]^2$, and every minimum
point is interior. Indeed, $\partial_uq(u,w)\to-\infty$ as $u\downarrow0$ and
$\partial_uq(1,w)=\alpha+\beta w>0$. The derivative $\partial_wq$ behaves in
the same way. Hence every minimum point solves
\begin{equation}\label{s6:crit}
\log u=-\alpha u-\beta w,\quad \log w=-\beta u-\alpha w .
\end{equation}
By \eqref{s4:eq:EL}, $(a,b)$ and $(b,a)$ solve \eqref{s6:crit}. Let
$c=\beta-\alpha>0$.

A solution with $u=w=s$ satisfies $\log s+(\alpha+\beta)s=0$. The left side
increases strictly from $-\infty$ to $\infty$ on $(0,\infty)$, so it has a
unique zero $s_0$. Hence $(s_0,s_0)$ is the only solution of \eqref{s6:crit}
with $u=w$.

Consider a solution with $u<w$. Let $s=(u+w)/2$ and
$y=\operatorname{artanh}((w-u)/(w+u))>0$. Then $u=s(1-\tanh y)$,
$w=s(1+\tanh y)$, $\log(w/u)=2y$ and $uw=s^2/\cosh^2y$. Subtracting the two
equations of \eqref{s6:crit} gives $2y=c(w-u)=2cs\tanh y$, that is,
$cs=y\coth y$. Adding them gives $2\log s-2\log\cosh y=-2(\alpha+\beta)s$.
Substituting $s=y\coth y/c$ shows that $y$ is a root of
\[
G(y)=\log\frac{y}{\sinh y}-\log c+\frac{\alpha+\beta}cy\coth y .
\]
The root $y$ determines $(u,w)$, since $s=y\coth y/c$. For $y>0$,
\[
G'(y)=\Big(\frac1y-\frac{y}{\sinh^2y}\Big)
+\frac{2\alpha}c\Big(\coth y-\frac{y}{\sinh^2y}\Big)>0,
\]
because $\sinh y>y$ and $\sinh y\cosh y>y$. Hence $G$ has at most one root.
The solution $(a,b)$ gives the root $y_0=\frac12\log(b/a)$. Thus $(a,b)$ and
$(b,a)$ are the only solutions of \eqref{s6:crit} with $u\ne w$.

The point $(s_0,s_0)$ is not a local minimum. Indeed, $G$ is strictly
increasing with $G(y_0)=0$, and $G(y)\to(\alpha+\beta)/c-\log c$ as
$y\downarrow0$. Hence $(\alpha+\beta)/c-\log c<0$. This number is the value at
$s=1/c$ of the increasing function $\log s+(\alpha+\beta)s$, so $s_0>1/c$. The
Hessian of $q$ at $(s_0,s_0)$ has diagonal entries $1/s_0+\alpha$ and off
diagonal entries $\beta$. Its quadratic form in the direction $(1,-1)$
therefore equals $2(1/s_0-c)<0$. We conclude that the minimum points of $q$ on
$[0,1]^2$ are exactly $(a,b)$ and $(b,a)$.

\noindent\textbf{Step 3: nondegeneracy.} The Hessian of $q$ at $(a,b)$ is
$M_B$. For the solution $(a,b)$, Step 2 gives the root $y_0=\frac12\log(b/a)$
of $G$, $c(a+b)/2=y_0\coth y_0$ and $ab=((a+b)/2)^2/\cosh^2y_0$. Hence
$c^2ab=y_0^2/\sinh^2y_0$. Since
$\beta^2-\alpha^2=c(\alpha+\beta)=c^2+2\alpha c$,
\begin{align*}
ab\det M_B&=(1+\alpha a)(1+\alpha b)-\beta^2ab
=1+\alpha(a+b)-(\beta^2-\alpha^2)ab\\
&=\Big(1-\frac{y_0^2}{\sinh^2y_0}\Big)
+\frac{2\alpha}c\Big(y_0\coth y_0-\frac{y_0^2}{\sinh^2y_0}\Big)>0,
\end{align*}
because $y_0\coth y_0>1>y_0^2/\sinh^2y_0$. Since $\det M_B>0$ and the
diagonal entries of $M_B$ are positive, $M_B$ is positive definite, which is
(a). Since $q(u,w)=q(w,u)$, the Hessian of $q$ at $(b,a)$ is $M_B$ with both
coordinates swapped, which is also positive definite.

\noindent\textbf{Step 4: quadratic growth.} Let $\omega_1=(a,b)$,
$\omega_2=(b,a)$ and $d(\omega)=\min_i|\omega-\omega_i|$ for
$\omega\in[0,1]^2$, where $|\cdot|$ is the Euclidean norm. Recall that
$\nabla q(\omega_i)=0$ and $q(\omega_i)=q(\omega_1)$. By Step 3, the
continuity of the Hessian of $q$ on $(0,1)^2$ and Taylor's formula at
$\omega_1$ and $\omega_2$, there are $\sigma,\lambda>0$ such that
$q(\omega)-q(\omega_1)\ge\lambda d(\omega)^2$ whenever $d(\omega)\le\sigma$. By
Step 2, the continuous function $q-q(\omega_1)$ is positive on the compact set
$\{\omega\in[0,1]^2:d(\omega)\ge\sigma\}$, hence at least some $\gamma>0$
there. Since $d^2\le2$, we obtain $q(\omega)-q(\omega_1)\ge\theta d(\omega)^2$
on $[0,1]^2$ with $\theta=\min\{\lambda,\gamma/2\}$.

Now let $v$, $u$, $w$ and $\bar v$ be as in Step 1. The function $v-\bar v$
has mean zero on $A_1$ and on $A_2$, and $\bar v-f_i$ is constant on both
sets. Hence $v-\bar v$ is orthogonal to $\bar v-f_i$. Since
$\|\bar v-f_i\|_2^2=\frac12|(u,w)-\omega_i|^2$, we get
$\min_i\|v-f_i\|_2^2=\|v-\bar v\|_2^2+\frac12d(u,w)^2$. Let
$c_B=\min\{1,2\theta\}$. By \eqref{s6:jensen}, $J_B(f_1)=q(\omega_1)$ and the
growth of $q$,
\[
J_B(v)-J_B(f_1)\ge\|v-\bar v\|_2^2+\theta d(u,w)^2
\ge c_B\min_{i=1,2}\|v-f_i\|_2^2 .
\]
This is (c). Part (b) follows from (c) and $J_B(f_1)=J_B(f_2)$.
\end{proof}

\subsection{Minimizers near \texorpdfstring{$B_z$}{B\_z} are two block SBM graphons}

Two block SBM graphons are defined in Section~\ref{sec:intro} and again in
Section~\ref{s4:sec:bipodal}.

\begin{theorem}\label{thm:rigidity}
Let $0<a<b<1$ and $\mu>0$ satisfy \eqref{s4:eq:EL}, and let $B=B_{a,b}$.
For each $n$ let $r_n\in(0,1)$ and let $U_n$ be a minimizer at $r_n$.
Suppose that $\delta_\cut(U_n,B)\to0$ and $\mu_{U_n}\to\mu$. Then $U_n$ is a
two block SBM graphon for all sufficiently large $n$.
\end{theorem}

\begin{proof}
By the definition of the cut distance, there are measure preserving
bijections $\varphi_n$ with $\|U_n^{\varphi_n}-B\|_\cut\to0$. Since $e$, $t$
and $F$ are invariant under relabeling, $U_n^{\varphi_n}$ is a minimizer at
$r_n$ with multiplier $\mu_{U_n}$. Since $\varphi_n^{-1}$ is measure
preserving, $U_n$ is a two block SBM graphon if and only if $U_n^{\varphi_n}$
is one. We may therefore assume that $\|U_n-B\|_\cut\to0$.
Let $\mu_n=\mu_{U_n}$, $\delta_n=\|U_n-B\|_\cut$ and $u_{n,x}=U_n(x,\cdot)$.
By Lemma~\ref{lem:insertion}, there is a measurable set $X_n\subset[0,1]$ of
full measure such that for every $x\in X_n$,
\begin{equation}\label{s6:rowprops}
J_{U_n}(u_{n,x})=\min J_{U_n},\quad
\log u_{n,x}+\mu_nT_{U_n}u_{n,x}=0\quad\text{a.e.},
\end{equation}
where $(T_{U_n}u_{n,x})(y)=C_{U_n}(x,y)$.

\noindent\textbf{Step 1: the functionals $J_{U_n}$ and $J_B$ are uniformly
close.} For measurable $v:[0,1]\to[0,1]$,
\[
J_{U_n}(v)-J_B(v)=\mu_n\langle v,T_{U_n-B}v\rangle+(\mu_n-\mu)\langle v,T_Bv\rangle .
\]
By Lemma~\ref{lem:cut}(a), $|\langle v,T_{U_n-B}v\rangle|\le\delta_n$.
Moreover, $0\le\langle v,T_Bv\rangle\le1$, since $v$ and $B$ take values in
$[0,1]$. Hence
\begin{equation}\label{s6:beta}
\sup_v|J_{U_n}(v)-J_B(v)|\le\beta_n:=\mu_n\delta_n+|\mu_n-\mu|\longrightarrow0 .
\end{equation}

\noindent\textbf{Step 2: every row is close to a row of $B$.} For $x\in X_n$,
\[
J_B(u_{n,x})\overset{\mathrm{(a)}}{\le}J_{U_n}(u_{n,x})+\beta_n
\overset{\mathrm{(b)}}{\le}J_{U_n}(f_1)+\beta_n
\overset{\mathrm{(a)}}{\le}J_B(f_1)+2\beta_n,
\]
where $\mathrm{(a)}$ follows from \eqref{s6:beta} and $\mathrm{(b)}$ from the
minimality in \eqref{s6:rowprops}. Theorem~\ref{thm:twowell}(c) gives
\begin{equation}\label{s6:rows}
\min_{i=1,2}\|u_{n,x}-f_i\|_2\le\varepsilon_n:=(2\beta_n/c_B)^{1/2}
\quad\text{for every }x\in X_n .
\end{equation}

\noindent\textbf{Step 3: uniform approximation by a two block SBM kernel.}
Note that $\|f_1-f_2\|_2=b-a$. For large $n$ we have $2\varepsilon_n<b-a$, so
the closed $L^2$ balls of radius $\varepsilon_n$ about $f_1$ and $f_2$ are
disjoint. Let
\[
P_{n,1}=\{x\in X_n:\|u_{n,x}-f_1\|_2\le\varepsilon_n\},\quad
P_{n,2}=[0,1]\setminus P_{n,1}.
\]
We call $P_{n,1}$ and $P_{n,2}$ the blocks. They are measurable by Fubini's
theorem. By \eqref{s6:rows}, $\|u_{n,x}-f_i\|_2\le\varepsilon_n$ for
$x\in P_{n,i}\cap X_n$. Let $\tilde B_n$ be the kernel equal to $a$ on
$P_{n,1}^2\cup P_{n,2}^2$ and to $b$ elsewhere. Let
$x\in P_{n,i}\cap X_n$ and $y\in P_{n,j}\cap X_n$. Since
$C_{U_n}(x,y)=\langle u_{n,x},u_{n,y}\rangle$ and all functions involved
take values in $[0,1]$,
\[
|C_{U_n}(x,y)-\langle f_i,f_j\rangle|
\le\|u_{n,x}-f_i\|_2\|u_{n,y}\|_2+\|f_i\|_2\|u_{n,y}-f_j\|_2\le2\varepsilon_n .
\]
Since $\langle f_i,f_i\rangle=(a^2+b^2)/2$ and $\langle f_1,f_2\rangle=ab$,
\eqref{s4:eq:EL} gives $\tilde B_n(x,y)=\exp(-\mu\langle f_i,f_j\rangle)$.
Moreover, $U_n=\exp(-\mu_nC_{U_n})$ a.e.\ by Proposition~\ref{prop:euler}.
Since $s\mapsto e^{-s}$ is $1$-Lipschitz on $[0,\infty)$ and
$0\le\langle f_i,f_j\rangle\le1$, we obtain
\begin{equation}\label{s6:unif}
|U_n-\tilde B_n|\le\zeta_n:=2\mu_n\varepsilon_n+|\mu_n-\mu|\quad\text{a.e.}
\end{equation}
Note that $\zeta_n\to0$, since $\mu_n\to\mu$ and $\varepsilon_n\to0$. Let
$p_n=|P_{n,1}|$. Then
\[
e(\tilde B_n)=a\big(p_n^2+(1-p_n)^2\big)+2bp_n(1-p_n)
=\frac{a+b}2-2(b-a)\Big(p_n-\frac12\Big)^2 .
\]
By Lemma~\ref{lem:cut}(b), $|e(U_n)-(a+b)/2|\le\delta_n$. By
\eqref{s6:unif}, $|e(U_n)-e(\tilde B_n)|\le\zeta_n$. Hence
$2(b-a)(p_n-\frac12)^2\le\delta_n+\zeta_n\to0$, so $p_n\to\frac12$. In
particular, both blocks have positive measure for large $n$.

\noindent\textbf{Step 4: uniform coercivity at the rows of $\tilde B_n$.} Let
$g_{n,1}=a\1_{P_{n,1}}+b\1_{P_{n,2}}$ and $g_{n,2}=b\1_{P_{n,1}}+a\1_{P_{n,2}}$
be the rows of $\tilde B_n$. For $w\in L^2([0,1])$, let
\[
Q_{n,i}(w)=\int_0^1\frac{w^2}{g_{n,i}}+\mu\langle w,T_{\tilde B_n}w\rangle .
\]
Let $e_j=|P_{n,j}|^{-1/2}\1_{P_{n,j}}$, $\alpha_j=\langle w,e_j\rangle$,
$\alpha=(\alpha_1,\alpha_2)$ and $w^\perp=w-\alpha_1e_1-\alpha_2e_2$. Then
$w^\perp$ has mean zero on each block. The function $g_{n,i}$ is constant on
each block, and $\tilde B_n$ is constant on each product of blocks. Hence the
mixed terms vanish and $T_{\tilde B_n}w^\perp=0$. Using also $g_{n,i}\le1$, we
obtain
\[
Q_{n,i}(w)=\alpha^{\top}M_{n,i}\alpha+\int_0^1\frac{(w^\perp)^2}{g_{n,i}}
\ge\lambda_{\min}(M_{n,i})|\alpha|^2+\|w^\perp\|_2^2,
\]
where $M_{n,1}=\operatorname{diag}(1/a,1/b)+\mu K_n$,
$M_{n,2}=\operatorname{diag}(1/b,1/a)+\mu K_n$ and
\[
K_n=\begin{pmatrix}ap_n&b\sqrt{p_n(1-p_n)}\\ b\sqrt{p_n(1-p_n)}&a(1-p_n)\end{pmatrix}.
\]
As $p_n\to\frac12$, we have $M_{n,1}\to M_B$, while $M_{n,2}$ converges to
$M_B$ with both coordinates swapped. By Theorem~\ref{thm:twowell}(a), both
limits have smallest eigenvalue $\lambda_B>0$. Since
$\|w\|_2^2=|\alpha|^2+\|w^\perp\|_2^2$, we obtain, for large $n$ and both $i$,
\begin{equation}\label{s6:Qcoercive}
Q_{n,i}(w)\ge\eta\|w\|_2^2\quad(w\in L^2([0,1])),\quad
\eta=\min\{\lambda_B/2,1\}.
\end{equation}

\noindent\textbf{Step 5: two rows in the same block coincide.} By
\eqref{s6:unif} and Fubini's theorem, there is a set $X_n'\subset X_n$ of
full measure such that $|u_{n,x}-g_{n,i}|\le\zeta_n$ a.e.\ for every
$x\in P_{n,i}\cap X_n'$. Fix $n_0$ such that for every $n\ge n_0$, both
blocks have positive measure, \eqref{s6:Qcoercive} holds and
$\zeta_n/a^2+\mu_n\zeta_n+|\mu_n-\mu|\le\eta/2$.

Let $n\ge n_0$, let
$x,x'\in P_{n,i}\cap X_n'$, and let $v_0=u_{n,x}$, $v_1=u_{n,x'}$ and
$w=v_1-v_0$. By \eqref{s6:rowprops}, $\log v_0$ and $\log v_1$ are finite
a.e. Thus, by the choice of $X_n'$, both $v_0$ and $v_1$ take values in
$(0,g_{n,i}+\zeta_n]$ a.e. By the mean value theorem, a.e.\ on $[0,1]$ we
have $w(\log v_1-\log v_0)=w^2/\xi$ for some $\xi$ between $v_0$ and $v_1$.
Hence $w(\log v_1-\log v_0)\ge w^2/(g_{n,i}+\zeta_n)$ a.e. Recall from
Section~\ref{s2:sub:notation} that $\|T_W\|_{\mathrm{op}}\le\|W\|_2$ for a
kernel $W$. Therefore \eqref{s6:unif} and $0\le\tilde B_n\le1$ give
$\|\mu_nT_{U_n}-\mu T_{\tilde B_n}\|_{\mathrm{op}}\le\mu_n\zeta_n+|\mu_n-\mu|$.
Subtracting the two row equations in \eqref{s6:rowprops} gives
$\log v_1-\log v_0+\mu_nT_{U_n}w=0$ a.e. Pairing with the bounded function
$w$, we obtain
\begin{align*}
0&=\int_0^1w(\log v_1-\log v_0)+\mu_n\langle w,T_{U_n}w\rangle
\overset{\mathrm{(a)}}{\ge}\int_0^1\frac{w^2}{g_{n,i}+\zeta_n}+\mu_n\langle w,T_{U_n}w\rangle\\
&\overset{\mathrm{(b)}}{\ge}Q_{n,i}(w)-\Big(\frac{\zeta_n}{a^2}+\mu_n\zeta_n+|\mu_n-\mu|\Big)\|w\|_2^2
\overset{\mathrm{(c)}}{\ge}\frac\eta2\|w\|_2^2,
\end{align*}
where $\mathrm{(a)}$ is the pointwise bound on $w(\log v_1-\log v_0)$;
$\mathrm{(b)}$ uses the bound on $\|\mu_nT_{U_n}-\mu T_{\tilde B_n}\|_{\mathrm{op}}$
and $1/(g_{n,i}+\zeta_n)\ge1/g_{n,i}-\zeta_n/a^2$, which holds since
$g_{n,i}\ge a$; and $\mathrm{(c)}$ uses \eqref{s6:Qcoercive} and the choice
of $n_0$. Hence $u_{n,x}=u_{n,x'}$ a.e.

Since $|P_{n,i}|>0$, we can fix $x_i\in P_{n,i}\cap X_n'$ for $i=1,2$. By
Fubini's theorem, $U_n(x,y)=u_{n,x_i}(y)$ for a.e.\
$(x,y)\in P_{n,i}\times[0,1]$. Hence Lemma~\ref{lem:rowsbipodal} shows that
$U_n$ is a two block SBM graphon for every $n\ge n_0$.
\end{proof}

\subsection{Uniqueness of the minimizer at nearby targets}

We now apply Theorem~\ref{thm:rigidity} to sequences of minimizers at targets
$r_n\to r_0$. Corollary~\ref{cor:open} treats an interior point $r_0<r_*$,
and Corollary~\ref{cor:seed} treats $r_0=r_*$. We prove both by a single
compactness argument.

\begin{corollary}\label{cor:open}
Let $r_0\in(0,r_*)$ be such that the minimizers at $r_0$ are exactly the
relabelings of $B_{r_0}$. Then there is $\delta>0$ with
$(r_0-\delta,r_0+\delta)\subset(0,r_*)$ such that for every
$r\in(r_0-\delta,r_0+\delta)$, the minimizers at $r$ are exactly the
relabelings of $B_r$.
\end{corollary}

\begin{corollary}\label{cor:seed}
There is $\delta>0$ such that for every $r\in(r_*-\delta,r_*)$, the minimizers
at $r$ are exactly the relabelings of $B_r$.
\end{corollary}

\begin{proof}[Proof of Corollaries~\ref{cor:open} and~\ref{cor:seed}]
For Corollary~\ref{cor:open}, let $r_0$ be as in its statement, and for
Corollary~\ref{cor:seed}, let $r_0=r_*$. Recall that minimizers exist
(Lemma~\ref{lem:existence}) and that relabeling preserves $F$ and $t$. Hence,
if every minimizer at $r$ is a relabeling of $B_r$, then the minimizers at $r$
are exactly the relabelings of $B_r$. Thus, if the conclusion fails, there are
$r_n\in(0,r_*)$ with $r_n\to r_0$ and minimizers $U_n$ at $r_n$ that are not
relabelings of $B_{r_n}$. By Lemma~\ref{lem:cut}(c), after passing to a subsequence, there is a graphon $W$
with $\delta_\cut(U_n,W)\to0$. By Lemma~\ref{lem:limits}, $W$ is a minimizer
at $r_0$ and $\mu_{U_n}\to\mu_W$.

We claim that $W$ is a relabeling of $B_{r_0}$. If $r_0<r_*$, this follows
from the hypothesis of Corollary~\ref{cor:open}, since $W$ is a minimizer at
$r_0$. If $r_0=r_*$, this is Lemma~\ref{lem:leftlimits}.

By Lemma~\ref{lem:cut}(d), $\delta_\cut(W,B_{r_0})=0$. Hence
$\delta_\cut(U_n,B_{r_0})\to0$ by the triangle inequality. Since $W$ has the
cost and the edge density of $B_{r_0}$, \eqref{eq:mult} and
Lemma~\ref{lem:Bderivative}(a) give $\mu_W=\mu_{r_0}$. By
Lemma~\ref{lem:branch}(a), $(a_{r_0},b_{r_0},\mu_{r_0})$ satisfies
\eqref{s4:eq:EL} with $0<a_{r_0}<b_{r_0}<1$ and $\mu_{r_0}>0$. Hence
Theorem~\ref{thm:rigidity} with $B=B_{r_0}$ shows that $U_n$ is a two block SBM
graphon for all large $n$. Since $r_n<r_*$, Corollary~\ref{cor:bipodalmin}
shows that $U_n$ is a relabeling of $B_{r_n}$. This contradicts the choice of
$U_n$.
\end{proof}
\section{Minimizer for small targets}\label{sec:smalltail}

In this section we classify the minimizers at all targets up to the constant $r_-=0.03965\ldots$ defined in
Section~\ref{s7:pen}.

\begin{theorem}\label{thm:smalltail}
For every $r\in(0,r_-]$, the minimizers at $r$ are exactly the relabelings of $B_r$.
\end{theorem}

Here $B_r$ is defined because $r_-<r_*$ (Lemma~\ref{lem:order}). We outline the proof. Let $U$ be a
minimizer at $r\le r_-$. By \cite[Theorem~1.6]{KP26} and Section~\ref{s7:pen}, $U$ minimizes
$F+\frac{\mu_r}3t$ over all graphons and has the edge density, cost, triangle density and multiplier of
$B_r$ (Proposition~\ref{prop:smalldata}). Hence $T_U$ is negative semidefinite on functions of mean zero
(Lemma~\ref{lem:penalized}). Since $t(U)=\tr(T_U^3)$ is small, we can use a spectral argument to split $[0,1]$ into
two blocks. For each block $P$, the graphon $U$ is exponentially small in $\mu_r$ on $P\times P$.
Comparing the Euler--Lagrange equations for two rows in the same block, we show that these rows coincide.
Hence $U$ is a two block SBM graphon (Section~\ref{s7:bip}), and Corollary~\ref{cor:bipodalmin} applies.

\subsection{The penalized problem over the graphons \texorpdfstring{$B_{a,b}$}{B\_\{a,b\}}}\label{s7:pen}

For $\gamma\ge0$ and $a,b\in[0,1]$ define
\[
\Phi_\gamma(a,b)=F(B_{a,b})+\gamma t(B_{a,b})=\frac{h(a)+h(b)}2+\gamma\frac{a^3+3ab^2}4,\quad
C(\gamma)=\min_{[0,1]^2}\Phi_\gamma .
\]
For $c>0$, Kwan and Pham study the stationary points of $\Phi_{c^2/3}$ in the discussion after
\cite[Theorem~1.6]{KP26}. They explain qualitatively that nonconstant stationary points appear at
$c\approx3.350$, that one of them is a local minimum, and that it becomes the global minimizer beyond
$c_*\approx4.341$. We prove
the quantitative form of this picture that we need later.

We use $a(z)$, $b(z)$, $\mu(z)$, $\rho(z)$, $B_z$ and $z_1$ from Section~\ref{sec:balanced}, with
$L=\log(1/z)$. The nonconstant stationary points of $\Phi_\gamma$ are the pairs $(a(z),b(z))$ with
$\mu(z)=3\gamma$ (Lemma~\ref{lem:stationary}). For $z\in(0,1)$ let $Q(z)=2z(1+z)L(2L-1+z)/(1-z)^3$,
which is positive since $L\ge1-z$. We next show that $Q$ governs the monotonicity of $\mu$ and
the determinant of the Hessian matrix of $\Phi_{\mu(z)/3}$ at $(a(z),b(z))$.

\begin{lemma}\label{lem:fold}
The following hold.
\begin{enumerate}[label=(\alph*)]
\item $Q$ is strictly increasing on $(0,1)$, and $Q(10^{-3})<1<Q(1/50)$. Hence there is a unique
$z_m\in(10^{-3},1/50)$ with $Q(z_m)=1$.
\item $(\log\mu)'(z)=(Q(z)-1)/(zL)$ for $z\in(0,1)$. Hence $\mu$ is strictly decreasing on $(0,z_m]$
and strictly increasing on $[z_m,1)$. Moreover, $\mu(z)\to\infty$ as $z\to0$.
\item For $z\in(0,1)$ and $\gamma=\mu(z)/3$, the Hessian matrix of $\Phi_\gamma$ at $(a(z),b(z))$ has
determinant $(1-Q(z))/(4zb(z)^2)$.
\end{enumerate}
\end{lemma}

\begin{proof}
(a) Using $L'=-1/z$, we compute $z(1-z^2)L(2L-1+z)(\log Q)'(z)=P(z,L)$ with
\[
P(z,L)=2(1+4z+z^2)L^2-(1-z)(5+7z)L+(1-z)^2(1+z).
\]
As a quadratic in $L$, $P(z,\cdot)$ has positive leading coefficient and discriminant $(1-z)^2S(z)$,
where $S(z)=17+30z+9z^2-8z^3>0$. Hence $P(z,L)>0$ for $L>L_+(z)=(1-z)(5+7z+S(z)^{1/2})/(4(1+4z+z^2))$.
The identities $5(1+4z+z^2)-(1-z)(5+7z)=18z+12z^2$ and
$17(1+4z+z^2)^2-(1-z)^2S(z)=4z(35+85z+33z^2-2z^3(1-z))$ give $L_+(z)\le(5+17^{1/2})/4<2.29<\log10$.
This settles the case $z\le1/10$. For $z\in[1/10,1)$ we have $\log(1/z)>2(1-z)/(1+z)$, since the
difference vanishes at $z=1$ and has derivative $-(1-z)^2/(z(1+z)^2)<0$. Moreover,
$4(1+z)(1+4z+z^2)(2(1-z)/(1+z)-L_+(z))=(1-z)(3+20z+z^2-(1+z)S(z)^{1/2})$ and
$(3+20z+z^2)^2-(1+z)^2S(z)=8(1+4z+z^2)(z^3-3z^2+11z-1)$. Here the cubic has positive derivative
$3z^2-6z+11$ and value $0.071$ at $z=1/10$. Hence $L>L_+(z)$ on $(0,1)$, and $(\log Q)'>0$. Finally,
$\log1000<6.908$ and $\log50>3.912$ give $Q(10^{-3})<2\cdot10^{-3}\cdot1.001\cdot6.908\cdot12.82/0.997<0.18$
and $Q(1/50)>0.04\cdot1.02\cdot3.912\cdot6.844>1.09$.

(b) Since $\log b(z)=-2zL/(1-z)^2$, we have $(\log b)'=(2(1-z)-2(1+z)L)/(1-z)^3$. Then
$\log\mu=\log2+\log L-2\log(1-z)-2\log b$ gives
\[
zL(\log\mu)'(z)=-1+\frac{2zL}{1-z}+zL\frac{4(1+z)L-4(1-z)}{(1-z)^3}=Q(z)-1,
\]
since $2(1-z)^2-4(1-z)=-2(1-z^2)$. The monotonicity follows from (a). Since $b(z)<1$, we have
$\mu(z)\ge2L\to\infty$ as $z\to0$.

(c) The Hessian matrix of $\Phi_\gamma$ at $(a,b)\in(0,1)^2$ has diagonal entries $\frac12(1/a+3\gamma a)$
and $\frac12(1/b+3\gamma a)$, and off diagonal entries $\frac32\gamma b$. Let $a=zb$, $3\gamma=\mu(z)$ and
$\omega=\mu(z)b^2=2L/(1-z)^2$. Then $4b^2$ times the determinant equals
$(1/z+\omega z)(1+\omega z)-\omega^2=1/z+\omega(1+z)(1-(1-z)\omega)$. Moreover,
$\omega(1+z)(1-(1-z)\omega)=2L(1+z)(1-z-2L)/(1-z)^3=-Q(z)/z$.
\end{proof}

Let $\gamma>0$. The derivative $\log s+3\gamma s^2$ of the strictly convex function $h(s)+\gamma s^3$
increases strictly from $-\infty$ to $3\gamma>0$ on $(0,1]$. Hence it has exactly one zero
$s_\gamma\in(0,1)$, which is the only minimizing point of $h(s)+\gamma s^3$ on $[0,1]$. Define
$C_0(\gamma)=h(s_\gamma)+\gamma s_\gamma^3=\Phi_\gamma(s_\gamma,s_\gamma)$. Since
$\partial_s(\log s+3\gamma s^2)=1/s+6\gamma s>0$, the implicit function theorem shows that
$\gamma\mapsto s_\gamma$ is continuously differentiable on $(0,\infty)$. Hence
$C_0'(\gamma)=(\log s_\gamma+3\gamma s_\gamma^2)s_\gamma'+s_\gamma^3=s_\gamma^3$. Since $\log s\ge1-1/s$,
we have $3\gamma s_\gamma^3=-s_\gamma\log s_\gamma\le1$, so $s_\gamma\to0$ as $\gamma\to\infty$.

\begin{lemma}\label{lem:stationary}
Let $\gamma>0$.
\begin{enumerate}[label=(\alph*)]
\item Every minimizer of $\Phi_\gamma$ on $[0,1]^2$ lies in $(0,1)^2$ and is a stationary point of
$\Phi_\gamma$ with positive semidefinite Hessian matrix.
\item The stationary points of $\Phi_\gamma$ in $(0,1)^2$ are $(s_\gamma,s_\gamma)$ and the pairs
$(a(z),b(z))$ with $z\in(0,1)$ and $\mu(z)=3\gamma$.
\item Every stationary point $(a,b)\in(0,1)^2$ of $\Phi_\gamma$ satisfies
$F(B_{a,b})+3\gamma t(B_{a,b})=1-e(B_{a,b})$.
\end{enumerate}
\end{lemma}

\begin{proof}
The partial derivatives of $\Phi_\gamma$ are $\partial_a\Phi_\gamma=\frac12\log a+\frac34\gamma(a^2+b^2)$ and
$\partial_b\Phi_\gamma=\frac12\log b+\frac32\gamma ab$.

(a) As $a\to0$, we have $\partial_a\Phi_\gamma\le\frac12\log a+\frac32\gamma\to-\infty$ uniformly in $b$,
and similarly for $\partial_b\Phi_\gamma$ as $b\to0$. Moreover, $\partial_a\Phi_\gamma(1,b)>0$, and
$\partial_b\Phi_\gamma(a,1)>0$ for $a>0$. Hence a minimizer, which exists by compactness, is not a
boundary point, and it is therefore a stationary point with positive semidefinite Hessian matrix.

(b) The gradient vanishes exactly when $(a,b,3\gamma)$ satisfies \eqref{s4:eq:EL}. For $a=b=s$, both
equations read $\log s+3\gamma s^2=0$, so $s=s_\gamma$. For $a\ne b$, Lemma~\ref{lem:branch} applies
because $(a^3+3ab^2)/4<1$, and it gives exactly the pairs $(a(z),b(z))$ with $\mu(z)=3\gamma$.

(c) At a stationary point, $0=a\partial_a\Phi_\gamma+b\partial_b\Phi_\gamma=\frac12(a\log a+b\log b)
+3\gamma t(B_{a,b})$. Moreover, $\frac12(a\log a+b\log b)=F(B_{a,b})+e(B_{a,b})-1$ because
$u\log u=h(u)+u-1$.
\end{proof}

For $z\in(0,1)$ define $\Delta(z)=\Phi_{\mu(z)/3}(a(z),b(z))-C_0(\mu(z)/3)$. This is the difference
between the values of $\Phi_{\mu(z)/3}$ at its stationary points $(a(z),b(z))$ and
$(s_{\mu(z)/3},s_{\mu(z)/3})$.

\begin{lemma}\label{lem:tie}
The function $\Delta$ is continuously differentiable on $(0,1)$, with
\begin{equation}\label{eq:s7Deltaprime}
\Delta'(z)=\frac{\mu'(z)}3\big(\rho(z)^3-s_{\mu(z)/3}^3\big).
\end{equation}
If $\Delta(z)=0$, then $\rho(z)<s_{\mu(z)/3}$. Moreover, $\Delta(z)\to-\frac12$ as $z\to0$, and
$\Delta(10^{-3})>0$.
\end{lemma}

\begin{proof}
Let $\gamma=\mu(z)/3$. By Lemma~\ref{lem:stationary}(b), the derivative of
$y\mapsto\Phi_\gamma(a(y),b(y))$ vanishes at $y=z$. Since also
$\partial_\gamma\Phi_\gamma(a(z),b(z))=t(B_z)=\rho(z)^3$ and $C_0'(\gamma)=s_\gamma^3$, the chain rule
gives \eqref{eq:s7Deltaprime}.

Now suppose that $\Delta(z)=0$, and let $s=s_\gamma$ and $B=B_z$. Lemma~\ref{lem:stationary}(c),
applied to $(a(z),b(z))$ and to $(s,s)$, gives $F(B)+3\gamma t(B)=1-e(B)$ and $h(s)+3\gamma s^3=1-s$.
Subtracting $\Delta(z)=0$, that is, $F(B)+\gamma t(B)=h(s)+\gamma s^3$, we get
$2\gamma(t(B)-s^3)=s-e(B)$. If $e(B)\le s$, then $t(B)\ge s^3$, and $F(B)>h(e(B))\ge h(s)$ since $h$ is
strictly convex and decreasing on $[0,1]$. Then $F(B)+\gamma t(B)>h(s)+\gamma s^3$, a contradiction.
Hence $e(B)>s$, and $\rho(z)^3=t(B)<s^3$.

As $z\to0$, we have $a(z)\to0$, $b(z)\to1$ and $\gamma=\mu(z)/3\to\infty$ (Lemma~\ref{lem:fold}(b)), so
$s_\gamma\to0$. Hence Lemma~\ref{lem:stationary}(c) gives
$\Phi_\gamma(a(z),b(z))=F(B_z)+\frac13(1-F(B_z)-e(B_z))\to\frac12$ and
$C_0(\gamma)=h(s_\gamma)+\frac13(1-h(s_\gamma)-s_\gamma)\to1$. Thus $\Delta(z)\to-\frac12$.

Finally, let $z=10^{-3}$ and $\gamma=\mu(z)/3$. Since $a(z)<z$ and $h$ is decreasing and nonnegative on
$[0,1]$, we have $\Phi_\gamma(a(z),b(z))\ge h(10^{-3})/2=(1-10^{-3}(1+\log1000))/2>0.496$. On the other
hand, $\log b(z)^2=4z\log z/(1-z)^2>-0.0277$, so $b(z)^2>0.972$. Hence
$\gamma=2L/(3(1-z)^2b(z)^2)<13.816/(3\cdot0.998\cdot0.972)<4.75$. Since $C_0(\gamma)$ is the minimum of
$h(s)+\gamma s^3$, the choice $s=0.3$ gives $C_0(\gamma)\le h(0.3)+0.027\gamma<0.3389+0.1283<0.468$. Hence
$\Delta(10^{-3})>0$.
\end{proof}

For $z\in(0,1)$, Lemma~\ref{lem:stationary}(a) and Lemma~\ref{lem:fold}(a),(c) show that
$(a(z),b(z))$ can minimize $\Phi_{\mu(z)/3}$ only if $Q(z)\le1$, that is, only if $z\le z_m$. We next
determine the sign of $\Delta$ on $(0,z_m]$.

\begin{lemma}\label{lem:zc}
There is exactly one $z_c\in(0,z_m]$ with $\Delta(z_c)=0$. Moreover, $z_c<10^{-3}$, $\Delta<0$ on
$(0,z_c)$ and $\Delta>0$ on $(z_c,z_m]$.
\end{lemma}

\begin{proof}
By Lemma~\ref{lem:fold}(a),(b), $\mu'<0$ on $(0,z_m)$. Let $z_0\in(0,z_m]$ be a zero of $\Delta$. By
Lemma~\ref{lem:tie} and continuity, $\rho<s_{\mu/3}$ near $z_0$. Hence, by \eqref{eq:s7Deltaprime},
$\Delta'>0$ on $(0,z_m)$ near $z_0$. Thus $\Delta<0$ just left of $z_0$, and $\Delta>0$ just right of
$z_0$ if $z_0<z_m$. In particular, the zeros of $\Delta$ in $(0,z_m]$ are isolated. Since
$\Delta\to-\frac12$ at $0$, they form a compact, hence finite, set. Between two consecutive zeros,
$\Delta$ would change sign from positive to negative and vanish once more. Hence $\Delta$ has at most one
zero in $(0,z_m]$. Since $\Delta(10^{-3})>0$ and $10^{-3}<z_m$, it has exactly one zero
$z_c\in(0,10^{-3})$, and $\Delta$ has the stated signs.
\end{proof}

Let $\gamma_c=\mu(z_c)/3$, $r_+=s_{\gamma_c}$ and $r_-=\rho(z_c)$. Then $t(B_{z_c})=r_-^3$. Since
$\Delta(z_c)=0$, the graphon $B_{z_c}$ and the constant graphon $r_+$ have the same value of
$F+\gamma_ct$.

\begin{proposition}\label{prop:pen}
We have $r_-<r_+$, and the following hold for $\gamma>0$.
\begin{enumerate}[label=(\alph*)]
\item If $\gamma<\gamma_c$, then $(s_\gamma,s_\gamma)$ is the unique minimizer of $\Phi_\gamma$, and
$C(\gamma)=C_0(\gamma)$.
\item If $\gamma=\gamma_c$, then the minimizers of $\Phi_\gamma$ are exactly $(r_+,r_+)$ and
$(a(z_c),b(z_c))$. Moreover, $C(\gamma_c)=C_0(\gamma_c)=\Phi_{\gamma_c}(a(z_c),b(z_c))$.
\item If $\gamma>\gamma_c$, then there is exactly one $z\in(0,z_c)$ with $\mu(z)=3\gamma$, and
$(a(z),b(z))$ is the unique minimizer of $\Phi_\gamma$.
\end{enumerate}
\end{proposition}

\begin{proof}
Since $\Delta(z_c)=0$, Lemma~\ref{lem:tie} gives $r_-=\rho(z_c)<s_{\gamma_c}=r_+$.

By Lemma~\ref{lem:stationary}(a),(b) and Lemma~\ref{lem:fold}(a),(c), every minimizer of
$\Phi_\gamma$ is $(s_\gamma,s_\gamma)$ or a pair $(a(z),b(z))$ with $\mu(z)=3\gamma$ and $z\le z_m$.
Since $\mu$ decreases strictly and continuously from $\infty$ on $(0,z_m]$, such a pair exists exactly
when $3\gamma\ge\mu(z_m)$, and then it is unique. Let $z(\gamma)$ be its ratio. Then $z(\gamma)<z_c$ if
and only if $\gamma>\gamma_c$, and $z(\gamma)=z_c$ if and only if $\gamma=\gamma_c$. Since
$\Phi_\gamma(a(z(\gamma)),b(z(\gamma)))-\Phi_\gamma(s_\gamma,s_\gamma)=\Delta(z(\gamma))$, the signs of
$\Delta$ in Lemma~\ref{lem:zc} give (a), (b) and (c).
\end{proof}

We obtain the digits of these constants from interval evaluations of explicit functions of one variable.

\begin{lemma}\label{lem:penconstants}
We have
\begin{align*}
z_c&\in(8.352460192\cdot10^{-5},8.352460193\cdot10^{-5}), & \gamma_c&\in[6.280969860,6.280969861],\\
C(\gamma_c)&\in[0.4999589178,0.4999589179], & r_-&\in[0.03965310391,0.03965310392],
\end{align*}
and $r_+\in[0.2653474271,0.2653474272]$.
\end{lemma}

\begin{proof}
We enclose $s_\gamma$ by the signs of $s\mapsto\log s+3\gamma s^2$. By Lemmas~\ref{lem:zc} and~\ref{lem:fold}, interval evaluation of $\Delta$ at
the two endpoints of the interval for $z_c$ then locates $z_c$. Evaluation of $\mu/3$,
$\rho$, $\Phi_{\mu/3}(a,b)$ and $s_\gamma$ on the resulting enclosures gives the other constants. See
Appendix~\ref{app:constants} for details.
\end{proof}

By Proposition~\ref{prop:pen} and Lemma~\ref{lem:penconstants}, $c=(3\gamma_c)^{1/2}=4.340\ldots$ is the
transition point $c_*\approx4.341$ of Kwan and Pham~\cite[Section~1.1.5]{KP26}. We next show that the
small targets lie below the threshold.

\begin{lemma}\label{lem:order}
$r_-<r_*$.
\end{lemma}

\begin{proof}
Recall that $r_+=s_{\gamma_c}$ is the only minimizing point of $h(s)+\gamma_cs^3$ on $[0,1]$, and
$r_-\ne r_+$. Since $\Delta(z_c)=0$ and $t(B_{z_c})=r_-^3$, we get
$F(B_{z_c})+\gamma_cr_-^3=C_0(\gamma_c)<h(r_-)+\gamma_cr_-^3$. Hence $\Vbal(r_-)\le F(B_{z_c})<h(r_-)$.
Since $\Vbal(r)=h(r)$ for $r\ge r_*$ by Corollary~\ref{cor:Vbal}, we conclude that $r_-<r_*$.
\end{proof}

\subsection{The penalized problem over all graphons}\label{s7:data}

Let $c>0$. For the strictly $2$-balanced graph $H=K_3$ and $z=0$, the quantity
$\Psi_{\mathrm{Bin}}(K_3,c,0)$ of Kwan and Pham~\cite[Corollary~1.2]{KP26} is the minimum over
graphons $W$ of $F(W)+\frac{c^2}3t(W)$. Indeed, the
penalty there is $\psi_0(t(W))=t(W)$, with coefficient $2c^{|E(K_3)|-1}/|\mathrm{Aut}(K_3)|=\frac{c^2}3$.
Kwan and Pham~\cite[Theorem~1.6]{KP26} showed that $\Psi_{\mathrm{Bin}}(K_3,c,0)=\min_{[0,1]^2}\Phi_{c^2/3}$. With
$c=(3\gamma)^{1/2}$ this reads
\begin{equation}\label{eq:KP}
\min\{F(W)+\gamma t(W):W\text{ a graphon}\}=C(\gamma)\quad\text{for every }\gamma>0 .
\end{equation}
Since every $B_{a,b}$ is a graphon, the content of \eqref{eq:KP} is the inequality
$F(W)+\gamma t(W)\ge C(\gamma)$ for all graphons $W$.

By \eqref{eq:KP} and Proposition~\ref{prop:pen}, every minimizer at a target $r\le r_-$ also minimizes
$F+\gamma t$ over all graphons for some $\gamma>0$ (Proposition~\ref{prop:smalldata}). For such a minimizer, we identify the multiplier by scaling and obtain a
spectral constraint from the second variation under the vertex reweighting
of Section~\ref{s3:sub:vertex}.

\begin{lemma}\label{lem:penalized}
Let $U$ be a minimizer at $r\in(0,1)$, and let $\gamma>0$. Suppose that
$F(U)+\gamma t(U)\le F(W)+\gamma t(W)$ for every graphon $W$. Then
$\mu_U=3\gamma$, and $\langle f,T_Uf\rangle\le0$ for every
$f\in L^2([0,1])$ with $\int_0^1f=0$.
\end{lemma}

\begin{proof}
By Proposition~\ref{prop:euler}(d), $\lambda U$ agrees a.e.\ with a graphon
for all $\lambda$ near $1$. By \eqref{eq:scale1}, Lemma~\ref{lem:active} and
\eqref{eq:mult}, the derivative of $F(\lambda U)+\gamma t(\lambda U)$ at
$\lambda=1$ is $F(U)+e(U)-1+3\gamma r^3=(3\gamma-\mu_U)r^3$. It vanishes,
since $\lambda=1$ is an interior minimum. Hence $\mu_U=3\gamma$.

Let $f$ be measurable with $|f|\le1$ and $\int_0^1f=0$, and take $\psi=f$
in the construction of $U^\varepsilon$ in Section~\ref{s3:sub:vertex}.
Since every $U^\varepsilon$ is a graphon, the polynomial
$\varepsilon\mapsto F(U^\varepsilon)+\gamma t(U^\varepsilon)$ of
\eqref{s3:expand} has a minimum at $\varepsilon=0$ on $[-\frac12,\frac12]$.
Hence its second derivative at $0$ is nonnegative. By \eqref{s3:expand},
this second derivative is
\[
2\iint\big(h(U)+3\gamma UC_U\big)(f\otimes f)
\overset{\mathrm{(a)}}{=}2\iint(1-U)(f\otimes f)
\overset{\mathrm{(b)}}{=}-2\langle f,T_Uf\rangle,
\]
where $\mathrm{(a)}$ uses $3\gamma=\mu_U$ and
$h(U)=U\log U-U+1=-\mu_UUC_U-U+1$ by \eqref{eq:euler}; $\mathrm{(b)}$
uses $\iint f\otimes f=(\int_0^1f)^2=0$. Hence $\langle f,T_Uf\rangle\le0$.
By scaling, this holds for every bounded $f$ with mean zero. These functions
are dense in the functions of mean zero in $L^2([0,1])$, and $T_U$ is
bounded on $L^2([0,1])$. Hence the inequality holds for every
$f\in L^2([0,1])$ with $\int_0^1f=0$.
\end{proof}

\begin{proposition}\label{prop:smalldata}
Let $r\in(0,r_-]$. Then $z_r\in(0,z_c]$, $\mu_r\ge3\gamma_c$ and $V(r)=F(B_r)$. Every minimizer $U$ at
$r$ satisfies $t(U)=r^3$, $e(U)=m_r$, $F(U)=F(B_r)$ and $\mu_U=\mu_r$. Moreover, every such $U$ minimizes
$F+\frac{\mu_r}3t$ over all graphons.
\end{proposition}

\begin{proof}
Since $z_c<10^{-3}<z_*<z_{\mathrm{sp}}$ (Lemma~\ref{lem:zc}, Lemma~\ref{lem:digits} and
Theorem~\ref{thm:crossing}(a)) and $\rho(z_c)=r_-$, we have $z_1(r_-)=z_c$. Since $z_1$ is increasing,
$z_r=z_1(r)\in(0,z_c]$. Let $\gamma=\mu_r/3$. Since $\mu$ is decreasing on $(0,z_m]$, we have
$\gamma\ge\gamma_c$. Hence Proposition~\ref{prop:pen}(b),(c) shows that $(a_r,b_r)$ minimizes
$\Phi_\gamma$. For $\gamma>\gamma_c$, we use here that $z_r$ is the only $z\in(0,z_c)$ with $\mu(z)=3\gamma$.
Thus $C(\gamma)=F(B_r)+\gamma r^3$.
By \eqref{eq:KP}, every graphon $U$ with $t(U)\le r^3$ satisfies
\begin{equation}\label{eq:s7dual}
F(U)\ge C(\gamma)-\gamma t(U)=F(B_r)+\gamma(r^3-t(U))\ge F(B_r).
\end{equation}
As $t(B_r)=r^3$, this gives $V(r)=F(B_r)$. Now let $U$ be a minimizer at $r$. Then equality holds in
\eqref{eq:s7dual}, so $t(U)=r^3$ and $F(U)+\gamma t(U)=C(\gamma)$. By \eqref{eq:KP}, $U$ minimizes
$F+\gamma t$ over all graphons. Hence Lemma~\ref{lem:penalized} gives $\mu_U=3\gamma=\mu_r$. Finally, by
\eqref{eq:mult} and Lemma~\ref{lem:stationary}(c), $e(U)=1-F(B_r)-3\gamma r^3=e(B_r)=m_r$.
\end{proof}

The penalized problem does not determine $V$ near $r_*$.

\begin{remark}\label{rem:convex}
For $r\in(0,1)$ let $V_{\mathrm{pen}}(r)=\sup_{\gamma\ge0}(C(\gamma)-\gamma r^3)$. By \eqref{eq:KP} and
$C(0)=0$, we have $V_{\mathrm{pen}}\le V\le\Vbal$. The proof of Proposition~\ref{prop:smalldata} gives
$V_{\mathrm{pen}}=\Vbal$ on $(0,r_-]$. However, $V_{\mathrm{pen}}<\Vbal$ on $(r_-,r_+)$. Indeed, let
$r\in(r_-,r_+)$, $t_1=r_-^3$, $t_2=r_+^3$ and $\theta=(t_2-r^3)/(t_2-t_1)$. For $\gamma\ge0$, the values
of $\Phi_\gamma$ at $(a(z_c),b(z_c))$ and $(r_+,r_+)$ give $C(\gamma)\le F(B_{z_c})+\gamma t_1$ and
$C(\gamma)\le h(r_+)+\gamma t_2$. By Proposition~\ref{prop:pen}(b), both are equalities for
$\gamma=\gamma_c$. Averaging with weights $\theta$ and $1-\theta$ gives
$C(\gamma)-\gamma r^3\le\theta F(B_{z_c})+(1-\theta)h(r_+)=C(\gamma_c)-\gamma_cr^3$. Hence
$V_{\mathrm{pen}}(r)=C(\gamma_c)-\gamma_cr^3$. Now let $(a,b)$ attain $\Vbal(r)$. Then
$\Vbal(r)+\gamma_cr^3\ge\Phi_{\gamma_c}(a,b)\ge C(\gamma_c)$. Equality would make $(a,b)$ a minimizer of
$\Phi_{\gamma_c}$ with $t(B_{a,b})=r^3\notin\{t_1,t_2\}$, contrary to Proposition~\ref{prop:pen}(b).
Hence $V_{\mathrm{pen}}(r)<\Vbal(r)$. By Lemmas~\ref{lem:digits} and~\ref{lem:penconstants},
$r_*\in(r_-,r_+)$ and $h(r_*)-V_{\mathrm{pen}}(r_*)>0.02$.

\end{remark}

\subsection{Minimizers at small targets are two block SBM graphons}\label{s7:bip}

The next proposition is the remaining step in the proof of Theorem~\ref{thm:smalltail}.

\begin{proposition}\label{prop:structural}
Let $r\in(0,r_-]$ and let $U$ be a minimizer at $r$. Then $U$ is a two block SBM graphon.
\end{proposition}

\begin{proof}[Proof of Theorem~\ref{thm:smalltail}]
Let $U$ be a minimizer at $r\in(0,r_-]$. By Proposition~\ref{prop:structural}, $U$ is a two block SBM
graphon. Since $r<r_*$ (Lemma~\ref{lem:order}), Corollary~\ref{cor:bipodalmin} shows that $U$ is a
relabeling of $B_r$. Conversely, $V(r)=F(B_r)$ by Proposition~\ref{prop:smalldata}, and $t(B_r)=r^3$.
Hence every relabeling of $B_r$ is a minimizer.
\end{proof}

For the rest of the section, fix $r\in(0,r_-]$ and a minimizer $U$ at $r$. Let $\mu=\mu_U$, $m=e(U)$,
$t=t(U)$, $d=d_U$ and $\bar k=\bar k_U$. By Proposition~\ref{prop:smalldata}, we have $z=z_r\in(0,z_c]$,
$\mu=\mu(z)$, $m=(a(z)+b(z))/2$ and $t=\rho(z)^3$. By \eqref{eq:mult} and Proposition~\ref{prop:vertex},
$\bar k=1+F(U)-m=2-2m-\mu t$. We use the following four facts.
\begin{enumerate}[label=(F\arabic*)]
\item $U=\exp(-\mu C_U)$ a.e.\ (Proposition~\ref{prop:euler}).
\item $\int_0^1k(U(x,y))\dd y=\bar k$ and $\mu\tau_U(x)+2d(x)=\mu t+2m$ for a.e.\ $x$
(Proposition~\ref{prop:vertex}).
\item For a.e.\ $x$, the row $U(x,\cdot)$ satisfies $\log U(x,\cdot)+\mu T_UU(x,\cdot)=0$ a.e.\
(Lemma~\ref{lem:insertion}).
\item $\langle f,T_Uf\rangle\le0$ for every $f\in L^2([0,1])$ with $\int_0^1f=0$
(Lemma~\ref{lem:penalized}, which applies by Proposition~\ref{prop:smalldata}).
\end{enumerate}
Let $D=d-m$, $\eta=\mu t/2$, $A=m+\mu t$ and $K(x,y)=d(x)+d(y)-m-U(x,y)$. Let $f\otimes g$ denote the
operator $\varphi\mapsto\langle g,\varphi\rangle f$. We control all error terms by $\mu t$, $\mu^2t$
and $\mu(\frac12-m)$. These three numbers tend to $0$ with $z$, although $\mu\to\infty$. We first show
that $T_K$ is positive and nearly of rank one.

\begin{lemma}\label{lem:spectral}
The operator $T_K$ is positive and trace class, $T_K\1=0$ and $\tr T_K\le m$. Let
$T_K=\sum_j\lambda_je_j\otimes e_j$ with $\lambda_1\ge\lambda_2\ge\dots>0$ and $(e_j)$ orthonormal. Define
$K(x,x)=\sum_j\lambda_je_j(x)^2$. Then, for almost every $x$,
\begin{equation}\label{eq:s7diag}
0\le K(x,x)\le2d(x)-m,
\end{equation}
and consequently
\begin{equation}\label{eq:s7degree}
\frac m2\le d\le m+\eta,\quad K(x,x)\le A,\quad \|D\|_1\le2\eta,\quad \|D\|_2^2\le\frac{\eta m}2 .
\end{equation}
Moreover, $t=m^3-\tr(T_K^3)+3\langle D,(mI-T_K)D\rangle$. Let $\lambda=\|T_K\|$ be the largest
eigenvalue, and let $s=m-\lambda$. Then
\begin{equation}\label{eq:s7tail}
m\lambda^2\ge m^3-t,\quad 3m\lambda s\le t .
\end{equation}
\end{lemma}

\begin{proof}
\emph{Positivity.} We have $T_K\1=d+m-m-d=0$. For $f_0\perp\1$ we have
$T_Kf_0=\langle d,f_0\rangle\1-T_Uf_0$. Hence $\langle f_0,T_Kf_0\rangle=-\langle f_0,T_Uf_0\rangle\ge0$ by
(F4). As $T_K$ is self-adjoint and $T_K\1=0$, it follows that $T_K\ge0$. Since $K$ is bounded, $T_K$ is
compact and has a decomposition as in the statement.

\emph{Trace and diagonal.} Fix $N$ and let $K_N(x,y)=\sum_{j\le N}\lambda_je_j(x)e_j(y)$. Let
$\mathcal D_n$ be the partition of $[0,1]$ into dyadic intervals of length $2^{-n}$. For $x$ in a cell
$I\in\mathcal D_n$, positivity of $T_K-T_{K_N}$ gives
\[
\sum_{j\le N}\lambda_j\Big(\frac1{|I|}\int_Ie_j\Big)^2=\frac1{|I|^2}\iint_{I\times I}K_N\le\frac1{|I|^2}\iint_{I\times
I}K\le\frac2{|I|}\int_Id-m,
\]
where the last step uses $U\ge0$. As $n\to\infty$, the Lebesgue differentiation theorem along dyadic
cells gives $\sum_{j\le N}\lambda_je_j(x)^2\le2d(x)-m$ for a.e.\ $x$. Integrating gives
$\sum_{j\le N}\lambda_j\le m$. Letting $N\to\infty$, we see that $T_K$ is trace class with $\tr T_K\le m$,
and hence $T_K\le mI$. The same limit proves \eqref{eq:s7diag}.

\emph{Consequences.} The bound $d\ge m/2$ follows from \eqref{eq:s7diag}. Since $\tau_U\ge0$, (F2) gives
$D\le\eta$. Hence $K(x,x)\le2d-m\le A$. As $\int D=0$, we have $\|D\|_1=2\int D_+\le2\eta$. Integrating
$(D+\frac m2)(\eta-D)\ge0$ gives $\|D\|_2^2\le\eta m/2$.

\emph{Eigenvalues.} The operator $T_U$ is Hilbert--Schmidt, and $T_U^2$ has kernel $C_U$. Hence
$\tr(T_U^3)=\iint UC_U=t$ by Lemma~\ref{lem:trilinear}. In the decomposition $L^2=\R\1\oplus\1^\perp$,
we have $T_U\1=m\1+D$ and $T_Uf=\langle D,f\rangle\1-T_Kf$ for $f\perp\1$. Thus
\[
T_U=\begin{pmatrix}m&\langle D,\cdot\rangle\\D&-T_K\end{pmatrix}.
\]
In $T_U^3$, the upper left entry is $m^3+2m\|D\|_2^2-\langle D,T_KD\rangle$. The lower right block is
$mD\otimes D-D\otimes T_KD-T_KD\otimes D-T_K^3$. Adding their traces gives the identity for $t$. Since
$T_K\le mI$, this identity gives $m^3-t\le\tr(T_K^3)\le\lambda^2\tr T_K\le m\lambda^2$. The other
eigenvalues have sum at most $\tr T_K-\lambda\le s$. Hence the sum of their cubes is at most $s^3$, and
$m^3-t\le\tr(T_K^3)\le\lambda^3+s^3=m^3-3m\lambda s$.
\end{proof}

Since $m^3-t=(b(z)-a(z))^3/8>0$ by \eqref{s4:eq:cube}, \eqref{eq:s7tail} gives $\lambda>0$. Let
$\psi=e_1$, so that $T_K\psi=\lambda\psi$, $\|\psi\|_2=1$ and $\int\psi=\langle\psi,\1\rangle=0$. Let
$T_R=T_K-\lambda\psi\otimes\psi$, with kernel $R(x,y)=K(x,y)-\lambda\psi(x)\psi(y)$, and define
$R(x,x)=\sum_{j\ge2}\lambda_je_j(x)^2$. Then $T_R\ge0$, $T_R\1=T_R\psi=0$ and $\|T_R\|\le\tr T_R\le s$.
For almost every $x$, we have $\lambda\psi(x)^2\le K(x,x)\le A$, $0\le R(x,x)\le A$ and
$\sum_j\lambda_j^2e_j(x)^2\le\lambda_1K(x,x)<\infty$. With $K_N$ as in the proof of
Lemma~\ref{lem:spectral}, we have $\|K-K_N\|_2^2=\int_0^1\|K(x,\cdot)-K_N(x,\cdot)\|_2^2\dd x\to0$. Hence,
for a.e.\ $x$, a subsequence of $K_N(x,\cdot)$ converges to $K(x,\cdot)$ in $L^2$. Thus, for a.e.\ $x$,
we have $K(x,\cdot)=\sum_j\lambda_je_j(x)e_j$ and $R(x,\cdot)=\sum_{j\ge2}\lambda_je_j(x)e_j$ in $L^2$, and
$(T_Rf)(x)=\sum_{j\ge2}\lambda_je_j(x)\langle e_j,f\rangle$ for $f\in L^2$. Let $M=(A/\lambda)^{1/2}\ge1$.
These facts give, for almost every $x$,
\begin{equation}\label{eq:s7Rrow}
|\psi(x)|\le M,\quad \|R(x,\cdot)\|_2^2=\sum_{j\ge2}\lambda_j^2e_j(x)^2\le sR(x,x)\le sA .
\end{equation}
Let $P_+=\{\psi\ge0\}$, $P_-=\{\psi<0\}$, $p_\pm=|P_\pm|$, $\sigma=\1_{P_+}-\1_{P_-}$ and
$p_0=1/(2M^2)$. We call $P_+$ and $P_-$ the blocks. Since $1=\|\psi\|_2^2\le M\|\psi\|_1$ and
$\int\psi=0$, we have $Mp_\pm\ge\int_{P_\pm}|\psi|=\frac12\|\psi\|_1\ge1/(2M)$, so
\begin{equation}\label{eq:s7mass}
p_0\le p_\pm\le1-p_0,\quad \|\sigma-\psi\|_2^2=2-2\|\psi\|_1\le2-\frac2M .
\end{equation}
For measurable $S$ with $|S|>0$ let $\E_S[g]=|S|^{-1}\int_Sg$. For almost every $x$, let
$u_x=\E_{P_+}[U(x,\cdot)]$ and $v_x=\E_{P_-}[U(x,\cdot)]$ be the block averages of the row $U(x,\cdot)$.
Then $d(x)=p_+u_x+p_-v_x$. Define
\begin{align}
E_m&=\eta(2M^2-1)+3\lambda M(M-1)+M^2(sA)^{1/2}(2-2/M)^{1/2},\notag\\
E_C&=\frac{\eta m}2+2\lambda M^2\eta+\Big(\frac{As\eta m}2\Big)^{1/2}+4M^2(m+\eta)\eta\notag\\
&\quad+2M^2s\Big(\frac{\eta m}2\Big)^{1/2}+5\lambda^2M(M-1)+2M^2s(sA)^{1/2},\notag\\
E_0&=|m-\tfrac12|(m+\eta)+\lambda|\lambda-\tfrac12|M,\quad
\varepsilon=\mu(E_C+E_0+E_m/2).\label{eq:s7errors}
\end{align}
We next show that the block averages $u_x$ and $v_x$ of almost every row satisfy an
approximate Euler--Lagrange system and an inequality from the vertex identity.

\begin{lemma}\label{lem:rowreduction}
With the constants \eqref{eq:s7errors}, for almost every $x$,
\begin{align}
|u_x-(d(x)-\lambda\psi(x))|&\le E_m,\quad |v_x-(d(x)+\lambda\psi(x))|\le E_m,\label{eq:s7rowmean}\\
|u_x-e^{-\mu v_x/2}|&\le\varepsilon,\quad |v_x-e^{-\mu u_x/2}|\le\varepsilon,\label{eq:s7approxeuler}\\
\frac{k(u_x)+k(v_x)}2&\le\bar k+1-M^{-2}.\label{eq:s7approxk}
\end{align}
\end{lemma}

\begin{proof}
\emph{Averages of $\psi$ and $D$.} By \eqref{eq:s7mass}, $\E_{P_+}[\psi]=\|\psi\|_1/(2p_+)$ lies between
$(1/M)/(2-2p_0)=M/(2M^2-1)$ and $M$. Moreover, $1-M/(2M^2-1)=(M-1)(2M+1)/(2M^2-1)\le3(M-1)$ for $M\ge1$.
On $P_+$, the bounds $0\le\psi\le M$ give $|\psi-1|\le1-\psi+2(M-1)$, so $\E_{P_+}[|\psi-1|]\le5(M-1)$.
Arguing in the same way for $-\psi$ on $P_-$, we get for $S=P_\pm$
\begin{equation}\label{eq:s7psiavg}
|\E_S[\psi-\sigma]|\le3(M-1),\quad \E_S[|\psi-\sigma|]\le5(M-1).
\end{equation}
Let $\beta=\langle D,\psi\rangle$. As $\int D=\int\psi=0$ and $D\le\eta$, we have
$|\beta|=|\int(D-\eta)\psi|\le M\int(\eta-D)=M\eta$. Let $S=P_\pm$. Then $\E_S[D]\le\eta$ and
$|S|\E_S[D]=-\int_{S^c}D\ge-(1-|S|)\eta$. Since $(1-p_0)/p_0=2M^2-1$, it follows that
\begin{equation}\label{eq:s7Davg}
|\E_S[D]|\le\eta(2M^2-1),\quad \E_S[|D|]\le\|D\|_1/p_0\le4M^2\eta .
\end{equation}

\emph{Row averages.} Since $T_R\1=T_R\psi=0$ and $\1_{P_\pm}=(\1\pm\sigma)/2$, we have
$T_R\1_{P_\pm}=\pm\frac12T_R(\sigma-\psi)$. Hence, by \eqref{eq:s7Rrow}, \eqref{eq:s7mass} and the
Cauchy--Schwarz inequality, $|(T_R\1_{P_\pm})(x)|\le\frac12(sA)^{1/2}(2-2/M)^{1/2}$. Averaging
$U(x,y)=d(x)+D(y)-\lambda\psi(x)\psi(y)-R(x,y)$ over $y\in P_+$ gives
\[
u_x-(d(x)-\lambda\psi(x))=\E_{P_+}[D]-\lambda\psi(x)\big(\E_{P_+}[\psi]-1\big)-(T_R\1_{P_+})(x)/p_+.
\]
By \eqref{eq:s7Davg}, \eqref{eq:s7psiavg} and $1/p_+\le2M^2$, the three terms on the right are bounded by
the three terms of $E_m$. The same computation on $P_-$ proves \eqref{eq:s7rowmean}.

\emph{The codegree.} Since $U(x,w)=d(x)+D(w)-K(x,w)$, $\int D=0$ and $T_K\1=0$, we have
\[
C_U(x,y)=d(x)d(y)+\|D\|_2^2-(T_KD)(x)-(T_KD)(y)+K^{(2)}(x,y),
\]
where $K^{(2)}$ is the kernel of $T_K^2$. Since $T_R\psi=0$, we have
$T_K^2=\lambda^2\psi\otimes\psi+T_R^2$ and $T_KD=\lambda\beta\psi+T_RD$. Let $R^{(2)}$ be the kernel of
$T_R^2$, and let $c_\pm(x)=md(x)\pm\lambda^2\psi(x)$. Writing $d(x)d(y)=md(x)+d(x)D(y)$, we get, for
almost every $(x,y)$ with $y\in P_\pm$,
\begin{align}
C_U(x,y)-c_\pm(x)&=\|D\|_2^2-\lambda\beta(\psi(x)+\psi(y))-(T_RD)(x)+d(x)D(y)\notag\\
&\quad-(T_RD)(y)+\lambda^2\psi(x)(\psi(y)-\sigma(y))+R^{(2)}(x,y).\label{eq:s7codeg}
\end{align}
Let $S=P_\pm$. The averages over $y\in S$ of the absolute values of the seven terms in
\eqref{eq:s7codeg} are, in order, at most the seven terms of $E_C$. For the first, second, fourth and
sixth terms, this follows from \eqref{eq:s7degree}, $|\beta|\le M\eta$, $|\psi|\le M$,
\eqref{eq:s7Davg} and \eqref{eq:s7psiavg}. For the third term, the Cauchy--Schwarz inequality in the
expansion of $T_R$ gives $(T_RD)(x)^2\le R(x,x)\langle D,T_RD\rangle\le As\|D\|_2^2$. For the fifth and
seventh terms, we use $\E_S[|g|]\le\|g\|_1/p_0\le2M^2\|g\|_2$, $\|T_RD\|_2\le s\|D\|_2$ and
$\|R^{(2)}(x,\cdot)\|_2=\|T_RR(x,\cdot)\|_2\le s(sA)^{1/2}$, where the last bound follows from
\eqref{eq:s7Rrow}. Hence
\begin{equation}\label{eq:s7EC}
\E_{P_\pm}\big[|C_U(x,\cdot)-c_\pm(x)|\big]\le E_C .
\end{equation}

\emph{Approximate Euler--Lagrange equations.} We have $C_U\ge0$ and $c_\pm(x)\ge0$. For the second
bound, note that $\lambda\psi(x)^2\le K(x,x)\le2d(x)-m$ by \eqref{eq:s7diag}. As $0<\lambda\le m$ and
$2|\psi(x)|\le1+\psi(x)^2$, we get
$md(x)\ge\frac12m(m+\lambda\psi(x)^2)\ge\frac12\lambda^2(1+\psi(x)^2)\ge\lambda^2|\psi(x)|$. Since
$y\mapsto e^{-\mu y}$ is $\mu$-Lipschitz on $[0,\infty)$, (F1) and \eqref{eq:s7EC} give
$|u_x-e^{-\mu c_+(x)}|\le\mu E_C$ and $|v_x-e^{-\mu c_-(x)}|\le\mu E_C$. Next,
$c_\pm(x)-\frac12(d(x)\pm\lambda\psi(x))=(m-\frac12)d(x)\pm\lambda(\lambda-\frac12)\psi(x)$. This number is
at most $E_0$ in absolute value. Hence, by \eqref{eq:s7rowmean}, $|c_+(x)-v_x/2|$ and $|c_-(x)-u_x/2|$ are
at most $E_0+E_m/2$. Applying the Lipschitz bound once more gives \eqref{eq:s7approxeuler}.

\emph{The inequality \eqref{eq:s7approxk}.} By convexity of $k$, Jensen's inequality and (F2),
\[
p_+k(u_x)+p_-k(v_x)\le\int_0^1k(U(x,y))\dd y=\bar k.
\]
Moreover, $\frac12(k(u_x)+k(v_x))=p_+k(u_x)+p_-k(v_x)+(\frac12-p_+)(k(u_x)-k(v_x))$. Since $k$ maps
$[0,1]$ into $[0,2]$ and $|\frac12-p_+|\le\frac12-p_0$, the last term is at most $1-2p_0=1-M^{-2}$.
This proves \eqref{eq:s7approxk}.
\end{proof}

We next bound the error terms explicitly, using $z_c<8.36\cdot10^{-5}$, which follows from
Lemma~\ref{lem:penconstants}.

\begin{lemma}\label{lem:smallparams}
Let $Z=8.36\cdot10^{-5}$.
\begin{enumerate}[label=(\alph*)]
\item Let $z\in(0,Z]$, $m=(a(z)+b(z))/2$, $\mu=\mu(z)$ and $t=\rho(z)^3$. Then $0.4992\le m\le\frac12$,
$\mu\ge18.6$, $t\le6.3\cdot10^{-5}$, $\mu t\le0.0012$, $\mu^2t\le0.023$ and $\mu(\frac12-m)\le0.0155$.
\item In the setting of Lemma~\ref{lem:rowreduction}, $\varepsilon<0.114$, $E_m<0.003$,
$\bar k+1-M^{-2}<1.005$, $p_0>0.498$ and $d\le0.5006$.
\end{enumerate}
\end{lemma}

\begin{proof}
(a) Let $\delta=2zL/(1-z)^2$, so that $b(z)=e^{-\delta}$. The functions $zL$ and $zL^2$ increase on
$(0,e^{-2})\supset(0,Z]$. Moreover, $9.3<\log(1/Z)<9.4$. Hence $\delta<2\cdot9.4Z/(1-Z)^2<0.0016$, so
$b(z)>0.9984$ and $m=(1+z)b(z)/2>0.4992$. Since $\log(1+z)\le z\le\delta$, we also have $m\le\frac12$.
Next, $\mu\ge2L>18.6$ and $t=b(z)^3z(z^2+3)/4<Z(3+Z^2)/4<6.3\cdot10^{-5}$. From the branch formulas,
\begin{align*}
\mu t&=\frac{zL(z^2+3)b(z)}{2(1-z)^2}\le\frac{9.4Z(3+Z^2)}{2(1-Z)^2}<0.0012,\\
\mu^2t&=\frac{zL^2(z^2+3)}{(1-z)^4b(z)}\le\frac{9.4^2Z(3+Z^2)}{0.9984(1-Z)^4}<0.023,\\
\mu(\tfrac12-m)&\le\frac{\mu\delta}2=\frac{2zL^2}{(1-z)^4b(z)^2}\le\frac{2\cdot9.4^2Z}{0.9984^2(1-Z)^4}<0.0155,
\end{align*}
where the last line uses $\frac12-m\le\frac12(1-b(z))\le\frac\delta2$.

(b) Part (a) applies since $z=z_r\le z_c<Z$. By \eqref{eq:s7tail},
$\lambda^2\ge m^2-t/m\ge0.4992^2-6.3\cdot10^{-5}/0.4992>0.4989^2$ and $s\le t/(3m\lambda)<1.34t$. As
$A-\lambda=\mu t+s$, we get $M^2-1=(\mu t+s)/\lambda\le(0.0012+1.34\cdot6.3\cdot10^{-5})/0.4989<0.002575$.
Hence $M<1.0013$, $2M^2<2.00515$ and $p_0>0.4987$. We also use $\eta\le0.0006$, $A\le0.5012$,
$d\le m+\eta\le0.5006$ (by \eqref{eq:s7degree}) and $\lambda\le m\le\frac12$. Moreover,
$\mu\eta=\mu^2t/2$ and $\mu s\le1.34\mu t\le0.001608$, so that $\nu=\mu^2t+\mu s\le0.024608$. Since
$M-1\le\frac12(M^2-1)$ and $2-2/M\le2(M-1)$, we have $\mu(M-1)\le\nu/(2\lambda)$ and
$\mu^2sA(2-2/M)\le A\mu s\nu/\lambda$. Moreover, $(\eta m/2)^{1/2}<0.01225$,
$(sA)^{1/2}\le(1.34\cdot6.3\cdot10^{-5}\cdot0.5012)^{1/2}<0.0066$,
$(A\mu s\mu^2tm/4)^{1/2}\le(0.5012\cdot0.001608\cdot0.002875)^{1/2}<0.001523$ and
$(A\mu s\nu/\lambda)^{1/2}\le(0.5012\cdot0.001608\cdot0.024608/0.4989)^{1/2}<0.006305$. Hence $\mu$ times
each term of $E_C$ and of $E_m$, in the order of \eqref{eq:s7errors}, satisfies the bound in the
following table.
\begin{center}
\begin{tabular}{lll}
\toprule
& bound for $\mu$ times the term & numerical value\\
\midrule
$E_C$ & $\mu^2tm/4$ & $0.023\cdot0.5/4=0.002875$\\
& $\lambda M^2\mu^2t$ & $0.5\cdot1.002575\cdot0.023<0.01153$\\
& $(A\mu s\mu^2tm/4)^{1/2}$ & $<0.001523$\\
& $2M^2(m+\eta)\mu^2t$ & $2.00515\cdot0.5006\cdot0.023<0.02309$\\
& $2M^2\mu s(\eta m/2)^{1/2}$ & $2.00515\cdot0.001608\cdot0.01225<0.00004$\\
& $\frac52\lambda M\nu$ & $1.25\cdot1.0013\cdot0.024608<0.0308$\\
& $2M^2\mu s(sA)^{1/2}$ & $2.00515\cdot0.001608\cdot0.0066<0.000022$\\
\midrule
$E_m$ & $\mu^2t(2M^2-1)/2$ & $0.0115\cdot1.00515<0.01156$\\
& $\frac32M\nu$ & $1.5\cdot1.0013\cdot0.024608<0.03696$\\
& $M^2(A\mu s\nu/\lambda)^{1/2}$ & $1.002575\cdot0.006305<0.006322$\\
\bottomrule
\end{tabular}
\end{center}
Summing, $\mu E_C<0.0699$ and $\mu E_m<0.0549$. Since $|m-\frac12|=\frac12-m$ and
$|\lambda-\frac12|=(\frac12-m)+s$, we also have
$\mu E_0\le0.0155\cdot0.5006+0.5\cdot1.0013\cdot(0.0155+0.001608)<0.01633$. Hence
$\varepsilon<0.0699+0.01633+0.02745<0.114$ and $E_m<0.0549/18.6<0.003$. Finally, $1-M^{-2}\le M^2-1$, so
$\bar k+1-M^{-2}\le2-2m+0.002575\le1.0016+0.002575<1.005$.
\end{proof}

We use the following scalar inequality to separate the two block averages of a row. Its hypotheses follow from
\eqref{eq:s7approxeuler} with $\mu/2\ge9.3$ and from \eqref{eq:s7approxk}. Recall that
$k(u)=2-u(2-\log u)$ is convex and strictly decreasing on $[0,1]$, with $k(0)=2$ and $k(1)=0$.

\begin{lemma}\label{lem:scalar}
If $0\le u\le v\le1$, $(k(u)+k(v))/2\le1.005$ and $v\le e^{-9.3u}+0.13$, then $v>0.775$.
\end{lemma}

\begin{proof}
If $v\le0.3$, then $(k(u)+k(v))/2\ge k(0.3)=2-0.3(2-\log0.3)>2-0.3\cdot3.204>1.038$, a contradiction.
Now let $v\in(0.3,0.775]$ and $w(v)=-\log(v-0.13)/9.3\in(0,0.2)$. The last hypothesis gives $u\le w(v)$.
Let $G(v)=k(w(v))+k(v)$. Since $k$ is decreasing, $k(u)+k(v)\ge G(v)$. Moreover,
$G'(v)=\log v-1+(1-\log w(v))/(9.3(v-0.13))$. As $w$ is decreasing, we have
$G'\le\log v_2-1+(1-\log w(v_2))/(9.3(v_1-0.13))$ on an interval $[v_1,v_2]$. We use $w(0.4)>0.139$,
$w(0.5)>0.106$ and $w(0.775)>0.047$. On the intervals $[0.3,0.4]$, $[0.4,0.5]$ and $[0.5,0.775]$, this
bound gives $G'<-1.91+2.98/1.581<0$, $G'<-1.69+3.25/2.511<0$ and $G'<-1.25+4.06/3.441<0$, respectively.
Hence $G(v)\ge G(0.775)$. Since $w(0.775)<0.04716$, $\log0.04716>-3.055$ and $\log0.775>-0.255$, we get
$G(0.775)>k(0.04716)+k(0.775)>(2-0.04716\cdot5.055)+(2-0.775\cdot2.255)>2.01$. This contradicts
$k(u)+k(v)\le2.01$. Hence $v>0.775$.
\end{proof}

\begin{proof}[Proof of Proposition~\ref{prop:structural}]
By Lemma~\ref{lem:rowreduction} and (F3), there is a null set $N$ such that for every $x\notin N$ the
conclusions of Lemma~\ref{lem:rowreduction} hold and $\log U(x,\cdot)+\mu T_UU(x,\cdot)=0$ a.e.

\noindent\textbf{Step 1: the block averages separate.} Let $x\notin N$, and let $u\le v$ be the numbers
$u_x,v_x$ in increasing order. In either case, \eqref{eq:s7approxeuler} and Lemma~\ref{lem:smallparams}(b)
give $v\le e^{-\mu u/2}+\varepsilon\le e^{-9.3u}+0.13$. Similarly, \eqref{eq:s7approxk} and
Lemma~\ref{lem:smallparams}(b) give $(k(u)+k(v))/2<1.005$. Hence Lemma~\ref{lem:scalar} yields $v>0.775$.
Let $q\in\{p_+,p_-\}$ be the measure of the block over which $v$ is the average. Since
$d(x)=qv+(1-q)u\le0.5006$ and $q,1-q\ge p_0>0.498$, we get
$u\le(d(x)-qv)/(1-q)<(0.5006-0.498\cdot0.775)/0.498<0.231$. Hence $|v_x-u_x|>0.544$.

\noindent\textbf{Step 2: the sign of $\psi$ determines the block.} By \eqref{eq:s7rowmean},
$|(v_x-u_x)-2\lambda\psi(x)|\le2E_m<0.006$ for $x\notin N$. Hence, by Step~1, $\psi(x)\ne0$ and $v_x-u_x$
has the sign of $\psi(x)$. Thus $v_x>0.775$ for almost every $x\in P_+$, and $u_x>0.775$ for almost every
$x\in P_-$.

\noindent\textbf{Step 3: $U$ is small on $P_+\times P_+$ and $P_-\times P_-$.} Let $c_0=0.27$ and
$a_\mu=e^{-c_0\mu}$. For almost every $x,y\in P_+$, the inequality $fg\ge f+g-1$ on $[0,1]$ gives
$C_U(x,y)\ge\int_{P_-}U(x,w)U(y,w)\dd w\ge p_-(v_x+v_y-1)>0.498\cdot0.55>c_0$. A similar bound holds on
$P_-\times P_-$. By (F1), $U\le a_\mu$ almost everywhere on $(P_+\times P_+)\cup(P_-\times P_-)$. After
enlarging $N$, we may assume that every $x\notin N$ satisfies $U(x,y)\le a_\mu$ for almost every $y$ in
the block of $x$.

\noindent\textbf{Step 4: rows in one block coincide.} Let $x_0,x_1\in P_+\setminus N$, $f_i=U(x_i,\cdot)$,
$f=f_1-f_0$, $X=\|f\1_{P_+}\|_2$ and $Y=\|f\1_{P_-}\|_2$. Since $x_0,x_1\notin N$ and $0\le T_Uf_i\le1$,
we have $\log f_i=-\mu T_Uf_i\in[-\mu,0]$ a.e. Since $(\log y)'=1/y\ge1/\omega$ on $(0,\omega]$, the mean
value theorem gives $(\alpha-\alpha')(\log\alpha-\log\alpha')\ge(\alpha-\alpha')^2/\omega$ for
$\alpha,\alpha'\in(0,\omega]$. By Step~3, we may use this bound with $\omega=a_\mu$ on $P_+$ and with
$\omega=1$ on $P_-$. Let $\mathcal Q=a_\mu^{-1}X^2+Y^2+\mu\langle f,T_Uf\rangle$. Then
\[
a_\mu^{-1}X^2+Y^2\overset{\mathrm{(a)}}{\le}\int_0^1f(\log f_1-\log f_0)\overset{\mathrm{(b)}}{=}-\mu\langle f,T_Uf\rangle,
\]
where $\mathrm{(a)}$ integrates the bound with $\alpha=f_1(y)$ and $\alpha'=f_0(y)$, and $\mathrm{(b)}$ uses
$\log f_i=-\mu T_Uf_i$. Hence $\mathcal Q\le0$. On the other hand, we split $\langle f,T_Uf\rangle$ over the
four products of blocks. On $P_\pm\times P_\pm$ we use $0\le U\le a_\mu$, and elsewhere we use $U\le1$.
Together with $\int_{P_\pm}|f|\le p_\pm^{1/2}\|f\1_{P_\pm}\|_2$, this gives
\[
\mathcal Q\ge(a_\mu^{-1}-\mu a_\mu p_+)X^2-2\mu(p_+p_-)^{1/2}XY+(1-\mu a_\mu p_-)Y^2 .
\]
Since $a_\mu\le1$, $p_++p_-=1$ and $p_+p_-\le\frac14$, the coefficient of $Y^2$ is at least
$1-\mu a_\mu$. The determinant of this form is $a_\mu^{-1}-\mu p_--\mu a_\mu p_++\mu^2p_+p_-(a_\mu^2-1)$.
By the same three facts, it is at least $a_\mu^{-1}-\mu-\mu^2/4$. Both lower bounds are positive because
$a_\mu^{-1}=e^{c_0\mu}>\mu+\mu^2/4$. Indeed, $\mu^{-2}e^{c_0\mu}$ increases for $\mu\ge2/c_0$, and
$\mu\ge18$ by Lemma~\ref{lem:smallparams}(a). Hence
$e^{c_0\mu}\ge\mu^2e^{18c_0}/324>0.37\mu^2\ge\mu+\mu^2/4$, using $e^{4.86}>120$ and $0.12\mu\ge1$. Thus
the form is positive definite. Since $\mathcal Q\le0$, it follows that $X=Y=0$, that is, $f_0=f_1$ almost
everywhere. The same argument applies to two rows indexed by $P_-\setminus N$.

\noindent\textbf{Step 5: $U$ is a two block SBM graphon.} By Step~4 and Fubini's theorem, there are
functions $f_\pm$ with $U(x,y)=f_\pm(y)$ for almost every $(x,y)\in P_\pm\times[0,1]$. Since $0<p_+<1$
by \eqref{eq:s7mass}, Lemma~\ref{lem:rowsbipodal} shows that $U$ is a two block SBM graphon.
\end{proof}

\section{Global continuation and minimizer below the threshold}\label{sec:intermediate}

Let $G$ be the set of $r\in(0,r_*)$ at which the minimizers are exactly the
relabelings of $B_r$. This set is open by Corollary~\ref{cor:open}. It
contains $(0,r_-]$ by Theorem~\ref{thm:smalltail} and an interval
$(r_*-\delta,r_*)$ by Corollary~\ref{cor:seed}. In this section we show that $G=(0,r_*)$. Let $(r_0,r_*)$ be the maximal interval ending at $r_*$ that
is contained in $G$, and suppose that $r_0\ge r_-$. Using a first order comparison, we show that every
minimizer $U$ at $r_0$ has the cost of $B_{r_0}$ and edge
density at least $m_{r_0}$ (Lemma~\ref{lem:boundary}). These two facts and
the Euler--Lagrange equation force $U$ to be a relabeling of $B_{r_0}$
(Lemma~\ref{lem:tieexclusion}). Hence $r_0\in G$, and the openness of $G$
contradicts the maximality of $(r_0,r_*)$.

In the proof of Lemma~\ref{lem:tieexclusion} we use an explicit function
$\Psi_z$ of the three edge values of a triangle, where $z=z_{r_0}$. With
computer assistance, we show in Proposition~\ref{prop:family} that $\Psi_z$ is
nonnegative and vanishes exactly at the edge value triples $(a,a,a)$ and
$(a,b,b)$ of $B_z=B_{a,b}$ and at their permutations.

\subsection{Certificates for all targets in \texorpdfstring{$[r_-,r_*]$}{[r-,r*]}}\label{s8:sub:cert}

For $z\in(0,1)$ let $L=\log(1/z)$, and let $a=a(z)$, $b=b(z)$ and
$\mu=\mu(z)$ be given by \eqref{eq:branch}. Let
$m=(a+b)/2=e(B_z)$ and $t_z=\rho(z)^3=t(B_z)=b^3\kappa(z)^3$, where
$\kappa(z)^3=(z^3+3z)/4$. Define
\[
\Gamma=6L^2z^2+3Lz^3-Lz^2-Lz-L-z^3+3z^2-3z+1.
\]
If $\Gamma\ne0$, define
\begin{align}
d_1&=\frac{z(z-1)(z+1)(L+z-1)(2L+z-1)}{2\Gamma},\quad
d_2=\frac{z(z-1)^2(2Lz+z^2-1)}{2\Gamma},\notag\\
d_3&=-\frac{z(z-1)^3(3Lz+L+2z^2-2)}{2(z+1)\Gamma},\quad
d_4=\frac{(z-1)^4(2Lz+z^2-1)}{2L(z+1)\Gamma},\label{eq:certcoef}
\end{align}
and $c_1=b^2d_1$, $c_2=b^2d_2$, $c_3=b^3d_3$, $c_4=b^2d_4$.
Let $c_0=-c_1m+c_2\mu t_z-t_z/3$, and let
$S(u)=c_0+c_1u+c_2u\log u+c_3\log u+c_4u(\log u)^2$ for $u>0$. Define
the symmetric function
\begin{equation}\label{eq:certpsi}
\Psi_z(x,y,w)=S(x)+S(y)+S(w)+\mu c_3(xy+yw+wx)+xyw\big(1+\mu c_4\log(xyw)\big)
\end{equation}
for $(x,y,w)\in(0,1]^3$.

We now explain the form of $\Psi_z$. Let $U$ be a graphon with values in
$(0,1]$, let $X_1,X_2,X_3$ be as in Definition~\ref{def:trilaw}, and let
$X=U(X_1,X_2)$, $Y=U(X_2,X_3)$ and $Z=U(X_1,X_3)$. Up to the order of the
coordinates, $(X,Y,Z)$ is the triple of Definition~\ref{def:trilaw}. Since
$\nu_U$ is invariant under permutations of the coordinates, $(X,Y,Z)$ also
has law $\nu_U$. Expanding $\log(XYZ)$ in \eqref{eq:certpsi} gives
\begin{align}
\Psi_z(X,Y,Z)&=3c_0+c_1\sum X+c_2\sum X\log X+XYZ\notag\\
&\quad+\sum(c_3+c_4X\log X)(\log X+\mu YZ),\label{s8:rewrite}
\end{align}
where each sum runs over the three cyclic shifts of $(X,Y,Z)$. Since
$\E[YZ\mid X_1,X_2]=C_U(X_1,X_2)$, the conditional mean of the factor
$\log X+\mu YZ$ given $(X_1,X_2)$ is $\log U(X_1,X_2)+\mu C_U(X_1,X_2)$,
which is the residual of the Euler--Lagrange equation with multiplier $\mu$.
The last sum in \eqref{s8:rewrite} tests this residual against the functions
$1$ and $u\log u$. In the next lemma we compute the mean of $\Psi_z(X,Y,Z)$ when
$U$ satisfies an Euler--Lagrange equation. In particular, the constant $c_0$
makes this mean vanish when $U=B_z$. Hence, if $\Psi_z$ is nonnegative, then
$\Psi_z$ and its gradient must vanish at the edge value triples $(a,a,a)$
and $(a,b,b)$ of $B_z$ and at their permutations. The coefficients
\eqref{eq:certcoef} are the unique solution of these tangency conditions
(Lemma~\ref{lem:certcoef}).

\begin{lemma}\label{s8:expect}
Let $z\in(0,1)$ with $\Gamma\ne0$. Let $U$ be a graphon with $t(U)=t_z$.
Assume that $U\ge\ell$ a.e.\ for some $\ell>0$ and that
$\log U=-\tilde\mu C_U$ a.e.\ for some $\tilde\mu\in\R$. Let $X,Y,Z$ be as
in \eqref{s8:rewrite}. Then
\begin{equation}\label{s8:expectid}
\E[\Psi_z(X,Y,Z)]=3c_1\big(e(U)-m\big)+3(\mu-\tilde\mu)\Big[c_2t_z+c_3\iint C_U+c_4\iint U\log U C_U\Big].
\end{equation}
In particular $\Psi_z(a,a,a)+3\Psi_z(a,b,b)=0$.
\end{lemma}

\begin{proof}
Since $\ell\le U\le1$ a.e., every term of \eqref{s8:rewrite} is bounded. We
have $\E[X]=e(U)$ and $\E[XYZ]=t_z$. By Lemma~\ref{lem:trilinear},
$t(U)=\iint UC_U$, so
$\E[X\log X]=\iint U\log U=-\tilde\mu\iint UC_U=-\tilde\mu t_z$. Let
$\varphi(u)=c_3+c_4u\log u$. Since $X$ is a function of $(X_1,X_2)$,
$\log X=-\tilde\mu C_U(X_1,X_2)$ a.s.\ and
$\E[YZ\mid X_1,X_2]=C_U(X_1,X_2)$, we get
\begin{align*}
\E[\varphi(X)(\log X+\mu YZ)]&=(\mu-\tilde\mu)\E[\varphi(X)C_U(X_1,X_2)]\\
&=(\mu-\tilde\mu)\Big[c_3\iint C_U+c_4\iint U\log U C_U\Big].
\end{align*}
Since the law of $(X,Y,Z)$ is invariant under permutations, the three cyclic
terms in \eqref{s8:rewrite} have the same expectation. Summing these
expectations and using $3c_0+t_z=-3c_1m+3c_2\mu t_z$ gives
\eqref{s8:expectid}.

For the last claim, let $U=B_z$. Then $t(B_z)=t_z$, $e(B_z)=m$, $B_z\ge a>0$
and $\log B_z=-\mu C_{B_z}$ by Lemma~\ref{lem:branch}. Hence the hypotheses
hold with $\tilde\mu=\mu$, so the right side of \eqref{s8:expectid}
vanishes. Since $\Psi_z$ is symmetric, Lemma~\ref{lem:patterns}(a) shows that
the left side equals $\frac14\Psi_z(a,a,a)+\frac34\Psi_z(a,b,b)$.
\end{proof}

We next verify the tangency conditions and determine the signs of
the coefficients.

\begin{lemma}\label{lem:certcoef}
Let $0<z\le1/40$. Then $\Gamma<0$, $c_1,c_2,c_4>0$, $c_3<0$, and
\[
\bar c:=c_1+c_2+\frac{c_3m^2}{t_z}<0 .
\]
Moreover $\Psi_z$ and its gradient vanish at $(a,a,a)$ and at the three
permutations of $(a,b,b)$.
\end{lemma}

\begin{proof}
\noindent\textbf{Tangency.} This part holds for every $z\in(0,1)$ with
$\Gamma\ne0$. Formula \eqref{eq:certpsi} defines a symmetric function on
$(0,\infty)^3$. Hence, by Lemma~\ref{lem:patterns}(b), its three partial
derivatives coincide at $(a,a,a)$, and $\partial_y\Psi_z=\partial_w\Psi_z$ at
$(a,b,b)$. The other permutations of $(a,b,b)$ follow by symmetry, and
$\Psi_z(a,a,a)+3\Psi_z(a,b,b)=0$ by Lemma~\ref{s8:expect}. It therefore
suffices to show that the four quantities $\Psi_z(a,a,a)$,
$\partial_x\Psi_z(a,a,a)$, $\partial_x\Psi_z(a,b,b)$ and
$\partial_y\Psi_z(a,b,b)$ vanish. Let $M=\mu b^2$. By \eqref{eq:branch},
\begin{equation}\label{s8:subs}
\log a=-\frac{L(1+z^2)}{(1-z)^2},\quad \log b=-\frac{2Lz}{(1-z)^2},\quad
M=\frac{2L}{(1-z)^2},
\end{equation}
and $a=zb$, $m=b(1+z)/2$, $t_z=b^3\kappa(z)^3$. We divide the first of the
four quantities by $3b^3$ and the other three by $b^2$. By
\eqref{eq:certpsi} and \eqref{s8:subs}, setting the results equal to zero
gives a linear system in $(d_1,d_2,d_3,d_4)$ whose coefficients and right
sides are rational functions of $z$ and $L$. Its determinant equals
$-L(z+1)\Gamma/(z(z-1)^2)$, which is nonzero. Its unique solution is
\eqref{eq:certcoef}. Both claims are identities between rational functions
of the independent variables $z$ and $L$, which we check by clearing
denominators.

\noindent\textbf{Signs.} Let $0<z\le1/40$. Then $L\ge\log40>3.68$. Since
$u\mapsto u\log(1/u)$ is increasing on $(0,1/e)$, we also have
$Lz\le(\log40)/40<0.09225$. Dropping the negative terms
$-Lz^2-Lz-z^3-3z$ of $\Gamma$ gives
\begin{align*}
\Gamma&\le1-L+6(Lz)^2+3(Lz)z^2+3z^2\\
&<1-3.68+6(0.09225)^2+\frac{3(0.09225)}{1600}+\frac3{1600}<-2.6 .
\end{align*}
Moreover $2Lz+z^2-1<0.19-1<0$. The factors $L+z-1$, $2L+z-1$ and
$3Lz+L+2z^2-2>L-2$ are positive. Reading off the signs of the factors in
\eqref{eq:certcoef} gives $d_1,d_2,d_4>0$ and $d_3<0$. Since $b>0$, the $c_i$
have the same signs.

For $\bar c$, note that $c_3m^2/t_z=b^2d_3(1+z)^2/(4\kappa(z)^3)$. A rational identity
in $z$ and $L$ gives
\[
\frac{\bar c}{b^2}=d_1+d_2+\frac{(1+z)^2d_3}{z^3+3z}=\frac{(z-1)P}{2(z^2+3)\Gamma},
\]
where
\begin{align*}
P&=2L^2z(z+1)(z^2+3)+L(5z^5-5z^4+14z^3-2z^2-11z-1)\\
&\quad+2z^6-4z^5+6z^4-10z^2+4z+2 .
\end{align*}
Since $(z-1)/\Gamma>0$, it suffices to show $P<0$. Let
$\chi(z)=(1+z)(1+z^2/3)$ and $g(z)=6Lz\chi(z)-11z+14z^3+5z^5$, with
$L=\log(1/z)$. Then $2L^2z(z+1)(z^2+3)=6L^2z\chi(z)$. Dropping the negative
terms $-L(5z^4+2z^2)-4z^5-10z^2$ of $P$ gives
\[
P\le L\big(g(z)-1\big)+2+4z+6z^4+2z^6 .
\]
Since $\dd L/\dd z=-1/z$, $\chi\ge1$ and $\chi'>0$,
\[
g'(z)=6\big((L-1)\chi(z)+Lz\chi'(z)\big)-11+42z^2+25z^4\ge6(L-1)-11>0.
\]
Thus $g(z)\le g(1/40)<\frac{6(3.69)}{40}\cdot\frac{41}{40}\big(1+\frac1{4800}\big)
-\frac{11}{40}+\frac{14}{40^3}+\frac5{40^5}<0.293$. Hence
$P<-0.707L+2.101<-0.707(3.68)+2.101<-0.5$. This proves $\bar c<0$.
\end{proof}

Using the next lemma, we reduce
Proposition~\ref{prop:family} to a compact domain in logarithmic
coordinates. For $q\in\R^3$ let
$x_i=e^{-q_i}$, $s=q_1+q_2+q_3$ and $A=\mu c_4$. Define
$\Phi_z(q)=\Psi_z(e^{-q_1},e^{-q_2},e^{-q_3})$, that is,
\begin{equation}\label{s8:phi}
\Phi_z(q)=\sum_{i=1}^3\big(c_0+c_1x_i-c_2q_ix_i-c_3q_i+c_4q_i^2x_i\big)+\mu c_3\sum_{i<j}x_ix_j+e^{-s}(1-As).
\end{equation}
Differentiation gives
\begin{equation}\label{s8:gradq}
\partial_{q_i}\Phi_z(q)=-c_3+x_i\big(-c_1+c_2(q_i-1)+c_4(2q_i-q_i^2)\big)-\mu c_3x_i\sum_{j\ne i}x_j+e^{-s}(As-1-A).
\end{equation}
The map $q\mapsto(e^{-q_1},e^{-q_2},e^{-q_3})$ is a bijection from
$[0,\infty)^3$ onto $(0,1]^3$ that commutes with permutations. By the chain
rule, $\partial_{q_i}\Phi_z=-x_i\partial_{x_i}\Psi_z$. Let $q_a=-\log a$ and
$q_b=-\log b$.

\begin{lemma}\label{lem:certtail}
Let $0<z\le1/40$, $q\in[0,\infty)^3$ and $i\in\{1,2,3\}$ with
$q_i\ge q_a+2$. Then $\partial_{q_i}\Phi_z(q)>0$.
\end{lemma}

\begin{proof}
By Lemma~\ref{lem:certcoef}, $c_1,c_2,c_4>0$, $c_3<0$ and $A>0$. If
$q_i\ge\max\{1,1+1/A,2+\sqrt2\}$, then in \eqref{s8:gradq} we have
$c_2x_i(q_i-1)\ge0$, $-\mu c_3x_i\sum_{j\ne i}x_j\ge0$, and
$As-1-A\ge Aq_i-1-A\ge0$. Hence
\[
\partial_{q_i}\Phi_z(q)\ge|c_3|-\gamma(q_i),\quad
\gamma(u)=e^{-u}\big(c_1+c_4(u^2-2u)\big).
\]
Moreover, $\gamma'(u)=-e^{-u}(c_1+c_4(u^2-4u+2))<0$ for $u\ge2+\sqrt2$. It
therefore suffices to prove (i) $q_a+2\ge\max\{1,1+1/A,2+\sqrt2\}$ and (ii)
$\gamma(q_a+2)<|c_3|$.

(i) By \eqref{s8:subs}, $q_a=L(1+z^2)/(1-z)^2>L>3.68$, so
$q_a+2>2+\sqrt2$. Since $A=\mu c_4=Md_4$, a rational identity in $z$ and $L$
gives
\[
A(q_a+1)-1=-\frac{2N}{(z+1)\Gamma},
\]
where $N=2L^2z^3+3L^2z^2-L^2z+Lz^4+Lz^2-2Lz-z^4+2z^3-2z+1$.
Since the function $z\mapsto zL^2$ is increasing when $L>2$, we have
$zL^2\le(\log40)^2/40<3.69^2/40$. Dropping positive terms,
\begin{align*}
N&\ge1-L^2z-2Lz-z^4-2z\\
&>1-\frac{3.69^2}{40}-2(0.09225)-\frac1{40^4}-\frac2{40}>0.425 .
\end{align*}
Since $\Gamma<0$, we get $A(q_a+1)>1$, that is, $q_a+2>1+1/A$.

(ii) Since $e^{-q_a}=a=zb$ and $(q_a+2)^2-2(q_a+2)=q_a(q_a+2)$,
\[
\gamma(q_a+2)=e^{-2}zb^3\big(d_1+d_4q_a(q_a+2)\big),\quad |c_3|=-b^3d_3 .
\]
A rational identity in $z$ and $L$ gives
\[
4(-d_3)-z\big(d_1+d_4q_a(q_a+2)\big)=-\frac{zT}{2(z+1)\Gamma},
\]
where
\begin{align*}
T&=2L^2z^5+2L^2z^4+6L^2z^3-2L^2z^2\\
&\quad+Lz^6+7Lz^5-19Lz^4+34Lz^3-33Lz^2+7Lz+3L\\
&\quad+3z^6-13z^5+24z^4-14z^3-17z^2+27z-10 .
\end{align*}
Dropping all positive terms except $3L$,
\begin{align*}
T&\ge3L-10-2(Lz)^2-19(Lz)z^3-33(Lz)z-13z^5-14z^3-17z^2\\
&>3(3.68)-10-2(0.09225)^2-\frac{19(0.09225)}{40^3}-\frac{33(0.09225)}{40}\\
&\quad-\frac{13}{40^5}-\frac{14}{40^3}-\frac{17}{40^2}>0.93 .
\end{align*}
Since $\Gamma<0$, it follows that $z(d_1+d_4q_a(q_a+2))<4(-d_3)$. Hence
$\gamma(q_a+2)<4e^{-2}b^3(-d_3)<|c_3|$, where the last step uses $e^2>4$ and
$-d_3>0$.
\end{proof}

\begin{proposition}[computer assisted]\label{prop:family}
Let $z\in[e^{-9.5},e^{-3.9}]$, that is, $L\in[3.9,9.5]$. Then every
$(x,y,w)\in(0,1]^3$ satisfies $\Psi_z(x,y,w)\ge0$. Equality holds if and
only if $(x,y,w)=(a,a,a)$ or $(x,y,w)$ is a permutation of $(a,b,b)$.
\end{proposition}

\begin{proof}
Since $e^{-3.9}<1/40$, Lemmas~\ref{lem:certcoef} and~\ref{lem:certtail}
apply. As $Le^{-L}$ and $(1-e^{-L})^{-2}$ are decreasing for $L>1$,
\eqref{s8:subs} gives
\[
q_a+2=L+2+\frac{2Le^{-L}}{(1-e^{-L})^2}\le9.5+2+\frac{2(3.9)e^{-3.9}}{(1-e^{-3.9})^2}<11.67<12 .
\]
Let $q\in[0,\infty)^3$, and suppose that $q_i>q_a+2$ for some $i$. By
Lemma~\ref{lem:certtail} applied along the segment, replacing $q_i$ by
$q_a+2$ strictly decreases $\Phi_z$. After at most three such replacements we
reach $q'\in[0,q_a+2]^3\subset[0,12]^3$ with $\Phi_z(q)>\Phi_z(q')$ unless
$q=q'$. It therefore suffices to show that $\Phi_z\ge0$ on $[0,12]^3$ with
equality exactly at the permutations of $(q_a,q_a,q_a)$ and $(q_a,q_b,q_b)$.
By symmetry it suffices to consider $12\ge q_1\ge q_2\ge q_3\ge0$. This set
contains both points since $q_a>q_b$. The statement on this set is
Proposition~\ref{appA3:compact}, which we verify in
Appendix~\ref{app:family}.
\end{proof}

By Lemmas~\ref{lem:penconstants} and~\ref{lem:digits},
$e^{-9.5}<7.5\cdot10^{-5}<z_c<z_*<0.0201<e^{-3.9}$. Since $r_-=\rho(z_c)$ and
$z_c<z_*\le z_{\mathrm{sp}}$, we have $z_1(r_-)=z_c$. Recall also that
$z_1(r_*)=z_*$. As $z_1$ is increasing and $r_-<r_*$
(Lemma~\ref{lem:order}), Proposition~\ref{prop:family} applies to $z=z_r$
for every $r\in[r_-,r_*]$.

\subsection{Minimizers with the cost of $B_r$}

In the next lemma we apply the certificate $\Psi_z$ with $z=z_r$.

\begin{lemma}\label{lem:tieexclusion}
Let $r\in[r_-,r_*]$, and let $U$ be a minimizer at $r$ with $F(U)=F(B_r)$ and
$e(U)\ge m_r$. Then $U$ is a relabeling of $B_r$.
\end{lemma}

\begin{proof}
Let $z=z_r$. In the notation of \eqref{eq:certpsi}, we then have $a=a_r$,
$b=b_r$, $\mu=\mu_r$, $m=m_r$ and $t_z=r^3$. Moreover, $L\in[3.9,9.5]$, as
shown at the end of Section~\ref{s8:sub:cert}. Hence $z<1/40$, and
Lemma~\ref{lem:certcoef} and Proposition~\ref{prop:family} apply to $z$. Let $\Delta=e(U)-m\ge0$. By Lemma~\ref{lem:active},
Proposition~\ref{prop:euler} and Lemma~\ref{lem:support}, $t(U)=t_z$,
$\log U=-\mu_U C_U$ a.e.\ and $U\ge\ell_U>0$ a.e. By \eqref{eq:mult} and
Lemma~\ref{lem:Bderivative}(a), $\mu_U t_z=1-F(U)-e(U)$ and
$\mu t_z=1-F(B_r)-m$. Since $F(U)=F(B_r)$, subtracting gives
$\mu-\mu_U=\Delta/t_z$. Lemma~\ref{s8:expect} with $\tilde\mu=\mu_U$ then
gives
\begin{equation}\label{s8:tie}
\E[\Psi_z(X,Y,Z)]=3\Delta\Big[c_1+c_2+\frac{c_3\iint C_U+c_4\iint U\log U C_U}{t_z}\Big].
\end{equation}
By Fubini's theorem and the Cauchy--Schwarz inequality,
$\iint C_U=\int_0^1d_U^2\ge e(U)^2\ge m^2$. Since $c_3<0$ by
Lemma~\ref{lem:certcoef}, it follows that $c_3\iint C_U\le c_3m^2$. Moreover
$U\log U C_U\le0$ and $c_4>0$. Hence the bracket in \eqref{s8:tie} is at most
$\bar c<0$. Since $\Delta\ge0$, we get $\E[\Psi_z(X,Y,Z)]\le3\Delta\bar c\le0$.
On the other hand, $(X,Y,Z)\in(0,1]^3$ a.s., so $\E[\Psi_z(X,Y,Z)]\ge0$ by
Proposition~\ref{prop:family}. Therefore $\Delta=0$ and $\Psi_z(X,Y,Z)=0$ a.s.

By Proposition~\ref{prop:family}, $(X,Y,Z)$ is a.s.\ equal to $(a,a,a)$ or
to a permutation of $(a,b,b)$. Hence its law $\nu_U$ is carried by these four
points. Since $e(U)=m=(a+b)/2$ and $a\ne b$, Lemma~\ref{lem:patterns}(c)
shows that $U$ is a relabeling of $B_{a,b}=B_r$.
\end{proof}

We next show that the hypotheses of Lemma~\ref{lem:tieexclusion} hold at the lower end of an interval on which $V(r)=F(B_r)$.

\begin{lemma}\label{lem:boundary}
Let $0<r_0<r_1\le r_*$ be such that $V(r)=F(B_r)$ for every
$r\in(r_0,r_1)$. Then every minimizer $U$ at $r_0$ satisfies
$F(U)=F(B_{r_0})$ and $e(U)\ge m_{r_0}$.
\end{lemma}

\begin{proof}
By Lemmas~\ref{lem:existence} and~\ref{lem:Bderivative}, $V$ and
$r\mapsto F(B_r)$ are continuous. Hence $F(U)=V(r_0)=F(B_{r_0})$.

For $s>r_0$ let $W_s=(s/r_0)U$. By Lemma~\ref{lem:active}, $t(W_s)=s^3$.
Lemma~\ref{lem:truncation} applies to every nonnegative kernel, so it gives
$F(W_s)\ge V(s)$. Hence $\phi(s)=F(W_s)-F(B_s)$ satisfies $\phi(s)\ge0$ for
$s\in(r_0,r_1)$ and $\phi(r_0)=0$. Thus the right derivative of $\phi$ at
$r_0$ is nonnegative. By \eqref{eq:scale1} with $\lambda=s/r_0$, we have
$F(W_s)=\lambda F(U)+(\lambda\log\lambda)e(U)+1-\lambda$. The derivative of
$F(W_s)$ in $s$ at $s=r_0$ is $(F(U)+e(U)-1)/r_0=-\mu_Ur_0^2$ by
\eqref{eq:mult}. By Lemma~\ref{lem:Bderivative}(b), the derivative of
$F(B_s)$ at $s=r_0$ is $-\mu_{r_0}r_0^2$. Therefore $\mu_{r_0}\ge\mu_U$. Finally, by
\eqref{eq:mult}, Lemma~\ref{lem:Bderivative}(a) and $F(U)=F(B_{r_0})$,
\[
\mu_{r_0}-\mu_U=\frac{(1-F(B_{r_0})-m_{r_0})-(1-F(U)-e(U))}{r_0^3}
=\frac{e(U)-m_{r_0}}{r_0^3},
\]
so $e(U)\ge m_{r_0}$.
\end{proof}

\subsection{Proof of Theorem~\ref{thm:main}}

\begin{proof}[Proof of Theorem~\ref{thm:main}]
By Corollary~\ref{cor:Vbal}, the number $r_*$ of Theorem~\ref{thm:main} is
the constant of Theorem~\ref{thm:crossing}. Lemma~\ref{lem:digits} gives its
digits and those of $a_*$ and $b_*$. Parts (a) and (c) are
Theorem~\ref{thm:threshold}(i) and (ii), respectively.

We next show that $G=(0,r_*)$. Recall from the beginning of
Section~\ref{sec:intermediate} that $G$ is open
(Corollary~\ref{cor:open}), that $(0,r_-]\subset G$
(Theorem~\ref{thm:smalltail}), and that $(r_*-\delta,r_*)\subset G$ for some
$\delta\in(0,r_*)$ (Corollary~\ref{cor:seed}). Let $r_0$ be the infimum of
the set $\mathcal S$ of all $s\in(0,r_*)$ with $(s,r_*)\subset G$. Then
$r_0\le r_*-\delta$ and $(r_0,r_*)=\bigcup_{s\in\mathcal S}(s,r_*)\subset G$.

Suppose that $r_0\ge r_-$. We show that the endpoint $r_0$ lies in $G$. For
$r\in(r_0,r_*)$ the graphon $B_r$ is a minimizer, so $V(r)=F(B_r)$. By
Lemma~\ref{lem:boundary} with $r_1=r_*$, every minimizer $U$ at $r_0$
satisfies $F(U)=F(B_{r_0})$ and $e(U)\ge m_{r_0}$. Lemma~\ref{lem:tieexclusion}
then shows that $U$ is a relabeling of $B_{r_0}$. Since minimizers exist and
relabeling preserves $F$ and $t$, we get $r_0\in G$. As $G$ is open,
$(r_0-\eta,r_*)\subset G$ for some $\eta\in(0,r_0)$. Then
$r_0-\eta\in\mathcal S$, which contradicts the definition of $r_0$. Hence
$r_0<r_-$, and $G\supset(0,r_-]\cup(r_0,r_*)=(0,r_*)$.

Now let $r\in(0,r_*)$. Since $r\in G$, the minimizers at $r$ are exactly the
relabelings of $B_r$, so $V(r)=F(B_r)$. By Corollary~\ref{cor:Vbal}(a),
$F(B_r)=\Vbal(r)$ and $(a_r,b_r)$ is the only minimizing pair at $r$.
Moreover, $0<a_r<r<b_r<1$ by \eqref{s4:eq:arbr}. This proves (b).
\end{proof}
\section{Open problems}\label{sec:concluding}

We conclude with three questions that our results leave open.

\medskip
\noindent\textbf{Sandwiching.} Theorem~\ref{thm:structure} describes the
conditioned graph only up to a small error in cut norm. For $\eta>\eta_*$,
can it be coupled so that, with high probability, it contains
$G(n,(1-\varepsilon)\eta^{1/3}p)$ and is contained in
$G(n,(1+\varepsilon)\eta^{1/3}p)$? Chin~\cite{Chi25} asked this question for
$\eta^{1/3}>0.4664\ldots$. For $\eta<\eta_*$ one can ask the same with two block SBMs in
place of the Erdős–Rényi graphs.

\medskip
\noindent\textbf{Other graphs.} For a fixed graph $H$ and
$n^{-1/m_2(H)}\ll p\ll1$, where $m_2(H)$ is the maximum $2$-density of $H$,
the lower tail for the number of copies of $H$ is governed by the analogue
$V_H$ of \eqref{eq:V}, with $t(U)$ replaced by the homomorphism density
$t(H,U)$~\cite{KS23,Chi25,Chi26}. If $H$ is bipartite and satisfies Sidorenko's
inequality, Jensen's inequality shows that the constant graphon is always
the unique minimizer. For nonbipartite $H$ there is a single
threshold~\cite{Chi26}, but its value is unknown, even for odd cycles. For $H$ of chromatic number
$k$, a natural guess is that below the threshold the minimizers are balanced
SBMs with $k-1$ blocks. This is consistent with the results of Janson and
Warnke~\cite{JW16} as $\eta\to0$, and for small targets it follows from the
work of Kwan and Pham~\cite{KP26} when $H$ is strictly $2$-balanced and edge
critical, as Chin~\cite{Chi26} observed.

\medskip
\noindent\textbf{The dense case.} For fixed $p$, the phase diagram of the
dense problem $V_p$ of Section~\ref{s1:sub:related} is open between the two
curves of Zhao~\cite{Zha17}, and Theorem~\ref{thm:main} suggests
that its boundary has slope $r_*$ at the origin. The regime $p\to1$ may be
especially rich. There $V_p$ approaches the problem of keeping as many
edges as possible under a bound on the triangle density, which
Razborov~\cite{Raz08} solved with complete multipartite graphons, and whose
extremal graphs were described by Pikhurko and Razborov~\cite{PR17} and by
Liu, Pikhurko and Staden~\cite{LPS20}. This
suggests many phase transitions in $\eta$.

\section*{Acknowledgements}
We thank Byron Chin and Sam van der Poel for a helpful
conversation about future directions and for sharing their perspective on
these problems.
This work used the Delta CPU system at NCSA through allocation CIS261417
from the ACCESS program, which is supported by U.S. National Science
Foundation grants \#2138259, \#2138286, \#2138307, \#2137603, and
\#2138296. We also thank the Yale Center for Research Computing for the use
of the Bouchet cluster.

\bibliographystyle{amsplain0}
\bibliography{main}

@article{Zha17,
  author  = {Zhao, Yufei},
  title   = {On the lower tail variational problem for random graphs},
  journal = {Combinatorics, Probability and Computing},
  volume  = {26},
  number  = {2},
  pages   = {301--320},
  year    = {2017}
}

@article{CV11,
  author  = {Chatterjee, Sourav and Varadhan, S. R. S.},
  title   = {The large deviation principle for the {E}rd{\H{o}}s--{R}{\'e}nyi random graph},
  journal = {European Journal of Combinatorics},
  volume  = {32},
  number  = {7},
  pages   = {1000--1017},
  year    = {2011}
}

@article{CD16,
  author  = {Chatterjee, Sourav and Dembo, Amir},
  title   = {Nonlinear large deviations},
  journal = {Advances in Mathematics},
  volume  = {299},
  pages   = {396--450},
  year    = {2016}
}

@article{CD20,
  author  = {Cook, Nicholas A. and Dembo, Amir},
  title   = {Large deviations of subgraph counts for sparse {E}rd{\H{o}}s--{R}{\'e}nyi graphs},
  journal = {Advances in Mathematics},
  volume  = {373},
  pages   = {107289},
  year    = {2020}
}

@article{CDP,
  author  = {Cook, Nicholas A. and Dembo, Amir and Pham, Huy Tuan},
  title   = {Regularity method and large deviation principles for the {E}rd{\H{o}}s--{R}{\'e}nyi hypergraph},
  journal = {Duke Mathematical Journal},
  volume  = {173},
  number  = {5},
  pages   = {873--946},
  year    = {2024}
}

@article{KS23,
  author  = {Kozma, Gady and Samotij, Wojciech},
  title   = {Lower tails via relative entropy},
  journal = {The Annals of Probability},
  volume  = {51},
  number  = {2},
  pages   = {665--698},
  year    = {2023}
}

@article{JW16,
  author  = {Janson, Svante and Warnke, Lutz},
  title   = {The lower tail: {P}oisson approximation revisited},
  journal = {Random Structures \& Algorithms},
  volume  = {48},
  number  = {2},
  pages   = {219--246},
  year    = {2016}
}

@article{Chi25,
  author  = {Chin, Byron},
  title   = {Structure of lower tails in sparse random graphs},
  journal = {Random Structures \& Algorithms},
  volume  = {67},
  number  = {1},
  pages   = {e70028},
  year    = {2025}
}

@misc{JPPS24,
  author  = {Jenssen, Matthew and Perkins, Will and Potukuchi, Aditya and Simkin, Michael},
  title   = {Lower tails for triangles inside the critical window},
  year    = {2024},
  note    = {Preprint, arXiv:2411.18563}
}

@misc{KP26,
  author  = {Kwan, Matthew and Pham, Huy Tuan},
  title   = {Nonlinear lower-tail large deviations at criticality},
  year    = {2026},
  note    = {Preprint, arXiv:2610.10233}
}

@article{LZ15,
  author  = {Lubetzky, Eyal and Zhao, Yufei},
  title   = {On replica symmetry of large deviations in random graphs},
  journal = {Random Structures \& Algorithms},
  volume  = {47},
  number  = {1},
  pages   = {109--146},
  year    = {2015}
}

@article{LZ17,
  author  = {Lubetzky, Eyal and Zhao, Yufei},
  title   = {On the variational problem for upper tails in sparse random graphs},
  journal = {Random Structures \& Algorithms},
  volume  = {50},
  number  = {3},
  pages   = {420--436},
  year    = {2017}
}

@article{BGLZ17,
  author  = {Bhattacharya, Bhaswar B. and Ganguly, Shirshendu and Lubetzky, Eyal and Zhao, Yufei},
  title   = {Upper tails and independence polynomials in random graphs},
  journal = {Advances in Mathematics},
  volume  = {319},
  pages   = {313--347},
  year    = {2017}
}

@article{RS13,
  author  = {Radin, Charles and Sadun, Lorenzo},
  title   = {Phase transitions in a complex network},
  journal = {Journal of Physics A: Mathematical and Theoretical},
  volume  = {46},
  number  = {30},
  pages   = {305002},
  year    = {2013}
}

@article{KRRS17,
  author  = {Kenyon, Richard and Radin, Charles and Ren, Kui and Sadun, Lorenzo},
  title   = {The phases of large networks with edge and triangle constraints},
  journal = {Journal of Physics A: Mathematical and Theoretical},
  volume  = {50},
  number  = {43},
  pages   = {435001},
  year    = {2017}
}

@article{NRS24,
  author  = {Neeman, Joe and Radin, Charles and Sadun, Lorenzo},
  title   = {Existence of a symmetric bipodal phase in the edge-triangle model},
  journal = {Journal of Physics A: Mathematical and Theoretical},
  volume  = {57},
  pages   = {095003},
  year    = {2024}
}

@book{Lov12,
  author    = {Lov{\'a}sz, L{\'a}szl{\'o}},
  title     = {Large Networks and Graph Limits},
  series    = {American Mathematical Society Colloquium Publications},
  volume    = {60},
  publisher = {American Mathematical Society},
  address   = {Providence, RI},
  year      = {2012}
}

@article{LS06,
  author  = {Lov{\'a}sz, L{\'a}szl{\'o} and Szegedy, Bal{\'a}zs},
  title   = {Limits of dense graph sequences},
  journal = {Journal of Combinatorial Theory, Series B},
  volume  = {96},
  number  = {6},
  pages   = {933--957},
  year    = {2006}
}

@article{LS07,
  author  = {Lov{\'a}sz, L{\'a}szl{\'o} and Szegedy, Bal{\'a}zs},
  title   = {Szemer{\'e}di's lemma for the analyst},
  journal = {Geometric and Functional Analysis},
  volume  = {17},
  number  = {1},
  pages   = {252--270},
  year    = {2007}
}

@article{Goo59,
  author  = {Goodman, A. W.},
  title   = {On sets of acquaintances and strangers at any party},
  journal = {The American Mathematical Monthly},
  volume  = {66},
  number  = {9},
  pages   = {778--783},
  year    = {1959}
}

@article{Raz07,
  author  = {Razborov, Alexander A.},
  title   = {Flag algebras},
  journal = {The Journal of Symbolic Logic},
  volume  = {72},
  number  = {4},
  pages   = {1239--1282},
  year    = {2007}
}

@article{Raz08,
  author  = {Razborov, Alexander A.},
  title   = {On the minimal density of triangles in graphs},
  journal = {Combinatorics, Probability and Computing},
  volume  = {17},
  number  = {4},
  pages   = {603--618},
  year    = {2008}
}

@book{MKC09,
  author    = {Moore, Ramon E. and Kearfott, R. Baker and Cloud, Michael J.},
  title     = {Introduction to Interval Analysis},
  publisher = {Society for Industrial and Applied Mathematics},
  address   = {Philadelphia, PA},
  year      = {2009}
}

@book{Tuc11,
  author    = {Tucker, Warwick},
  title     = {Validated Numerics: A Short Introduction to Rigorous Computations},
  publisher = {Princeton University Press},
  address   = {Princeton, NJ},
  year      = {2011}
}

@article{Joh17,
  author  = {Johansson, Fredrik},
  title   = {Arb: efficient arbitrary-precision midpoint-radius interval arithmetic},
  journal = {IEEE Transactions on Computers},
  volume  = {66},
  number  = {8},
  pages   = {1281--1292},
  year    = {2017}
}

@article{BCLSV08,
  author  = {Borgs, Christian and Chayes, Jennifer T. and Lov{\'a}sz, L{\'a}szl{\'o} and S{\'o}s, Vera T. and Vesztergombi, Katalin},
  title   = {Convergent sequences of dense graphs {I}: {S}ubgraph frequencies, metric properties and testing},
  journal = {Advances in Mathematics},
  volume  = {219},
  number  = {6},
  pages   = {1801--1851},
  year    = {2008}
}

@book{Cha17,
  author    = {Chatterjee, Sourav},
  title     = {Large Deviations for Random Graphs},
  series    = {Lecture Notes in Mathematics},
  volume    = {2197},
  note      = {Lecture notes from the 45th Probability Summer School held in Saint-Flour, June 2015},
  publisher = {Springer, Cham},
  year      = {2017}
}

@article{HMS22,
  author  = {Harel, Matan and Mousset, Frank and Samotij, Wojciech},
  title   = {Upper tails via high moments and entropic stability},
  journal = {Duke Mathematical Journal},
  volume  = {171},
  number  = {10},
  pages   = {2089--2192},
  year    = {2022}
}

@misc{CD24,
  author  = {Cook, Nicholas A. and Dembo, Amir},
  title   = {Typical structure of sparse exponential random graph models},
  year    = {2024},
  note    = {Preprint, arXiv:2208.06397}
}

@article{Eld18,
  author  = {Eldan, Ronen},
  title   = {Gaussian-width gradient complexity, reverse log-{S}obolev inequalities and nonlinear large deviations},
  journal = {Geometric and Functional Analysis},
  volume  = {28},
  pages   = {1548--1596},
  year    = {2018}
}

@article{Aug20,
  author  = {Augeri, Fanny},
  title   = {Nonlinear large deviation bounds with applications to {W}igner matrices and sparse {E}rd{\H{o}}s--{R}{\'e}nyi graphs},
  journal = {The Annals of Probability},
  volume  = {48},
  pages   = {2404--2448},
  year    = {2020}
}

@article{BB23,
  author  = {Basak, Anirban and Basu, Riddhipratim},
  title   = {Upper tail large deviations of regular subgraph counts in {E}rd{\H{o}}s--{R}{\'e}nyi graphs in the full localized regime},
  journal = {Communications on Pure and Applied Mathematics},
  volume  = {76},
  pages   = {3--72},
  year    = {2023}
}

@article{Har60,
  author  = {Harris, T. E.},
  title   = {A lower bound for the critical probability in a certain percolation process},
  journal = {Proceedings of the Cambridge Philosophical Society},
  volume  = {56},
  pages   = {13--20},
  year    = {1960}
}

@incollection{JLR90,
  author    = {Janson, Svante and {\L}uczak, Tomasz and Ruci{\'n}ski, Andrzej},
  title     = {An exponential bound for the probability of nonexistence of a specified subgraph in a random graph},
  booktitle = {Random Graphs '87 ({P}ozna{\'n}, 1987)},
  editor    = {Karo{\'n}ski, Micha{\l} and Jaworski, Jerzy and Ruci{\'n}ski, Andrzej},
  publisher = {Wiley},
  address   = {Chichester},
  pages     = {73--87},
  year      = {1990}
}

@article{Jan90,
  author  = {Janson, Svante},
  title   = {{P}oisson approximation for large deviations},
  journal = {Random Structures \& Algorithms},
  volume  = {1},
  number  = {2},
  pages   = {221--229},
  year    = {1990}
}

@article{JR02,
  author  = {Janson, Svante and Ruci{\'n}ski, Andrzej},
  title   = {The infamous upper tail},
  journal = {Random Structures \& Algorithms},
  volume  = {20},
  number  = {3},
  pages   = {317--342},
  year    = {2002}
}

@article{JOR04,
  author  = {Janson, Svante and Oleszkiewicz, Krzysztof and Ruci{\'n}ski, Andrzej},
  title   = {Upper tails for subgraph counts in random graphs},
  journal = {Israel Journal of Mathematics},
  volume  = {142},
  number  = {1},
  pages   = {61--92},
  year    = {2004}
}

@article{KV04,
  author  = {Kim, Jeong Han and Vu, Van H.},
  title   = {Divide and conquer martingales and the number of triangles in a random graph},
  journal = {Random Structures \& Algorithms},
  volume  = {24},
  number  = {2},
  pages   = {166--174},
  year    = {2004}
}

@article{Cha12,
  author  = {Chatterjee, Sourav},
  title   = {The missing log in large deviations for triangle counts},
  journal = {Random Structures \& Algorithms},
  volume  = {40},
  number  = {4},
  pages   = {437--451},
  year    = {2012}
}

@article{DK12a,
  author  = {DeMarco, Bobby and Kahn, Jeff},
  title   = {Upper tails for triangles},
  journal = {Random Structures \& Algorithms},
  volume  = {40},
  number  = {4},
  pages   = {452--459},
  year    = {2012}
}

@article{DK12b,
  author  = {DeMarco, Bobby and Kahn, Jeff},
  title   = {Tight upper tail bounds for cliques},
  journal = {Random Structures \& Algorithms},
  volume  = {41},
  number  = {4},
  pages   = {469--487},
  year    = {2012}
}

@article{Cha16,
  author  = {Chatterjee, Sourav},
  title   = {An introduction to large deviations for random graphs},
  journal = {Bulletin of the American Mathematical Society},
  volume  = {53},
  number  = {4},
  pages   = {617--642},
  year    = {2016}
}

@article{MNPS20,
  author  = {Mousset, Frank and Noever, Andreas and Panagiotou, Konstantinos and Samotij, Wojciech},
  title   = {On the probability of nonexistence in binomial subsets},
  journal = {The Annals of Probability},
  volume  = {48},
  number  = {1},
  pages   = {493--525},
  year    = {2020}
}

@book{DZ10,
  author    = {Dembo, Amir and Zeitouni, Ofer},
  title     = {Large Deviations Techniques and Applications},
  series    = {Stochastic Modelling and Applied Probability},
  volume    = {38},
  publisher = {Springer-Verlag},
  address   = {Berlin},
  year      = {2010},
  note      = {Corrected reprint of the second (1998) edition}
}

@article{Var66,
  author  = {Varadhan, S. R. S.},
  title   = {Asymptotic probabilities and differential equations},
  journal = {Communications on Pure and Applied Mathematics},
  volume  = {19},
  number  = {3},
  pages   = {261--286},
  year    = {1966}
}

@article{Cra38,
  author  = {Cram{\'e}r, Harald},
  title   = {Sur un nouveau th{\'e}or{\`e}me-limite de la th{\'e}orie des probabilit{\'e}s},
  journal = {Actualit{\'e}s Scientifiques et Industrielles},
  volume  = {736},
  pages   = {5--23},
  year    = {1938}
}

@article{CDey10,
  author  = {Chatterjee, Sourav and Dey, Partha S.},
  title   = {Applications of {S}tein's method for concentration inequalities},
  journal = {The Annals of Probability},
  volume  = {38},
  number  = {6},
  pages   = {2443--2485},
  year    = {2010}
}

@article{FK99,
  author  = {Frieze, Alan and Kannan, Ravi},
  title   = {Quick approximation to matrices and applications},
  journal = {Combinatorica},
  volume  = {19},
  number  = {2},
  pages   = {175--220},
  year    = {1999}
}

@article{BCCZ19,
  author  = {Borgs, Christian and Chayes, Jennifer T. and Cohn, Henry and Zhao, Yufei},
  title   = {An {$L^p$} theory of sparse graph convergence {I}: {Limits}, sparse random graph models, and power law distributions},
  journal = {Transactions of the American Mathematical Society},
  volume  = {372},
  number  = {5},
  pages   = {3019--3062},
  year    = {2019}
}

@article{BCCZ18,
  author  = {Borgs, Christian and Chayes, Jennifer T. and Cohn, Henry and Zhao, Yufei},
  title   = {An {$L^p$} theory of sparse graph convergence {II}: {LD} convergence, quotients and right convergence},
  journal = {The Annals of Probability},
  volume  = {46},
  number  = {1},
  pages   = {337--396},
  year    = {2018}
}

@misc{JPPS26a,
  author  = {Jenssen, Matthew and Perkins, Will and Potukuchi, Aditya and Simkin, Michael},
  title   = {Non-existence probabilities and lower tails in the critical regime via {B}elief {P}ropagation},
  year    = {2026},
  note    = {Preprint, arXiv:2604.04652}
}

@misc{JPPS26b,
  author  = {Jenssen, Matthew and Perkins, Will and Potukuchi, Aditya and Simkin, Michael},
  title   = {The replica symmetric solution for hypergraph independent sets in the critical regime},
  year    = {2026},
  note    = {Preprint, arXiv:2610.10433}
}

@misc{Chi26,
  author  = {Chin, Byron},
  title   = {A first order phase transition for lower tails in sparse random graphs},
  year    = {2026},
  note    = {Preprint, arXiv:2610.27457}
}

@incollection{EKR76,
  author    = {Erd{\H{o}}s, Paul and Kleitman, Daniel J. and Rothschild, Bruce L.},
  title     = {Asymptotic enumeration of {$K_n$}-free graphs},
  booktitle = {Colloquio Internazionale sulle Teorie Combinatorie (Rome, 1973), Tomo II},
  series    = {Atti dei Convegni Lincei},
  number    = {17},
  publisher = {Accademia Nazionale dei Lincei},
  address   = {Rome},
  pages     = {19--27},
  year      = {1976}
}

@article{KPR87,
  author  = {Kolaitis, Phokion G. and Pr{\"o}mel, Hans J{\"u}rgen and Rothschild, Bruce L.},
  title   = {{$K_{\ell+1}$}-free graphs: asymptotic structure and a 0-1 law},
  journal = {Transactions of the American Mathematical Society},
  volume  = {303},
  number  = {2},
  pages   = {637--671},
  year    = {1987}
}

@article{PS96,
  author  = {Pr{\"o}mel, Hans J{\"u}rgen and Steger, Angelika},
  title   = {On the asymptotic structure of sparse triangle free graphs},
  journal = {Journal of Graph Theory},
  volume  = {21},
  number  = {2},
  pages   = {137--151},
  year    = {1996}
}

@article{Luc00,
  author  = {{\L}uczak, Tomasz},
  title   = {On triangle-free random graphs},
  journal = {Random Structures \& Algorithms},
  volume  = {16},
  number  = {3},
  pages   = {260--276},
  year    = {2000}
}

@article{OPT03,
  author  = {Osthus, Deryk and Pr{\"o}mel, Hans J{\"u}rgen and Taraz, Anusch},
  title   = {For which densities are random triangle-free graphs almost surely bipartite?},
  journal = {Combinatorica},
  volume  = {23},
  number  = {1},
  pages   = {105--150},
  year    = {2003}
}

@article{Ste05,
  author  = {Steger, Angelika},
  title   = {On the evolution of triangle-free graphs},
  journal = {Combinatorics, Probability and Computing},
  volume  = {14},
  number  = {1-2},
  pages   = {211--224},
  year    = {2005}
}

@article{BMSW16,
  author  = {Balogh, J{\'o}zsef and Morris, Robert and Samotij, Wojciech and Warnke, Lutz},
  title   = {The typical structure of sparse {$K_{r+1}$}-free graphs},
  journal = {Transactions of the American Mathematical Society},
  volume  = {368},
  number  = {9},
  pages   = {6439--6485},
  year    = {2016}
}

@article{JPP25,
  author  = {Jenssen, Matthew and Perkins, Will and Potukuchi, Aditya},
  title   = {On the evolution of structure in triangle-free graphs},
  journal = {Advances in Mathematics},
  volume  = {480},
  pages   = {Paper No. 110499},
  year    = {2025}
}

@article{MPW26,
  author  = {Mizgerd, Clayton and Perkins, Will and Wang, Yuzhou},
  title   = {On the chromatic number of random triangle-free graphs},
  journal = {Random Structures \& Algorithms},
  volume  = {68},
  number  = {4},
  pages   = {Paper No. e70081},
  year    = {2026}
}

@inproceedings{JPPS24S,
  author    = {Jenssen, Matthew and Perkins, Will and Potukuchi, Aditya and Simkin, Michael},
  title     = {Sampling, counting, and large deviations for triangle-free graphs near the critical density},
  booktitle = {2024 IEEE 65th Annual Symposium on Foundations of Computer Science (FOCS)},
  pages     = {151--165},
  year      = {2024}
}

@misc{LPY26,
  author  = {Lintott, Jade and Perkins, Will and Yap, Corrine},
  title   = {Local large deviations for triangles in sparse random graphs},
  year    = {2026},
  note    = {Preprint, arXiv:2609.21890}
}

@article{DK15,
  author  = {DeMarco, Bobby and Kahn, Jeff},
  title   = {Mantel's theorem for random graphs},
  journal = {Random Structures \& Algorithms},
  volume  = {47},
  number  = {1},
  pages   = {59--72},
  year    = {2015}
}

@article{CG16,
  author  = {Conlon, David and Gowers, W. Timothy},
  title   = {Combinatorial theorems in sparse random sets},
  journal = {Annals of Mathematics},
  volume  = {184},
  number  = {2},
  pages   = {367--454},
  year    = {2016}
}

@article{Sch16,
  author  = {Schacht, Mathias},
  title   = {Extremal results for random discrete structures},
  journal = {Annals of Mathematics},
  volume  = {184},
  number  = {2},
  pages   = {333--365},
  year    = {2016}
}

@article{BMS15,
  author  = {Balogh, J{\'o}zsef and Morris, Robert and Samotij, Wojciech},
  title   = {Independent sets in hypergraphs},
  journal = {Journal of the American Mathematical Society},
  volume  = {28},
  number  = {3},
  pages   = {669--709},
  year    = {2015}
}

@article{ST15,
  author  = {Saxton, David and Thomason, Andrew},
  title   = {Hypergraph containers},
  journal = {Inventiones Mathematicae},
  volume  = {201},
  number  = {3},
  pages   = {925--992},
  year    = {2015}
}

@article{Wor96,
  author  = {Wormald, Nicholas C.},
  title   = {The perturbation method and triangle-free random graphs},
  journal = {Random Structures \& Algorithms},
  volume  = {9},
  number  = {1-2},
  pages   = {253--269},
  year    = {1996}
}

@article{SW18,
  author  = {Stark, Dudley and Wormald, Nick},
  title   = {The probability of non-existence of a subgraph in a moderately sparse random graph},
  journal = {Combinatorics, Probability and Computing},
  volume  = {27},
  number  = {4},
  pages   = {672--715},
  year    = {2018}
}

@article{GGS20,
  author  = {Goldschmidt, Christina and Griffiths, Simon and Scott, Alex},
  title   = {Moderate deviations of subgraph counts in the {E}rd{\H{o}}s--{R}{\'e}nyi random graphs {$G(n,m)$} and {$G(n,p)$}},
  journal = {Transactions of the American Mathematical Society},
  volume  = {373},
  number  = {8},
  pages   = {5517--5585},
  year    = {2020}
}

@article{NRSm23,
  author  = {Neeman, Joe and Radin, Charles and Sadun, Lorenzo},
  title   = {Moderate deviations in cycle count},
  journal = {Random Structures \& Algorithms},
  volume  = {63},
  number  = {3},
  pages   = {779--820},
  year    = {2023}
}

@article{ADG25,
  author  = {Alvarado, Jos{\'e} D. and Dias do Couto, Gabriel and Griffiths, Simon},
  title   = {Moderate deviations of triangle counts in the {E}rd{\H{o}}s--{R}{\'e}nyi random graph {$G(n,m)$}: the lower tail},
  journal = {European Journal of Combinatorics},
  volume  = {129},
  pages   = {Paper No. 104189},
  year    = {2025}
}

@misc{PvdP25,
  author  = {Perkins, Will and van der Poel, Sam},
  title   = {The typical structure of dense claw-free graphs},
  year    = {2025},
  note    = {Preprint, arXiv:2501.17816}
}

@misc{DS25,
  author  = {DiCarlo, William and Sadun, Lorenzo},
  title   = {Tripodal structure in undersaturated random graphs},
  year    = {2025},
  note    = {Preprint, arXiv:2508.16797}
}

@article{PR17,
  author  = {Pikhurko, Oleg and Razborov, Alexander},
  title   = {Asymptotic structure of graphs with the minimum number of triangles},
  journal = {Combinatorics, Probability and Computing},
  volume  = {26},
  number  = {1},
  pages   = {138--160},
  year    = {2017}
}

@article{LPS20,
  author  = {Liu, Hong and Pikhurko, Oleg and Staden, Katherine},
  title   = {The exact minimum number of triangles in graphs of given order and size},
  journal = {Forum of Mathematics, Pi},
  volume  = {8},
  pages   = {e8},
  year    = {2020}
}

\appendix
\section{From the conditioned graph to the variational problem}\label{sec:reduction}


In this appendix we give a self-contained proof, for triangles, of the reduction described in Section~\ref{s1:sub:reduction}. Most of this reduction is already contained in previous work.
For fixed $p$, Chatterjee and Varadhan~\cite{CV11} reduced both questions to the dense problem $V_p$ of Section~\ref{s1:sub:related}. Zhao~\cite[Section~2.3]{Zha17} obtained
\eqref{eq:V} as the limit of $p^{-1}V_p(rp)$ as $p\to0$, after
$n\to\infty$. For $n^{-1/2}\ll p\ll1$, Kozma and Samotij~\cite[Theorem~2]{KS23} reduced
the rate to the discrete problem $\Phi_{n,p}$. Their proof reveals a small set of edges, after which the three edges of a
typical triangle are nearly independent.
Chin~\cite{Chi25} combined this method with a quantitative stability estimate for
the constant graphon. This gave the structure of the conditioned graph for
$r>0.4664\ldots$, where the constant graphon is the unique minimizer.
More recently, Chin~\cite[Corollary~2.6 and Theorem~1.4]{Chi26} proved the
reduction
for every nonbipartite graph $H$ in place of the triangle. We obtained the
proof below for triangles independently, along the same lines, and include
it for completeness.

We compare graphs with graphons through step
functions and the cut distance. For a real symmetric $n\times n$ matrix
$M$, let $W_M$ be the function on $[0,1]^2$ equal to $M_{ij}$ on
$I_i\times I_j$. Here $I_i=[\frac{i-1}n,\frac in)$ for $i<n$ and
$I_n=[\frac{n-1}n,1]$. For a bounded measurable function $D$ on $[0,1]^2$,
let $\|D\|_\cut=\sup_{S,T}|\int_{S\times T}D|$ be its cut norm, where the
supremum is over measurable $S,T\subset[0,1]$. By Lemma~\ref{lem:stepcut},
$\|W_M\|_\cut=\|M\|_\cut$. For bounded symmetric measurable functions $U$
and $W$ on $[0,1]^2$, define the cut distance
\[
\delta_\cut(U,W)=\inf_\varphi\|U^\varphi-W\|_\cut ,
\]
where $U^\varphi(x,y)=U(\varphi(x),\varphi(y))$ and the infimum is over
measure preserving bijections $\varphi$ of $[0,1]$. Then $p^{-1}W_A$ is the
step function of $G(n,p)$ normalized by its expected edge density.

\begin{theorem}\label{thm:reduction}
Let $n^{-1/2}\ll p\ll1$, $\eta\in(0,1)$ and $r=\eta^{1/3}$.
\begin{enumerate}[label=(\alph*)]
\item $\mathrm{LT}_{n,p}(\eta)=(1+o(1))\frac{n^2p}2V(r)$.
\item Let $\varepsilon>0$. Conditionally on $T_n\le\eta\E[T_n]$, with
probability $1-o(1)$ there is a minimizer $W$ at $r$ with
$\delta_\cut(p^{-1}W_A,W)\le\varepsilon$.
\end{enumerate}
\end{theorem}

We deduce part (a) from
\cite[Theorem~2]{KS23} and a comparison of $\Phi_{n,p}$ with \eqref{eq:V}
when $n\to\infty$ and $p\to0$ together (Lemma~\ref{lem:discrete}). Part (b) compares the conditioned graph with the random graph
with edge probabilities $pW$ for a minimizer $W$. After a relabeling of the
vertices, the two graphs have about the same edge counts between vertex
sets, when these counts are compared at the scale of the edge density.

For part (b) we use only the existence of minimizers, as does the proof
of~\cite[Theorem~1.4]{Chi26}. Following
\cite{Chi25,KS23}, we reveal a small set of edges. The other edges are then
nearly independent, hence close in cut norm to their conditional
probabilities. We divide the step function of these probabilities by $p$
and set the revealed pairs to $0$. By the Harris inequality, these
probabilities are at most $p$, so the result is a graphon. This graphon is
nearly optimal in \eqref{eq:V}, so by the compactness of the space of
graphons it is close in cut distance to a minimizer. 
When the constant graphon $r$ is the
unique minimizer, we have $\delta_\cut(p^{-1}W_A,r)=\|p^{-1}W_A-r\|_\cut$.
Hence Theorem~\ref{thm:reduction}(b) gives the conclusion of
\cite[Corollary~1.2]{Chi25} for triangles when $p\to0$ and
$n^{1/2}p\to\infty$.

Finally, in Appendix~\ref{s3:sub:cor} we combine Theorem~\ref{thm:reduction}
with Theorem~\ref{thm:main}. We round a relabeling of $B_r$ to a partition of $[n]$ into two blocks
(Lemma~\ref{lem:rounding}), which gives Theorems~\ref{thm:rate}
and~\ref{thm:structure}.

\subsection{The lower tail rate via discretization}\label{s3:sub:value}

We first write the minimization $\Phi_{n,p}$ of
Section~\ref{s1:sub:previous} in coordinates, with $0\log0=0$ in $i_p$. Let
$E_n$ be the set of two element subsets of $[n]=\{1,\dots,n\}$. For
$q=(q_e)_{e\in E_n}\in[0,1]^{E_n}$ write $q_{ij}=q_{\{i,j\}}$ and let
$N(q)=\sum_{\{i,j,l\}}q_{ij}q_{jl}q_{il}$, where the sum is over the three
element subsets of $[n]$. Then $N(q)$ is the expected number of triangles
in $G_q$. For $\eta\ge0$ define
\begin{equation}\label{s3:Phi}
\Phi_{n,p}(\eta)=\min\Big\{\sum_{e\in E_n}i_p(q_e):q\in[0,1]^{E_n}\text{
and }N(q)\le\eta\binom n3p^3\Big\}.
\end{equation}
The minimum exists, since $i_p$ is continuous on $[0,1]$ and the feasible
set is compact and contains the constant vector $\min(\eta,1)^{1/3}p$.
For triangles, Kozma and Samotij~\cite[Theorem~2]{KS23} proved the following
result, which bounds $-\log\P[T_n\le\eta\E[T_n]]$ above and below by $\Phi_{n,p}$
evaluated near $\eta$, up to factors close to $1$.

\begin{lemma}[{\cite[Theorem~2]{KS23}}]\label{lem:KS}
For every $p_0\in(0,1)$ and $\varepsilon>0$ there is a constant $L$ with the
following property. If $Ln^{-1/2}\le p\le p_0$, then every $\eta\in[0,1]$
satisfies
\[
(1-\varepsilon)\Phi_{n,p}(\eta+\varepsilon)\le-\log\P\big[T_n\le\eta\E[T_n]\big]
\le(1+\varepsilon)\Phi_{n,p}((1-\varepsilon)\eta).
\]
\end{lemma}

Using Lemma~\ref{lem:KS}, we reduce part (a) of
Theorem~\ref{thm:reduction} to the behavior of $\Phi_{n,p}$ as $p\to0$. We next compare $i_p$ at a small mean with its sparse limit $h$, uniformly
on $[0,1]$.

\begin{lemma}\label{lem:fp}
For $p\in(0,1)$ and $w\in[0,1]$,
\[
0\le\frac{i_p(pw)}p-h(w)\le\frac{p(1-w)^2}{2(1-p)}\le\frac p{2(1-p)}.
\]
\end{lemma}

\begin{proof}
Let $g(w)=i_p(pw)/p-h(w)=(\frac1p-w)\log\frac{1-pw}{1-p}+w-1$, a smooth
function on $[0,1]$. Then $g(1)=0$ and $g'(w)=-\log\frac{1-pw}{1-p}$, so
$g'(1)=0$. Moreover, $g''(w)=p/(1-pw)\in[p,p/(1-p)]$. Hence Taylor's
formula $g(w)=\int_w^1(v-w)g''(v)\dd v$ gives
$0\le g(w)\le\frac p{1-p}\cdot\frac{(1-w)^2}2$.
\end{proof}

In the next lemma we identify the normalized discrete minimum with $V$ as
$p\to0$. Part (a) is a lower bound for every finite $n$, which we use again in Appendix~\ref{s3:sub:proof}. In the proof we pass to step graphons.
For $q\in[0,p]^{E_n}$ let $Q$, the \emph{matrix of $q$}, be the symmetric
matrix with entries $q_e$ off the diagonal and $0$ on it. Recall the
intervals $I_i$ and the step functions $W_M$ of
Section~\ref{s1:sub:reduction}. Then $W_{Q/p}$ is a graphon that vanishes
on the diagonal blocks $I_i\times I_i$. Only
triples of distinct indices contribute to $t(W_{Q/p})$, and each three
element subset of $[n]$ gives $6$ such triples. Since $h(0)=1$, it follows that
\begin{equation}\label{s3:step}
t(W_{Q/p})=\frac{6N(q)}{n^3p^3},\quad
F(W_{Q/p})=\frac2{n^2}\sum_{e\in E_n}h\Big(\frac{q_e}p\Big)+\frac1n.
\end{equation}

\begin{lemma}\label{lem:discrete}
Let $\eta\in(0,1)$.
\begin{enumerate}[label=(\alph*)]
\item For all $n\ge3$ and $p\in(0,1)$,
$\Phi_{n,p}(\eta)\ge\frac{n^2p}2\big(V(\eta^{1/3})-\frac1n\big)$.
\item For all sequences $n_j\to\infty$ and $p_j\in(0,1)$ with $p_j\to0$,
\[
\limsup_{j\to\infty}\frac{\Phi_{n_j,p_j}(\eta)}{n_j^2p_j/2}\le V(\eta^{1/3}).
\]
\end{enumerate}
\end{lemma}

\begin{proof}
(a) Let $q$ attain the minimum in \eqref{s3:Phi} and $q'_e=\min(q_e,p)$.
Since $i_p\ge0=i_p(p)$, we have $i_p(q'_e)\le i_p(q_e)$ for every $e$.
Since $N$ has nonnegative coefficients, $N(q')\le N(q)\le\eta\binom n3p^3$.
Let $W=W_{Q'/p}$, where $Q'$ is the matrix of $q'$. By \eqref{s3:step},
$t(W)\le\frac{6\eta}{n^3}\binom n3\le\eta$, so $F(W)\ge V(\eta^{1/3})$.
Since $i_p(q'_e)\ge ph(q'_e/p)$ by Lemma~\ref{lem:fp}, \eqref{s3:step}
gives
\[
\Phi_{n,p}(\eta)\ge\sum_ei_p(q'_e)\ge p\sum_eh\Big(\frac{q'_e}p\Big)
=\frac{n^2p}2\Big(F(W)-\frac1n\Big)\ge\frac{n^2p}2\Big(V(\eta^{1/3})
-\frac1n\Big).
\]

(b) Fix $\delta>0$. By Lemma~\ref{lem:existence} there are
$s\in(0,\eta^{1/3})$ with $V(s)\le V(\eta^{1/3})+\delta$ and a minimizer
$U$ at $s$. By Lemmas~\ref{lem:lsc}(a) and~\ref{lem:cut}(b), the graphons
$A_kU$ satisfy $F(A_kU)\le F(U)$ and $|t(A_kU)-t(U)|\le3\|A_kU-U\|_1\to0$.
Since $t(U)\le s^3<\eta$, we can fix $k$ with $t(A_kU)<\eta$. Let $K=2^k$,
and let $u_{ab}$ be the value of $A_kU$ on $J_a\times J_b$.

For $i\in[n]$ let $a(i)\in[K]$ be the index with
$\frac{2i-1}{2n}\in J_{a(i)}$. Let $q_{ij}=pu_{a(i)a(j)}$ for $i\ne j$, and
let $U_n=W_{Q/p}$ for the matrix $Q$ of $q$. Then $U_n$ does not depend on
$p$. Moreover, $U_n=A_kU$ outside the diagonal blocks and outside the sets
$(I_i\times[0,1])\cup([0,1]\times I_i)$ for the indices $i$ with
$I_i\not\subset J_{a(i)}$. Each such $I_i$ contains a point $\frac aK$ with
$a\in[K-1]$ in its interior, so there are at most $K-1$ of them. Hence
$|\{U_n\ne A_kU\}|\le\frac{2K}n$. Since $0\le h\le1$ on $[0,1]$, this
gives $|F(U_n)-F(A_kU)|\le\frac{2K}n$. Since $\|\cdot\|_\cut\le\|\cdot\|_1$,
Lemma~\ref{lem:cut}(b) gives $|t(U_n)-t(A_kU)|\le\frac{6K}n$. By
\eqref{s3:step}, $N(q)=\frac{n^3p^3}6t(U_n)$. Since $t(A_kU)<\eta$ and
$\frac6{n^3}\binom n3\to1$, the vector $q$ is feasible in \eqref{s3:Phi}
for all large $n$, uniformly in $p$. For such $n$,
\[
\Phi_{n,p}(\eta)\le\sum_{e\in E_n}i_p(q_e)
\overset{\mathrm{(a)}}{\le}p\sum_{e\in E_n}h\Big(\frac{q_e}p\Big)
+\frac{n^2p}2\cdot\frac p{2(1-p)}
\overset{\mathrm{(b)}}{\le}\frac{n^2p}2\Big(F(U_n)+\frac p{2(1-p)}\Big),
\]
where $\mathrm{(a)}$ uses Lemma~\ref{lem:fp} and $|E_n|\le\frac{n^2}2$,
and $\mathrm{(b)}$ is \eqref{s3:step}. Along $n_j\to\infty$ and $p_j\to0$,
the right side divided by $n_j^2p_j/2$ tends to $F(A_kU)\le F(U)=V(s)\le
V(\eta^{1/3})+\delta$. Since $\delta>0$ was arbitrary, (b) follows.
\end{proof}

\begin{proof}[Proof of Theorem~\ref{thm:reduction}(a)]
Fix $\eta\in(0,1)$ and $\varepsilon\in(0,1-\eta)$, and let
$R_n=\mathrm{LT}_{n,p}(\eta)/(n^2p/2)$. Let $L$ be given by
Lemma~\ref{lem:KS} with $p_0=\frac12$ and this $\varepsilon$. Since $p\to0$
and $n^{1/2}p\to\infty$, we have $Ln^{-1/2}\le p\le\frac12$ for all large
$n$. As $\eta+\varepsilon$ and $(1-\varepsilon)\eta$ lie in $(0,1)$,
\[
(1-\varepsilon)V\big((\eta+\varepsilon)^{1/3}\big)
\overset{\mathrm{(a)}}{\le}\liminf_{n\to\infty}R_n
\le\limsup_{n\to\infty}R_n\overset{\mathrm{(b)}}{\le}
(1+\varepsilon)V\big(((1-\varepsilon)\eta)^{1/3}\big),
\]
where $\mathrm{(a)}$ combines the lower bound of Lemma~\ref{lem:KS} with
Lemma~\ref{lem:discrete}(a) at $\eta+\varepsilon$, and $\mathrm{(b)}$
combines the upper bound of Lemma~\ref{lem:KS} with
Lemma~\ref{lem:discrete}(b) at $(1-\varepsilon)\eta$. Both are applied
along the sequence $(n,p(n))$. Letting $\varepsilon\to0$ and using the continuity of
$V$ (Lemma~\ref{lem:existence}(b)) gives $R_n\to V(\eta^{1/3})$. Since
$V(\eta^{1/3})>0$ by Lemma~\ref{lem:existence}(a), this is equivalent to
part (a).
\end{proof}

\subsection{Near minimizers and step functions}\label{s3:sub:stability}

By compactness, a graphon that nearly satisfies the constraint and nearly
attains $V(r)$ is, after relabeling, close to a minimizer at $r$. In
Appendix~\ref{s3:sub:proof} we apply this to the step graphon of the
conditional edge probabilities.

\begin{lemma}\label{lem:stability}
Let $r\in(0,1)$ and $\varepsilon>0$. There is $\delta>0$ with the following
property. If $U$ is a graphon with $t(U)\le r^3+\delta$ and
$F(U)\le V(r)+\delta$, then there are a minimizer $W_0$ at $r$ and a measure
preserving bijection $\varphi$ of $[0,1]$ with
$\|U^\varphi-W_0\|_\cut\le\varepsilon$.
\end{lemma}

\begin{proof}
Otherwise there are graphons $U_k$ with $t(U_k)\le r^3+1/k$ and
$F(U_k)\le V(r)+1/k$ for which no such $W_0$ and $\varphi$ exist. By
Lemma~\ref{lem:cut}(c), along a subsequence there are measure preserving
bijections $\varphi_k$ and a graphon $U$ with
$\|U_k^{\varphi_k}-U\|_\cut\to0$. Since relabeling preserves $t$ and $F$,
Lemmas~\ref{lem:cut}(b) and~\ref{lem:lsc}(b) give $t(U)\le r^3$ and
$F(U)\le V(r)$. Thus $U$ is a minimizer at $r$, which contradicts the
choice of $U_k$ for large $k$.
\end{proof}

We next pass from graphons to $n\times n$ matrices. Recall the normalized
cut norm $\|M\|_\cut$ of Section~\ref{sec:intro}. In the next lemma we transfer cut norm bounds between
step functions and matrices.

\begin{lemma}\label{lem:stepcut}
Let $M$ be a real symmetric $n\times n$ matrix and let $\varphi$ be measure
preserving. Then $\|W_M^\varphi\|_\cut=\|W_M\|_\cut=\|M\|_\cut$. Moreover,
the maximum defining $\|M\|_\cut$ is attained at some $x,y\in\{0,1\}^n$.
\end{lemma}

\begin{proof}
The sets $I_i^\varphi=\varphi^{-1}(I_i)$ partition $[0,1]$ into sets of
measure $\frac1n$, and $W_M^\varphi$ equals $M_{ij}$ on
$I_i^\varphi\times I_j^\varphi$. Hence $\iint W_M^\varphi fg=n^{-2}x^{\top}My$
for measurable $f,g:[0,1]\to[0,1]$, where $x_i=n\int_{I_i^\varphi}f$ and
$y_j=n\int_{I_j^\varphi}g$ lie in $[0,1]$, and every $x,y\in[0,1]^n$ arises
from $f=\sum_ix_i\1_{I_i^\varphi}$ and $g=\sum_jy_j\1_{I_j^\varphi}$.
Since indicator functions are such $f,g$, Lemma~\ref{lem:cut}(a) shows that
$\|W_M^\varphi\|_\cut$ is the supremum of $|\iint W_M^\varphi fg|$ over
such $f,g$. Hence $\|W_M^\varphi\|_\cut=\|M\|_\cut$ for every measure
preserving $\varphi$, including $\varphi=\mathrm{id}$. Finally, $x^{\top}My$ is linear in $x$ for fixed $y$ and in $y$ for fixed
$x$. Optimizing first in $x$ and then in $y$ shows that its maximum and
minimum over $[0,1]^n\times[0,1]^n$ are attained at some
$x,y\in\{0,1\}^n$. Hence so is the maximum of $|x^{\top}My|$.
\end{proof}

\subsection{Decorrelation and concentration under the conditioned measure}\label{s3:sub:measure}

Let $Y=(Y_e)_{e\in E_n}$ be the edge indicators of $G(n,p)$, and let $\P$
be their law. Fix $\eta\in(0,1)$. Recall the adjacency matrix $A$ and the
event $\mathcal L=\{T_n\le\eta\E[T_n]\}$ of Section~\ref{sec:intro}, and
let $\mathbb Q=\P[\cdot\mid\mathcal L]$. For $W\subset E_n$ write
$Y_W=(Y_e)_{e\in W}$. For a set $\kappa$ of edges let
$\chi_\kappa=\prod_{e\in\kappa}Y_e$. We view a triangle of $K_n$ as a set of
three edges, and we let $\mathcal T_W$ be the set of triangles of $K_n$
that are disjoint from $W$. We apply the results of \cite{KS23} to the $3$
uniform hypergraph $\mathcal H$ on the vertex set $E_n$ whose edges are the
triangles of $K_n$, all with unit weight. A set of edges of $K_n$ contains
as many edges of $\mathcal H$ as the corresponding graph has triangles, and
$\mathcal H$ has $\binom n3=\E[T_n]/p^3$ edges. Hence $\mathcal L$ is the
lower tail event of \cite[Section~5]{KS23} for $\mathcal H$, and $\mathbb Q$
is the conditioned law studied there.

An event is \emph{increasing} if it is preserved by adding edges, and
\emph{decreasing} if it is preserved by removing edges. The event
$\mathcal L$ is decreasing. We next record two consequences of the Harris inequality. Part (a) controls the edges of the conditioned graph
inside a small set, and part (b) bounds the conditional edge probabilities
by $p$.

\begin{lemma}\label{lem:harris}
\begin{enumerate}[label=(\alph*)]
\item For every increasing event $\mathcal D$,
$\mathbb Q[\mathcal D]\le\P[\mathcal D]$.
\item Let $W\subset E_n$ and $\kappa\subset E_n\setminus W$, and let $y$ be a
value of $Y_W$ with $\mathbb Q[Y_W=y]>0$. Then $\mathbb
Q[\chi_\kappa=1\mid Y_W=y]\le p^{|\kappa|}$.
\end{enumerate}
\end{lemma}

\begin{proof}
(a) Since $\mathcal L$ is decreasing, the Harris inequality (see, e.g.,
\cite{Har60}) for the product law $\P$ gives $\P[\mathcal
D\cap\mathcal L]\le\P[\mathcal D]\P[\mathcal L]$. Dividing by
$\P[\mathcal L]>0$ gives (a).

(b) This is \cite[Claim~16]{KS23} for $\mathcal H$. See also
\cite[Lemma~3.1]{Chi25}. Under $\P[\cdot\mid Y_W=y]$, the edges outside
$W$ are independent with mean $p$. The event $\{\chi_\kappa=1\}$ is
increasing, and the section $\{z:(y,z)\in\mathcal L\}$ is decreasing with
positive probability, since $\mathbb Q[Y_W=y]>0$. Hence the Harris
inequality applies as in (a).
\end{proof}

The next lemma is the decorrelation estimate of Kozma and
Samotij~\cite{KS23}, which is also the starting point of the work of
Chin~\cite{Chi25}.

\begin{lemma}\label{lem:decorrelation}
For all $\alpha,\beta\in(0,1)$ there is $L\ge4$ with the following
property. If $n\ge L$ and $Ln^{-1/2}\le p\le\frac12$, then there is a set
$W\subset E_n$ with $|W|\le\alpha\binom n2$ and
\[
\E_{\mathbb Q}\Big[\sum_{\kappa\in\mathcal T_W}\Big|\mathbb Q[\chi_\kappa=1\mid Y_W]
-\prod_{e\in\kappa}\mathbb Q[Y_e=1\mid Y_W]\Big|\Big]\le\beta p^3\binom n3.
\]
Here $W$ may depend on $n$, $p$, $\eta$, $\alpha$ and $\beta$.
\end{lemma}

\begin{proof}
In the notation of \cite{KS23}, $v(\mathcal H)=\binom n2$ and
$e(\mathcal H)=\binom n3$, so $e(\mathcal H)/v(\mathcal H)=\frac{n-2}3$. A
pair lies in $n-2$ triangles, and two or three pairs lie in at most one
triangle. Thus $\Delta_1(\mathcal H)=n-2$ and $\Delta_2(\mathcal
H)=\Delta_3(\mathcal H)=1$. Hence the codegree condition
$\Delta_s(\mathcal H)\le3(\lambda p)^{s-1}e(\mathcal H)/v(\mathcal H)$ of
\cite{KS23} holds for $s=1,2,3$ as soon as $n\ge4$, $\lambda pn\ge2$ and
$\lambda^2p^2n\ge2$. In particular, it holds for $n\ge4$ and
$p\ge2\lambda^{-1}n^{-1/2}$.

We apply \cite[Lemma~18]{KS23} to $\mathcal H$, with uniformity $3$,
codegree constant $3$, $p_0=\frac12$, the given $\alpha$, and $\beta^2/2$ in
place of $\beta$. This gives $\lambda>0$ and $V_0$. Let $\binom n2\ge
V_0$, and let the codegree condition hold with this $\lambda$. Then there is
a set $W$ with $|W|\le\alpha\binom n2$ and
$\mathcal E(W)\le\frac{\beta^2}2e(\mathcal H)$. Here $\mathcal E(W)$ is the
expected weighted sum of squared conditional covariances given $Y_W$
defined in \cite[p.~686]{KS23}. Since $e(\mathcal H-W)\le e(\mathcal H)$,
the Cauchy--Schwarz bound on the same page shows that the left side of the
inequality in Lemma~\ref{lem:decorrelation} is at most
$p^3(2e(\mathcal H))^{1/2}\mathcal E(W)^{1/2}\le\beta p^3\binom n3$. See
also \cite[Lemmas~3.2 and~3.3]{Chi25}. It suffices to take
$L\ge\max\{2\lambda^{-1},n_0,4\}$, where $n_0$ is such that $\binom n2\ge
V_0$ for all $n\ge n_0$.
\end{proof}

We also use two facts about the relative entropy $D(\cdot\|\cdot)$ of
Section~\ref{s1:sub:previous}. Let $W\subset E_n$, and let $y$ be a value
of $Y_W$ with $\mathbb Q[Y_W=y]>0$. Let $\mathbb Q^y$ be the law of
$Y_{E_n\setminus W}$ under $\mathbb Q[\cdot\mid Y_W=y]$. For $e\notin W$ let
$q^y_e=\mathbb Q[Y_e=1\mid Y_W=y]$, and let $\widehat{\mathbb Q}^y$ be the
product of the Bernoulli laws $\widehat{\mathbb Q}^y_e$ with means $q^y_e$.
Part (a) of the next lemma bounds by $-\log\P[\mathcal L]$ the sum of two
expected costs, namely the relative entropy of $\mathbb Q^{Y_W}$ with
respect to the product of its marginals and the cost of the conditional
edge probabilities. Part (b) transfers small probabilities from
$\widehat{\mathbb Q}^y$ to $\mathbb Q^y$.

\begin{lemma}\label{lem:entropy}
Let $W\subset E_n$.
\begin{enumerate}[label=(\alph*)]
\item $\E_{\mathbb Q}\big[D(\mathbb Q^{Y_W}\|\widehat{\mathbb Q}^{Y_W})\big]
+\E_{\mathbb Q}\big[\sum_{e\notin W}i_p(q^{Y_W}_e)\big]\le-\log\P[\mathcal
L]$.
\item Let $y$ be a value of $Y_W$ with $\mathbb Q[Y_W=y]>0$, and let
$\mathcal E$ be an event. If $0<\widehat{\mathbb Q}^y[\mathcal E]<1$, then
\[
\mathbb Q^y[\mathcal E]\le\frac{D(\mathbb Q^y\|\widehat{\mathbb Q}^y)+\log2}
{\log(1/\widehat{\mathbb Q}^y[\mathcal E])}.
\]
If $\widehat{\mathbb Q}^y[\mathcal E]=0$, then $\mathbb Q^y[\mathcal E]=0$.
\end{enumerate}
\end{lemma}

\begin{proof}
(a) By \cite[Proposition~10]{KS23}, $D(\mathbb Q\|\P)=-\log\P[\mathcal L]$.
Write an outcome as $(y,z)$, where $y$ is the value of $Y_W$. Let $\mathbb
Q_W$ and $\P_W$ be the laws of $Y_W$ under $\mathbb Q$ and $\P$, and let
$\P'$ be the law of $Y_{E_n\setminus W}$ under $\P$, so that
$\P=\P_W\otimes\P'$. Taking expectations under $\mathbb Q$ in
$\log\frac{\mathbb Q(y,z)}{\P(y,z)}=\log\frac{\mathbb Q_W(y)}{\P_W(y)}
+\log\frac{\mathbb Q^y(z)}{\P'(z)}$ gives
\[
-\log\P[\mathcal L]=D(\mathbb Q_W\|\P_W)+\E_{\mathbb Q}\big[D(\mathbb
Q^{Y_W}\|\P')\big]\ge\E_{\mathbb Q}\big[D(\mathbb Q^{Y_W}\|\P')\big].
\]
If $\mathbb Q^y(z)>0$, then every coordinate of $z$ has positive
probability under $\mathbb Q^y$, so $\widehat{\mathbb Q}^y(z)>0$. Thus
$\mathbb Q^y\ll\widehat{\mathbb Q}^y$. Let $\P'_e$ be the Bernoulli law
with mean $p$. Under $\mathbb Q^y$, each $z_e$ has law $\widehat{\mathbb
Q}^y_e$. Hence taking expectations under $\mathbb Q^y$ in
$\log\frac{\mathbb Q^y(z)}{\P'(z)}=\log\frac{\mathbb
Q^y(z)}{\widehat{\mathbb Q}^y(z)}+\sum_{e\notin
W}\log\frac{\widehat{\mathbb Q}^y_e(z_e)}{\P'_e(z_e)}$ gives
\[
D(\mathbb Q^y\|\P')=D(\mathbb Q^y\|\widehat{\mathbb Q}^y)+\sum_{e\notin
W}i_p(q^y_e).
\]
Averaging this identity over $Y_W$ under $\mathbb Q$ and inserting it into
the lower bound for $-\log\P[\mathcal L]$ proves (a).

(b) The last claim follows from $\mathbb Q^y\ll\widehat{\mathbb Q}^y$. Let
$x=\mathbb Q^y[\mathcal E]$ and $u=\widehat{\mathbb Q}^y[\mathcal E]$. By the
concavity of the logarithm, grouping the outcomes in $\mathcal E$ and in its
complement can only decrease the relative entropy. Hence
\[
D(\mathbb Q^y\|\widehat{\mathbb Q}^y)\ge x\log\frac xu+(1-x)\log\frac{1-x}{1-u}
\ge x\log\frac1u-\log2,
\]
where the second inequality uses $x\log x+(1-x)\log(1-x)\ge-\log2$ and
$\log\frac1{1-u}\ge0$. Since $0<u<1$, dividing by $\log\frac1u>0$ gives the
first claim.
\end{proof}

Finally, cut norm deviations are exponentially unlikely under a product
measure with small means.

\begin{lemma}\label{lem:cutconc}
Let $E'\subset E_n$, and let $(Z_e)_{e\in E'}$ be independent Bernoulli
variables with means $q_e\le p$. Let $Z$ be the symmetric $n\times n$
matrix with entries $Z_e-q_e$ at the two positions of each $e\in E'$ and $0$
elsewhere. Then for every $\varepsilon\in(0,1]$,
\[
\P\big[\|Z\|_\cut>\varepsilon p\big]\le2\cdot4^n\exp\Big(-\frac{\varepsilon^2pn^2}8\Big).
\]
\end{lemma}

\begin{proof}
By Lemma~\ref{lem:stepcut} it suffices to consider the $4^n$ pairs
$x,y\in\{0,1\}^n$. For such $x,y$ we have $x^{\top}Zy=\sum_{e\in
E'}c_e(Z_e-q_e)$ with $c_e=x_iy_j+x_jy_i\in\{0,1,2\}$ for $e=\{i,j\}$. Let
$\lambda=\varepsilon/4$. For $|u|\le1$ we have $e^u-1-u\le u^2$, so for
$\tau\in\{1,-1\}$,
\[
\E\big[e^{\tau\lambda c_e(Z_e-q_e)}\big]=e^{-\tau\lambda c_eq_e}\big(1-q_e
+q_ee^{\tau\lambda c_e}\big)\le\exp\big(q_e(e^{\tau\lambda c_e}-1-\tau\lambda
c_e)\big)\le e^{4\lambda^2q_e}.
\]
Since $\sum_eq_e\le p\binom n2\le\frac{pn^2}2$, independence and Markov's
inequality for $e^{\tau\lambda x^{\top}Zy}$ give
$\P[\tau x^{\top}Zy>\varepsilon pn^2]\le\exp(2\lambda^2pn^2-\lambda\varepsilon
pn^2)=\exp(-\varepsilon^2pn^2/8)$. A union bound over $x$, $y$ and $\tau$
completes the proof.
\end{proof}

\subsection{Proof of Theorem~\ref{thm:reduction}(b)}\label{s3:sub:proof}

Fix $\eta\in(0,1)$ and $\varepsilon>0$, and let $r=\eta^{1/3}$. The
conclusion weakens as $\varepsilon$ grows, so we may assume $\varepsilon<1$.
Let $\mathcal B$ be the event that $\delta_\cut(p^{-1}W_A,W_0)>\varepsilon$
for every minimizer $W_0$ at $r$. Fix $\vartheta\in(0,1)$. We consider four
events. On $\mathcal G_1$, the conditional edge probabilities nearly satisfy
the triangle constraint. On $\mathcal G_2$, their cost is nearly $V(r)$. On
$\mathcal E^c$, the unrevealed edges are close in cut norm to their
conditional probabilities. The fourth event says that few revealed edges
are present. We show that each of these events fails with $\mathbb Q$
probability at most $\vartheta/8$ for all large $n$, and that
$\mathcal B$ fails on their intersection. Hence $\limsup_n\mathbb
Q[\mathcal B]\le\vartheta/2$. Since $\vartheta$ is arbitrary, $\mathbb
Q[\mathcal B]\to0$.

We choose, in this order, the stability tolerance, the allowed excess
cost, the entropy error, the triangle slack, the decorrelation error and
the size of the revealed set. Each choice depends only on the earlier ones. Let $\delta>0$ be the
tolerance given by Lemma~\ref{lem:stability} for $\varepsilon/4$ in place
of $\varepsilon$, and let $\gamma=\delta/4$ be the allowed excess cost.
Let $\omega_0=\frac\vartheta8\min\{\gamma,\varepsilon^2/400\}$. It bounds
the average excess cost and the average relative entropy to the product
measure. Choose the slack $\zeta\in(0,\min\{\delta,\omega_0/8\}]$ in the
triangle constraint with $\eta+\zeta<1$ and $1-r/(\eta+\zeta)^{1/3}\le\omega_0/8$.
Such $\zeta$ exists because $(\eta+\zeta)^{1/3}\to r$ as $\zeta\to0$. Let
$\beta=\zeta\min\{\omega_0,\vartheta\}/8$ be the allowed decorrelation
error. Finally, let $\alpha=\min\{\varepsilon/16,\delta/16,\omega_0/8\}$
bound the size of the revealed set, and let $L$ be given by
Lemma~\ref{lem:decorrelation} for $\alpha$ and $\beta$. Let
$\Sigma=\frac{n^2p}2$. For all large $n$ we have
\begin{align}
&Ln^{-1/2}\le p\le\tfrac12,\quad-\log\P[\mathcal L]\le\Sigma\big(V(r)+\zeta\big),
\quad\frac1n\le\frac{\omega_0}8,\notag\\
&2\cdot4^ne^{-\varepsilon^2pn^2/128}\le e^{-\varepsilon^2pn^2/256},\quad
\log2\le\omega_0n^2p,\quad e^{-\alpha pn^2}\le\frac\vartheta8.
\label{s3:large}
\end{align}
Indeed, Theorem~\ref{thm:reduction}(a) gives the bound on
$-\log\P[\mathcal L]$. The other bounds follow from $p\to0$ and
$n^{1/2}p\to\infty$. We consider such $n$, and note that $n\ge L$ since
$Ln^{-1/2}\le\frac12$.

Next we define the revealed set and the conditional edge probabilities.
Let $W$ be given by Lemma~\ref{lem:decorrelation}, so $|W|\le\alpha\binom
n2$. Let $y$ be a value of $Y_W$ with positive probability, and let
$q^y_e$, $\mathbb Q^y$ and $\widehat{\mathbb Q}^y$ be as defined before
Lemma~\ref{lem:entropy}. Let $\tilde q^y_e=q^y_e$ for $e\notin W$ and
$\tilde q^y_e=0$ for $e\in W$. By Lemma~\ref{lem:harris}(b),
$0\le\tilde q^y_e\le p$. Let $\tilde Q^y$ be the matrix of $\tilde q^y$,
and let $U^y=W_{\tilde Q^y/p}$. From now on, $y$ is the
random value of $Y_W$, and we drop it from the notation.

We start with the triangle constraint. For a triangle $\kappa\in\mathcal
T_W$ we have $\prod_{e\in\kappa}\tilde q_e=\prod_{e\in\kappa}q_e$. Every
other triangle contains an edge of $W$, so the product of the numbers
$\tilde q_e$ over its edges vanishes. Hence $N(\tilde q)=\sum_{\kappa\in
\mathcal T_W}\prod_{e\in\kappa}q_e$. On the other hand, $\sum_{\kappa\in
\mathcal T_W}\mathbb Q[\chi_\kappa=1\mid Y_W]$ is the conditional
expectation of the number of triangles in $\mathcal T_W$. Since $\mathbb
Q[\mathcal L]=1$, it is at most $\E_{\mathbb Q}[T_n\mid Y_W]\le\eta\binom
n3p^3$. Let $\mathcal G_1=\{N(\tilde q)\le(\eta+\zeta)\binom n3p^3\}$. On
$\mathcal G_1^c$, the sum inside the expectation in
Lemma~\ref{lem:decorrelation} therefore exceeds $\zeta\binom n3p^3$. Hence
Lemma~\ref{lem:decorrelation} and Markov's inequality give $\mathbb
Q[\mathcal G_1^c]\le\beta/\zeta\le\vartheta/8$. On $\mathcal G_1$ we have
$t(U)=6N(\tilde q)/(n^3p^3)\le\eta+\zeta$ by \eqref{s3:step}.

Next we show that the cost of the conditional edge probabilities is
nearly $V(r)$. Let $X=\sum_{e\notin W}i_p(q_e)$ and
$v_0=V((\eta+\zeta)^{1/3})-\frac1n-2\alpha$. Since $p\le\frac12$ by
\eqref{s3:large}, we have $i_p(0)=\log\frac1{1-p}\le2p$. Since also
$|W|\le\alpha n^2/2$, it follows that $\sum_{e\in E_n}i_p(\tilde q_e)\le
X+\alpha n^2p$. On $\mathcal G_1$ the vector $\tilde q$ is feasible in
\eqref{s3:Phi} with $\eta+\zeta$ in place of $\eta$. Hence
Lemma~\ref{lem:discrete}(a) and the bound on $\sum_{e\in E_n}i_p(\tilde
q_e)$ give $X\ge\Sigma v_0$ on $\mathcal G_1$. To compare $v_0$ with
$V(r)$, we apply Lemma~\ref{lem:existence}(b) with
$s=r<(\eta+\zeta)^{1/3}$. Together with our choices of $\zeta$, $\alpha$
and $n$, this gives
\begin{equation}\label{s3:v0}
0\le V(r)-v_0\le1-\frac r{(\eta+\zeta)^{1/3}}+\frac1n+2\alpha\le\frac{\omega_0}2.
\end{equation}
By Lemma~\ref{lem:entropy}(a) and \eqref{s3:large},
$\E_{\mathbb Q}[X]\le-\log\P[\mathcal L]\le\Sigma(V(r)+\zeta)$. Since $X\ge0$,
$v_0\le1$ and $X\ge\Sigma v_0$ on $\mathcal G_1$,
\[
\E_{\mathbb Q}\big[(X/\Sigma-v_0)^+\big]\le\E_{\mathbb Q}[X]/\Sigma-v_0
+\mathbb Q[\mathcal G_1^c]\le\frac{\omega_0}2+\zeta+\frac\beta\zeta\le\omega_0 .
\]
Let $\mathcal G_2=\{X\le\Sigma(v_0+\gamma)\}$. On $\mathcal G_2^c$ we have
$(X/\Sigma-v_0)^+>\gamma$, so Markov's inequality gives $\mathbb
Q[\mathcal G_2^c]\le\omega_0/\gamma\le\vartheta/8$. On $\mathcal
G_1\cap\mathcal G_2$,
\[
F(U)\overset{\mathrm{(a)}}{=}\frac2{n^2}\sum_{e\in E_n}h\Big(\frac{\tilde
q_e}p\Big)+\frac1n\overset{\mathrm{(b)}}{\le}\frac{X+\alpha n^2p}\Sigma+\frac1n
\overset{\mathrm{(c)}}{\le}V(r)+\gamma+2\alpha+\frac1n\le V(r)+\delta,
\]
where $\mathrm{(a)}$ is \eqref{s3:step}; $\mathrm{(b)}$ uses
Lemma~\ref{lem:fp} and $\sum_{e\in E_n}i_p(\tilde q_e)\le X+\alpha n^2p$;
$\mathrm{(c)}$ uses $\mathcal G_2$ and $v_0\le V(r)$. Moreover
$t(U)\le\eta+\zeta\le r^3+\delta$ on $\mathcal G_1$. Hence, on $\mathcal
G_1\cap\mathcal G_2$, Lemma~\ref{lem:stability} gives a minimizer $W_0$ at
$r$ and a measure preserving bijection $\varphi$ with
$\|U^\varphi-W_0\|_\cut\le\varepsilon/4$.

We now compare the unrevealed edges with their conditional probabilities
in cut norm. Let $A'$ be obtained from $A$ by replacing the entries at pairs
in $W$ by $0$. Let $\mathcal E=\{\|A'-\tilde Q\|_\cut>\varepsilon p/4\}$.
Conditionally on $Y_W=y$, the event $\mathcal E$ depends only on
$Y_{E_n\setminus W}$. If $Y_{E_n\setminus W}$ has law $\widehat{\mathbb
Q}^y$, then $A'-\tilde Q^y$ is the matrix $Z$ of Lemma~\ref{lem:cutconc}
with $E'=E_n\setminus W$ and means $q^y_e\le p$. By
Lemma~\ref{lem:cutconc} with $\varepsilon/4$ in place of $\varepsilon$ and
by \eqref{s3:large}, we get $\widehat{\mathbb Q}^y[\mathcal E]\le
e^{-\varepsilon^2pn^2/256}$. Hence Lemma~\ref{lem:entropy}(b) gives
\[
\mathbb Q[\mathcal E]=\E_{\mathbb Q}\big[\mathbb Q^{Y_W}[\mathcal E]\big]\le
\frac{256\big(\E_{\mathbb Q}[D(\mathbb Q^{Y_W}\|\widehat{\mathbb
Q}^{Y_W})]+\log2\big)}{\varepsilon^2pn^2}.
\]
Using $\E_{\mathbb Q}[X]\ge\Sigma v_0\mathbb Q[\mathcal G_1]\ge\Sigma(v_0-\mathbb
Q[\mathcal G_1^c])$, we have
\[
\E_{\mathbb Q}\big[D(\mathbb Q^{Y_W}\|\widehat{\mathbb Q}^{Y_W})\big]
\overset{\mathrm{(a)}}{\le}-\log\P[\mathcal L]-\E_{\mathbb Q}[X]
\overset{\mathrm{(b)}}{\le}\Sigma\Big(V(r)-v_0+\zeta+\frac\beta\zeta\Big)
\overset{\mathrm{(c)}}{\le}\Sigma\omega_0,
\]
where $\mathrm{(a)}$ is Lemma~\ref{lem:entropy}(a); $\mathrm{(b)}$ uses
\eqref{s3:large} and $\mathbb Q[\mathcal G_1^c]\le\beta/\zeta$; $\mathrm{(c)}$
uses \eqref{s3:v0} and the choices of $\zeta$ and $\beta$. Since
$\log2\le\omega_0n^2p$ by \eqref{s3:large}, we get $\mathbb Q[\mathcal
E]\le384\omega_0/\varepsilon^2\le\vartheta/8$.

Finally, we show that few revealed edges are present. Let
$N_W=\sum_{e\in W}Y_e$. The event $\{N_W\ge2\alpha pn^2\}$ is increasing.
Under $\P$ we have $\E[e^{N_W}]=(1+p(e-1))^{|W|}\le e^{2p|W|}\le e^{\alpha
pn^2}$. Hence Markov's inequality gives $\P[N_W\ge2\alpha pn^2]\le
e^{-\alpha pn^2}$. By
Lemma~\ref{lem:harris}(a) and \eqref{s3:large}, $\mathbb Q[N_W\ge2\alpha
pn^2]\le\vartheta/8$. If $N_W<2\alpha pn^2$, then
$\|A-A'\|_\cut\le2N_W/n^2<4\alpha p$.

We can now complete the proof. The events $\mathcal G_1$, $\mathcal G_2$,
$\mathcal E^c$ and $\{N_W<2\alpha pn^2\}$ all hold outside an event of
$\mathbb Q$ probability at most $\vartheta/2$. On their intersection,
\[
\|A-\tilde Q\|_\cut\le\|A-A'\|_\cut+\|A'-\tilde Q\|_\cut\le4\alpha
p+\frac{\varepsilon p}4\le\frac{\varepsilon p}2.
\]
Also on this intersection, the minimizer $W_0$ at $r$ and the bijection
$\varphi$ given by Lemma~\ref{lem:stability} satisfy
$\|U^\varphi-W_0\|_\cut\le\varepsilon/4$. Since
$p^{-1}W_A-U=W_{(A-\tilde Q)/p}$, Lemma~\ref{lem:stepcut} gives
\[
\|(p^{-1}W_A)^\varphi-W_0\|_\cut\le\|(W_{(A-\tilde Q)/p})^\varphi\|_\cut
+\|U^\varphi-W_0\|_\cut\le\frac\varepsilon2+\frac\varepsilon4<\varepsilon.
\]
Hence $\delta_\cut(p^{-1}W_A,W_0)<\varepsilon$, so $\mathcal B$ fails. This
completes the proof. \qed

\subsection{Proofs of Theorems~\ref{thm:rate} and~\ref{thm:structure}}\label{s3:sub:cor}

Recall the matrices $M^S_{a,b}$ of Section~\ref{sec:intro} and the
graphons $B^{S_0}_{a,b}$ of Section~\ref{sec:prelim}, defined for measurable
$S_0\subset[0,1]$ with $|S_0|=\frac12$. In the next lemma we round a step
function $W_M$ to an explicit vertex partition when $W_M$ is close to a
relabeling of a graphon $B_{a,b}$. The lemma applies to any real symmetric matrix
$M$, in particular to $A/p$.

\begin{lemma}\label{lem:rounding}
Let $M$ be a real symmetric $n\times n$ matrix, $0\le a<b\le1$ and
$\delta>0$. Suppose that there are a measure preserving $\varphi$ and a
measurable $S_0\subset[0,1]$ with $|S_0|=\frac12$ such that
$\|W_M^\varphi-B^{S_0}_{a,b}\|_\cut\le\delta$. Then there is $S\subset[n]$
with $\big||S|-\frac n2\big|\le\frac{16\delta n}{b-a}$ and
$\|M-M^S_{a,b}\|_\cut\le33\delta$.
\end{lemma}

\begin{proof}
Let $\lambda_1=\frac{a+b}2$, $\lambda_2=\frac{a-b}2<0$ and
$g=\1_{S_0}-\1_{[0,1]\setminus S_0}$. Then $\|g\|_2=1$ and
$B^{S_0}_{a,b}=\lambda_1+\lambda_2g\otimes g$, where $(g\otimes
g)(x,y)=g(x)g(y)$. Let $I_i^\varphi=\varphi^{-1}(I_i)$ for $i\in[n]$. As in
the proof of Lemma~\ref{lem:stepcut}, these sets partition $[0,1]$ into
sets of measure $\frac1n$, and $W_M^\varphi$ is constant on each
$I_i^\varphi\times I_j^\varphi$. For an integrable function
$f$ on $[0,1]$ let $\Pi f$ be the function equal to $n\int_{I_i^\varphi}f$
on $I_i^\varphi$ for each $i$. For an integrable function $D$ on $[0,1]^2$
let $\Pi D$ be the function equal to $n^2\iint_{I_i^\varphi\times
I_j^\varphi}D$ on $I_i^\varphi\times I_j^\varphi$ for all $i,j$.

We first show that $\Pi$ does not increase the cut norm. For measurable
$f,f':[0,1]\to[0,1]$, Fubini's theorem gives
\[
\iint(\Pi D)(x,y)f(x)f'(y)\dd x\dd y=\iint D(x,y)(\Pi f)(x)(\Pi
f')(y)\dd x\dd y.
\]
Since $\Pi f$ and $\Pi f'$ take values in $[0,1]$, Lemma~\ref{lem:cut}(a)
gives $\|\Pi D\|_\cut\le\|D\|_\cut$. Let $\tilde g=\Pi g$. We apply this
bound to $D=W_M^\varphi-B^{S_0}_{a,b}$. Since $\Pi W_M^\varphi=W_M^\varphi$,
$\Pi1=1$ and $\Pi(g\otimes g)=\tilde g\otimes\tilde g$, we get
$\|W_M^\varphi-\lambda_1-\lambda_2\tilde g\otimes\tilde g\|_\cut\le\delta$.
By the triangle inequality and the hypothesis
$\|W_M^\varphi-B^{S_0}_{a,b}\|_\cut\le\delta$, it follows that
\[
|\lambda_2|\|g\otimes g-\tilde g\otimes\tilde g\|_\cut\le2\delta .
\]

Next we show that $\tilde g$ is close to $g$ in $L^2([0,1])$. Let $f,f'$
take values in $[-1,1]$. Write each of them as a difference of two
functions with values in $[0,1]$. Applying Lemma~\ref{lem:cut}(a) to the
four resulting terms gives
$|\iint D(x,y)f(x)f'(y)\dd x\dd y|\le4\|D\|_\cut$. The map $\Pi$ is the
orthogonal projection of $L^2([0,1])$ onto the functions that are constant
on each $I_i^\varphi$. Hence $\langle\tilde g,g\rangle=\|\tilde
g\|_2^2\le1$. Testing with $f=f'=g$ gives
\[
|\lambda_2|\big(1-\|\tilde g\|_2^4\big)=|\lambda_2|\Big|\iint(g\otimes
g-\tilde g\otimes\tilde g)(x,y)g(x)g(y)\dd x\dd y\Big|\le8\delta .
\]
Let $\theta=8\delta/|\lambda_2|=16\delta/(b-a)$. Since $\|\tilde g\|_2\le1$,
it follows that
$\|g-\tilde g\|_2^2=1-\|\tilde g\|_2^2\le1-\|\tilde g\|_2^4\le\theta$.

It remains to round $\tilde g$ to a vertex partition. Let $c_i$ be the
value of $\tilde g$ on $I_i^\varphi$, and let $S=\{i:c_i>0\}$. Let $s_i=1$
for $i\in S$ and $s_i=-1$ otherwise, and let
$\sigma=\sum_is_i\1_{I_i^\varphi}$. At every point where $\sigma\ne g$ we have
$|g-\tilde g|\ge1$, since $g=\pm1$ and $\tilde g$ has the sign of $\sigma$
or vanishes. Hence $|\{\sigma\ne g\}|\le\|g-\tilde g\|_2^2\le\theta$. Since
$|\{\sigma=1\}|=|S|/n$ and $|\{g=1\}|=|S_0|=\frac12$, this gives
$\big||S|-\frac n2\big|\le\theta n$. Moreover,
$W_{M^S_{a,b}}^\varphi=\lambda_1+\lambda_2\sigma\otimes\sigma$. Since
$|\sigma(x)\sigma(y)-g(x)g(y)|\le2(\1_{\{\sigma\ne g\}}(x)+\1_{\{\sigma\ne
g\}}(y))$, we get $\|W_{M^S_{a,b}}^\varphi-B^{S_0}_{a,b}\|_1\le4|\lambda_2|\theta
=32\delta$. By Lemma~\ref{lem:stepcut} and $\|\cdot\|_\cut\le\|\cdot\|_1$,
\[
\|M-M^S_{a,b}\|_\cut=\|W_M^\varphi-W_{M^S_{a,b}}^\varphi\|_\cut
\le\|W_M^\varphi-B^{S_0}_{a,b}\|_\cut+\|B^{S_0}_{a,b}-W_{M^S_{a,b}}^\varphi\|_1
\le33\delta .
\]
\end{proof}

\begin{proof}[Proof of Theorem~\ref{thm:rate}]
By Theorem~\ref{thm:reduction}(a),
$\mathrm{LT}_{n,p}(\eta)=(1+o(1))\frac{n^2p}2V(r)$, and $r\ge r_*$ if and
only if $\eta\ge\eta_*$. If $r\ge r_*$, then the constant graphon $r$ is a
minimizer by Theorem~\ref{thm:main}(a) and~(c), so $V(r)=h(r)$. If $r<r_*$,
then $B_{a_r,b_r}$ is a minimizer and $(a_r,b_r)$ attains the minimum in
\eqref{eq:Vbal} by Theorem~\ref{thm:main}(b), so
$V(r)=F(B_{a_r,b_r})=\Vbal(r)$ by \eqref{eq:Bdens}.
\end{proof}

\begin{proof}[Proof of Theorem~\ref{thm:structure}]
Recall that
$(a_{r_*},b_{r_*})=(a_*,b_*)$, and fix $\varepsilon\in(0,1)$. Let
$\varepsilon'=\varepsilon(b_r-a_r)/66$ if $r\le r_*$, and
$\varepsilon'=\varepsilon/2$ if $r>r_*$. Condition on
$T_n\le\eta\E[T_n]$. By Theorem~\ref{thm:reduction}(b) with $\varepsilon'$
in place of $\varepsilon$, the following holds with probability $1-o(1)$.
There are a minimizer $W$ at $r$ and a measure preserving bijection
$\varphi$ with $\|(p^{-1}W_A)^\varphi-W\|_\cut\le2\varepsilon'$. Consider
such an outcome. By Theorem~\ref{thm:main}, $W$ is either the constant
graphon $r$ or a relabeling of $B_{a_r,b_r}$. The first case is possible
only if $r\ge r_*$, and the second only if $r\le r_*$.

In the first case, $p^{-1}W_A-r=W_{A/p-r\1\1^{\top}}$. Hence
Lemma~\ref{lem:stepcut} gives $\|A/p-r\1\1^{\top}\|_\cut=
\|(p^{-1}W_A)^\varphi-r\|_\cut\le2\varepsilon'\le\varepsilon$. That is,
$\|A-rp\1\1^{\top}\|_\cut\le\varepsilon p$. In the second case,
Lemma~\ref{lem:relabelB} gives a measurable $S_0\subset[0,1]$ with
$|S_0|=\frac12$ and $W=B^{S_0}_{a_r,b_r}$ a.e. We apply
Lemma~\ref{lem:rounding} to $M=A/p$ with $2\varepsilon'$ in place of
$\delta$. It gives $S\subset[n]$ with $\big||S|-\frac
n2\big|\le32\varepsilon'n/(b_r-a_r)\le\varepsilon n$ and
$\|A/p-M^S_{a_r,b_r}\|_\cut\le66\varepsilon'\le\varepsilon$. That is,
$\|A-pM^S_{a_r,b_r}\|_\cut\le\varepsilon p$. For $\eta>\eta_*$ only the
first case occurs, and for $\eta<\eta_*$ only the second. For $\eta=\eta_*$
one of the two cases occurs. This proves (a), (b) and (c).
\end{proof}

\section{Validated numerics}\label{app:numerics}

Four statements rest on finite computations, namely Lemmas~\ref{lem:digits}
and~\ref{lem:penconstants} and Propositions~\ref{prop:cover}
and~\ref{prop:family}. We reduce each of them to finitely many checks of three kinds. Checks of the first kind concern the signs of explicit real numbers and
validated solutions of linear systems of size four. Checks of the second kind
give lower bounds for explicit functions of a few variables on boxes, with
local Taylor bounds near prescribed zeros. Checks of the third kind use exact
rational geometry in the plane. We describe the arithmetic in Appendix~\ref{appA:arith}, verify the two
propositions in Appendices~\ref{app:cover} and~\ref{app:family}, and treat
the constants in Appendix~\ref{app:constants}. We found the
multipliers of the certificates by numerical optimization, and the checks
take them as exact inputs. The programs performing the checks are available as supplementary material and at
\url{https://github.com/fifalsp/sparse-triangle-lower-tail-verification}.

\subsection{Interval arithmetic and lower bounds on boxes}\label{appA:arith}

Interval arithmetic (see, e.g., \cite{MKC09,Tuc11}) works with closed
intervals, which we write in bold. Each real operation is replaced by an
interval operation whose result contains every value of the real operation
at points of the arguments. By this inclusion property, the interval
evaluation of an expression on a box encloses its range on the box.
Johansson's Arb library~\cite{Joh17} implements interval arithmetic with this property in midpoint
radius (ball) form at arbitrary precision. We use ball arithmetic at a
precision of at least $100$ bits. The branch and bound of
Lemma~\ref{lem:comp-bnb} uses intervals with binary64 endpoints and outward rounding, with
an a priori bound for the rounding error of long sums. Statements about
polygons are checked in exact rational arithmetic, in which binary64 numbers
enter exactly as dyadic rationals. A strict inequality is accepted only when its
enclosure has the required sign. For an interval $\mathbf x$ let
$|\mathbf x|=\max\{|x|:x\in\mathbf x\}$, and let $\inf\mathbf x$ be its lower
end.

We solve linear systems by a standard validated method (see, e.g.,
\cite{MKC09}). Let $\mathbf M$ be a square matrix and $\mathbf r$ a vector of
balls, let $R$ be an approximate inverse of the midpoint of $\mathbf M$, and
let $x_0$ be an approximate solution. The method succeeds if ball arithmetic
gives $\|I-RM\|_\infty\le\beta<1$ for all $M\in\mathbf M$. Then every
$M\in\mathbf M$ is invertible, since $Mv=0$ implies $v=(I-RM)v$. The identity
$x-x_0=(I-RM)(x-x_0)+R(r-Mx_0)$ gives
$|x-x_0|_\infty\le|R(r-Mx_0)|_\infty/(1-\beta)$ for the solution of $Mx=r$.
Ball arithmetic bounds the right side uniformly for $M\in\mathbf M$ and
$r\in\mathbf r$.

We next turn enclosures of a function and its derivatives into lower
bounds on a box. Let $\Phi:\R^3\to\R$ be twice continuously differentiable,
let $B$ be a box, let $c\in B$, and let $r\in\R^3$ satisfy
$|q_i-c_i|\le r_i$ for all $q\in B$ and $i=1,2,3$. Let the intervals
$\mathbf V$ and $\mathbf g_i$ contain $\Phi(c)$ and $\partial_i\Phi(c)$, and
let $\mathbf V_B$, $\mathbf G_i$ and $\mathbf K_{ij}=\mathbf K_{ji}$ contain
the ranges of $\Phi$, $\partial_i\Phi$ and $\partial_i\partial_j\Phi$ on $B$.
Let $L_0=\inf\mathbf V_B$ and
\begin{align}
L_1&=\inf\mathbf V-\sum_{i=1}^3r_i|\mathbf G_i|,\notag\\
L_2&=\inf\mathbf V-\sum_{i=1}^3r_i|\mathbf g_i|+\frac12\sum_{i=1}^3\min\{0,\inf\mathbf K_{ii}\}r_i^2
-\sum_{1\le i<j\le3}|\mathbf K_{ij}|r_ir_j .\label{appA:L12}
\end{align}
If every symmetric matrix $H$ with $H_{ij}\in\mathbf K_{ij}$ for all $i,j$ is
positive definite, let $y^*\in\R^3$ and let $L_3$ be a number such that
\begin{equation}\label{appA:L3}
L_3\le\inf\mathbf V+g^{\top}y^*+\frac12y^{*\top}Hy^*+\sum_{i=1}^3\min_{s\in\{-1,1\}}(g+Hy^*)_i(sr_i-y^*_i)
\end{equation}
for every such $H$ and every $g\in\R^3$ with $g_i\in\mathbf g_i$. Otherwise let
$L_3=-\infty$.

\begin{lemma}\label{lem:comp-box}
In the setting of \eqref{appA:L12} and \eqref{appA:L3},
$\Phi(q)\ge\max\{L_0,L_1,L_2,L_3\}$ for every $q\in B$.
\end{lemma}

\begin{proof}
The bound $L_0$ holds by definition. Let $q\in B$ and $y=q-c$, so that
$|y_i|\le r_i$. Since $B$ is convex, Taylor's theorem on the segment from $c$
to $q$ gives points $\xi,\xi'\in B$ with
$\Phi(q)=\Phi(c)+\nabla\Phi(\xi)^{\top}y=\Phi(c)+g^{\top}y+\frac12y^{\top}Hy$,
where $g=\nabla\Phi(c)$ and $H=\nabla^2\Phi(\xi')$. Since
$\Phi(c)\in\mathbf V$, $\partial_i\Phi(\xi)\in\mathbf G_i$,
$g_i\in\mathbf g_i$ and $H_{ij}\in\mathbf K_{ij}$, bounding the terms one by
one gives
$\Phi(q)\ge L_1$ and $\Phi(q)\ge L_2$. If $L_3>-\infty$, then $H$ is positive
definite, so $\Theta(y)=g^{\top}y+\frac12y^{\top}Hy$ is convex and
$\Theta(y)\ge\Theta(y^*)+(g+Hy^*)^{\top}(y-y^*)$. The right side is affine in
each $y_i\in[-r_i,r_i]$, so it is at least
$\Theta(y^*)+\sum_i\min_{s\in\{-1,1\}}(g+Hy^*)_i(sr_i-y^*_i)$. Since
$\Phi(q)=\Phi(c)+\Theta(y)$, \eqref{appA:L3} gives $\Phi(q)\ge L_3$.
\end{proof}

We certify the positive definiteness required for $L_3$ and for
Lemma~\ref{lem:taylor} by the following criterion.

\begin{lemma}\label{lem:comp-chol}
Let $\mathbf A$ be a symmetric $3\times3$ array of intervals. Run the
Cholesky recursion $L_{jj}=(A_{jj}-\sum_{k<j}L_{jk}^2)^{1/2}$ and
$L_{ij}=(A_{ij}-\sum_{k<j}L_{ik}L_{jk})/L_{jj}$ for $i>j$ in interval
arithmetic, with $\mathbf A_{ij}$ in place of $A_{ij}$. If every radicand has
a strictly positive lower end, then every real symmetric matrix $A$ with
$A_{ij}\in\mathbf A_{ij}$ for all $i,j$ is positive definite.
\end{lemma}

\begin{proof}
By the inclusion property, the exact recursion for $A$ has positive
radicands, so $A=LL^{\top}$ with $L$ lower triangular and invertible, and $A$
is positive definite.
\end{proof}

\medskip
\noindent\textbf{Branch and bound.} For $q\in\R^3$ let $q^\downarrow$ be the
vector of the coordinates of $q$ in decreasing order, and recall the ordered
set $\mathcal D_Q$ of Section~\ref{s5:relax}. The procedure takes as input a
permutation invariant, twice continuously differentiable function
$\Phi:\R^3\to\R$, together with enclosures of its value and of its
partial derivatives of order at most two at points and on boxes. It also
takes a dyadic number $Q>0$, a required lower bound $\varepsilon$, and
acceptance boxes $A_1,\dots,A_s$ with dyadic corners. It starts from the
boxes $I_i\times I_j\times I_k$ with $i\ge j\ge k$, where $I_1,\dots,I_N$
subdivide $[0,Q]$ into consecutive intervals with dyadic endpoints. These
boxes cover $\mathcal D_Q$. A box $B=\prod_i[l_i,u_i]$ is accepted if it lies
in some $A_j$, or if $\max\{L_0,L_1,L_2,L_3\}\ge\varepsilon$ for the bounds of
Lemma~\ref{lem:comp-box}, computed with the midpoint $c$ of $B$ and radii
$r_i\ge\max\{c_i-l_i,u_i-c_i\}$. The positive definiteness needed for $L_3$
is certified by Lemma~\ref{lem:comp-chol}. A box that is not accepted is
bisected at a dyadic point of a longest edge. A subbox with $u_1<l_2$ or
$u_2<l_3$ is discarded, since it is disjoint from $\mathcal D_Q$. Hence, at
every stage, the accepted boxes and the remaining boxes together cover
$\mathcal D_Q$ (the covering invariant). The procedure succeeds when no box
remains.

\begin{lemma}\label{lem:comp-bnb}
Let $\Phi$, $Q$, $\varepsilon$ and $A_1,\dots,A_s$ be inputs for which the
branch and bound succeeds. Then $\Phi(q)\ge\varepsilon$ for every
$q\in[0,Q]^3$ with $q^\downarrow\notin A_1\cup\dots\cup A_s$.
\end{lemma}

\begin{proof}
Let $q\in[0,Q]^3$ with $q^\downarrow\notin A_1\cup\dots\cup A_s$. Since
$q^\downarrow\in\mathcal D_Q$ and no box remains, the covering invariant
puts $q^\downarrow$ in an accepted box $B$. This box lies in no $A_j$, since it contains
$q^\downarrow$. Hence it was accepted by its bounds, and
Lemma~\ref{lem:comp-box} and permutation invariance give
$\Phi(q)=\Phi(q^\downarrow)\ge\varepsilon$.
\end{proof}

\subsection{Verification of Proposition~\ref{prop:cover}}\label{app:cover}

We use the notation of Section~\ref{sec:threshold}. The certificate data are
listed at the end of this subsection.

\medskip
\noindent\textbf{Constants.} By Appendix~\ref{app:constants}, $z_*$ lies in a
short interval with rational endpoints. Ball evaluation of $a$, $b$, $\mu$
and $\rho$ on this interval encloses $a_*$, $b_*$, $\mu_*$ and $r_*$. From these we
obtain enclosures of $t_*$, $m_*$, $F_*$, $\bar k_*$, $\ell_*$, $q_a$,
$q_b$ and $\omega_*$. Let $E$ be the
box with dyadic corners given by the enclosure of $\pi_*=(m_*,F_*)$.

\medskip
\noindent\textbf{The certificate functions.} Let $\theta\in\R^{12}$,
$\lambda\ge0$ and $\pi=(m,F)$ be given through enclosures, and let
$\Phi=\Phi_{\theta,\lambda,\pi}$. With $x_i=e^{-q_i}$ we have $u_i=x_i$,
$\log u_i=-q_i$, $h(u_i)=1-(1+q_i)x_i$, $k(u_i)=2-2x_i-q_ix_i$,
$u_1u_2u_3/u_i=x_jx_l$ for $\{i,j,l\}=\{1,2,3\}$, and
$g_*(u_i)=x_i+\omega_*h(u_i)$. By \eqref{eq:features} and \eqref{s5:Psi},
$\Phi$ is a polynomial in $q_1,q_2,q_3,x_1,x_2,x_3$ of degree at most $2$
in each $q_i$ and at most $4$ in each $x_i$. Its coefficients are polynomials
in the entries of $\theta$ and in $\lambda$, $\omega_*$, $m$, $F$, $\mu(\pi)$,
$\bar k(\pi)$ and $t_*$. Since $\partial_{q_i}x_i=-x_i$, every partial
derivative of $\Phi$ is a polynomial of the same kind. We enclose the
coefficients in ball arithmetic, using \eqref{s5:params} at a rational vertex
and $\mu(\pi_*)=\mu_*$ and $\bar k(\pi_*)=\bar k_*$ (Lemma~\ref{lem:apex}(a))
at $\pi_*$. Let $\alpha,\beta\ge0$ be
integers. If $\alpha\beta=0$, the function $q^\alpha e^{-\beta q}$ is
monotone on $[0,\infty)$. Otherwise it is increasing on $[0,\alpha/\beta]$ and
decreasing on $[\alpha/\beta,\infty)$. Hence its range on an interval
$[l,u]\subset[0,\infty)$ is enclosed by its values at $l$ and $u$, together
with its maximum $(\alpha/(\beta e))^\alpha$ if $\alpha/\beta\in[l,u]$.
Products of these ranges over the three coordinates enclose the monomials.
Summation with the coefficients then encloses the value and the partial
derivatives of $\Phi$ on a box or at a point.

\medskip
\noindent\textbf{The tangency equations.} Let $K$ be one of the triangles
$\mathrm{F1},\dots,\mathrm{F6}$. The tangency equations \eqref{s5:tangency}
form a linear system in $\theta_{K,0},\theta_{K,2},\theta_{K,3},\theta_{K,4}$.
Its coefficients and right sides are values and first derivatives of
explicit functions at $P_1$ and $P_2$. In the verification we write these
equations in the coordinates $q$, that is, at $\hat P_1$ and $\hat P_2$. The
resulting system is equivalent, since $\partial_{q_i}\Phi_{\theta,\lambda,\pi}(q)$
equals $-u_i\partial_{u_i}\Psi_{\theta,\lambda}(u;\pi)$ with
$u_i=e^{-q_i}>0$. The validated method for linear systems of
Appendix~\ref{appA:arith} succeeds for each of the six systems. Hence each
system has exactly one solution, and the returned balls enclose
it. All later evaluations of $\Psi_K$ use these balls and therefore
enclose the quantities of the exact certificate.

\medskip
\noindent\textbf{Item~(i), the region.} Since the region $\Omega$ has a curved
lower boundary, we first enclose it in a union $\Omega_+$ of two
rectangles. Let $m_0=0.2182$, $m_1=0.2186$, $m_2=0.3301$, $m_3=0.5095$, $L_1=0.391$,
$L_2=0.3355$ and $F^+=0.4627799$, as exact decimal numbers, and let
\[
\Omega_+=\big([m_1,m_2]\times[L_1,F^+]\big)\cup\big([m_2,m_3]\times[L_2,F^+]\big).
\]
For a rational interval $I\subset(0,1)$ and a ball $\mathbf I$ containing
it, let $\zeta=\max\{0,\inf(\mathbf I^3-t_*)\}$. Since $\beta>0$ and
$y\mapsto y_+^{2/3}$ is nondecreasing, the number
$\inf h(\mathbf I)+\max\{0,\inf\beta(\mathbf I)\}\zeta^{2/3}$ is a lower bound for
$F_{\mathrm{B1}}$ on $I$. Using these bounds on intervals obtained by repeated
bisection at rational midpoints, we prove in ball arithmetic that $0<m_0\le M_1$, that $F_{\mathrm{B1}}>F^+$
on $[m_0,m_1]$, that $F_{\mathrm{B1}}\ge L_1$ on $[m_1,m_2]$ and
$F_{\mathrm{B1}}\ge L_2$ on $[m_2,m_3]$, that $m_G\le m_3$, and that $F_*<F^+$. Let
$(m,F)\in\Omega$. Then $m\ge M_1\ge m_0$. Moreover, $m\le m_1$ is impossible,
since it would give $F_{\mathrm{B1}}(m)>F^+>F_*\ge F$. Hence
$m_1<m\le m_G\le m_3$. For the index $i$ with $m\in[m_i,m_{i+1}]$ we then
have $L_i\le F_{\mathrm{B1}}(m)\le F\le F_*<F^+$. Therefore
$\Omega\subset\Omega_+\cap\{F\le F_*\}$.

\medskip
\noindent\textbf{Item~(i), the triangles $\mathrm{F1},\dots,\mathrm{F6}$.}
Since the point $\pi_*$ is known only through an enclosure, we handle the
triangles $\mathrm{F1},\dots,\mathrm{F6}$ by the following elementary lemma.
For $p\in\R^2$ we write $p_2$ for its second coordinate.

\begin{lemma}\label{lem:fan}
Let $n\ge2$, $W\subset\R^2$ a convex polygon with vertices
$w_0,\dots,w_n$ in counterclockwise order, and $A\in W$. If $(w_0)_2>A_2$ and
$(w_n)_2>A_2$, then every $p\in W$ with $p_2\le A_2$ lies in
$\bigcup_{j=0}^{n-1}\mathrm{conv}\{A,w_j,w_{j+1}\}$.
\end{lemma}

\begin{proof}
Let $w_{n+1}=w_0$. Since $A\in W$ and $W$ is convex, the triangles
$\mathrm{conv}\{A,w_j,w_{j+1}\}$, $j=0,\dots,n$, cover $W$. Indeed, the ray
from $A$ through a point $p\in W\setminus\{A\}$ leaves $W$ through an edge
$[w_j,w_{j+1}]$, and $p$ lies between $A$ and the exit point. Now let $p\in W$ with
$p_2\le A_2$, and suppose that $p=s_0A+s_1w_n+s_2w_0$ with nonnegative
weights of sum one. Then $0\ge p_2-A_2=s_1((w_n)_2-A_2)+s_2((w_0)_2-A_2)$, so
$s_1=s_2=0$ by the hypotheses. Thus $p=A\in\mathrm{conv}\{A,w_0,w_1\}$.
\end{proof}

Let $W=\mathrm{conv}\{w_0,\dots,w_6\}$ for the points $w_j$ of
Section~\ref{s5:cover}, with indices taken modulo $7$. In
exact rational arithmetic we check the following.
\begin{enumerate}[label=(\arabic*)]
\item For every $j$ and every $k\notin\{j,j+1\}$, the point $w_k$ lies strictly
to the left of the directed line from $w_j$ to $w_{j+1}$.
\item The four corners of $E$ lie strictly to the left of each of these seven
directed lines.
\item The second coordinates of $w_0$ and $w_6$ exceed the upper end of the
enclosure of $F_*$.
\end{enumerate}
By (1), the points $w_0,\dots,w_6$ are in strictly convex position in
counterclockwise order. Hence they are the vertices of $W$ in this order, and
$W$ is the intersection of the closed half planes to the left of the seven
directed lines. The conditions in (2) are affine in the point, so they hold on
all of $E$. Hence $E\subset W$, and in particular $\pi_*\in W$. By (3) and
Lemma~\ref{lem:fan} with $A=\pi_*$,
$W\cap\{F\le F_*\}\subset\mathrm{F1}\cup\dots\cup\mathrm{F6}$.

\medskip
\noindent\textbf{Item~(i), the area argument.} Let $\mathcal K_0$ be the set
of the six cells other than $\mathrm{F1},\dots,\mathrm{F6}$, and let
$C=W\cup\bigcup_{K\in\mathcal K_0}K$. As in (1), we check that
the listed vertices of each cell of $\mathcal K_0$ are in strictly convex
counterclockwise position. Hence each cell of $\mathcal K_0$, like $W$, is the
intersection of the closed half planes to the left of its directed edges.
For a convex polygon $S$ and closed half planes $H_1,\dots,H_k$, the set
$S\setminus(H_1\cap\dots\cap H_k)$ is contained in the union of the closed
convex polygons $S\cap H_1\cap\dots\cap H_{l-1}\cap\overline{\R^2\setminus H_l}$,
$l=1,\dots,k$. Starting from each of the two rectangles of $\Omega_+$, we
subtract $W$ and the cells of $\mathcal K_0$ one at a time in this way from
every piece of positive area obtained so far. Pieces of area
zero are set aside. In exact rational arithmetic, no piece of positive area
remains after the last subtraction. Thus, for each rectangle $R$ of
$\Omega_+$, the set $R\setminus C$ is covered by finitely many closed polygons
of area zero. Since $C$ is closed, $R\setminus C$ is relatively open in $R$,
and a nonempty relatively open subset of a rectangle has positive area. Hence
$R\setminus C$ is empty, so $\Omega_+\subset C$ and
\[
\Omega\subset\Omega_+\cap\{F\le F_*\}\subset\big(W\cap\{F\le F_*\}\big)\cup\bigcup_{K\in\mathcal K_0}K
\subset\bigcup_{K\in\mathcal K}K,
\]
which is item~(i).

\medskip
\noindent\textbf{Item~(ii).} Let $K$ be a cell, and let $m'_K$ and $F'_K$ be the
minima of the first and of the second coordinates over the rational vertices
of $K$. If $K$ is one of the triangles $\mathrm{F1},\dots,\mathrm{F6}$, its
vertex $\pi_*$ is known only through an enclosure. The minima are then also
taken over the lower ends of the enclosures of $m_*$ and $F_*$. We check that
$m'_K>0$ and that some rational vertex of $K$ has second coordinate less
than $1$. We also check in ball arithmetic that
$((1-F'_K)^2-m_K'^2)/(2t_*)\le Q_K$. Then $0<m'_K\le m_K$
and $F'_K\le F_K<1$, so $(1-F_K)^2\le(1-F'_K)^2$ and $m_K^2\ge m_K'^2$. Hence
$((1-F_K)^2-m_K^2)/(2t_*)\le((1-F'_K)^2-m_K'^2)/(2t_*)\le Q_K$.

\medskip
\noindent\textbf{Item~(iii).} Let $K$ be a cell and let $v\ne\pi_*$ be a
vertex of $K$. The branch and bound succeeds for
$\Phi_{\theta_K,\lambda_K,v}$ with $Q=Q_K$, no acceptance box and the
required lower bound $\varepsilon_{K,v}$ of Table~\ref{tab:cellparams}. By
Lemma~\ref{lem:comp-bnb},
$\Psi_K(u;v)\ge\varepsilon_{K,v}>0$ for all $u\in[e^{-Q_K},1]^3$.

\medskip
\noindent\textbf{Item~(iv), the cells with vertex $\pi_*$.} Exact comparisons
show that no vertex of a cell of $\mathcal K_0$ lies in $E$. Since
$\pi_*\in E$, the cells with vertex $\pi_*$ are the six triangles
$\mathrm{F1},\dots,\mathrm{F6}$.

\medskip
\noindent\textbf{Item~(iv), local positivity.} Let $K$ be one of the triangles
$\mathrm{F1},\dots,\mathrm{F6}$, and let $\Phi=\Phi_{\theta_K,\lambda_K,\pi_*}$.
By the tangency equations \eqref{s5:tangency} and
Lemmas~\ref{lem:apex}(b) and~\ref{lem:patterns}(b), $\Phi(\hat P_i)=0$ and $\nabla\Phi(\hat P_i)=0$ for
$i=1,2$. Let $\delta_1=0.03$, let $\delta_2$ be the radius of
Table~\ref{tab:cellparams}, and let
$\mathcal C_{\delta_i}(\hat P_i)=\{q:|q-\hat P_i|_\infty\le\delta_i\}$. To check the hypotheses of
Lemma~\ref{lem:taylor}, we choose a
binary64 number $\lambda_0>0$ and subtract it from the diagonal entries of
the enclosure of $\nabla^2\Phi$ on the enclosure of $\hat P_i$. The resulting
array satisfies the hypothesis of Lemma~\ref{lem:comp-chol}. Hence
$\nabla^2\Phi(\hat P_i)-\lambda_0I$ is positive definite. The third
partial derivatives are enclosed on finitely many boxes whose union contains
$\mathcal C_{\delta_i}(\hat P_i)$ for every $\hat P_i$ in its enclosure. By symmetry of
the derivatives, the sum in Lemma~\ref{lem:taylor} is a sum over the ten
triples $j_1\le j_2\le j_3$ with multiplicities $1$, $3$ or $6$. Replacing each
maximum by the largest bound of its enclosures gives an admissible $S_3$.
Finally, we check that $\delta_iS_3<3\lambda_0$. Then
Lemma~\ref{lem:taylor} gives $\Phi>0$ on
$\mathcal C_{\delta_i}(\hat P_i)\setminus\{\hat P_i\}$.

\medskip
\noindent\textbf{Item~(iv), away from $\hat P_1$ and $\hat P_2$.} For $i=1,2$
let $A_i$ be a box with binary64 corners such that $A_i\subset \mathcal C_{\delta_i}(\hat P_i)$.
This inclusion is checked by exact comparison of the corners of $A_i$ with
the enclosure of $\hat P_i$. The branch and bound succeeds for $\Phi$ with
$Q=Q_K$, the acceptance boxes $A_1,A_2$ and the required lower bound
$\varepsilon_0=10^{-12}$. Hence $\Phi\ge\varepsilon_0$ on
$\mathcal D_{Q_K}\setminus(A_1\cup A_2)$ by Lemma~\ref{lem:comp-bnb}.

\medskip
\noindent\textbf{Item~(iv), conclusion.} By \eqref{s5:patterns},
$\hat P_1,\hat P_2\in\mathcal D_{Q_K}$. A point of $\mathcal D_{Q_K}$ outside
$A_1\cup A_2$ has $\Phi\ge\varepsilon_0$, and a point of $A_i$ lies in
$\mathcal C_{\delta_i}(\hat P_i)$, where $\Phi>0$ except at $\hat P_i$. Hence
$\Phi\ge0$ on $\mathcal D_{Q_K}$, with equality
exactly at $\hat P_1$ and $\hat P_2$. Now let $u\in[e^{-Q_K},1]^3$, and let
$q=-\log u$, taken coordinatewise. Then $q^\downarrow\in\mathcal D_{Q_K}$. By
symmetry, $\Psi_K(u;\pi_*)=\Phi(q^\downarrow)\ge0$. Equality holds if and only
if $q^\downarrow\in\{\hat P_1,\hat P_2\}$, that is,
$u\in\mathcal O(P_1)\cup\mathcal O(P_2)$. This completes the verification of
Proposition~\ref{prop:cover}.

\medskip
\noindent\textbf{Conventions for the data.} Each decimal in
Tables~\ref{tab:cells} to~\ref{tab:coef3} denotes the binary64 number nearest
to it, with ties resolved to the even significand. The same convention
applies to the numbers $\delta_1=0.03$ and $\varepsilon_0=10^{-12}$. A
significand followed by \texttt{e-5}, for example, is multiplied by $10^{-5}$
before this conversion. These decimal strings are therefore exact inputs. The
numbers $Q_K=N_K/1024$, the constants of $\Omega_+$ and the points $y_\pm$ of
the proof of Lemma~\ref{lem:digits} are exact rational numbers. Irrational constants
are given by enclosures. For reference, we give some of them rounded to ten
significant digits, namely $m_*=0.4331904817$, $F_*=0.4627797519$, $t_*=0.009217366870$,
$\bar k_*=1.029589270$, $\ell_*=0.004185359336$, $q_a=4.072542812$,
$q_b=0.1632875256$, $\omega_*=0.9231213775$ and
$g_*(a_*)=g_*(b_*)=0.8603923637$.

\medskip
\noindent\textbf{The cell data.} For each cell $K$,
Table~\ref{tab:cellparams} gives the weight $\lambda_K$, the number
$Q_K=N_K/1024$ and the required lower bounds $\varepsilon_{K,v}$ at the
vertices $v\ne\pi_*$. These vertices are taken in the order of
Table~\ref{tab:cells}. For the triangles $\mathrm{F1},\dots,\mathrm{F6}$, the
table also gives the radius $\delta_2$ of the cube
$\mathcal C_{\delta_2}(\hat P_2)$ in the local positivity test.
Tables~\ref{tab:coef1}, \ref{tab:coef2} and~\ref{tab:coef3} give the
multipliers, indexed as in \eqref{eq:features}. For the triangles
$\mathrm{F1},\dots,\mathrm{F6}$, the entries marked $*$ form
the unique solution of the tangency equations \eqref{s5:tangency}, given
$\lambda_K$ and the other entries.

\begin{table}[ht]
\centering
\small
\caption{The weights $\lambda_K$, the numbers $Q_K=N_K/1024$, the required
lower bounds $\varepsilon_{K,v}$ at the vertices $v\ne\pi_*$, and the radii
$\delta_2$ at $\hat P_2$ for the triangles $\mathrm{F1},\dots,\mathrm{F6}$.}
\label{tab:cellparams}
\begin{tabular}{lrrlr}
\toprule
Cell & $\lambda_K$ & $N_K$ & $\varepsilon_{K,v}$ & $\delta_2$\\
\midrule
F1 & \texttt{0.863} & 6574 & \texttt{27e-6}, \texttt{27e-6} & \texttt{0.0005}\\
F2 & \texttt{0.884} & 7377 & \texttt{42e-6}, \texttt{33e-6} & \texttt{0.00025}\\
F3 & \texttt{0.907} & 8084 & \texttt{64e-6}, \texttt{45e-4} & \texttt{0.001}\\
F4 & \texttt{0.975} & 7856 & \texttt{14e-5}, \texttt{14e-5} & \texttt{0.0005}\\
F5 & \texttt{0.596} & 6916 & \texttt{25e-6}, \texttt{91e-7} & \texttt{0.001}\\
F6 & \texttt{0} & 5786 & \texttt{31e-5}, \texttt{49e-5} & \texttt{0.002}\\
S1 & \texttt{0.96} & 13486 & \texttt{22e-6}, \texttt{36e-5}, \texttt{13e-6}, \texttt{44e-6} & \\
S2 & \texttt{0.931} & 11092 & \texttt{62e-5}, \texttt{51e-4}, \texttt{85e-5}, \texttt{52e-5} & \\
S3 & \texttt{0.857} & 9120 & \texttt{35e-5}, \texttt{45e-4}, \texttt{22e-5}, \texttt{30e-5} & \\
R1 & \texttt{0.164} & 6428 & \texttt{30e-4}, \texttt{39e-4}, \texttt{42e-4}, \texttt{38e-4} & \\
B1 & \texttt{1} & 18018 & \texttt{48e-5}, \texttt{14e-3}, \texttt{48e-5}, \texttt{56e-5} & \\
B2 & \texttt{0.268} & 18516 & \texttt{49e-6}, \texttt{18e-3}, \texttt{66e-6}, \texttt{35e-6} & \\
\bottomrule
\end{tabular}
\end{table}

\begin{table}[ht]
\centering
\scriptsize
\caption{The multipliers $\theta_0,\dots,\theta_3$ of the edge density, the
triangle density, the cost and the Euler--Lagrange equation tested with $1$.}
\label{tab:coef1}
\begin{tabular}{lrrrr}
\toprule
Cell & $\theta_{0}$ & $\theta_{1}$ & $\theta_{2}$ & $\theta_{3}$\\
\midrule
F1 & $*$ & $-3.847289099$ & $*$ & $*$\\
F2 & $*$ & $-3.807984177$ & $*$ & $*$\\
F3 & $*$ & $-4.086655072$ & $*$ & $*$\\
F4 & $*$ & $-5.874814709$ & $*$ & $*$\\
F5 & $*$ & $0.8832292362$ & $*$ & $*$\\
F6 & $*$ & $5.974396099$ & $*$ & $*$\\
S1 & $0.9591238521$ & $-4.720264495$ & $0.9999999996$ & $-0.0002044256$\\
S2 & $1$ & $-4.885203951$ & $0.944880088$ & $-0.001059666$\\
S3 & $1$ & $-3.727050662$ & $0.8955832267$ & $-0.0020680644$\\
R1 & $-0.6193316787$ & $-0.9944176312$ & $1$ & $0.1451388316$\\
B1 & $0.9184812124$ & $-6.023851127$ & $0.973384689$ & $-3.75629\cdot10^{-5}$\\
B2 & $0.15101055$ & $-3.004348059$ & $0.0883321776$ & $-7.28673\cdot10^{-5}$\\
\bottomrule
\end{tabular}
\end{table}

\begin{table}[ht]
\centering
\scriptsize
\caption{The multipliers $\theta_4,\dots,\theta_7$ of the Euler--Lagrange
equation tested with $u,u^2,u^3,u\log u$.}
\label{tab:coef2}
\begin{tabular}{lrrrr}
\toprule
Cell & $\theta_{4}$ & $\theta_{5}$ & $\theta_{6}$ & $\theta_{7}$\\
\midrule
F1 & $*$ & $0.6360661654$ & $-0.6671015204$ & $0.1527685923$\\
F2 & $*$ & $0.7230151623$ & $-0.7353524276$ & $0.1399870154$\\
F3 & $*$ & $0.4804572605$ & $-0.6501312432$ & $0.1736701008$\\
F4 & $*$ & $0.5791540482$ & $-0.7206926774$ & $0.2455163716$\\
F5 & $*$ & $0.4115910168$ & $-0.3925509204$ & $0.0376207471$\\
F6 & $*$ & $-0.9999999995$ & $0.1826854517$ & $0.1184025582$\\
S1 & $-0.1194056551$ & $0.9055043484$ & $-0.7993075377$ & $0.1217645018$\\
S2 & $-0.0584493425$ & $0.7437038004$ & $-0.6030087567$ & $0.1048504918$\\
S3 & $0.011889639$ & $0.5688704906$ & $-0.5362219892$ & $0.1080609339$\\
R1 & $-0.8217036627$ & $-1$ & $-1$ & $-1$\\
B1 & $-0.0761922072$ & $1$ & $-0.9443873005$ & $0.1875254028$\\
B2 & $-0.0951384338$ & $0.4226885325$ & $-0.266968611$ & $0.0360012488$\\
\bottomrule
\end{tabular}
\end{table}

\begin{table}[ht]
\centering
\scriptsize
\caption{The multipliers $\theta_8,\dots,\theta_{11}$ of the vertex identity
tested with $u,u^3,u^4,u\log u$.}
\label{tab:coef3}
\begin{tabular}{lrrrr}
\toprule
Cell & $\theta_{8}$ & $\theta_{9}$ & $\theta_{10}$ & $\theta_{11}$\\
\midrule
F1 & $0.0112250185$ & $-0.9674784838$ & $0.444859491$ & $0.5126250828$\\
F2 & $-0.143086236$ & $-0.689495481$ & $0.2912439931$ & $0.5279458108$\\
F3 & $-0.0962137491$ & $-0.9999999823$ & $0.5313007707$ & $0.6193744493$\\
F4 & $-0.08794637$ & $-0.6540880492$ & $0.1181991075$ & $0.8602008183$\\
F5 & $-0.0487809445$ & $-0.2132777015$ & $-0.0402843975$ & $0.2540583174$\\
F6 & $0.4056336334$ & $-0.39555763$ & $0.0974899047$ & $0.513893443$\\
S1 & $-0.1399729896$ & $0.3940162703$ & $-0.764047086$ & $0.7448317535$\\
S2 & $-0.3381151534$ & $0.4398461833$ & $-0.6093917014$ & $0.4999374123$\\
S3 & $-0.1816877103$ & $-0.2634340813$ & $-0.027956505$ & $0.406967943$\\
R1 & $-0.9999999997$ & $-1$ & $-0.099647709$ & $-0.2329168544$\\
B1 & $-0.0736682908$ & $0.2174062631$ & $-0.7236908707$ & $1$\\
B2 & $-0.4512883381$ & $0.8862667776$ & $-0.6829202941$ & $0.2220880755$\\
\bottomrule
\end{tabular}
\end{table}

\subsection{Verification of Proposition~\ref{prop:family}}\label{app:family}

By the reduction in the proof of Proposition~\ref{prop:family}, it remains to
verify a statement on a compact domain. For $L\in[3.9,9.5]$ let $z=e^{-L}$.
Then $\mu$, the coefficients $c_0,\dots,c_4$ of \eqref{eq:certcoef}
and \eqref{eq:certpsi}, $A=\mu c_4$, $q_a=L(1+z^2)/(1-z)^2$ and
$q_b=2Lz/(1-z)^2$ are explicit functions of $L$. Let $x_i=e^{-q_i}$
and $s=q_1+q_2+q_3$. The function $\Phi(L,q)=\Phi_z(q)$ of \eqref{s8:phi} is
\begin{equation}\label{appA3:phi}
\Phi(L,q)=\sum_{k=1}^7\alpha_k(L)f_k(q)+e^{-s},
\end{equation}
where $\alpha=(3c_0,c_1,-c_2,-c_3,c_4,\mu c_3,-A)$ and
\[
f=\Big(1,\sum_ix_i,\sum_iq_ix_i,\sum_iq_i,\sum_iq_i^2x_i,\sum_{i<j}x_ix_j,se^{-s}\Big).
\]
Let $P_1(L)=(q_a,q_a,q_a)$, $P_2(L)=(q_a,q_b,q_b)$ and
\[
\mathcal N=\{(L,q)\in\R^4:3.9\le L\le9.5\text{ and }0\le q_3\le q_2\le q_1\le12\}.
\]
For every $L\in[3.9,9.5]$, both $P_1(L)$ and $P_2(L)$ satisfy the conditions
on $q$ in the definition of $\mathcal N$, since $q_a-q_b=L>0$ and $q_a+2<12$
by the proof of Proposition~\ref{prop:family}. Let $\delta_0$ be the
binary64 number nearest to $1/10$, which is slightly larger than $1/10$.

\begin{proposition}[computer assisted]\label{appA3:compact}
Let $(L,q)\in\mathcal N$. Then $\Phi(L,q)\ge0$, with equality if and only if
$q=P_1(L)$ or $q=P_2(L)$.
\end{proposition}

We prove the proposition in two parts. First, we show that the Hessian in
$q$ is positive definite on
the cubes of radius $\delta_0$ about the zeros $P_1(L)$ and $P_2(L)$. Second, we
cover $\mathcal N$ by boxes that either lie in one of these cubes or
carry a positive lower bound for $\Phi$. All enclosures use ball arithmetic
on boxes with binary64 corners.

\medskip
\noindent\textbf{Enclosures in $L$.} For a function $\psi$ of $L$ let
$\psi^{[n]}=\psi^{(n)}/n!$. Taylor arithmetic (see, e.g., \cite{MKC09,Tuc11})
propagates $\psi^{[0]},\dots,\psi^{[3]}$ through the explicit formulas
\eqref{eq:branch} and \eqref{eq:certcoef} by the standard recurrences for
sums, products, quotients and the exponential. Let $J$ be an interval with
binary64 center $L_0$ and radius at most $w$. Started at the point $L_0$,
Taylor arithmetic encloses $\alpha_k^{[n]}(L_0)$. Started at the interval
$J$, it encloses $\alpha_k^{[n]}(\xi)$ for all $\xi\in J$. This gives
enclosures over $J$ of $\alpha_k$, $\mu$, $c_i$, $A$, $q_a$ and $q_b$, by
Taylor's theorem with the Lagrange remainder or by direct evaluation on $J$.
The denominators in \eqref{eq:certcoef} do not vanish, since $\Gamma<-2.6$
by the proof of Lemma~\ref{lem:certcoef}.

\medskip
\noindent\textbf{The Hessian.} Differentiating \eqref{s8:gradq} shows that the
Hessian $H=\nabla_q^2\Phi(L,q)$ is given by
\begin{equation}\label{appA3:hessian}
H_{ii}=W_i+\mu c_3x_i\sum_{j\ne i}x_j+G,\quad
H_{ij}=\mu c_3x_ix_j+G\quad(i\ne j),
\end{equation}
where $W_i=x_i(c_1+c_2(2-q_i)+c_4(q_i^2-4q_i+2))$ and $G=e^{-s}(1+2A-As)$.
Let $X=\mathrm{diag}(x_1,x_2,x_3)$. We certify positive definiteness for the
scaled matrix $\tilde H=X^{-1/2}HX^{-1/2}$, which is positive definite
if and only if $H$ is. By \eqref{appA3:hessian}, for
$\{i,j,k\}=\{1,2,3\}$,
\begin{align*}
\tilde H_{ii}&=\frac{W_i}{x_i}+\mu c_3(x_j+x_k)+x_jx_k(1+2A-As),\\
\tilde H_{ij}&=e^{-(q_i+q_j)/2}\big(\mu c_3+x_k(1+2A-As)\big).
\end{align*}
The test of positive definiteness on a box uses a real $3\times3$ matrix $T$,
computed from a numerical eigendecomposition at the center of the box. The test
succeeds for a symmetric
array $\tilde{\mathbf H}$ of intervals if the ball evaluation of
$\mathbf M=T^{\top}\tilde{\mathbf H}T$ satisfies
$\inf\mathbf M_{ii}-\sum_{j\ne i}|\mathbf M_{ij}|>0$ for $i=1,2,3$. Then, for
every symmetric $\tilde H$ with entries in $\tilde{\mathbf H}$, the symmetric
matrix $N=T^{\top}\tilde HT$ satisfies $N_{ii}>\sum_{j\ne i}|N_{ij}|$ for
$i=1,2,3$. By the Gershgorin theorem, the eigenvalues of $N$ are positive.
Hence $N$ is positive definite, so $T$ is invertible and $\tilde H$ is
positive definite.

\begin{lemma}[computer assisted]\label{appA3:hess}
Let $j\in\{1,2\}$, $L\in[3.9,9.5]$ and $q\in\R^3$ with
$|q-P_j(L)|_\infty\le\delta_0$. Then the matrix $\nabla_q^2\Phi(L,q)$ is
positive definite.
\end{lemma}

Since $\Phi(L,\cdot)$ is symmetric, $\nabla_q^2\Phi(L,q^\sigma)$ arises from
$\nabla_q^2\Phi(L,q)$ by permuting rows and columns. The point $P_1(L)$ is
fixed by all permutations, and $P_2(L)$ by the transposition of the last two
coordinates. Hence it suffices to treat $q=P_j(L)+v$ with $v_1\ge v_2\ge v_3$
for $j=1$ and with $v_2\ge v_3$ for $j=2$. We cover the set of these $(L,v)$ with
$L\in[3.9,9.5]$ and $|v|_\infty\le\delta_0$ by boxes $J\times V$ obtained by
bisection. On each box, $q_i$ is enclosed by $\mathbf P_{j,i}+V_i$,
where $\mathbf P_{j,i}$ encloses the $i$th coordinate of $P_j(L)$ for all
$L\in J$. The coefficients are enclosed over $J$. The box is accepted if the
test succeeds for the resulting enclosure of $\tilde H$. Every box is
accepted, which proves Lemma~\ref{appA3:hess}.

\medskip
\noindent\textbf{The cover of $\mathcal N$.} We maintain a
stack of closed boxes $J\times B$ in the coordinates $(L,q_1,q_2,q_3)$.
Initially the stack holds boxes whose union contains
$[3.9,9.5]\times[0,12]^3$. The procedure keeps the invariant that every
point of $\mathcal N$ lies in an accepted box or in a box on the stack. Before
a box $J\times\prod_i[l_i,u_i]$ is tested, its bounds are tightened by
$l_2\gets\max\{l_2,l_3\}$, $l_1\gets\max\{l_1,l_2\}$,
$u_2\gets\min\{u_2,u_1\}$ and $u_3\gets\min\{u_3,u_2\}$, in this order. This
removes no point of $\mathcal N$. A box in which some interval becomes empty
is discarded. A box is accepted if it passes one of the following tests.
Otherwise one coordinate interval is bisected at its binary64 midpoint, and
the two closed subboxes are pushed onto the stack.
\begin{enumerate}[label=(C\arabic*)]
\item \emph{Membership.} For some $j\in\{1,2\}$, the enclosure of
$[l_i,u_i]-\mathbf P_{j,i}$ lies in $(-\delta_0,\delta_0)$ for $i=1,2,3$. Then
$|q-P_j(L)|_\infty<\delta_0$ for every $(L,q)$ in the box.
\item \emph{Positivity.} Let $L_0$ be the center of $J$, and let $w$ be an
upper bound for its radius. Let $q^0$ be the binary64 midpoint of $B$, and
let $r_i\ge\max\{q^0_i-l_i,u_i-q^0_i\}$. For
$(L,q)\in J\times B$ let $\epsilon=L-L_0$ and $d=q-q^0$. Taylor's theorem
gives $\xi\in B$ and $L',L''\in J$ with
\begin{align*}
\Phi(L,q)&=\Phi(L,q^0)+\nabla_q\Phi(L,q^0)^{\top}d+\frac12d^{\top}\nabla_q^2\Phi(L,\xi)d,\\
\Phi(L,q^0)&=\phi_0(L_0)+\phi_1(L_0)\epsilon+\phi_2(L_0)\epsilon^2+\phi_3(L')\epsilon^3,\\
\partial_{q_i}\Phi(L,q^0)&=\gamma_{i,0}(L_0)+\gamma_{i,1}(L_0)\epsilon+\gamma_{i,2}(L'')\epsilon^2,
\end{align*}
where $\phi_n$ and $\gamma_{i,n}$ are the $n$th Taylor coefficients in $L$ of
$\Phi(\cdot,q^0)$ and $\partial_{q_i}\Phi(\cdot,q^0)$. By
\eqref{appA3:phi}, they are combinations of the $\alpha_k^{[n]}$ with the
weights $f_k(q^0)$ and $\partial_{q_i}f_k(q^0)$. The terms of $e^{-s}$ enter
only for $n=0$. Let $e_i$ be an upper bound for
$w|\gamma_{i,1}(L_0)|+w^2\sup_J|\gamma_{i,2}|$. Let $G_i^+$ be an upper bound
for $|\gamma_{i,0}(L_0)|+e_i$, and let $G_i^-$ be a lower bound for
$|\gamma_{i,0}(L_0)|-e_i$. Then $G_i^-\le|\partial_{q_i}\Phi(L,q^0)|\le G_i^+$
for $L\in J$. Let
$H_i^-=\inf_{J\times B}\partial_{q_i}^2\Phi$, and let $\beta_i$ be the largest of
the numbers $-G_i^+r_i+\frac12\min\{0,H_i^-\}r_i^2$, $-(G_i^+)^2/(2H_i^-)$ and
$-G_i^+r_i+\frac12H_i^-r_i^2$. Here the second number is used only if
$H_i^->0$, and the third only if $H_i^->0$ and $G_i^-\ge H_i^-r_i$. Then every
$g$ with $G_i^-\le|g|\le G_i^+$, every $H\ge H_i^-$ and every $|y|\le r_i$
satisfy $gy+\frac12Hy^2\ge\beta_i$. For the third number, note that if
$|g|<Hr_i$, the minimum over $y$ is
$-g^2/(2H)>-|g|r_i/2\ge-|g|r_i+\frac12H_i^-r_i^2$. The box is accepted if the
following expression has a positive lower end:
\begin{align*}
&\phi_0(L_0)-w|\phi_1(L_0)|+\min\{0,\inf\phi_2(L_0)\}w^2-w^3\sup_J|\phi_3|\\
&\quad+\sum_{i=1}^3\beta_i-\frac12\sum_{i\ne j}r_ir_j\sup_{J\times B}|\partial_{q_i}\partial_{q_j}\Phi|.
\end{align*}
This expression is a lower bound for $\Phi$ on the box. The suprema over
$J\times B$ are bounded by \eqref{appA3:hessian} with the coefficients
enclosed over $J$.
\end{enumerate}
The procedure ends with an empty stack.

\begin{proof}[Proof of Proposition~\ref{appA3:compact}]
Let $(L,q)\in\mathcal N$. By the invariant and the empty final stack, $(L,q)$
lies in an accepted box. If the box was accepted by (C2), then
$\Phi(L,q)>0$. Otherwise it was accepted by (C1) for some $j$. Then
$v=q-P_j(L)$ satisfies $|v|_\infty<\delta_0$. By Lemma~\ref{lem:certcoef} and
$\partial_{q_i}\Phi_z=-x_i\partial_{x_i}\Psi_z$, we have $\Phi(L,P_j(L))=0$ and
$\nabla_q\Phi(L,P_j(L))=0$. Hence
\[
\Phi(L,q)=\int_0^1(1-\tau)v^{\top}\nabla_q^2\Phi(L,P_j(L)+\tau v)v\dd\tau .
\]
The segment from $P_j(L)$ to $q$ lies in the cube
$|q'-P_j(L)|_\infty\le\delta_0$. Hence, if $v\ne0$, Lemma~\ref{appA3:hess}
shows that the integrand is positive for $\tau\in[0,1)$. Thus $\Phi(L,q)>0$ unless
$q=P_j(L)$, and $\Phi(L,P_j(L))=0$.
\end{proof}

\subsection{Enclosures of the numerical constants}\label{app:constants}

We enclose the constants of Sections~\ref{sec:balanced}
and~\ref{sec:smalltail} using the signs of explicit functions of one variable,
combined with monotonicity statements proved in those sections.

\medskip
\noindent\textbf{Enclosure of $z_*$.} Recall the function
$D(z)=h(\rho(z))-F(B_z)$ of Section~\ref{sec:balanced}, where
$F(B_z)=(h(a(z))+h(b(z)))/2$ by \eqref{eq:Bdens}. It is a composition of
rational operations, $\exp$, $\log$ and the cube root, and by
Theorem~\ref{thm:crossing}(c), $D>0$ on $(0,z_*)$ and $D<0$ on $(z_*,1)$.
For the exact decimals $y_\pm$ of the proof of Lemma~\ref{lem:digits}, ball
evaluation gives $D(y_-)>0>D(y_+)$, so $z_*\in(y_-,y_+)$. The evaluation of
$a(z)$, $b(z)$, $\mu(z)$ and $\rho(z)$ on $[y_-,y_+]$ gives intervals inside
the four intervals of Lemma~\ref{lem:digits}. For
Proposition~\ref{prop:cover}, we bisect the interval $[y_-,y_+]$ repeatedly
at rational midpoints, each time keeping the subinterval whose endpoints have
certified signs $D>0$ and $D<0$. This gives a shorter interval with rational
endpoints that contains $z_*$. Evaluation on this interval yields the
enclosures used in Appendix~\ref{app:cover}.

\medskip
\noindent\textbf{Enclosure of $z_c$.} In the proof of
Lemma~\ref{lem:penconstants} we enclose $z_c$ between
$z_{\mathrm{lo}}=8.352460192\cdot10^{-5}$ and
$z_{\mathrm{hi}}=8.352460193\cdot10^{-5}$, using the signs of the function
$\Delta$ of Lemma~\ref{lem:tie},
\[
\Delta(z)=\Phi_{\mu(z)/3}(a(z),b(z))-C_0(\mu(z)/3),\quad C_0(\gamma)=h(s_\gamma)+\gamma s_\gamma^3,
\]
where $s_\gamma\in(0,1)$ is the zero of $f(s,\gamma)=\log s+3\gamma s^2$. By
Lemmas~\ref{lem:zc} and~\ref{lem:fold}, $\Delta<0$ on $(0,z_c)$
and $\Delta>0$ on $(z_c,10^{-3}]$. Hence
$\Delta(z_{\mathrm{lo}})<0<\Delta(z_{\mathrm{hi}})$ implies
$z_c\in(z_{\mathrm{lo}},z_{\mathrm{hi}})$. The function $f$ is increasing in
$s\in(0,1)$ and in $\gamma>0$. Let $[\gamma_1,\gamma_2]\subset(0,\infty)$
contain $\gamma$, and let $s_1<s_2$ be rational numbers with
$f(s_1,\gamma_2)<0<f(s_2,\gamma_1)$. Then $s_\gamma\in(s_1,s_2)$ by this
monotonicity. Hence interval evaluation of $h(s)+\gamma s^3$ on
$[s_1,s_2]\times[\gamma_1,\gamma_2]$ encloses $C_0(\gamma)$. We apply this
with $[\gamma_1,\gamma_2]$ the ball of $\mu(z)/3$ at $z_{\mathrm{lo}}$, and
then at $z_{\mathrm{hi}}$. This gives balls for $\Delta(z_{\mathrm{lo}})$ and
$\Delta(z_{\mathrm{hi}})$ with the stated signs. Evaluation of $\mu(z)/3$,
$\rho(z)$ and $\Phi_{\mu(z)/3}(a(z),b(z))$ on
$[z_{\mathrm{lo}},z_{\mathrm{hi}}]$ encloses $\gamma_c$, $r_-=\rho(z_c)$ and
$C(\gamma_c)$ (Proposition~\ref{prop:pen}(b)). The same bracketing of
$s_\gamma$, applied with the enclosure of $\gamma_c$, encloses
$r_+=s_{\gamma_c}$.

\end{document}